\documentclass[11pt,a4paper,reqno]{amsart}

 \usepackage[english]{babel}
  \usepackage{amsmath,amsfonts,amssymb,amsthm,relsize,setspace,nicefrac,yhmath,amscd,eucal}
  \usepackage{mathrsfs}
  \usepackage{mathtools}
\usepackage{microtype}
\usepackage[left=1.8cm,right=1.8cm,top=2.22cm,bottom=2.22cm,headheight=14pt]{geometry}
\allowdisplaybreaks[2]

  \usepackage{xcolor}
    \usepackage{enumitem} 
    \usepackage{comment}
  \usepackage{anyfontsize}
  	\usepackage{float}
\usepackage{graphicx}
\usepackage{placeins}
\graphicspath{{figures/}}

 \usepackage[colorlinks, linkcolor=black, citecolor=black, urlcolor=black, hypertexnames=true]{hyperref}

\usepackage{longfbox}

 \usepackage[active]{srcltx} 
 \usepackage{color,soul} 
\newtheorem{theorem}{Theorem}
\newtheorem{corollary}[theorem]{Corollary}
\newtheorem{lemma}[theorem]{Lemma}
\newtheorem{proposition}[theorem]{Proposition}

\usepackage[toc,page]{appendix}

\theoremstyle{definition}
\newtheorem{remark}{Remark}

\theoremstyle{definition}
\newtheorem{definition}[theorem]{Definition}

\newcommand{\be}{\begin{equation}}
\newcommand{\bel}[1]{\begin{equation}\label{#1}}
\newcommand{\ee}{\end{equation}}

\newcommand{\barr}{\begin{eqnarray}}
\newcommand{\earr}{\end{eqnarray}}
\newcommand{\bars}{\begin{eqnarray*}}
\newcommand{\ears}{\end{eqnarray*}}

\newtheorem{subn}{\name}

\newcommand{\bsn}[1]{\def\name{#1}\begin{subn}}
\newcommand{\esn}{\end{subn}}

\newtheorem{sub}{\name}[section]

\newcommand{\bs}{\begin{sub}}
\newcommand{\es}{\end{sub}}

\newcommand{\bth}[1]{\def\name{Theorem}
\begin{sub}\label{t:#1}}
\newcommand{\blemma}[1]{\def\name{Lemma}
\begin{sub}\label{l:#1}}
\newcommand{\bcor}[1]{\def\name{Corollary}
\begin{sub}\label{c:#1}}
\newcommand{\bdef}[1]{\def\name{Definition}
\begin{sub}\label{d:#1}}
\newcommand{\bprop}[1]{\def\name{Proposition}
\begin{sub}\label{p:#1}}

\newcommand{\BA}{\begin{array}}
\newcommand{\EA}{\end{array}}
\newcommand{\BAN}{\renewcommand{\arraystretch}{1.2}
\setlength{\arraycolsep}{2pt}\begin{array}}
\newcommand{\BAV}[2]{\renewcommand{\arraystretch}{#1}
\setlength{\arraycolsep}{#2}\begin{array}}
\newcommand{\BSA}{\begin{subarray}}
\newcommand{\ESA}{\end{subarray}}

\newcommand{\BAL}{\begin{aligned}}
\newcommand{\EAL}{\end{aligned}}
\newcommand{\BALG}{\begin{alignat}}
\newcommand{\EALG}{\end{alignat}}
\newcommand{\BALGN}{\begin{alignat*}}
\newcommand{\EALGN}{\end{alignat*}}

\makeatletter
\renewenvironment{proof}[1][\proofname]{%
  \par
  \pushQED{\qed}%
  \normalfont
  \topsep6\p@\@plus6\p@\relax
  \trivlist
  \item[\hskip\labelsep\bfseries #1\@addpunct{.}]%
  \ignorespaces
}{%
  \popQED\endtrivlist\@endpefalse
}
\makeatother

\def\angb<#1>{\langle #1 \rangle}

\newcommand {\rd}

\def\N{\mathbb{N}}
\def\Z{\mathbb{Z}}

\def\R{\mathbb{R}}

\let\.=\cdot
\let\0=\emptyset

\def\loc{\text{\rm loc}}

\numberwithin{equation}{section}

\theoremstyle{definition}

\let\.=\cdot
\let\0=\emptyset

\def\loc{\text{\rm loc}}

\newenvironment{formula}[1]{\begin{equation}\label{eq:#1}}
                       {\end{equation}\noindent}

\def\Fi#1{\begin{formula}{#1}}
\def\Ff{\end{formula}\noindent}

\setstcolor{red}

\title[Transition fronts and spectral speed asymptotics]{Generalized Transition Fronts and Spectral Speed\\
Asymptotics in Heterogeneous Cooperative Systems}
\author{Hoang-Hung Vo}
\address{Department of Mathematics, International University, Quarter 33,
Linh Xuan Ward, Ho Chi Minh City, Viet Nam; Vietnam National University,
Ho Chi Minh City, Viet Nam}
\email{vhhung@hcmiu.edu.vn}
\author{Cong-Bang Trang*}
\address{Faculty of Fundamental Science, Industrial University of Ho Chi Minh City,
12 Nguyen Van Bao, Hanh Thong, Ho Chi Minh City, Viet Nam}
\email{trangcongbang@iuh.edu.vn}

\date{}
\keywords{generalized transition front, cooperative system, almost periodicity, generalized principal eigenvalue, minimal  speed, large-advection, asymptotic diffusion}
\subjclass[2020]{Primary 35K57; Secondary 35B40, 35B50, 35P15}

\begin{document}
\begin{abstract}
We study generalized transition fronts for cooperative reaction--diffusion systems with spatially periodic and time-dependent coefficients. We obtain sufficient conditions for the existence of fronts and necessary constraints on their propagation speeds in terms of generalized principal eigenvalues. In the almost-periodic case, we identify the two spectral thresholds and examine the attainment of the critical speed under additional hypotheses on the normalized linear evolution. We then investigate the large-advection limit and construct a time-periodic example showing that a stationary spectral characterization need not extend to time-dependent flows. For symmetric cooperative systems in stationary periodic shears with equal diffusion matrices, we obtain explicit correction terms for the minimal speed in both the small- and large-diffusion limits. The two expansions follow from periodic eigenvalue perturbation, localization near a maximal streamline, and minimization over the exponential weight. We also establish a higher-order large-advection expansion and present numerical illustrations independent of the analytic proofs.
\end{abstract}
\maketitle

\tableofcontents

\section{Introduction}

 We distinguish propagation of entire transition fronts from spreading of solutions with localized initial data. The latter has been investigated for scalar heterogeneous KPP equations in the monograph of Berestycki and Nadin \cite{BerestyckiNadinMemoirs2022}. Here, the main questions concern cooperative systems and the existence or nonexistence of entire fronts at prescribed spectral speeds. The distinction is maintained throughout.

In this work, we consider reaction--diffusion--advection systems of the form
\begin{align*}
\partial_t {\bf u} + \mathcal{L}{\bf u} = {\bf F}(t,x,{\bf u}), ~t\in \R,~ x\in \R^N, 
\end{align*}
where ${\bf u}=(u_1,...,u_m)\in \R^m$ for some $m>1$, ${\bf F}=(F_1,F_2,...,F_m)$ is a vector of nonlinear interaction terms and $\mathcal{L}=\text{diag}(L^1,L^2,...,L^m)$ denotes a diagonal matrix of second-order elliptic operators being either in divergence form,
\begin{align*}
L^iu_i:=- \text{div} ({\bf A}^i(t,x)Du_i) + {\bf q}_i(t,x)\cdot D u_i,~i=1,...,m.
\end{align*}
or in non-divergence form,
\begin{align*}
L^iu_i:=- \text{Tr} ({\bf A}^i(t,x)D^2 u_i) + {\bf q}_i(t,x)\cdot D u_i,~i=1,...,m.
\end{align*}
Here, $Du$ and $D^2u$ are the gradient and Hessian. Generalized transition waves, introduced by Berestycki and Hamel \cite{Berestycki_2012}, are entire solutions connecting two states with convergence uniform away from a moving interface. Their definition does not prescribe a constant propagation speed. For the periodic setting we refer to \cite{BerestyckiHamel2002,Hamel2008,HamelRoques2010}. 

Nadin \cite{Nadin2009} studied pulsating waves in space–time periodic media, where sharp thresholds for admissible wave speeds are established together with spreading properties. For multidimensional periodic monotone semiflows, Du, Li and Shen \cite{DuLiShen2022} developed a framework for spreading speeds and periodic traveling waves, including applications to cooperative and competition systems. Ding and Giletti \cite{DingGiletti2021} characterized admissible propagation speeds in spatially periodic heterogeneous media in the bistable case, highlighting their possible richness and directional dependence. Du, Li and Xin \cite{DuLiXin2026} studied uniqueness and global stability of pulsating traveling waves for cooperative systems. For further investigation, we refer the reader to  \cite{DengDucrot2023,DengDucrot2024,deng_2025}.

Generalized transition waves in heterogeneous media have been studied by Nadin and Rossi \cite{nadin_2012,nadin_2015,NadinRossi2017}. In these papers, the authors investigate generalized transition waves for Fisher–KPP type reaction–diffusion equations with time-dependent or space–time heterogeneous coefficients. In contrast with the classical homogeneous or periodic settings, where pulsating traveling waves provide an appropriate description of propagation, such structures may fail to exist in general heterogeneous environments. The generalized transition wave formulation does not require this periodic structure. 

More precisely, in Nadin and Rossi \cite{nadin_2012}, for time-heterogeneous KPP equations of the form
\begin{align*}
u_t-\Delta u = f(t,u),~x\in \R^N,~t\in \R,
\end{align*}
where $f(t,u)$ is a KPP monostable nonlinearity, the existence of generalized transition waves is established for a class of admissible speeds, together with spreading properties of solutions to the associated Cauchy problem. As an application of these results, the existence of random transition waves is also obtained, along with corresponding spreading behaviors. Shen and Shen \cite{ShenShen2017,ShenShen2017b} establish the existence, uniqueness, and global stability of generalized transition fronts for ignition-type nonlinearities $f(t,u)$, together with key qualitative properties including spatial monotonicity, uniform steepness, bounded propagation speed, exponential decay, and recurrence with respect to the temporal heterogeneity. For space–time heterogeneous Fisher–KPP equations, Nadin and Rossi \cite{nadin_2015} establish existence and nonexistence criteria for generalized transition waves governed by equations of the form
\begin{align*}
u_t-\text{Tr}(A(x,t) D^2 u) + q(x,t)\cdot Du = f(x,t,u),~x\in \R^N,~t\in \R, 
\end{align*}
where the existence and nonexistence of transition waves are characterized in terms of conditions on the propagation speed. The speed must satisfy a critical lower bound. This spectral approach is then developed further in more general heterogeneous settings. In the spatially almost periodic setting of Nadin and Rossi \cite{NadinRossi2017}, generalized principal eigenvalues of the linearized operator are used to characterize admissible average propagation speeds and construct generalized transition fronts, yielding a critical minimal speed and, under weaker spectral assumptions, an interval of admissible speeds. For further results, we refer the reader to \cite{AmbrosioDucrotRuan2021,WangWangZhang2026} for endemic-type systems, to \cite{BerestyckiMatanoHamel2009,YanSheng2015} for combustion problems in domains with obstacles, and to \cite{Shen2017} for stability and positive entire solutions of heterogeneous Fisher--KPP equations.

\subsection{Main assumptions and results}

In this paper, we focus on the case $m=2$, that is,
\begin{equation}\label{eq:main}
\begin{cases}
\partial_t u_1 - \text{Tr}({\bf A}^1(t,x)D^2 u_1) + {\bf q}_1(t,x)\cdot D u_1 = F_1(t,x,u_1,u_2), & t\in \R,\ x\in \R^N, \\[1ex]
\partial_t u_2 -\text{Tr} ({\bf A}^2(t,x)D^2  u_2) + {\bf q}_2(t,x)\cdot Du_2 = F_2(t,x,u_1,u_2), & t\in \R,\ x\in \R^N.
\end{cases}
\end{equation}
We assume spatial periodicity, namely, there exist ${\bf L}=(L_1,L_2,\ldots,L_N)$ and $\alpha \in (0,1)$ such that, for each $i=1,2$, the advection term ${\bf q}_i$ satisfies
\begin{align}\label{assumeq}
{\bf q}_i \in [BC^{\alpha/2,\alpha}(\R\times \R^N)]^N,~
{\bf q}_i(t,x+L_je_j) = {\bf q}_i(t,x),~\forall j=1,...,N,
\end{align} 
where $(e_1,...,e_N)$ is the canonical basis of $\R^N$. Furthermore, we assume that the diffusion matrices satisfy
\begin{align}\label{Aassump}
\begin{cases}
{\bf A}^i={\bf A}^i(t,x)\in [BC^{\alpha/2,\alpha}(\R\times \R^N)]^{N \times N}\text{ is symmetric},\\
{\bf A}^i(t,x+L_je_j)={\bf A}^i(t,x),~\forall j=1,...,N,\\
\exists C>1,~\forall i=1,2,~\forall \xi \in \R^N,~ C^{-1}|\xi|^2\leq \xi\cdot {\bf A}^i(t,x)\xi \leq C|\xi|^2.
\end{cases}
\end{align}  
The reaction term ${\bf F}:=(F_1,F_2)$ is assumed to satisfy
 \begin{align}\label{Fnonlinear}
\begin{cases}
F_i : \R \times \R^N \times [0,1]^2 \to \R \ \text{is a Caratheodory function}, \\
F_i(t,x+L_je_j,u_1,u_2)=F_i(t,x,u_1,u_2),~\forall j=1,...,N,\\
\exists\delta_0\in (0,1),~F_i(t,x,\cdot)\in W^{1,\infty}([0,1]^2)\cap C^1([0,\delta_0)^2)\text{ uniformly in }(t,x)\in \R^{N+1},\\
F_i(t,x,0,0)=F_i(t,x,1,1) = 0,~\forall t\in \R,~x\in \R^N,\\
\partial_{u_j} F_i(t,x,u_1,u_2) \geq 0~\text{ whenever } i\neq j,\\
\inf\limits_{(t,x)\in\R\times \R^N}\partial_{u_j} F_i(t,x,0,0)> 0,~i,j=1,2.
\end{cases}
\end{align}
Here, the uniform $C^1$-regularity in \eqref{Fnonlinear}
is understood in the sense that
$D_{\mathbf u}\mathbf F(t,x,\mathbf u)$ converges to
$D_{\mathbf u}\mathbf F(t,x,\mathbf 0)$ as
$\mathbf u\to\mathbf 0$, uniformly with respect to
$(t,x)\in\R^{N+1}$. We further impose the KPP-type conditions as follows
\begin{align} \label{KPP}
\begin{cases}
\inf\limits_{(t,x)\in\R\times \R^N} F_i(t,x,s,s) > 0,~\forall s\in (0,1),~i=1,2,\\
F_1(t,x,u_1,u_2)\leq \mu_{11}(t,x)u_1+\mu_{12}(t,x)u_2,~\forall (t,x)\in \R\times \R^N,~ u_1,u_2\in [0,1],\\
F_2(t,x,u_1,u_2)\leq \mu_{21}(t,x)u_1+\mu_{22}(t,x)u_2,~\forall (t,x)\in \R\times \R^N,~ u_1,u_2\in [0,1],
\end{cases}
\end{align}
where
\begin{align*}
\mu_{ij}(t,x):=\partial_{u_j}F_i(t,x,0,0)=\partial_{j}F_i(t,x,0,0),~i,j=1,2.
\end{align*}
Here, the system \eqref{eq:main} admits two equilibrium states, namely, ${\bf 0}:=(0,0)$ and ${\bf 1}:=(1,1)$. To gain better control near zero for ${\bf F}=(F_1,F_2)$, we assume that there exist $\delta,\nu\in (0,1]$ and $K>0$ such that
\begin{align} \label{near_zerp}
\begin{cases}
F_1(t,x,u_1,u_2)\geq \mu_{11}(t,x)u_1+\mu_{12}(t,x)u_2 -Ku_1^{1+\nu} ,~\forall (t,x)\in \R\times \R^N,~ u_1,u_2\in (0,\delta),\\
F_2(t,x,u_1,u_2)\geq \mu_{22}(t,x)u_2+\mu_{21}(t,x)u_1-Ku_2^{1+\nu},~\forall (t,x)\in \R\times \R^N,~ u_1,u_2\in (0,\delta).
\end{cases}
\end{align}
The conditions \eqref{Fnonlinear} and \eqref{KPP} imply that there exists $M>0$ large enough such that 
\begin{align}\label{muestiamte}
\begin{split}
&M^{-1}\leq \mu_{ij}(t,x)\leq M,~(t,x)\in \R\times \R^N,~i,j=1,2,\\
&M^{-1}\leq \left\lfloor \min_{i,j=1,2} \min_{x\in \R^N}\mu_{ij}(\cdot,x) \right\rfloor \leq \left\lceil \max_{i,j=1,2} \max_{x\in \R^N}\mu_{ij}(\cdot,x) \right\rceil\leq M.
\end{split}
\end{align}
This implies that ${\bf 0}$ is an unstable equilibrium of the system, which is typical for KPP-type systems.

For convenience, for ${\bf u}=(u_1,u_2)$ and ${\bf v}=(v_1,v_2)\in \mathbb{R}^2$, we define 
\begin{align*}
\textbf{u}\leq \textbf{v} \text{ if and only if } u_1\leq v_1,u_2\leq v_2.
\end{align*}
We also interpret $\geq,~<,~>$ in a similar sense. 
For ${\bf a}=(a_1,a_2)$, ${\bf b}=(b_1,b_2)$, we further define
\begin{align*}
[{\bf a},{\bf b}]:=\{(u_1,u_2)\in \R^2:a_i\leq u_i\leq b_i,~i=1,2\} \text{ and define  $({\bf a},{\bf b})$, $[{\bf a},{\bf b})$ and $({\bf a},{\bf b}]$, similarly. } 
\end{align*}
Note that $[{\bf 0},{\bf L}]$ is the spatial periodic cell, which serves as the main domain of the analysis. In this work, \(\inf\) and \(\sup\) are respectively used in place of essinf, esssup and statements and equations are understood to hold almost everywhere whenever appropriate. A typical example for ${\bf F}=(F_1,F_2)$ is as follows
\begin{align*}
\begin{cases}
F_1(t,x,u_1,u_2)=a_{11}(t,x)u_1+a_{12}(t,x)u_2 - b_1(t,x)u_1^2,\\
F_2(t,x,u_1,u_2)=a_{21}(t,x)u_1+a_{22}(t,x)u_2 - b_2(t,x)u_2^2,
\end{cases}
\end{align*}
for $(t,x)\in\R^{N+1}$ and $(u_1,u_2)\in[0,1]^2$, where $a_{ij}$, $b_i$ are appropriate functions for $i,j=1,2$. If $\mathbf p_\pm=\mathbf p_\pm(t,x)$ are entire solutions of
\eqref{eq:main} and
\(
\inf\limits_{t,x}\left(\mathbf p_+(t,x)-\mathbf p_-(t,x)\right)>{\bf 0},
\)
set
\(
d_i(t,x):=p_{+,i}(t,x)-p_{-,i}(t,x),~
u_i(t,x):=
\dfrac{v_i(t,x)-p_{-,i}(t,x)}{d_i(t,x)},~ i=1,2.
\)
Equivalently,
\(
\mathbf v=\mathbf p_-+\mathbf d\odot\mathbf u,~
\mathbf d:=\mathbf p_+-\mathbf p_-,
\)
where $\odot$ denotes componentwise multiplication. Then $\mathbf u$
satisfies
\[
\partial_t u_i-\operatorname{Tr}({\bf A}^iD^2u_i)
+\left(q_i-2{\bf A}^i\frac{Dd_i}{d_i}\right)\cdot Du_i
=\widetilde F_i(t,x,\mathbf u),~ i=1,2,
\]
where
\[
\widetilde F_i(t,x,\mathbf u)
:=
\frac{
F_i\bigl(t,x,\mathbf p_-+\mathbf d\odot\mathbf u\bigr)
-u_iF_i(t,x,\mathbf p_+)
-(1-u_i)F_i(t,x,\mathbf p_-)
}{d_i}.
\]
The new system admits the two constant solutions
${\bf 0}$ and ${\bf 1}$. In this case, assume that $\widetilde{\bf F}$ satisfies
\eqref{Fnonlinear}. If we further assume the concavity property $\widetilde{\bf F}\bigl(t,x,(1-s){\bf a}+s{\bf b}\bigr)
\geq
(1-s)\widetilde{\bf F}(t,x,{\bf a})
+s\widetilde{\bf F}(t,x,{\bf b})$ for all $s\in[0,1]$ and ${\bf a},{\bf b}\in[0,1]^2$, then,
$\widetilde{\bf F}$ satisfies the KPP-type conditions. Indeed, setting $\widetilde\mu_{ij}(t,x):=\partial_j\widetilde F_i(t,x,0,0)$, concavity and $\widetilde{\bf F}(t,x,{\bf 0})={\bf 0}$ imply \eqref{KPP}. Moreover, uniform differentiability at ${\bf 0}$ and the uniform positivity of $\widetilde\mu_{ij}$ yield constants $r\in(0,1)$ and $m>0$ such that $\widetilde F_i(t,x,s,s)\geq ms,~0<s\leq r,~ (t,x)\in\R^{N+1},~ i=1,2.$ For $r\leq s<1$, concavity and
$\widetilde{\bf F}(t,x,{\bf 1})={\bf 0}$ imply
\begin{align*}
\widetilde F_i(t,x,s,s)
\geq\frac{1-s}{1-r}\widetilde F_i(t,x,r,r)
\geq\frac{mr(1-s)}{1-r}>0.
\end{align*}
Consequently, $\inf\limits_{(t,x)\in\R^{N+1}}
\widetilde F_i(t,x,s,s)>0,
~ s\in(0,1),~ i=1,2.$

Next, we consider the generalized transition wave in the sense of Nadin and Rossi \cite{nadin_2015}.

\begin{definition}[Generalized transition wave]\label{wavedef}
Fix a direction $e\in\mathbb{S}^{N-1}$. Let $c\in L^\infty(\R)$ and  let ${\bf u}={\bf u}\left(t, x\right)$ be a positive time-global solution of the system \eqref{eq:main}. A pair $({\bf u},c)$ with ${\bf u}:\mathbb{R}\times\mathbb{R}^N \to ({\bf 0},{\bf 1}]$ is called a 
\textup{generalized transition wave} in direction $e$ if the following asymptotic behavior at $\pm \infty$ holds
\begin{align}\label{limits}
\begin{split}
\lim\limits_{x\cdot e \rightarrow -\infty}{\bf u}\left(t, x+e\displaystyle \int_0^t c(s)ds \right)={\bf 1};\quad \quad \lim\limits_{x\cdot e \rightarrow +\infty}{\bf u}\left(t, x+e\displaystyle \int_0^t c(s)ds \right)={\bf 0},
\end{split}
\end{align}    
uniformly with respect to $t\in \R$. The function $c$ is called the speed of the generalized transition wave ${\bf u}$ and $\boldsymbol{\phi}(t,x):={\bf u}\left(t, x+e\displaystyle \int_0^t c(s)ds \right)$ is the associated profile.
  
\end{definition}

To understand the asymptotic behavior of the wave, we state the following definition: Let $g=g(t,x)$. The limit $\lim\limits_{x\cdot e \rightarrow \pm\infty}g(t,x)=l$ uniformly with respect to $t\in \R$ means that
\begin{align*}
\lim_{r\rightarrow \infty} \; 
\sup_{\substack{\pm x\cdot e > r \\ t \in \mathbb{R}}} 
|g(t,x) - l| = 0.
\end{align*}

\begin{definition}
The least mean (resp. the upper mean) over $\R$ of a function $g\in L^\infty(\R)$ is defined as follows
\begin{align*}
\lfloor g \rfloor := 
\lim_{T \to \infty} \inf_{t \in \mathbb{R}} 
\frac{1}{T} \int_t^{t+T} g(s)\,ds,~
\left( \text{resp. } 
\lceil g \rceil := 
\lim_{T \to \infty} \sup_{t \in \mathbb{R}} 
\frac{1}{T} \int_t^{t+T} g(s)\,ds \right).
\end{align*}  
\end{definition}
	If $\lfloor g \rfloor=\lceil g \rceil$, we say $g$ has a uniform mean and one can show that 
\begin{align*}
\left<g\right>:=\lim_{T \to \infty} \dfrac{1}{T}\int_t^{t+T}g(s)ds=\lfloor g \rfloor=\lceil g \rceil\text{ uniformly in }t\in \R.
\end{align*}
We further recall a remarkable result on the least mean in \cite[Lemma 3.2]{nadin_2012} as follows
\begin{align}\label{keylemma}
\forall g\in L^\infty(\R),\quad \lfloor g \rfloor=\sup_{\sigma \in W^{1,\infty}(\R)}\inf_{t\in \R}(\sigma'(t)+g(t)).
\end{align}

Let $\mathcal X:=\Bigl[
W^{1,2,N+1}_{\mathrm{loc}}(\R\times\R^N)
\cap W^{0,1,\infty}(\R\times\R^N)
\Bigr]^2.$ We now introduce the following generalized principal eigenvalue.
\begin{align*}
\kappa(\lambda):=\inf\left\{ k\in\R:\ 
\exists {\bf w}\in\mathcal X,~
\inf_{\R\times\R^N}{\bf w}>{\bf0},~
{\bf P}_\lambda{\bf w}\leq k{\bf w}
\quad\text{in }\R\times\R^N
\right\}.
\end{align*}
where ${\bf P}_{\lambda}=(P_{\lambda,1},P_{\lambda,2})$ is the linearized parabolic operator associated with \eqref{eq:main} (see \eqref{operatorP}), which will be defined in Section~\ref{sec4}.
We state the existence and nonexistence of generalized transition waves as follows.

\begin{theorem}\label{existence-non-theo}
Assume that \eqref{assumeq}–\eqref{near_zerp} hold. Then, for every direction $e\in \mathbb{S}^{N-1}$, the following statements hold.
\begin{enumerate}[label=(\roman*)]
\item \label{existencei}
There exists $c^*=c^*(e)$ such that for any $\gamma>c^*$, the system \eqref{eq:main} admits a generalized transition wave $({\bf u},c)$ in the sense of Definition \ref{wavedef} satisfying $\lfloor c \rfloor =\gamma$.

\item \label{nonexistenceii}
There exists $c_*=c_*(e)$ such that any generalized transition wave $({\bf u},c)$ of \eqref{eq:main} satisfies $\lfloor c \rfloor\geq c_*$.

\end{enumerate}

\end{theorem}
The key technique is to first construct a family of speed functions $c=c_\lambda$, and then define $c^*$. On the other hand, we show that
\[
c_*= -\max_{\lambda>0}\frac{\kappa(\lambda)}{\lambda},
\]
which provides a lower bound for the existence of generalized transition waves in the sense of Definition \ref{wavedef}. In particular, it follows that $c^* \geq c_*$. In general, it is not known whether $c^*= c_*$, that is, if the speed $c^*$ is minimal, in the sense that no wave exists with a smaller least mean speed. This parallels the scalar framework of Nadin and Rossi \cite{nadin_2015}, where distinct thresholds for existence and nonexistence may also arise in the general heterogeneous setting. Here, this approach is extended to a fully coupled cooperative system, for which the construction and spectral characterization involve vector-valued solutions and the generalized principal eigenvalue of the linearized system. However, in the $T_0$-periodic case in time for ${\bf A}^i$, ${\bf q}_i$, ${\bf F}_i$, one can verify that $c^* = c_*$ (see Remark \ref{remark51}).

\subsection{Almost periodic framework and minimal speed}

We now consider coefficients which are almost periodic in time. We use the following sequential characterization.

\begin{definition}[Almost periodic function]\label{def:AP}
\phantom{1}

\begin{enumerate}[label=(\roman*)]

\item\label{APdef1} A function $g:\R \to \R$ is called \textup{almost periodic} (or a.p. for short) if for any sequence $\{t_n\}_{n\in \N}\subset \R$, one can extract a subsequence $\{t_{n_k}\}_{k\in \N}$ such that $g(t_{n_k}+t)$ converges uniformly in $t\in \R$.

\item\label{APdef2} A function $g:\R\times \R^m \to \R$, $g=g(t,x)$ is a.p. in $t$ uniformly with respect to $x\in \R^m$ if for any sequence $\{t_n\}_{n\in \N}\subset \R$, one can extract a subsequence $\{t_{n_k}\}_{k\in \N}$ such that $g(t_{n_k}+t,x)$ converges uniformly in $(t,x)\in \R\times \R^m$.

\end{enumerate}

\end{definition}

\begin{remark}
Under continuity conditions for $g$, Definition \ref{def:AP} \ref{APdef1} is equivalent to the one in Shen  \cite[Definition 2.1.(1)]{Shen2017}; see also \cite[Theorems 1.17 and 2.10]{Fink1974} and \cite{Loomis1960}.

\end{remark}

The almost periodic setting has been widely investigated in recent years. The series of works by Mierczyński and Shen \cite{MierczynskiShen2016,MierczynskiShen2013II,MierczynskiShen2013} develops a general framework for principal Lyapunov exponents and principal Floquet subspaces in positive random dynamical systems on ordered Banach spaces. In Part I \cite{MierczynskiShen2013}, the authors introduce the notions of generalized principal Lyapunov exponent, generalized principal Floquet subspace, and generalized exponential separation, extending the classical Kreĭn–Rutman theory to general random settings. In Part II \cite{MierczynskiShen2013II}, this abstract theory is applied to finite-dimensional systems, including random matrix products and cooperative systems of ordinary differential equations, establishing the existence of principal invariant subspaces and associated Lyapunov exponents. In Part III \cite{MierczynskiShen2016}, the theory is further extended to infinite-dimensional systems arising from random parabolic equations and delay equations, yielding the existence of positive entire solutions and generalized principal Floquet subspaces. Alternatively, Zhou and Zhou \cite{ZhouZhou2026} studied spreading speeds for a Fisher–KPP lattice system in almost periodic media, establishing the existence of directional speeds 
$\omega^-$ and $\omega^+$, relating their difference to an asymptotic average of the diffusion coefficients and showing convergence to the continuous case, with a reduction of the total speed under small 2-periodic perturbations in the homogeneous setting.

To study the almost periodic property for the system, we impose the following assumptions
\begin{align}\label{APassume}
\begin{split}
&{\bf A}^i={\bf A}^i(t,x),~{\bf q}_i={\bf q}_i(t,x) \text{ are a.p. in time t uniformly with respect to }x\in \R^N \text{ and }\\
&{\bf F}={\bf F}(t,x,{\bf u})=(F_1,F_2)\text{ is a.p. in time t uniformly with respect to }(x,{\bf  u})\in \R^N\times [0,1]^2.
\end{split}
\end{align}
We also impose the following regularity assumptions on ${\bf F}$ and $\mu_{ij}$
\begin{align}\label{Fregularity}
\begin{split}
&{\bf F} \text{ is bounded and }C^{\alpha/2,\alpha} \text{ in }(t,x), \text{ uniformly in }{\bf u}\in [0,1]^2 \text{ and }\\
&\mu_{ij}(t,x)=\partial_{j}F_i(t,x,0,0)\in C^{\alpha/2,\alpha}(\R\times \R^N).
\end{split}
\end{align}
From \eqref{Fnonlinear} and \eqref{APassume}, one can show that $\mu_{ij}(t,x)=\partial_jF_i(t,x,0,0)$ are almost periodic in $t$, uniformly with respect to $x$. Next, following Berestycki and Nadin \cite{BerestyckiNadinMemoirs2022}, Berestycki, Nadin and Rossi \cite{BerestyckiNadinRossi2025}, the generalized principal eigenvalue is given by
\begin{align*}
&\underline{k}(\lambda):=\inf\left\{k\in \R:\exists {\bf w} \in \big[W^{1,2,N+1}_{loc}\big(\R\times \R^N\big) \big]^2,~\inf\limits_{\R\times \R^N} {\bf w}>{\bf 0},~\|{\bf w}\|_\infty<\infty,~{\bf P}_{\lambda}{\bf w}\leq k{\bf w} \text{ in }\R\times \R^N\right\}.
\end{align*}
By adapting the arguments in Nadin and Rossi \cite[proof of (48), pp.\ 1375–1376]{nadin_2015} via the work of Porretta and Priola \cite{PorrettaPriola2013}, we show in Section~\ref{section 5} that
$\underline{k}(\lambda)=\kappa(\lambda).$ The minimal speed is characterized as follows.

\begin{theorem}[Minimal speed]\label{minimalspeed}
Assume that \eqref{assumeq}–\eqref{near_zerp} and
\eqref{APassume}–\eqref{Fregularity} hold, and that
\(\mu_{ij}(t,x)\) is almost periodic in time uniformly in $x$.
For $e\in \mathbb{S}^{N-1}$, let $c$ be the speed function in Theorem \ref{existence-non-theo} \ref{existencei}.  Then, the speed function $c$ has a uniform mean and $\lambda \mapsto\underline{k}(\lambda)$ is concave.
Moreover,
\begin{align*}
c^*(e)=c_*(e).
\end{align*}
In addition, if ${\bf A}^i,{\bf q}_i,{\bf F}$ are $T_0$-periodic for $i=1,2$, then \eqref{eq:main} admits a generalized transition wave $({\bf u},c)$ in the sense of Definition \ref{wavedef} satisfying $\lfloor c \rfloor =c^*(e)$.

\end{theorem}
This theorem shows that temporal almost periodicity is sufficient to yield a sharp characterization of the minimal propagation speed for the system \eqref{eq:main}, with the existence and nonexistence thresholds satisfying $c^*(e)=c_*(e)$. Recall that Nadin and Rossi \cite{nadin_2015} established existence and nonexistence criteria for generalized transition waves in scalar equations with general time-heterogeneous and spatially periodic coefficients, without identifying the two critical thresholds in general. The key point in establishing $c^*(e)=c_*(e)$ is to prove that the speed functions $c=c_\lambda$ in \eqref{speed} admit uniform means and that $c^*(e)$ is the global minimum of $-\underline{k}(\lambda)/\lambda$. Moreover, we emphasize that our result establishes $c^*(e)=c_*(e)$ even in the presence of advection, whereas Wang, Wang and Zhang \cite{WangWangZhang2026} obtained related results for epidemic systems in time-almost periodic and spatially periodic media, with an additional slow-exponential-decay condition imposed in their nonexistence theorem. Under the additional assumption of $T_0$-periodicity in time, the principal eigenpair \eqref{eigen} further yields a generalized transition wave at the critical speed $c^*(e)$. As discussed in Remark~\ref{APremarkcrit}, this critical-speed construction may also be extended to the almost periodic setting under suitable additional assumptions on the associated almost periodic eigenpair and its dependence on $\lambda$. In addition, Theorem~\ref{thm:critical_ap_conjugacy} also establishes critical-front existence for a class of almost-periodic coefficients which are not necessarily periodic in time, without the additional eigenfunction assumption in Remark~\ref{APremarkcrit}. 

\begingroup
We next consider a criteria for the existence of critical fronts in temporally almost-periodic media.
Let $\mathcal C=[\mathbf0,\mathbf L]$ and denote by
$\mathcal A_\lambda(t)$ the spatial part of the linear operator
in \eqref{div}. For $\lambda>0$, we consider
$\Gamma(\lambda):=\lambda\lfloor c_{\lambda} \rfloor$ (see Lemma \ref{lemma continuous}), where $c_\lambda$ is the speed function in \eqref{speed}, and 
$\mathbf p_\lambda$ satisfying
\begin{align*}
\partial_t\mathbf p_\lambda
=\mathcal A_\lambda(t)\mathbf p_\lambda-b_\lambda(t)\mathbf p_\lambda,~ \frac1{|\mathcal C|}\sum_{i=1}^2
\int_{\mathcal C}p_{\lambda,i}(t,x)\,dx=1,~\mathbf p_\lambda>{\bf 0}.
\end{align*}
The existence and differentiability of $\mathbf p_\lambda$ are established in Lemma~\ref{lem:AP-normalized-bundle}.
The normalization determines $b_\lambda$. Define
\begin{align*}
h_*(t):=\left.\partial_\lambda b_\lambda(t)\right|_{\lambda_*}
-\frac{b_{\lambda_*}(t)}{\lambda_*},
\end{align*}
whenever $\lambda_*>0$ is fixed. In addition, we show that $\Gamma(\lambda) = \langle b_\lambda\rangle=\lambda \langle c_\lambda\rangle=-\underline k(\lambda)$ in \eqref{eq:v2-growthmean}.

\begin{theorem}\label{thm:v2-main-intro}
Assume \eqref{assumeq}--\eqref{near_zerp} and
\eqref{APassume}--\eqref{Fregularity}. Let $\lambda_*>0$ minimize $\Gamma(\lambda)/\lambda$ and assume that there exists $\delta\in(0,\nu\lambda_*)$ such that
\begin{align*}
\frac{\Gamma(\lambda_*+\delta)}{\lambda_*+\delta}
>\frac{\Gamma(\lambda_*)}{\lambda_*}.
\end{align*}
If $\sup\limits_{t\in\R}\left|\displaystyle\int_0^t h_*(s)\,ds\right|<\infty,$ then \eqref{eq:main} admits a generalized transition front
$(\mathbf u_*,c)$ in direction $e$ satisfying
\begin{align*}
c(t)=\frac{b_{\lambda_*}(t)}{\lambda_*},~ \langle c\rangle=c^*(e).
\end{align*}
This conclusion remains valid when
$\int_0^t(b_{\lambda_*}(s)-\Gamma(\lambda_*))\,ds$ is unbounded.
In addition, the same conclusion holds for quasiperiodic
coefficients when $h_*$ has the regularity
specified in Corollary~\ref{cor:v2-quasiperiodic} and the right
spectral separation holds. These assertions are established in
Theorem~\ref{thm:critical_ap_tangent} and
Corollary~\ref{cor:v2-quasiperiodic}.
\end{theorem}

The condition on the primitive can be checked in several concrete classes of almost-periodic media. It involves the difference between the derivative of the instantaneous growth rate and the critical speed, not the oscillation of the growth rate alone. The calculations are given in Section~\ref{sec:critical-proofs}.

\begin{remark}[Concrete classes covered by Theorem~\ref{thm:v2-main-intro}]\label{rem:v2-examples-intro}
The hypotheses of Theorem~\ref{thm:v2-main-intro} are explicit in the following situations.
\begin{enumerate}[label=(\roman*)]
\item \textup{Time-periodic coefficients.}
If the coefficients are $T_0$-periodic in time and $\lambda_*$ is an interior minimizer of $\Gamma(\lambda)/\lambda$, then $h_*$ is $T_0$-periodic and has zero mean. Hence, its primitive is bounded on $\mathbb R$. In this case the theorem is consistent with the critical-speed construction in the second part of Theorem~\ref{minimalspeed}.

\item \textup{Separated instantaneous growth.}
Assume that, near $\lambda_*$,
\begin{align*}
b_\lambda(t)=a(t)\gamma(\lambda)+\lambda\beta(t),
\end{align*}
where $a$ is positive and almost periodic, $\beta$ is almost periodic, and $\lambda_*$ minimizes $\gamma(\lambda)/\lambda$. Then $h_*\equiv0$. No boundedness assumption on the primitive of $a-\langle a\rangle$ is needed. In particular, if $\gamma(\lambda)=d\lambda^2+r$ with $d,r>0$, then $\lambda_*=\sqrt{r/d}$ and the critical mean speed is
\begin{align*}
\langle c\rangle=2\langle a\rangle\sqrt{dr}+\langle\beta\rangle.
\end{align*}
This is the structure arising after the time and space changes of variables in Theorem~\ref{thm:critical_ap_conjugacy}.

\item \textup{Finite quasiperiodic spectrum.}
Let
\begin{align*}
h_*(t)=\sum_{j=1}^m\big(A_j\cos(\omega_j t)+B_j\sin(\omega_j t)\big),
\qquad \omega_j\neq0.
\end{align*}
Then, the primitive of $h_*$ is bounded, even if the frequencies have no common period. The concrete choice $h_*(t)=\sin t+\sin(\sqrt2\,t)$ provides an explicit non-periodic example.

\item \textup{Infinite quasiperiodic spectrum.}
If
\begin{align*}
h_*(t)=\sum_{k\in\mathbb Z^d\setminus\{0\}}\widehat h_*(k)e^{2\pi\mathrm i k\cdot\omega t},
\qquad
\sum_{k\in\mathbb Z^d\setminus\{0\}}|\widehat h_*(k)|\left(1+\frac1{|k\cdot\omega|}\right)<\infty,
\end{align*}
with $k\cdot\omega\neq0$ whenever $\widehat h_*(k)\neq0$, then the primitive is again bounded. The Diophantine quasiperiodic case from Corollary~\ref{cor:v2-quasiperiodic} is contained in this criterion.

\item \textup{A counterexample.}
The vanishing of the mean is not sufficient. The smooth Bohr almost-periodic function
\begin{align*}
h(t)=\sum_{n=1}^{\infty}2^{-n}\cos(4^{-n}t)
\end{align*}
has zero uniform mean, while its primitive on $\mathbb R$ is unbounded. Thus almost periodicity alone does not imply the bounded-primitive condition. Conversely, failure of this sufficient condition does not exclude critical-front existence; it only shows that Theorem~\ref{thm:v2-main-intro} cannot be invoked directly.
\end{enumerate}
\end{remark}
\endgroup

Theorem~\ref{minimalspeed} identifies the two thresholds under temporal almost periodicity. The scalar problem was studied by Nadin and Rossi \cite{nadin_2015}; related results for epidemic systems were obtained by Wang, Wang and Zhang \cite{WangWangZhang2026} under an additional assumption on slow exponential decay. Here the equality $c^*(e)=c_*(e)$ follows from the uniform means of the speeds $c_\lambda$ in \eqref{speed} and the spectral characterization through $-\underline{k}(\lambda)/\lambda$. In the time-periodic case, \eqref{eigen} also leads to a front with the critical mean speed. The almost-periodic critical-front result above uses the bounded primitive of $h_*$ in place of the eigenfunction condition in Remark~\ref{APremarkcrit}. The proof relies on a normalized positive solution of the parabolic evolution equation. The quasiperiodic case follows under a Fourier summability condition. Section~\ref{sec:critical-proofs} also treats coefficients obtained by a change of variables from a stationary system, for which the required bounded constant-growth eigenfunction need not exist.

\subsection{Asymptotic speeds for small and large diffusion}
We consider a stationary periodic shear for which the two species have the same constant diffusion matrix. In longitudinal and transverse coordinates $(s,z)$, with $z\in\mathbb T=\mathbb R/\mathbb Z$, the diffusion operator is $d(a\partial_{ss}+D\partial_{zz})$ and the longitudinal transport field is $-h\sin(2\pi z)\partial_s$. The positive constants $a,D,r,\gamma,h$ are fixed while $d>0$ varies. The cooperative reaction terms are $(r-\gamma)u_i+\gamma u_{3-i}-ru_i^2$ for $i=1,2$. The system is written in \eqref{eq:memoirs-diffusion-model}. This is an independent constant-coefficient-in-reaction setting: the assertion below permits $r-\gamma\leq0$, whereas the assumptions on the diagonal entries of the general linearization in \eqref{Fnonlinear} impose their positivity. The spectral reduction and linear determinacy used here are established directly for the symmetric system, without invoking the general front-existence theorem outside its hypotheses. We denote by $c^*(d)$ the minimal spectral speed in the positive $s$-direction, characterized by \eqref{eq:memoirs-diffusion-speed}. The following result is proved in Section~\ref{sec:memoirs-diffusion-regimes}.

\begin{theorem}[Two-term corrections in the two diffusion limits]\label{thm:memoirs-diffusion-two-regimes}
Assume \eqref{eq:memoirs-diffusion-model} with fixed $a,D,r,\gamma,h>0$ and put $k=2\pi$. Let $c^*(d)$ be given by \eqref{eq:memoirs-diffusion-speed}. Then, as $d\to\infty$,
\begin{equation}\label{eq:memoirs-diffusion-large}
c^*(d)=2\sqrt{rad}+C_{\infty,1}d^{-3/2}-C_{\infty,2}d^{-7/2}+O(d^{-9/2}).
\end{equation}
Here, $C_{\infty,1}=\sqrt r\,h^2/(2Dk^2\sqrt a)$ and
$C_{\infty,2}=\sqrt r\,h^4/(16D^2k^4a^{3/2})$.
As $d\downarrow0$ with $h>0$ fixed,
\begin{equation}\label{eq:memoirs-diffusion-small}
c^*(d)=h+2\sqrt{rad}-C_{0,1}d^{3/4}-C_{0,2}d+O(d^{5/4}).
\end{equation}
Here, $C_{0,1}=k\sqrt{Dh}(a/r)^{1/4}/\sqrt2$ and
$C_{0,2}=k^2Dh/(32r)$.
In particular,
\begin{equation}\label{eq:memoirs-diffusion-limits}
\lim_{d\to\infty}\frac{c^*(d)}{\sqrt d}=2\sqrt{ra},
\qquad
\lim_{d\downarrow0}c^*(d)=h.
\end{equation}
The large-diffusion expansion is uniform when $a,D,r,h$ range over a compact subset of $(0,\infty)^4$; the small-diffusion expansion is a fixed-positive-shear limit and is not claimed to be uniform as $h\downarrow0$.
\end{theorem}
The two limits in Theorem~\ref{thm:memoirs-diffusion-two-regimes} have different origins. For $d\to\infty$, the principal transverse eigenfunction approaches a constant, and the first correction to the speed is determined by the zero-mean shear through its quadratic spectral response. The term of order $d^{-7/2}$ follows from minimization of the resulting spectral quotient. For $d\downarrow0$, the eigenfunction concentrates around the point where the projected velocity is maximal. Its leading localization cost is the ground-state energy of a harmonic oscillator. The noninteger power $d^{3/4}$ in \eqref{eq:memoirs-diffusion-small} arises from this localization, whereas the term of order $d$ also accounts for the displacement of the minimizing exponential weight. The two limits cannot be interchanged with the zero-shear limit, for which $c^*(d)=2\sqrt{rad}$ holds exactly.

The variational characterization of KPP speeds for periodic media and shears was established in \cite{BerestyckiHamelNadirashviliJEMS2005,NolenXin2005}. The dependence of the speed on diffusion parameters and homogenization regimes was investigated by El Smaily \cite{ElSmaily2008}. Zlato\v{s} \cite{Zlatos2010} obtained sharp large-flow limits and related them to effective diffusivity. Theorem~\ref{thm:memoirs-diffusion-two-regimes} concerns a different parameter regime: the shear amplitude is fixed, both diffusion rates change together, and the correction coefficients are calculated on both sides of the diffusion scale. In particular, the small-diffusion calculation uses a localized eigenfunction and the perturbation of its optimal exponential weight, while the large-diffusion calculation combines analytic eigenvalue perturbation with an explicit minimization. The reduction of the symmetric two-component linearization to its positive scalar mode is exact. Consequently, the diffusion expansions are established through a scalar periodic eigenvalue problem. The formulas do not assert a new effect arising from genuinely unequal diffusion matrices or nonsymmetric coupling, and the present proof does not extend their higher-order remainders to either setting. The estimates of Nolen and Xin concern speed selection for periodic shears, whereas El Smaily investigates diffusion and homogenization regimes; neither reference alone establishes the two correction formulas written here. This comparison is not a claim of priority over all published higher-order expansions. These conclusions are proved only for the stated stationary shears with diagonal diffusion in $(s,z)$; they are not asserted for general heterogeneous cooperative systems.

For the two examples of Section~\ref{section 6}, the common physical diffusion matrix is $dI_2$. The horizontal and diagonal shears correspond respectively to $(a,D,r,\gamma,h)=(1,1,3,1,B)$ and $(1,2,5,2,\sqrt2 B)$; see Corollary~\ref{cor:memoirs-diffusion-cases}. Thus the leading small-diffusion speeds are $B$ and $\sqrt2 B$, while the leading large-diffusion speeds are $2\sqrt{3d}$ and $2\sqrt{5d}$. The dependence on $\gamma$ is absent from the spectral speed because the Perron mode of the symmetric coupling matrix is $(1,1)$.

\subsection{Large-drift asymptotics}

We study the behavior of $c^*(e,A)$ as $A\to\infty$, with ${\bf q}_i=A{\bf v}_i$, $i=1,2$, where the fields ${\bf v}_i$ satisfy \eqref{assumeq}. To connect with classical results in the divergence-form setting, we further assume that
\begin{align}\label{regularityofA}
{\bf A}^i\in [C^{0,1}(\R\times \R^N)\cap W^{0,1,\infty}(\R\times \R^N)]^{N\times N} \text{ and }(\nabla_x\cdot {\bf A}^i)_k =0,~k=1,...,N,~i=1,2,
\end{align}
and
\begin{align}\label{driftterm}
\text{${\bf v}_1={\bf v}_2 ={\bf v}(t,x)$, $\nabla_x \cdot {\bf v}=0$ and }\int_{[{\bf 0},{\bf L}]} {\bf v} dx = {\bf 0}.
\end{align}

We now introduce several important definitions for the subsequent results.
\begin{align*}
&\mathfrak{I}:=\left\{ \boldsymbol{\phi}\in [H_{per}^1(\R^N)]^2:{\bf v}\cdot D \phi_i=0 \text{ on }\R\times \R^N,~\forall i=1,2 \right\}, \\
&\mathfrak{I}_{{\bf A}}^{\mu}:=\left\{\boldsymbol{\phi}=(\phi_1,\phi_2)\in \mathfrak{I}:\sum_{i=1,2}\left<\int_{[{\bf 0},{\bf L}]} D\phi_i\cdot   {\bf A}^iD\phi_i dx\right>\leq \sum_{i,j=1,2}\left<\int_{[{\bf 0},{\bf L}]} {\mu}_{ij} \phi_i\phi_j dx\right>  \right\},\\
& \mathfrak{I}_{\text{time},p}
:=\Big\{\boldsymbol{\phi}\in [L^2_{per}(\R,H^1_{per}(\R^N))]^2 :
\mathbf{v}\cdot D\phi_i=0,\ \forall i=1,2
\Big\},
\end{align*}
and, \textup{only in the time-independent setting} \(\mathbf A^i=\mathbf A^i(x),\ \mathbf v=\mathbf v(x),\ \mu_{ij}=\mu_{ij}(x)\), define
\begin{align*}
\kappa_e(\lambda):=\sup_{{\bf w}=(w_1,w_2)\in \mathfrak{I}\setminus\{\bf 0\}}
\left\{
\dfrac{
\lambda \displaystyle\int_{[{\bf 0},{\bf L}]}
({\bf v}\cdot e)\,|{\bf w}|^2 \, dx
+\sum_{\substack{i,j=1,2}}
\int_{[{\bf 0},{\bf L}]}
\mu_{ij}\, w_i\, w_j  \, dx -  \left\| D {\bf w} \right\|_{2,{\bf A}^1,{\bf A}^2}^2
}{
\|{\bf w}\|_2^2
}
\right\},
\end{align*}
where $\left\| D {\bf w} \right\|_{2,{\bf A}^1,{\bf A}^2}^2=\left\| D w_1 \right\|_{2,{\bf A}^1}^2+\left\| D w_2 \right\|_{2,{\bf A}^2}^2=\displaystyle\int_{[{\bf 0},{\bf L}]}Dw_1 \cdot {\bf A}^1Dw_1 dx + \displaystyle\int_{[{\bf 0},{\bf L}]}Dw_2\cdot {\bf A}^2Dw_2 dx.$

We have the following large-advection estimates.

\begingroup
\begin{theorem}\label{limitspeed1}
Assume \eqref{assumeq}--\eqref{near_zerp} and \eqref{APassume}--\eqref{Fregularity}. The following assertions hold.

\smallskip\noindent
\textup{(a)} General space--time heterogeneous drift.
Let ${\bf q}_i=A{\bf v}_i+{\bf r}_i$, $i=1,2$, where ${\bf v}_i$ and ${\bf r}_i$ are bounded, spatially periodic and almost periodic in time, uniformly in space, and ${\bf r}_i$ are independent of $A$.
For each $e\in\mathbb S^{N-1}$,
\begin{equation}\label{eq:v5-general}
\max_{i=1,2}\left\langle\inf_{x\in\mathbb R^N}({\bf v}_i\cdot e)\right\rangle
\le \liminf_{A\to\infty}\frac{c^*(e,A)}A
\le \limsup_{A\to\infty}\frac{c^*(e,A)}A
\le\left\langle\max_{i=1,2}\sup_{x\in\mathbb R^N}({\bf v}_i\cdot e)\right\rangle.
\end{equation}

\smallskip\noindent
\textup{(b)} First-integral estimates and stationary limit.
In the setting of \textup{(a)}, assume ${\bf r}_1={\bf r}_2=0$ and \eqref{regularityofA}--\eqref{driftterm}, with the notation $\mathfrak I$, $\mathfrak I_{\bf A}^\mu$, $\mathfrak I_{\mathrm{time},p}$ and $\kappa_e$ introduced above.
If $\mu_{12}=\mu_{21}$, then
\begin{equation}\label{eq:v5-lower}
\liminf_{A\to\infty}\frac{c^*(e,A)}A\ge
\sup_{\boldsymbol\phi\in\mathfrak I_{\bf A}^\mu\setminus\{\bf0\}}
\frac{\left\langle\int_{[\bf0,\bf L]}({\bf v}\cdot e)|\boldsymbol\phi|^2dx\right\rangle}{\int_{[\bf0,\bf L]}|\boldsymbol\phi|^2dx}\ge0.
\end{equation}
If the coefficients are $T_0$-periodic in time, then
\begin{equation}\label{eq:v5-upper}
\limsup_{A\to\infty}\frac{c^*(e,A)}A\le
\sup_{\boldsymbol\phi\in\mathfrak I_{\mathrm{time},p}\setminus\{\bf0\}}
\frac{\left\langle\int_{[\bf0,\bf L]}({\bf v}\cdot e)|\boldsymbol\phi|^2dx\right\rangle}
{\left\langle\int_{[\bf0,\bf L]}|\boldsymbol\phi|^2dx\right\rangle}.
\end{equation}
If all the coefficients are independent of time and $\mu_{12}=\mu_{21}$, then
\begin{equation}\label{eq:v5-stationary}
\lim_{A\to\infty}\frac{c^*(e,A)}A
=\inf_{\lambda>0}\frac{\kappa_e(\lambda)}\lambda
=\sup_{\boldsymbol\phi\in\mathfrak I_{\bf A}^\mu\setminus\{\bf0\}}
\frac{\int_{[\bf0,\bf L]}({\bf v}\cdot e)|\boldsymbol\phi|^2dx}{\int_{[\bf0,\bf L]}|\boldsymbol\phi|^2dx}.
\end{equation}
If, moreover, $\mathfrak I\ne\mathbb R^2$, then
\begin{equation}\label{eq:v5-small}
\begin{aligned}
&\lim_{\varepsilon\downarrow0}\lim_{A\to\infty}
\frac{c^*(e,A,\varepsilon{\bf F})}{2A\sqrt\varepsilon}=\frac1V\sup_{\substack{\xi\in\mathbb S^1\\\xi_1,\xi_2\ge0}}
\left\{\sqrt{\xi\cdot M_0\xi}\,
\sup_{\substack{\boldsymbol\phi\in\mathfrak I\\
\int_{[{\bf0},{\bf L}]}\boldsymbol\phi\,dx={\bf0}\\
\|D\boldsymbol\phi\|_{2,{\bf A}^1,{\bf A}^2}>0}}
\frac{\int_{[{\bf0},{\bf L}]}({\bf v}\cdot e)
(\xi\cdot\boldsymbol\phi)\,dx}
{\|D\boldsymbol\phi\|_{2,{\bf A}^1,{\bf A}^2}}\right\}.
\end{aligned}
\end{equation}
Here $V=|[{\bf0},{\bf L}]|$ and
$M_0=(\int_{[{\bf0},{\bf L}]}\mu_{ij}\,dx)_{i,j=1,2}$.
The inner supremum is taken with $\xi$ fixed. The mean-zero
normalization does not change the quotient. By the symmetry of
the first-integral space and the positivity of the off-diagonal
entries of $M_0$, it is sufficient to consider $\xi_i\ge0$.

\end{theorem}
\endgroup

\begin{remark}[Two explicit first-integral classes]\label{rem:first-integral-explicit}
The abstract spaces in Theorem~\ref{limitspeed1}\textup{(b)} become completely explicit for the two shear flows used in Section~\ref{section 6}.
\begin{enumerate}[label=(\roman*)]
\item If ${\bf v}(x,y)=(\sin(2\pi y),0)$ on $[0,1]^2$, then every first integral of ${\bf v}$ is a function of $y$ only. Hence,
\begin{align*}
\mathfrak I=\big\{\boldsymbol\phi=(\phi_1(y),\phi_2(y)):\ \phi_i\in H^1_{per}(0,1)\big\},
\end{align*}
and the variational formula in \eqref{eq:v5-stationary} reduces to the one-dimensional periodic spectral problem in \eqref{eq:app-case1-1d}.

\item If ${\bf v}(x,y)=\big(\sin(2\pi(x-y)),\sin(2\pi(x-y))\big)$, then every first integral is a function of $z=x-y$ only. Therefore,
\begin{align*}
\mathfrak I=\big\{\boldsymbol\phi=(\phi_1(z),\phi_2(z)):\ \phi_i\in H^1_{per}(0,1)\big\},
\end{align*}
and \eqref{eq:v5-stationary} reduces to the one-dimensional problem in \eqref{eq:app-case2-1d}.
\end{enumerate}
In both cases the abstract first-integral characterization is thus converted into an explicit one-dimensional periodic eigenvalue problem, which explains the spectral reductions carried out in Section~\ref{section 6}.
\end{remark}

The bounds in Theorem~\ref{limitspeed1} extend the scalar large-drift estimates of Berestycki, Hamel and Nadirashvili \cite[Theorem 4.1]{BerestyckiHamelNadirashvili2005} to cooperative systems. In the stationary, symmetric case the limit admits an exact first-integral formula, as in the scalar setting studied by Zlato\v{s} \cite{Zlatos2010}. The linearization now involves the matrix $\mu_{ij}=\partial_j F_i(t,x,0,0)$, and the estimate must account for the coupling between components. Section~\ref{counter} shows why the corresponding energy-constrained equality fails in general for time-periodic coefficients. The case $\mathfrak I=\R^2$ is excluded from the final assertion because its limit is trivial.

\begin{remark}\label{remarktimecounter} Even if $\mu_{12}=\mu_{21}$, in general, the inequality in
the upper estimate \eqref{eq:v5-upper} cannot be improved to
an equality by introducing a natural energy constraint. Indeed, define
\begin{align*}
\mathfrak{I}_{{\bf A},time}^\mu
:=
\left\{
\boldsymbol{\phi}\in\mathfrak{I}_{\text{time},p}:
\sum_{i=1,2}
\left<
\int_{[{\bf 0},{\bf L}]}
D\phi_i\cdot{\bf A}^iD\phi_i\,dx
\right>
\leq
\sum_{i,j=1,2}
\left<
\int_{[{\bf 0},{\bf L}]}
\mu_{ij}\phi_i\phi_j\,dx
\right>
\right\}.
\end{align*}
Then, in general,
\begin{align*}
\lim_{A\rightarrow\infty}
\dfrac{c^*(e,A)}{A}
\neq
\sup_{\boldsymbol{\phi}\in
\mathfrak{I}_{{\bf A},time}^\mu\setminus\{{\bf 0}\}}
\dfrac{
\left<
\displaystyle\int_{[{\bf 0},{\bf L}]}
({\bf v}\cdot e)|\boldsymbol{\phi}|^2\,dx
\right>
}{
\left<
\displaystyle\int_{[{\bf 0},{\bf L}]}
|\boldsymbol{\phi}|^2\,dx
\right>
},~
\lim_{A\rightarrow\infty}
\dfrac{c^*(e,A)}{A}
\neq
\sup_{\boldsymbol{\phi}\in
\mathfrak{I}_{{\bf A}}^\mu\setminus\{{\bf 0}\}}
\dfrac{
\left<
\displaystyle\int_{[{\bf 0},{\bf L}]}
({\bf v}\cdot e)|\boldsymbol{\phi}|^2\,dx
\right>
}{
\displaystyle\int_{[{\bf 0},{\bf L}]}
|\boldsymbol{\phi}|^2\,dx
}.
\end{align*}
Thus, neither the time-periodic energy-constrained class nor its time-independent counterpart yields the exact large-drift limit. A counterexample is given in Section \ref{counter}.
\end{remark}

\subsection{Organization of the paper}
We first derive the existence and nonexistence thresholds from the tilted linear evolution. The two statements are proved separately, since the existence threshold need not agree with the spectral lower bound in a general heterogeneous medium. We then identify the minimal speed in the almost-periodic framework.

The proof of Theorem~\ref{thm:memoirs-diffusion-two-regimes} follows the minimal-speed analysis. It uses a scalar periodic spectral reduction, perturbation about the constant eigenfunction for large diffusion, and localization near a maximal streamline for small diffusion. We next establish the critical-front results under their additional hypotheses. Large-advection behavior and the time-dependent counterexample are considered thereafter. The explicit periodic models and numerical illustrations are collected at the end; the computations concern finite-time solutions and are independent of the analytic proofs.

The sections are arranged in this order. Section~\ref{sec2} contains the preliminaries. Sections~\ref{section:ex2} and~\ref{sec4} prove the existence and nonexistence results. Section~\ref{section 5} establishes the almost-periodic minimal-speed characterization, and Section~\ref{sec:memoirs-diffusion-regimes} proves the two diffusion-rate expansions. Section~\ref{sec:critical-proofs} concerns critical-front attainment, Section~\ref{sec:large-advection-proofs} deals with large-advection limits and Section~\ref{section 6} collects explicit models and their numerical illustrations.

\section{Preliminary Results}\label{sec2}

For ${\bf u}=(u_1,u_2)\in BC(\R^{N+1})$, define $\|{\bf u}(t,\cdot)\|_{\infty}=\max\limits_{i=1,2}\max\limits_{x\in \R^N}|u_i(t,x)|.$ Furthermore, for $T_0>0$, we consider the following norm for $[L^2\bigl(0,T_0; H^1_{\mathrm{per}}([{\bf 0},{\bf L}])]^2$ 
\begin{align*}
\begin{aligned}
\|{\bf u}\|^2 &=\|{\bf u}\|_2^2+\|D{\bf u}\|_{2,{\bf A}^1,{\bf A}^2}^2,~
{\bf u}=(u_1,u_2)\in[L^2(0,T_0;H^1_{\mathrm{per}}([{\bf0},{\bf L}]))]^2.
\end{aligned}
\end{align*}
where
\begin{equation*}
\|D{\bf u}\|_{2,{\bf A}^1,{\bf A}^2}^2
=\sum_{i=1}^2\int_{[0,T_0]\times[{\bf0},{\bf L}]}
Du_i\cdot {\bf A}^i Du_i\,dx\,dt.
\end{equation*} Similarly, for $H^1_{\mathrm{per}}([{\bf 0},{\bf L}])$, we consider 
\begin{align*}
\left\| D {\bf u} \right\|_{2,{\bf A}^1,{\bf A}^2}^2=\displaystyle\int_{[{\bf 0},{\bf L}]}Du_1 \cdot {\bf A}^1Du_1\, dx + \displaystyle\int_{[{\bf 0},{\bf L}]}Du_2\cdot {\bf A}^2Du_2\, dx.
\end{align*}

Next, we investigate the following ODE system
\begin{align}\label{super1}
\begin{cases}
\partial_t v_1
= v_1\left[C\lambda^2 +\lambda\max\limits_{ \R^{N+1}}|{\bf q}_1|+\max\limits_{x\in\R^{N}}\mu_{11}(t,x)\right]+v_2\max\limits_{x\in\R^{N}}\mu_{12}(t,x),\\
\partial_t v_2
= v_2\left[C\lambda^2 +\lambda\max\limits_{ \R^{N+1}}|{\bf q}_2|+\max\limits_{x\in\R^{N}}\mu_{22}(t,x)\right]+v_1\max\limits_{x\in\R^{N}}\mu_{21}(t,x).
\end{cases}
\end{align}
It is well-known that the system admits a unique evolution operator $\Phi(t,s)$ such that
\begin{align*}
{\bf v}(t)=(v_1(t),v_2(t))=\Phi(t,s){\bf v}_0,\text{ for }t,s\in \R,~t>s,
\end{align*}
for some initial condition ${\bf v}_0=(v_{0,1},v_{0,2})$. On the other hand, since $\max\limits_{x\in \R^N}\mu_{12}(t,x),\max\limits_{x\in \R^N}\mu_{21}(t,x)>0$, it follows from Metzler theory that ${\bf v}(t)=\Phi(t,s){\bf v}_0>{\bf 0},~\forall t>s,$ for some ${\bf v}_0\geq {\bf 0}$, ${\bf v}_0\neq {\bf 0}$. Moreover, define 
\begin{align*}
\mu_{\infty}(t)=\max\limits_{i=1,2}\left\{\max\limits_{x\in\R^{N}}\mu_{11}(t,x)+\max\limits_{x\in\R^{N}}\mu_{12}(t,x);\max\limits_{x\in\R^{N}}\mu_{22}(t,x)+\max\limits_{x\in\R^{N}}\mu_{21}(t,x)\right\}.
\end{align*}
Grönwall's inequality implies that
\begin{align}\label{upperexp}
\|\Phi(t,s)\|\leq \exp\left((C\lambda^2 +\lambda\max_{i=1,2}\max\limits_{ \R^{N+1}}|{\bf q}_i|)(t-s)+\int_s^t\mu_{\infty}(k)dk\right).
\end{align}

Similarly, 
\begin{align}\label{sub121}
\begin{cases}
\partial_t v_1
= v_1\left[C^{-1}\lambda^2 -\lambda\max\limits_{ \R^{N+1}}|{\bf q}_1|+\min\limits_{x\in\R^{N}}\mu_{11}(t,x)\right]+v_2\min\limits_{x\in\R^{N}}\mu_{12}(t,x),\\
\partial_t v_2
= v_2\left[C^{-1}\lambda^2 -\lambda\max\limits_{ \R^{N+1}}|{\bf q}_2|+\min\limits_{x\in\R^{N}}\mu_{22}(t,x)\right]+v_1\min\limits_{x\in\R^{N}}\mu_{21}(t,x),
\end{cases}
\end{align}
admits a unique evolution operator $\Psi(t,s)$ such that ${\bf v}(t)=\Psi(t,s){\bf v}_0>0,~\forall t>s,$ for some ${\bf v}_0\geq {\bf 0}$, ${\bf v}_0\neq {\bf 0}$. On the other hand, since $\mu_{ij}\geq 0$ for $i\neq j$, one has
\begin{align*}
\begin{cases}
\partial_t v_1
\geq v_1\left[C^{-1}\lambda^2 -\lambda\max\limits_{ \R^{N+1}}|{\bf q}_1|+\min\limits_{x\in\R^{N}}\mu_{11}(t,x)\right],\\
\partial_t v_2
\geq v_2\left[C^{-1}\lambda^2 -\lambda\max\limits_{ \R^{N+1}}|{\bf q}_2|+\min\limits_{x\in\R^{N}}\mu_{22}(t,x)\right].
\end{cases}
\end{align*}
This implies that 
\begin{align*}
\begin{cases}
v_1(t)\geq v_{0,1}\exp{\left((C^{-1}\lambda^2 -\lambda\max\limits_{ \R^{N+1}}|{\bf q}_1|)(t-s)+\displaystyle\int_{s}^t\min\limits_{x\in\R^{N}}\mu_{11}(k,x)dk\right)},\\
v_2(t)\geq v_{0,2}\exp{\left((C^{-1}\lambda^2 -\lambda\max\limits_{ \R^{N+1}}|{\bf q}_2|)(t-s)+\displaystyle\int_{s}^t\min\limits_{x\in\R^{N}}\mu_{22}(k,x)dk\right)},
\end{cases}
\end{align*}
for $t>s$. It then follows that, for $i=1,2$,
\begin{align}\label{estimate123}
v_i(t)\geq v_{0,i}\exp{\left((C^{-1}\lambda^2 -\lambda\max\limits_{i=1,2}\max\limits_{ \R^{N+1}}|{\bf q}_i|)(t-s)+\displaystyle\int_{s}^t\min_{i=1,2}\min\limits_{x\in\R^{N}}\mu_{ii}(k,x)dk\right)},~\forall t>s.
\end{align}

\section{Existence of Generalized Transition Waves}\label{section:ex2}

\subsection{Estimates for the linearized equation}

First, we consider the linearized system around ${\bf 0}=(0,0)$ as follows
\begin{equation}\label{eq:main1}
\begin{cases}
\partial_t u_1  - \text{Tr}({\bf A}^1(t,x)D^2 u_1) + {\bf q}_1(t,x)\cdot D u_1 = \mu_{11}(t,x)u_1+\mu_{12}(t,x)u_2, & t\in \R,\ x\in \R^N, \\[1ex]
\partial_t u_2  - \text{Tr} ({\bf A}^2(t,x)D^2 u_2) + {\bf q}_2(t,x)\cdot D u_2 = \mu_{21}(t,x)u_1+\mu_{22}(t,x)u_2, & t\in \R,\ x\in \R^N,
\end{cases}
\end{equation}
where $\mu_{ij}=\partial_{j}F_i(t,x,0,0)$ for $i,j=1,2$.
The coefficients $\mu_{ij}$ are bounded but need not be H\"older continuous. In this section, the linearized equations
are treated in the local parabolic strong $W^{1,2}_p$ class, with $p>N+2$.
The diffusion matrices are uniformly elliptic and globally H\"older continuous;
the zero-order terms are bounded. Local parabolic $L^p$ estimates therefore
apply to bounded solutions. The comparison principle and the cooperative
Harnack estimates below are used for such strong solutions. No spatial or temporal
Schauder regularity of $\mu_{ij}$ is claimed at this stage.
It is standard to expect the tail of ${\bf u} = (u_1,u_2)$ to decay exponentially. Thus, we are looking for the following type of solution
\begin{align*}
{\bf u}(t,x)=e^{-\lambda x\cdot e}\boldsymbol{\eta}_{\lambda}(t,x),
\end{align*}
where $\boldsymbol{\eta}_{\lambda}(t,x)=(\eta_{\lambda,1}(t,x),\eta_{\lambda,2}(t,x))$ is {\bf L}-periodic in $x$ and $\lambda$ is a threshold constant used to determine the existence of the wavefront solution. Direct calculations yield that
\begin{align}\label{div}
\begin{cases}
\partial_t \eta_{\lambda,1}
&= \text{Tr} \big( {\bf A}^1(t,x) D^2 \eta_{\lambda,1} \big)
- \big( {\bf q}_1(t,x) + 2\lambda {\bf A}^1(t,x)e \big) \cdot D \eta_{\lambda,1}+ K_1^\lambda\eta_{\lambda,1}+\mu_{12}(t,x)\eta_{\lambda,2},\\
\partial_t \eta_{\lambda,2}
&= \text{Tr} \big( {\bf A}^2(t,x) D^2 \eta_{\lambda,2} \big)
- \big( {\bf q}_2(t,x) + 2\lambda {\bf A}^2(t,x)e \big) \cdot D \eta_{\lambda,2}+K_2^\lambda\eta_{\lambda,2}+\mu_{21}(t,x)\eta_{\lambda,1},\\
\end{cases}
\end{align}
for $t\in\mathbb{R},~x\in\mathbb{R}^N$, where
\begin{align*}
K_{i}^\lambda:=\mu_{ii}(t,x) + \lambda^2\, e \cdot {\bf A}^i(t,x)e + \lambda\, {\bf q}_i(t,x)\cdot e,~i=1,2.
\end{align*}
Thanks to the estimates obtained in \eqref{super1} and \eqref{sub121}, Harnack-type results for $\boldsymbol{\eta}_{\lambda}$ are studied as follows

\begin{lemma}\label{lm:Haranack-type}
For any $\lambda >0$, \eqref{div} admits a positive {\bf L}-periodic entire solution $\boldsymbol{\eta}_{\lambda}$ on $(t,x)\in \R\times \R^N$. Moreover, $\boldsymbol{\eta}_{\lambda}$ is unique up to a multiplicative constant and enjoys the following estimates, for any $(t,x)\in \R\times \R^N$, $T\geq 0$,
\begin{align}
&\max_{\substack{i=1,2\\x\in \R^N }}\eta_{\lambda,i}(t+T,x)\leq  \max_{\substack{i=1,2\\x\in \R^N }}\eta_{\lambda,i}(t,x)\exp\left((C\lambda +\max_{i=1,2}\|{\bf q}_i\|_{\infty})\lambda T+\int_t^{t+T}\mu_{\infty}(k)dk\right),\label{upper}\\
&\min\limits_{\substack{i=1,2\\x\in \R^N }}\eta_{\lambda,i}(t+T,x)\geq C_0\max_{\substack{i=1,2\\x\in \R^N }}\eta_{\lambda,i}(t,x) \exp{\left((C^{-1}\lambda -\max_{i=1,2}\|{\bf q}_i\|_{\infty})\lambda T+\displaystyle\int_{t}^{t+T}\min\limits_{\substack{i=1,2\\x\in \R^N }}\mu_{ii}(k,x)dk\right)},\label{below}
\end{align}
where $C_0>0$ depends on $\lambda$, the periodic cell, the ellipticity and coefficient regularity bounds, the uniform bounds of $\mathbf q_i,\mu_{ij}$, and the positive lower bounds of the off-diagonal couplings $\mu_{12},\mu_{21}$. Here, for simplicity, we use $\max\limits_{\substack{i=1,2\\x\in \R^N }}:=\max\limits_{i=1,2}\max\limits_{x\in \R^N}$ and $\min\limits_{\substack{i=1,2\\x\in \R^N }}:=\min\limits_{i=1,2}\min\limits_{x\in \R^N}$.

\end{lemma}

\begin{proof}

For each integer $n>1$, solve \eqref{div} on $[-n,\infty)$ with
strictly positive constant initial data at $-n$, and multiply the
solution by a positive constant so that
$\|\boldsymbol\eta^n_\lambda(0,\cdot)\|_\infty=1$.
The cooperative Harnack inequality on the periodic cell, applied on
fixed positive-length cylinders around time $0$, has constants
independent of $n$. Hence, there is $\alpha_\lambda>0$ such that
$\min\limits_{i,x}\eta^n_{\lambda,i}(0,x)\ge \alpha_\lambda$ for all
sufficiently large $n$. Comparison with spatially constant linear
systems controls the solutions on each fixed time interval $[-T,T]$,
while the local parabolic $W^{1,2}_p$ estimates ($p>N+2$) and
periodicity in $x$ yield compactness in the uniform topology on
$[-T,T]\times[\mathbf0,\mathbf L]$. A diagonal subsequence
converges to an entire solution $\boldsymbol\eta_\lambda$ of
\eqref{div}. The normalization and Harnack lower bound pass to
the limit at time $0$; strong cooperativity and comparison then
show that both components are strictly positive for every time.
No global-in-time upper bound for this entire solution is asserted.

Next, we prove \eqref{upper}. Choose ${\bf v}_0=\left(\max\limits_{x\in \R^N}\eta_{\lambda,1}(t,x),\max\limits_{x\in \R^N}\eta_{\lambda,2}(t,x)\right)$ for some $t\in \R$ in the system \eqref{super1}. By the parabolic comparison principle, for any $T>0$, one has
\begin{align*}
(\eta_{\lambda,1}(t+T,x),\eta_{\lambda,2}(t+T,x))\leq  \Phi(t+T,t)\left(\max\limits_{x\in \R^N}\eta_{\lambda,1}(t,x),\max\limits_{x\in \R^N}\eta_{\lambda,2}(t,x)\right).
\end{align*}
It follows from \eqref{upperexp} that 
\begin{align*}
\max_{i=1,2}\max_{x\in \R^N}\eta_{\lambda,i}(t+T,x)\leq  \max_{i=1,2}\max_{x\in \R^N}\eta_{\lambda,i}(t,x)\exp\left((C\lambda^2 +\lambda\max_{i=1,2}\max\limits_{ \R^{N+1}}|{\bf q}_i|)T+\int_t^{t+T}\mu_{\infty}(k)dk\right).
\end{align*}

Now, we show \eqref{below}. It follows from $\Psi$ of  \eqref{sub121} that
\begin{align*}
(\eta_{\lambda,1}(t+T,x),\eta_{\lambda,2}(t+T,x))&\geq  \Psi(t+T,t)\left(\min\limits_{x\in \R^N}\eta_{\lambda,1}(t,x),\min\limits_{x\in \R^N}\eta_{\lambda,2}(t,x)\right)\\
&\geq \min_{i=1,2}\min\limits_{x\in \R^N}\eta_{\lambda,i}(t,x) \Psi(t+T,t){\bf 1}.
\end{align*}

We complete the lower estimate by making the time shift in the
cooperative Harnack inequality explicit. Let
$\mathcal C=[\mathbf0,\mathbf L]$ and put
\(
M(t)=\max\limits_{i,x\in\mathcal C}\eta_{\lambda,i}(t,x),~ m(t)=\min\limits_{i,x\in\mathcal C}\eta_{\lambda,i}(t,x),~ Q=\max\limits_{i=1,2}\|\mathbf q_i\|_\infty.
\)
The uniform cooperative Harnack inequality for irreducible parabolic
systems (see \cite[Theorem 3.9]{FoldesPolacik2009} and
\cite[Proposition 2.4]{Girardin_2024}) implies that for some
$\tau>0$ and $C_1\ge1$ independent of $t$,
\[
M(t-\tau)\le C_1m(t),~ t\in\mathbb R.
\]
The comparison estimate \eqref{upper} over the interval
$[t-\tau,t]$ also implies
\[
M(t)\le \exp\!\left((C\lambda^2+\lambda Q)\tau
+\int_{t-\tau}^{t}\mu_\infty(s)\,ds\right)M(t-\tau).
\]
We therefore have
$m(t)\ge C_0M(t)$ for every $t$, where
\(
C_0=C_1^{-1}\exp\!\left(-[C\lambda^2+\lambda Q+\|\mu_\infty\|_\infty]\tau\right)>0.
\)
Finally, \eqref{sub121} and \eqref{estimate123}, applied from
$t$ to $t+T$, show that
\begin{align*}
m(t+T)&\ge m(t)\exp\!\left(
(C^{-1}\lambda^2-\lambda Q)T
+\int_t^{t+T}\min_{i=1,2}\inf_{x\in\mathcal C}
\mu_{ii}(s,x)\,ds\right)\\
&\ge C_0M(t)\exp\!\left(
(C^{-1}\lambda^2-\lambda Q)T
+\int_t^{t+T}\min_{i=1,2}\inf_{x\in\mathcal C}
\mu_{ii}(s,x)\,ds\right).
\end{align*}
This proves \eqref{below}, including the case $T=0$. It remains to show the uniqueness. Fix $\lambda>0$ and suppose that $\boldsymbol{\eta}^1=(\eta^1_1,\eta^1_2)$ and $\boldsymbol{\eta}^2=(\eta^2_1,\eta^2_2)$ are two positive, {\bf L}-periodic solutions of \eqref{div}. Let  $h>0$ be a constant such that $\boldsymbol{\eta}^1(0,\cdot)\leq h\boldsymbol{\eta}^2(0,\cdot)$. For $t\leq 0$, one has
\begin{align*}
\min_{i=1,2}\min_{x\in \R^N}\eta^1_i(t,\cdot)\leq  h\max_{i=1,2}\max_{x\in \R^N}\eta_i^2(t,\cdot).
\end{align*}
Indeed, suppose otherwise, $\eta^1_i(t,\cdot)>h\eta^2_i(t,\cdot)$ for $i =1,2$. Then, the strong maximum principle implies that $\boldsymbol{\eta}^1(0,\cdot)> h\boldsymbol{\eta}^2(0,\cdot)$. Contradiction. Hence,  \eqref{Fnonlinear} and \eqref{below} imply that there exists $K>0$ large enough so that
\begin{align*}
\max_{i=1,2}\max_{x\in \R^N}\eta^1_i(t,\cdot)\leq  K\min_{i=1,2}\min_{x\in \R^N}\eta_i^2(t,\cdot).
\end{align*}
By interchanging the roles of $\eta^1$ and $\eta^2$ and choosing a larger $K>1$ if necessary, we deduce that
\begin{align}\label{ine1}
K^{-1}\eta^2_i(t,\cdot)\leq \eta^1_i(t,\cdot)\leq K\eta^2_i(t,\cdot)\text{ on }\R^N,~\forall t\in \R,~ i=1,2.
\end{align}

Next, define $k:=\limsup\limits_{t\rightarrow -\infty}\max\limits_{i=1,2} \max\limits_{x\in \R^N}\dfrac{\eta^1_i(t,x)}{\eta^2_i(t,x)}.$ One has $k\in [K^{-1},K]$. Consider  $\{t_n\}_{n\in \N}$ such that
\begin{align*}
\lim_{n\rightarrow \infty} t_n=-\infty,~ k=\lim_{n\rightarrow \infty}\max_{i=1,2} \max_{x\in \R^N}\dfrac{\eta^1_i(t_n,x)}{\eta^2_i(t_n,x)}.
\end{align*}
From here, we define
\begin{align*}
\eta^i_{n,j}(t,x):=\dfrac{\eta^i_j(t+t_n,x)}{\max\limits_{k=1,2}\max\limits_{y\in \R^N}\eta^1_k(t_n,y)},~i,j=1,2.
\end{align*}
Without loss of generality, assume $t_n<0$. By \eqref{upper} and \eqref{below}, for $j=1,2$,  $\eta^1_{n,j}$ are locally uniformly bounded in $\R^{N+1}$, say, for example, $[-2,2]\times \R^N$. This also holds for $\eta^2_{n,j}$ with $j=1,2$, by \eqref{ine1}. Thanks to parabolic estimates and periodicity in $x$, $\eta^i_{n,j}$, $D\eta^i_{n,j}$, $D^2\eta^i_{n,j}$ converge, up to subsequences, weakly in $L^p_{\loc}(\R^{N+1})$ for any $p>1$. By applying Sobolev embedding, $\eta^i_{n,j}$ converge locally uniformly to some functions $\hat{\eta}^i_j$ for $i,j=1,2$. On the other hand, we consider
\begin{align*}
{\bf A}^i_n(t,x)={\bf A}^i(t+t_n,x),~{\bf q}_{i}^n(t,x)={\bf q}_i(t+t_n,x),~\mu_{ij}^n(t,x)=\mu_{ij}(t+t_n,x),~i,j=1,2.
\end{align*}
Since ${\bf A}^i$ and ${\bf q}_i$ are uniformly continuous, we deduce that ${\bf A}^i_n$ and ${\bf q}_{i}^n$ converge, up to subsequences, to some $\hat{{\bf A}}^i$ and $\hat{{\bf q}}_i$ in $L_{loc}^{\infty}(\R^{N+1})$, respectively, while $\mu_{ij}^n$ converges to $\hat{\mu}_{ij}$ in the weak-$\star$ topology of $L^\infty(\R^{N+1})$ by Banach–Alaoglu–Bourbaki theorem for $i,j=1,2$. Passing to the limit as $n\rightarrow \infty$ weakly in $L^p_{\loc}(\R^{N+1})$, one has
\begin{align}\label{div12}
\begin{cases}
\partial_t \hat{\eta}^i_{1}
&= \text{Tr} \big( \hat{{\bf A}}^1(t,x) D^2 \hat{\eta}^i_{1} \big)
- \big( \hat{{\bf q}}_1(t,x) + 2\lambda \hat{{\bf A}}^1(t,x)e \big) \cdot D \hat{\eta}^i_{1}
+ \hat{K}_1^\lambda \hat{\eta}^i_{1}
+ \hat{\mu}_{12}(t,x)\hat{\eta}^i_{2},\\
\partial_t \hat{\eta}^i_{2}
&= \text{Tr} \big( \hat{{\bf A}}^2(t,x) D^2 \hat{\eta}^i_{2} \big)
- \big( \hat{{\bf q}}_2(t,x) + 2\lambda \hat{{\bf A}}^2(t,x)e \big) \cdot D \hat{\eta}^i_{2}
+ \hat{K}_2^\lambda \hat{\eta}^i_{2}
+ \hat{\mu}_{21}(t,x)\hat{\eta}^i_{1},
\end{cases}\text{ on } \R\times\R^N,
\end{align}
in the weak sense, where $\hat{K}_{i}^\lambda := \hat{\mu}_{ii}(t,x)
+ \lambda^2\, e \cdot \hat{{\bf A}}^i(t,x)e
+ \lambda\, \hat{{\bf q}}_i(t,x)\cdot e,~i=1,2.$ This implies that $\boldsymbol{\hat{\eta}}^1$ and 
$\boldsymbol{\hat{\eta}}^2$ are strong positive solutions of \eqref{div12}. Moreover, we get
\begin{align*}
\boldsymbol{\hat{\eta}}^1\leq k\boldsymbol{\hat{\eta}}^2,~\max_{i=1,2}\max_{x\in \R^{N}}\dfrac{\hat{\eta}^1_i(0,x)}{\hat{\eta}^2_i(0,x)}=k.
\end{align*}

We justify the last comparison globally. Write
$R(t)=\max\limits_{i,x\in\mathcal C}\eta_i^1(t,x)/\eta_i^2(t,x)$,
where $\mathcal C=[\mathbf0,\mathbf L]$. By \eqref{ine1},
$K^{-1}\le R(t)\le K$. The definition of $k$ implies that,
for every $\epsilon>0$, one has $R(t)\le k+\epsilon$
whenever $t$ is sufficiently negative. For every fixed
$s\in\mathbb R$, $t_n+s\to-\infty$, and therefore
$\widehat{\boldsymbol\eta}^{\,1}(s,x)\le
k\widehat{\boldsymbol\eta}^{\,2}(s,x)$
for every $(s,x)$, by local uniform convergence and the
positivity of the limiting second solution. Moreover, a maximum
point for $R(t_n)$ can be chosen on the compact periodic cell;
after passing to a subsequence, it converges to a point at which
$\widehat\eta_j^1(0,x)=k\widehat\eta_j^2(0,x)$
for some component $j$. The strong maximum principle for fully
coupled systems, applied to
$k\widehat{\boldsymbol\eta}^{\,2}
-\widehat{\boldsymbol\eta}^{\,1}\ge\mathbf0$,
shows that this difference vanishes for all $s\le0$.
Uniqueness for the linear Cauchy problem at $s=0$ extends the
identity to all $s\in\mathbb R$.

In particular,
$\widehat{\boldsymbol\eta}^{\,1}(0,\cdot)
=k\widehat{\boldsymbol\eta}^{\,2}(0,\cdot)$
uniformly on $\mathcal C$. The cooperative Harnack estimate
and \eqref{ine1}, together with the normalization of
$\boldsymbol\eta_n^1(0,\cdot)$, ensure that
$\min_{i,x\in\mathcal C}\eta_{n,i}^2(0,x)$
is bounded below by a positive constant independent of $n$.
Consequently, for every $\epsilon>0$ and all sufficiently
large $n$,
\begin{align*}
(k-\epsilon)\boldsymbol\eta^2(t_n,x)
\le\boldsymbol\eta^1(t_n,x)
\le(k+\epsilon)\boldsymbol\eta^2(t_n,x),~ x\in\mathbb R^N.
\end{align*}
Spatial periodicity is used here to pass from local uniform
convergence on the cell to an inequality for every $x\in\mathbb R^N$.
The comparison principle, starting at the single time $t_n$,
now proves the same inequalities for every $t\ge t_n$ and
$x\in\mathbb R^N$. Given any $t\in\mathbb R$, choose
$n$ large enough that $t_n\le t$, and then let
$\epsilon\downarrow0$. We obtain
$\boldsymbol\eta^1(t,x)=k\boldsymbol\eta^2(t,x)$
for every $(t,x)\in\mathbb R^{N+1}$.
This proves uniqueness up to multiplication by a positive constant.
\end{proof}

The conditions \eqref{upper} and \eqref{below}, together with \eqref{upperexp} and \eqref{estimate123}, yield $C_2>1$ large enough such that
\begin{align}\label{boundedtwosides}
C_2^{-1}e^{-C_2(\lambda+1) T}e^{C_2^{-1}\lambda^2 T}\|\boldsymbol{\eta}_{\lambda}(t,\cdot)\|_{\infty}\leq \eta_{\lambda,i}(t+T,x)\leq \|\boldsymbol{\eta}_{\lambda}(t,\cdot)\|_{\infty} C_2 e^{C_2(\lambda^2+\lambda+1) T},
\end{align}
for $(t,x)\in \R^{N+1},~T>0.$ This implies the local boundedness of $\eta$ on $\R\times \R^N$. On the other hand, the global boundedness of $\boldsymbol{\eta}_{\lambda}$ is not guaranteed on $\mathbb{R}\times\mathbb{R}^{N}$. In addition, without loss of generality, suppose that $\|\boldsymbol{\eta}_{\lambda}(0,\cdot)\|_{\infty}=1.$

\begin{remark}\label{periodicity}
If we further assume the periodicity in time $t$ for ${\bf A}^i$, ${\bf q}_i$ and ${\bf F}$, say, $T_0$-periodicity, the function $\boldsymbol{\eta}_{\lambda}$ can be chosen as $\boldsymbol{\eta}_{\lambda}(t,x)=e^{k(\lambda) t}\boldsymbol{\phi}_\lambda(t,x)$, where $(k(\lambda),\boldsymbol{\phi}_\lambda)$ is a principal eigen-pair of the following problem
\begin{align}\label{eigen}
\begin{cases}
\begin{aligned}
\partial_t \phi_{\lambda,1}
- \operatorname{Tr} \big( {\bf A}^1(t,x) D^2 \phi_{\lambda,1} \big)
+\big({\bf q}_1+2\lambda{\bf A}^1e\big)\cdot D\phi_{\lambda,1}
-K_1^\lambda\phi_{\lambda,1}
-\mu_{12}\phi_{\lambda,2}+k(\lambda)\phi_{\lambda,1}&=0,
\end{aligned}\\
\begin{aligned}
\partial_t \phi_{\lambda,2}
- \operatorname{Tr} \big( {\bf A}^2(t,x) D^2 \phi_{\lambda,2} \big)
+\big({\bf q}_2+2\lambda{\bf A}^2e\big)\cdot D\phi_{\lambda,2}
-K_2^\lambda\phi_{\lambda,2}
-\mu_{21}\phi_{\lambda,1}+k(\lambda)\phi_{\lambda,2}&=0,
\end{aligned}\\
K_{i}^\lambda:=\mu_{ii}(t,x) + \lambda^2\, e \cdot {\bf A}^i(t,x)e + \lambda\, {\bf q}_i(t,x)\cdot e,~i=1,2,\\
\boldsymbol{\phi}_{\lambda}>{\bf 0},\\
\boldsymbol{\phi}_{\lambda} \text{ is $T_0$-periodic in time and {\bf L}-periodic in space $x$},
\end{cases}
\end{align}
for $(t,x)\in\mathbb{R}^{N+1}$. Furthermore, $(k(\lambda),\boldsymbol{\phi}_\lambda)$ exists for any $\lambda \in \R$. Existence of the positive periodic eigen-pair follows from the strongly positive period map on the spatial periodic cell and the Krein--Rutman theorem. The off-diagonal coefficients are uniformly positive, so the period map is strongly positive; its leading eigenvector determines the spatially and temporally periodic principal mode after Floquet normalization \cite{Nadin20091,Nadin2009}.
\end{remark}

\subsection{Wave-speed thresholds}

In this section, we study the speed of the waves. First, the estimate \eqref{boundedtwosides} leads to 
\begin{align}\label{globalLi}
\big|\ln  \|\boldsymbol{\eta}_{\lambda}(t+T,\cdot)\|_{\infty}-\ln  \|\boldsymbol{\eta}_{\lambda}(t,\cdot)\|_{\infty} \big| \leq C_3( \lambda^2 +\lambda+1) (1+T),
\end{align}
for some $C_3>0$ large enough. From here, we define $S_{\lambda}:\R\rightarrow \R$ as follows
\begin{align*}
S_{\lambda}(n)=\dfrac{1}{\lambda}\ln \dfrac{\|\boldsymbol{\eta}_{\lambda}(n,\cdot)\|_{\infty}}{\|\boldsymbol{\eta}_{\lambda}(0,\cdot)\|_{\infty}};\quad S_{\lambda}(n+1)=\dfrac{1}{\lambda}\ln \dfrac{\|\boldsymbol{\eta}_{\lambda}(n+1,\cdot)\|_{\infty}}{\|\boldsymbol{\eta}_{\lambda}(0,\cdot)\|_{\infty}},
\end{align*}
and $S_{\lambda}$ is an affine function on $[n,n+1]$ for $n\in \Z$. Thus, one has $S_\lambda(0)=0$ and 
\begin{align*}
|S'_{\lambda}(t)|\leq 	C_3\dfrac{1+\lambda^2}{\lambda},~\forall t\in (n,n+1),~n\in \Z,
\end{align*}
which implies that $S_{\lambda}$ is uniformly Lipschitz-continuous on $\R$. Moreover, for $t\in [n,n+1]$, one has
\begin{align}\label{estimate}
\left|S_{\lambda}(t)-\dfrac{1}{\lambda}\ln \|\boldsymbol{\eta}_{\lambda}(t,\cdot)\|_{\infty}\right|\leq C_3\dfrac{1+\lambda^2}{\lambda}.
\end{align}
This implies that $t\mapsto S_{\lambda}(t)-\dfrac{1}{\lambda}\ln \|\boldsymbol{\eta}_{\lambda}(t,\cdot)\|_{\infty}$ is uniformly bounded on $\R$. At this point, we define 
\begin{align}\label{speed}
c_{\lambda}:=S'_{\lambda}(t),~t\in \R.
\end{align}
Then, $c_{\lambda} \in L^{\infty}(\R)$. From the estimate \eqref{estimate}, we get
\begin{align}\label{clambdaest}
\lfloor c_{\lambda} \rfloor = \dfrac{1}{\lambda}
\lim_{T \to +\infty} \inf_{t \in \mathbb{R}} 
\frac{1}{T} \ln\dfrac{ \|\boldsymbol{\eta}_{\lambda}(t+T,\cdot)\|_{\infty}}{ \|\boldsymbol{\eta}_{\lambda}(t,\cdot)\|_{\infty}}.
\end{align}
Furthermore, \eqref{upper} and \eqref{below} imply that there exists $K_1>0$ such that
\begin{align}\label{estiupper}
K_1\left(C^{-1}\lambda -\max_{i=1,2}\|{\bf q}_i\|_{\infty}+\dfrac{1}{\lambda}\left\lfloor \min_{i=1,2}\min\limits_{x\in\R^{N}}\mu_{ii}\right\rfloor\right)\leq \lfloor c_{\lambda} \rfloor \leq C\lambda +\max_{i=1,2}\|{\bf q}_i\|_{\infty}+\dfrac{1}{\lambda}\left\lfloor \mu_{\infty}\right\rfloor.
\end{align}
An analogous estimate holds for the upper mean. There exists $K_2>0$ such that
\begin{align}\label{estilow}
K_2\left(C^{-1}\lambda -\max_{i=1,2}\|{\bf q}_i\|_{\infty}+\dfrac{1}{\lambda}\left\lceil \min_{i=1,2}\min\limits_{x\in\R^{N}}\mu_{ii}\right\rceil\right)\leq \lceil c_{\lambda} \rceil \leq C\lambda +\max_{i=1,2}\|{\bf q}_i\|_{\infty}+\dfrac{1}{\lambda}\left\lceil \mu_{\infty}\right\rceil,
\end{align}

\begin{remark}\label{c_lambdaper}
In the $T_0$-periodic case, as in Nadin and Rossi \cite{nadin_2015}, one can choose $S_\lambda(t)=(k(\lambda)/\lambda) t$, where $k(\lambda)$ is the eigenvalue in Remark \ref{periodicity}. This implies that $c_\lambda = k(\lambda)/\lambda$. 
\end{remark}

\begin{lemma}\label{lemma continuous}
The functions $\lambda \mapsto \lfloor c_{\lambda} \rfloor$ and $\lambda \mapsto \lceil c_{\lambda} \rceil$ are locally Lipschitz continuous on $(0,\infty)$.
\end{lemma}

\begin{proof}
Fix $0<\epsilon<\Lambda_0$. We show that both functions are Lipschitz
continuous on $[\epsilon,\Lambda_0]$. Let
$\lambda_1,\lambda\in[\epsilon,\Lambda_0]$. If $\lambda=\lambda_1$,
there is nothing to prove. Otherwise, set
\[
\lambda_2:=
\begin{cases}
\lambda_+=2\Lambda_0, & \lambda>\lambda_1,\\
\lambda_-=\epsilon/2, & \lambda<\lambda_1.
\end{cases}
\]
Then, there exists $\theta\in(0,1)$ such that
\(
\lambda=(1-\theta)\lambda_1+\theta\lambda_2
\) and
\(
0<\lambda_-\leq \lambda_1,\lambda_2\leq \lambda_+,~
|\lambda_2-\lambda_1|\geq \epsilon/2.
\)
Define
\begin{align*}
\boldsymbol{\eta}^j(t,x)
&:=\boldsymbol{\eta}_{\lambda_j}(t,x),~
{\bf v}^j(t,x)
:=e^{-\lambda_j x\cdot e}\boldsymbol{\eta}^j(t,x),~j=1,2.
\end{align*}
Then, ${\bf v}^j=(v_1^j,v_2^j)$ solves \eqref{eq:main1} for $j=1,2$. Set ${\bf w}^j:=\ln {\bf v}^j=\left(\ln v^j_1,\ln v^j_2 \right)$. Direct differentiation yields
\begin{equation}\label{nonlinear1}
\begin{cases}
\begin{aligned}
\partial_t w_1^j-\operatorname{Tr}({\bf A}^1D^2w_1^j)
-\operatorname{Tr}({\bf A}^1(Dw_1^j\otimes Dw_1^j))
+q_1\cdot Dw_1^j
&=\mu_{11}+\mu_{12}e^{w_2^j-w_1^j},
\end{aligned}\\
\begin{aligned}
\partial_t w_2^j-\operatorname{Tr}({\bf A}^2D^2w_2^j)
-\operatorname{Tr}({\bf A}^2(Dw_2^j\otimes Dw_2^j))
+q_2\cdot Dw_2^j
&=\mu_{22}+\mu_{21}e^{w_1^j-w_2^j}.
\end{aligned}
\end{cases}
\end{equation}
Here and in the next displayed identity, coefficients without displayed arguments are evaluated at $(t,x)$.

With $\theta \in (0,1)$, setting ${\bf w}:=(1-\theta){\bf w}^1+\theta{\bf w}^2$, direct calculation yields that for $x\in \R^N$, $t\in \R$, 
\begin{align*}
\begin{cases}
\partial_t w_1
-\text{Tr}\!\big({\bf A}^1(t,x) D^2 w_1\big)
-\text{Tr}\!\big({\bf A}^1(t,x)\big(Dw_1 \otimes Dw_1\big)\big)
+ q_1(t,x)\!\cdot Dw_1
\geq \mu_{11}(t,x) + \mu_{12}(t,x)\,e^{\,w_2-w_1},\\
\partial_t w_2
-\text{Tr}\!\big({\bf A}^2(t,x) D^2 w_2\big)
-\text{Tr}\!\big({\bf A}^2(t,x)\big(Dw_2 \otimes Dw_2\big)\big)
+ q_2(t,x)\!\cdot Dw_2
\geq \mu_{22}(t,x) + \mu_{21}(t,x)\,e^{\,w_1-w_2},
\end{cases}
\end{align*}
which implies that {\bf w} is a supersolution of \eqref{nonlinear1}. Hence, $e^{{\bf w}(t,x)}:=\left(e^{w_1(t,x)},e^{w_2(t,x)}\right)$ is a supersolution of \eqref{eq:main1}. It then follows from
\begin{align*}
e^{{\bf w}(t,x)}=e^{-[(1-\theta)\lambda_1+\theta \lambda_2] x\cdot e}\left((\eta^1_1)^{1-\theta}(t,x)(\eta^2_1)^{\theta}(t,x),(\eta^1_2)^{1-\theta}(t,x)(\eta^2_2)^{\theta}(t,x)\right)
\end{align*}
that $(\boldsymbol{\eta}^1)^{1-\theta}(\boldsymbol{\eta}^2)^{\theta}(t,x):=\left((\eta^1_1)^{1-\theta}(t,x)(\eta^2_1)^{\theta}(t,x),(\eta^1_2)^{1-\theta}(t,x)(\eta^2_2)^{\theta}(t,x)\right)$ is a supersolution of \eqref{div} for $\lambda=\lambda_{\theta}:=(1-\theta)\lambda_1+\theta \lambda_2$. The parabolic comparison principle between $(\boldsymbol{\eta}^1)^{1-\theta}(\boldsymbol{\eta}^2)^{\theta}$ and $\boldsymbol{\eta}_{\lambda_\theta}$ leads to 
\begin{align*}
\dfrac{\|\boldsymbol{\eta}_{\lambda_\theta}(t+T,\cdot)\|_{\infty}}{\|\boldsymbol{\eta}_{\lambda_\theta}(t,\cdot)\|_{\infty}}
&\leq \dfrac{\max\limits_{i=1,2}\max\limits_{x\in \R^N}(\eta^1_i)^{1-\theta}(\eta^2_i)^{\theta}(t+T,x)}{\min\limits_{i=1,2}\min\limits_{x\in \R^N}(\eta^1_i)^{1-\theta}(\eta^2_i)^{\theta}(t,x)}
\leq \left(\dfrac{\|\boldsymbol{\eta}^1(t+T,\cdot)\|_{\infty}}{\min\limits_{i=1,2}\min\limits_{x\in \R^N}\eta^1_i(t,x)}\right)^{1-\theta}\left(\dfrac{\|\boldsymbol{\eta}^2(t+T,\cdot)\|_{\infty}}{\min\limits_{i=1,2}\min\limits_{x\in \R^N}\eta^2_i(t,x)}\right)^\theta.
\end{align*}
Thanks to \eqref{below} for $\boldsymbol{\eta}^1$ and $\boldsymbol{\eta}^2$ with the same $C_0>0$, depending only on $\epsilon,\Lambda_0$, we derive that
\begin{align}\label{estimate1}
\dfrac{\|\boldsymbol{\eta}_{\lambda_\theta}(t+T,\cdot)\|_{\infty}}{\|\boldsymbol{\eta}_{\lambda_\theta}(t,\cdot)\|_{\infty}}
&\leq C_0^{-1} \left(\dfrac{\|\boldsymbol{\eta}^1(t+T,\cdot)\|_{\infty}}{\|\boldsymbol{\eta}^1(t,\cdot)\|_{\infty}}\right)^{1-\theta}\left(\dfrac{\|\boldsymbol{\eta}^2(t+T,\cdot)\|_{\infty}}{\|\boldsymbol{\eta}^2(t,\cdot)\|_{\infty}}\right)^\theta.
\end{align}
Set $\Gamma(\lambda):=\lambda\lfloor c_{\lambda} \rfloor$. Equation~\eqref{estimate1} implies
\begin{align*}
\Gamma(\lambda_{\theta})\leq \lim_{T\rightarrow \infty} \inf_{t\in \R}\dfrac{1}{T}\left[(1-\theta)\ln\dfrac{\|\boldsymbol{\eta}^1(t+T,\cdot)\|_{\infty}}{\|\boldsymbol{\eta}^1(t,\cdot)\|_{\infty}}+\theta\ln\dfrac{\|\boldsymbol{\eta}^2(t+T,\cdot)\|_{\infty}}{\|\boldsymbol{\eta}^2(t,\cdot)\|_{\infty}}\right].
\end{align*}
We infer that $\Gamma(\lambda_{\theta})\leq (1-\theta)\Gamma(\lambda_1)+\theta \lambda_2 \lceil c_{\lambda_2} \rceil,~\forall \theta \in (0,1),$ which leads to 
\begin{align*}
\Gamma(\lambda_{\theta}) - \Gamma(\lambda_1)\leq \theta(\lambda_2 \lceil c_{\lambda_2} \rceil - \lambda_1\lfloor c_{\lambda_1} \rfloor).
\end{align*}
Since $\lambda_1,\lambda\in[\epsilon,\Lambda_0]$ and
$\lambda_2\in\{\lambda_-,\lambda_+\}$, the estimates \eqref{estiupper} and \eqref{estilow} yield a constant $K_0>0$, depending only on $\epsilon$, $\Lambda_0$ and the coefficients, such that
\begin{align*}
\Gamma(\lambda)-\Gamma(\lambda_1)=\Gamma(\lambda_\theta)-\Gamma(\lambda_1)
&\leq K_0\theta=\frac{K_0}{|\lambda_2-\lambda_1|}
|\lambda_\theta-\lambda_1|\leq \frac{2K_0}{\epsilon}
|\lambda-\lambda_1|,
\end{align*}
provided that $|\lambda_2-\lambda_1|\geq \epsilon/2$. Since $\lambda,\lambda_1\in[\epsilon,\Lambda_0]$ are arbitrary, interchanging the roles of $\lambda$ and $\lambda_1$ shows that $\Gamma$ is Lipschitz continuous on $[\epsilon,\Lambda_0]$. Finally, since
$\Gamma$ is bounded on $[\epsilon,\Lambda_0]$ by
\eqref{estiupper} and \eqref{estilow}, one has
\begin{align*}
|\lfloor c_\lambda\rfloor-\lfloor c_{\lambda_1}\rfloor|
\leq
\frac{1}{\epsilon}
|\Gamma(\lambda)-\Gamma(\lambda_1)|
+
|\Gamma(\lambda_1)|
\left|\frac{1}{\lambda}-\frac{1}{\lambda_1}\right|
\leq C_{\epsilon,\Lambda_0}|\lambda-\lambda_1|,
\end{align*} 
for some $C_{\epsilon,\Lambda_0}>0.$ The same argument applies to $\lceil c_{\lambda} \rceil$, which completes the proof.
\end{proof} 
 
Next, to define $c^*$, we first consider the following set
\begin{align*}
\Lambda:=\left\{\lambda>0:\exists \overline{k}>0,~\forall 0<k<\overline{k},~\lfloor c_{\lambda} -c_{\lambda+k}\rfloor>0\right\}.
\end{align*}

Following the key idea of \cite{nadin_2015}, we obtain the following result.

\begin{lemma}\label{monotonec}
There exists $\lambda^*>0$ such that $\Lambda=(0,\lambda^*)$. Furthermore, $\lambda \mapsto \lfloor c_{\lambda} \rfloor$ is decreasing on $\Lambda$.

\end{lemma}

\begin{proof}
Fix $\lambda_0, \lambda_1>0$. For $\theta \in (0,1)$, we consider $\lambda_\theta := (1-\theta)\lambda_0 + \theta \lambda_1$. Taking the natural log of \eqref{estimate1} and recalling that $c_\lambda = S'_\lambda$ with $S_\lambda$ satisfying \eqref{estimate} yields
\[
\int_t^{t+T} \big[ (1-\theta)\lambda_0 c_{\lambda_0} + \theta \lambda_1 c_{\lambda_1} - \lambda_\theta c_{\lambda_\theta} \big] ds \geq \ln C_0 - 4\lambda_\theta \alpha,
\]
where $\alpha=C_3\dfrac{1+\lambda^2}{\lambda}$. Hence,
\[
\lambda_\theta \int_t^{t+T} (c_{\lambda_0} - c_{\lambda_\theta})\,ds \geq 
\theta \lambda_1 \int_t^{t+T} (c_{\lambda_0} - c_{\lambda_1})\,ds + \ln C_0 - 4\lambda_\theta \alpha.
\]
Dividing both sides by $T$, taking the infimum over $t \in \R$ and then passing to the limit as $T \to +\infty$, for any $\theta \in (0,1)$, we derive
\begin{align}\label{3.11}
\lfloor c_{\lambda_0} - c_{\lambda_\theta}\rfloor \geq \theta \frac{\lambda_1}{\lambda_\theta} \lfloor c_{\lambda_0} - c_{\lambda_1}\rfloor.
\end{align}
On the other hand, if dividing by $-T$, for any $\theta \in (0,1)$, one has
\begin{align}\label{3.12}
\lfloor c_{\lambda_\theta} - c_{\lambda_0}\rfloor \leq \theta \frac{\lambda_1}{\lambda_\theta} \lfloor c_{\lambda_1} - c_{\lambda_0}\rfloor.
\end{align}
Analogous estimates hold for the upper mean $ \lceil c_{\lambda} \rceil$. We derive the characterization of $\Lambda$ from these inequalities.

We first prove that $\Lambda\neq\emptyset$.  
The estimate \eqref{estiupper}, together with \eqref{muestiamte}, yields  
\[
\lim_{\lambda \to 0^+} \lfloor c_\lambda - c_1\rfloor \geq \lim_{\lambda \to 0^+} \lfloor c_\lambda\rfloor - \lceil c_1\rceil = +\infty.
\]
Then there exists $0 < \lambda_3 < 1$ such that $\lfloor c_{\lambda_3} - c_1\rfloor > 0$. It follows from \eqref{3.11} with $\lambda_0 = \lambda_3$ and $\lambda_1 = 1$ that $\lfloor c_{\lambda_3} - c_{\lambda_3+k}\rfloor > 0$, for all $0 < k < 1 - \lambda_3$. This implies that $\lambda_3 \in \Lambda$.

Next, we show that $\Lambda$ is bounded from above.   By \eqref{estiupper}, we get   
\[
\lim_{\lambda \to +\infty} \lfloor c_1 - c_\lambda\rfloor\leq \lfloor c_1 \rfloor - \lim_{\lambda \to +\infty} \lfloor c_\lambda \rfloor  = -\infty.
\]
Then, there exists $\lambda' > 1$ such that, for any $\lambda > \lambda'$, we have 
$\lfloor c_1 - c_\lambda \rfloor < 0$. 
Hence, for any $k > 0$, by applying \eqref{3.12} with $\lambda_0 = \lambda + k$, $\lambda_1 = 1$ and 
$\theta = k/(k + \lambda - 1)$, we deduce that
\[
\lfloor c_\lambda - c_{\lambda + k} \rfloor 
\le \frac{k}{(k + \lambda - 1)\lambda} 
\lfloor c_1 - c_{\lambda + k} \rfloor < 0.
\]
Namely, $\lambda \notin \Lambda$ and thus $\Lambda$ is bounded from above by $\lambda'$.

We now prove that $\lambda\in\Lambda$ implies $(0,\lambda]\subset\Lambda$. 
Let $0 < \lambda' < \lambda$ and $k > 0$. Using first \eqref{3.11} and then \eqref{3.12} with appropriate $\theta$, we get
\[
\lfloor c_{\lambda'} - c_{\lambda' + k} \rfloor 
\geq \left( \frac{k}{k + \lambda - \lambda'} \right)
\left( \frac{\lambda + k}{\lambda' + k} \right)
\lfloor c_{\lambda'} - c_{\lambda + k} \rfloor
\geq \left( \frac{\lambda + k}{\lambda' + k} \right)
\frac{\lambda}{\lambda'}
\lfloor c_\lambda - c_{\lambda + k} \rfloor.
\]
Thus, $\lambda \in \Lambda$ implies $\lambda' \in \Lambda$.

Finally, we show that $\lambda^*:=\sup\Lambda$ does not belong to $\Lambda$. 
Let $k > 0$. 
For all $n \in \mathbb{N}$, there exists $0 < k_n < 1/n$ such that 
\begin{align*}
\lfloor c_{\lambda^* + 1/n} - c_{\lambda^* + 1/n + k_n} \rfloor \le 0.
\end{align*}
For $n$ large enough, one has $1/n + k_n < k$ and then \eqref{3.11} implies
\[
0 \geq \lfloor c_{\lambda^* + 1/n} - c_{\lambda^* + 1/n + k_n} \rfloor 
\geq \left( \frac{k_n}{k - 1/n} \right)
\left( \frac{\lambda^* + k}{\lambda^* + 1/n + k_n} \right)
\lfloor c_{\lambda^* + 1/n} - c_{\lambda^* + k} \rfloor.
\]
Hence,
\[
\lfloor c_{\lambda^*} - c_{\lambda^* + k} \rfloor 
\leq \lfloor c_{\lambda^* + 1/n} - c_{\lambda^* + k} \rfloor 
+ \lceil c_{\lambda^*} - c_{\lambda^* + 1/n} \rceil 
\leq \lceil c_{\lambda^*} - c_{\lambda^* + 1/n} \rceil.
\]
Using an analogous estimate to \eqref{3.12} for the upper mean, we can obtain that
\[
\lceil c_{\lambda^*} - c_{\lambda^* + 1/n}\rceil  \le \frac{1/n}{\lambda^* + 2/n} \lceil c_{\lambda^*/2} - c_{\lambda^* + 1/n}\rceil \leq \frac{1/n}{\lambda^* + 2/n} \left(\lceil c_{\lambda^*/2}\rceil - \lfloor c_{\lambda^* + 1/n}\rfloor\right),
\]
which goes to $0$ as $n \to \infty$ (recall that $\lambda \mapsto \lfloor c_{\lambda}\rfloor$ is continuous by Lemma \ref{lemma continuous}). Hence, we infer that $ \lfloor c_{\lambda^*} - c_{\lambda^* + k} \rfloor \le 0$; that is, $\lambda^* \notin \Lambda$. Consequently, $\Lambda=(0,\lambda^*)$.

To end the proof, we need to show that $\lambda \mapsto \lfloor c_\lambda\rfloor$ is decreasing on $\Lambda$.  
Assume by way of contradiction that there are $0 < \lambda_1 < \lambda_2<\lambda^*$ such that $\lfloor c_{\lambda_1	}\rfloor \leq \lfloor c_{\lambda_2}\rfloor$.  
The function $\lambda \mapsto \lfloor c_\lambda\rfloor$, being continuous, attains its minimum on $[\lambda_1, \lambda_2]$ at some $\lambda_0$.  
Since $\lfloor c_{\lambda_1	}\rfloor \leq \lfloor c_{\lambda_2}\rfloor$, one can assume that $\lambda_0 \in [\lambda_1, \lambda_2)$.  
The definition of $\Lambda$ implies that there exists $\lambda'_0 \in (\lambda_0, \lambda_2)$ such that $\lfloor c_{\lambda_0} - c_{\lambda'}\rfloor > 0$.  
Hence, we obtain the following contradiction:
\[
\left\lfloor c_{\lambda'_0}\right\rfloor \leq \left\lfloor c_{\lambda_0}\right\rfloor + \left\lceil c_{\lambda'_0} - c_{\lambda_0}\right\rceil  = \left\lfloor c_{\lambda_0}\right\rfloor - \left\lfloor c_{\lambda_0} - c_{\lambda'_0}\right\rfloor < \left\lfloor c_{\lambda_0}\right\rfloor.
\]
This completes the proof.
\end{proof}

We are now in a position to define the critical speed as follows
\begin{align}\label{critical speed}
\lambda^*:=\sup \Lambda,~c^*:=\left\lfloor c_{\lambda^*}\right\rfloor.
\end{align}

\begin{remark}
In the $T_0$-periodic case, the Floquet variational formula and the uniform ellipticity imply the convexity of the tilted principal eigenvalue $\lambda\mapsto k(\lambda)$, as in \cite{Nadin2009,nadin_2015}. The growth of $k(\lambda)/\lambda$ at both endpoints gives
\begin{align*}
c^*=\min_{\lambda>0} \dfrac{k(\lambda)}{\lambda}.
\end{align*}
\end{remark}
\subsection{Existence of transition fronts}

We establish the existence of generalized transition waves by constructing sub- and supersolutions. We first define a family $\left\{\boldsymbol{\varphi}_{\lambda}\right\}_{\lambda>0}$ as follows
\begin{align*}
\boldsymbol{\varphi}_{\lambda}(t,x):=e^{-\lambda S_{\lambda}(t)}\boldsymbol{\eta}_{\lambda}(t,x),~(t,x)\in \R^{N+1}.
\end{align*}
By choosing $T=0$ in \eqref{below}, one obtains a “no–delay” Harnack-type inequality as follows
\begin{align}\label{Harnack-zero}
\min_{i=1,2}\min\limits_{x\in \R^{N}}\eta_{\lambda,i}(t,x)\geq C_0\max_{i=1,2}\max_{x\in \R^{N}}\eta_{\lambda,i}(t,x),~\forall t\in \R.
\end{align}
This, together with \eqref{estimate}, yields constants $C_{\lambda},C'_{\lambda}>0$ such that
\begin{align}\label{estimatevarphi}
C_{\lambda}\leq \varphi_{\lambda,i}(t,x) \leq e^{C'_{\lambda}\lambda },~\forall (t,x)\in \R^{N+1},~i=1,2.
\end{align} 
\begin{theorem}\label{existence}
Assume that \eqref{assumeq} to  \eqref{near_zerp} hold. For $e\in \mathbb{S}^{N-1}$, there exists $c^*=c^*(e)\in \R$ such that for any $\gamma> c^*$, the system \eqref{eq:main} admits a generalized transition wave $({\bf u},c)$ in the sense of Definition \ref{wavedef} and $\left\lfloor c\right\rfloor=\gamma$.
\end{theorem}

\begin{proof}
For $\gamma > c^*$, recall that $\lambda \mapsto \lfloor c_{\lambda} \rfloor$ is continuous on $(0,\infty)$ and is decreasing on $\Lambda=(0,\lambda^*)$. Thanks to the fact that $\lfloor c_{\lambda} \rfloor\rightarrow \infty$ as $\lambda \rightarrow 0^+$, one can find $0<\lambda<\lambda^*$ such that $\lfloor c_{\lambda} \rfloor=\gamma$. Define
\begin{align*}
{\bf \overline{u}}(t,x):=\min\left\{e^{-\lambda x\cdot e}\boldsymbol{\eta}_{\lambda}(t,x),{\bf 1}\right\}=\left(\min \left\{e^{-\lambda x\cdot e}\eta_{\lambda,1}(t,x),1\right\},\min\left\{e^{-\lambda x\cdot e}\eta_{\lambda,2}(t,x),1\right\}\right).
\end{align*}
The componentwise minimum is a generalized supersolution. Indeed, where its $i$-th component equals $1$, cooperativity gives
$F_i(t,x,1,\overline u_j)\le F_i(t,x,1,1)=0$.
Where it equals $e^{-\lambda x\cdot e}\eta_{\lambda,i}<1$,
the other component does not exceed its corresponding linearized profile. The KPP bound \eqref{Fnonlinear}, \eqref{near_zerp} and $\mu_{ij}>0$ then yield
\begin{align*}
F_i(t,x,\overline{\mathbf u})
&\le\sum_{j=1}^2\mu_{ij}(t,x)\overline u_j
\le\sum_{j=1}^2\mu_{ij}(t,x)e^{-\lambda x\cdot e}\eta_{\lambda,j}.
\end{align*}
The last expression equals the diffusion--transport operator applied to the active linearized branch. Across a switching set, write the scalar minimum as a smooth concave approximation of $\min\{1,v_i\}$, where $v_i=e^{-\lambda x\cdot e}\eta_{\lambda,i}$. The chain rule for $\partial_t-\operatorname{Tr}({\bf A}^iD^2)+{\bf q}_i\cdot D$ has the additional term $-\beta''(v_i)Dv_i\cdot {\bf A}^iDv_i\ge0$, since $\beta''\le0$. Passing to the limit gives the supersolution inequality in distributions. Thus, no later result is used to justify the clipping and $F_i$ is evaluated only on $[0,1]^2$.

Next, let us construct a generalized sub-solution. We recall the constant $\nu\in (0,1]$ in \eqref{near_zerp}. Thanks to the definition of $\Lambda$, for the above $\lambda$, there exists $\lambda'\in(\lambda,(1+\nu)\lambda)$ such that $\left\lfloor c_{\lambda}-c_{\lambda'}\right\rfloor>0$. Set
\begin{align*}
\boldsymbol{\psi}(t,x)=(\psi_1(t,x),\psi_2(t,x)):=e^{\sigma(t)-\lambda'\left(x\cdot e-S_{\lambda}(t)+S_{\lambda'}(t)\right)}\boldsymbol{\eta}_{\lambda'}(t,x)=e^{\sigma(t)-\lambda'\left(x\cdot e-S_{\lambda}(t)\right)}\boldsymbol{\varphi}_{\lambda'}(t,x),
\end{align*} 
where $\sigma \in W^{1,\infty}(\R)$ will be determined later. Direct calculation yields that
\begin{align*}
\begin{cases}
\partial_t \psi_1  - \text{Tr}({\bf A}^1(t,x)D^2 \psi_1) + {\bf q}_1(t,x)\cdot D \psi_1 - \mu_{11}(t,x)\psi_1-\mu_{12}(t,x)\psi_2=[\sigma'(t)+\lambda'(c_{\lambda}(t)-c_{\lambda'}(t))]\psi_1, \\[1ex]
\partial_t \psi_2  - \text{Tr} ({\bf A}^2(t,x)D^2 \psi_2 ) + {\bf q}_2(t,x)\cdot D \psi_2  - \mu_{21}(t,x)\psi_1-\mu_{22}(t,x)\psi_2 =[\sigma'(t)+\lambda'(c_{\lambda}(t)-c_{\lambda'}(t))]\psi_2.
\end{cases}
\end{align*}
Thanks to $\left\lfloor \lambda' (c_{\lambda}-c_{\lambda'})\right\rfloor= \lambda'\left\lfloor c_{\lambda}-c_{\lambda'}\right\rfloor>0$, it follows from \eqref{keylemma} that there exists $\sigma\in W^{1,\infty}(\R)$ such that $K_0:=\inf\limits_{t\in \R}\left[\sigma'(t)+\lambda'(c_{\lambda}(t)-c_{\lambda'}(t))\right]>0.$ This implies that
\begin{align*}
\begin{cases}
\partial_t \psi_1  - \text{Tr}({\bf A}^1(t,x)D^2 \psi_1) + {\bf q}_1(t,x)\cdot D \psi_1 - \mu_{11}(t,x)\psi_1-\mu_{12}(t,x)\psi_2\geq K_0\psi_1, \\[1ex]
\partial_t \psi_2  - \text{Tr} ({\bf A}^2(t,x)D^2 \psi_2 ) + {\bf q}_2(t,x)\cdot D \psi_2  - \mu_{21}(t,x)\psi_1-\mu_{22}(t,x)\psi_2 \geq K_0\psi_2.
\end{cases}
\end{align*}

Consider ${\bf w}(t,x)=(w_1(t,x),w_2(t,x)):=e^{-\lambda x\cdot e}\boldsymbol{\eta}_{\lambda}(t,x)	-k\boldsymbol{\psi}(t,x)$ for some $k>0$ that will be chosen later. A direct calculation shows that
\begin{align}\label{propertyw}
\begin{split}
&{\bf w}(t,x)\leq e^{-\lambda x\cdot e}\boldsymbol{\eta}_{\lambda}(t,x),
\\
&{\bf w}(t,x) = e^{-\lambda\left(x \cdot e - S_{\lambda}(t)\right)}
\left[
  \boldsymbol{\varphi}_{\lambda}(t,x)
  - k\, \boldsymbol{\varphi}_{\lambda'}(t,x)\,
  e^{\sigma(t) - (\lambda' - \lambda)\left(x \cdot e - S_{\lambda}(t)\right)}
\right].
\end{split}
\end{align}
Since $\varphi_{\lambda}$, $\varphi_{\lambda'}$ satisfy \eqref{estimatevarphi} for $\lambda$, $\lambda'$, respectively and $\sigma\in W^{1,\infty}(\R)$, from \eqref{propertyw},  choose $k$ large enough so that 
\begin{align*}
\begin{cases}
{\bf w} (t,x)\leq {\bf 0} & \text{ if }x\cdot e-S_{\lambda}(t)\leq 0,\\
{\bf w} (t,x)\leq \delta \bf{1} & \text{ if }x\cdot e-S_{\lambda}(t)> 0,
\end{cases}
\end{align*}
where $\delta\in (0,1]$ is in \eqref{near_zerp}. It follows that
\begin{align*}
\begin{cases}
\partial_t w_1  - \text{Tr}({\bf A}^1(t,x)D^2 w_1) + {\bf q}_1(t,x)\cdot D w_1 - \mu_{11}(t,x)w_1-\mu_{12}(t,x)w_2\leq -kK_0\psi_1, \\[1ex]
\partial_t w_2  - \text{Tr} ({\bf A}^2(t,x)D^2 w_2 ) + {\bf q}_2(t,x)\cdot D w_2  - \mu_{21}(t,x)w_1-\mu_{22}(t,x)w_2 \leq -kK_0\psi_2.
\end{cases}
\end{align*}
Furthermore, in the case ${\bf w}(t,x)>\bf{0}$, this leads to $x\cdot e-S_{\lambda}(t)> 0$ and, for $i=1,2$,
\begin{align*}
-kK_0\psi_i\leq   -kK_0\psi_i \dfrac{w_i^{1+\nu}}{e^{-(1+\nu)\lambda x\cdot e}\eta_{\lambda,i}^{1+\nu}}
&=-kK_0w_i^{1+\nu} \dfrac{\varphi_{\lambda',i}}{\varphi^{1+\nu}_{\lambda,i}}e^{\sigma(t)-[\lambda'-(1+\nu)\lambda]\left(x\cdot e-S_{\lambda}(t)\right)}\\
&\leq -kK_0w_i^{1+\nu}C_{\lambda'}e^{-(1+\nu)\lambda C'_{\lambda}}\inf_{s\in \R}e^{\sigma(s)}\leq -Kw_i^{1+\nu},
\end{align*}
where $K>0$ is the constant in \eqref{near_zerp}. Here, we use \eqref{estimatevarphi} for $\lambda$ and $\lambda'$ to estimate $\varphi_{\lambda',i}/\varphi^{1+\nu}_{\lambda,i}$ and choose $k$ large enough so that $kK_0C_{\lambda'}e^{-(1+\nu)\lambda C'_{\lambda}}\inf\limits_{s\in \R}e^{\sigma(s)}>K$. Hence, ${\bf w}$ is a subsolution of \eqref{eq:main} whenever ${\bf w}>\bf{0}$. Direct calculation yields that
\begin{align*}
{\bf w}(t,x+S_\lambda(t)e)&=e^{-\lambda x \cdot e}
\left(
  \boldsymbol{\varphi}_{\lambda}(t,x+S_\lambda(t)e)
  - k\, \boldsymbol{\varphi}_{\lambda'}(t,x+S_\lambda(t)e)\,
  e^{\sigma(t) - (\lambda' - \lambda)x \cdot e}
\right)\\
&\geq e^{-\lambda x \cdot e}\left(C_\lambda-ke^{\lambda'C'_{\lambda'}+\|\sigma\|_{\infty}- (\lambda' - \lambda)x \cdot e}\right){\bf 1}.
\end{align*}
Hence, for $R_0>0$ large enough, it follows that
\begin{align*}
\inf_{\substack{x \cdot e = R_0 \\ t \in \mathbb{R}}} {\bf w}\bigl(t,x + S_{\lambda}(t)e\bigr)\geq e^{-\lambda R_0}\left(C_\lambda-ke^{\lambda'C'_{\lambda'}+\|\sigma\|_{\infty}- (\lambda' - \lambda)R_0}\right){\bf 1}=:\omega{\bf 1}.
\end{align*}
This implies that $\omega\in (0,1)$ since the right-hand side tends to zero as $R_0\rightarrow \infty$. It then follows from \eqref{KPP} that $F_i(t,x,\omega,\omega)>0$ for $i=1,2$. By decreasing $\omega$ if necessary, we assume that $\omega \in (0,\delta),~
\inf\limits_{\substack{t\in\R\\x\cdot e=S_{\lambda}(t)+R_0}}
\min\limits_{i=1,2}w_i(t,x)>\omega,$ where $\delta$ is in \eqref{near_zerp}. Consequently, the function
\begin{align*}
\underline{\bf u}(t,x):=\begin{cases}
{\bf w}(t,x) &\text{ if }x\cdot e\geq S_{\lambda}(t)+R_0,\\
\max\{\omega{\bf 1},{\bf w}(t,x)\}
&\text{ if }x\cdot e<S_{\lambda}(t)+R_0
\end{cases}
\end{align*}
is continuous on $\R^{N+1}$ and coincides locally with ${\bf w}$
near the interface $x\cdot e=S_{\lambda}(t)+R_0$.
Moreover, the preceding construction implies
$0<\underline{\bf u}\leq\delta{\bf 1}$ and
$\underline{\bf u}\geq{\bf w}$. Fix $i,j\in\{1,2\}$ with $i\neq j$.
The preceding estimate for $-kK_0\psi_i$ holds whenever $w_i>0$,
independently of the sign of $w_j$. Hence,
\begin{align*}
\partial_t w_i-\text{Tr}({\bf A}^iD^2w_i)
+{\bf q}_i\cdot Dw_i
\leq\mu_{ii}w_i+\mu_{ij}w_j-Kw_i^{1+\nu}
\text{ where }w_i>0.
\end{align*}
On each region where $\underline u_i=w_i>0$, the lower estimate
near ${\bf 0}$, together with $\underline u_j\geq w_j$ and
$\mu_{ij}\geq0$, yields
\begin{align*}
F_i(t,x,\underline u_1,\underline u_2)
&\geq\mu_{ii}\underline u_i+\mu_{ij}\underline u_j
-K\underline u_i^{1+\nu}\geq\mu_{ii}w_i+\mu_{ij}w_j-Kw_i^{1+\nu}\geq\partial_t w_i-\text{Tr}({\bf A}^iD^2w_i)
+{\bf q}_i\cdot Dw_i.
\end{align*}
On the region where $x\cdot e<S_{\lambda}(t)+R_0$ and
$\underline u_i=\omega$, one has $\underline u_j\geq\omega$.
By cooperativity,
\begin{align*}
\partial_t\underline u_i
-\text{Tr}({\bf A}^iD^2\underline u_i)
+{\bf q}_i\cdot D\underline u_i
=0<F_i(t,x,\omega,\omega)
\leq F_i(t,x,\underline u_1,\underline u_2).
\end{align*}
The componentwise maximum construction preserves these inequalities
in the generalized sense. Thus, $\underline{\bf u}$ is a generalized
subsolution of \eqref{eq:main}. Moreover, the definitions imply that ${\bf w}\leq {\bf \overline{u}}$ and ${\bf \overline{u}}(t,x+S_{\lambda}(t)e)\geq \min\left\{1,C_{\lambda}e^{-\lambda R_0}\right\}
>\omega$ if $x\cdot e< R_0$. It then follows that $\underline{\bf u}\leq \overline{\bf u}$. Now, the following system
\begin{equation*}
\begin{cases}
\partial_t u_1 - \text{Tr}({\bf A}^1(t,x)D^2 u_1) + {\bf q}_1(t,x)\cdot D u_1 = F_1(t,x,u_1,u_2), & t>-n,\ x\in \R^N, \\[1ex]
\partial_t u_2 -\text{Tr} ({\bf A}^2(t,x)D^2  u_2) + {\bf q}_2(t,x)\cdot Du_2 = F_2(t,x,u_1,u_2), & t> -n,\ x\in \R^N,\\
u_1(-n,x)=\underline u_1(-n,x),\quad u_2(-n,x)=\underline u_2(-n,x),
\end{cases}
\end{equation*}
admits a unique solution ${\bf u}_n$ such that, thanks to the parabolic comparison principle, $\underline{\bf u}\leq {\bf u}_n\leq \overline{\bf u}.$ By applying the standard parabolic estimates, there exists a solution ${\bf u}$ of \eqref{eq:main} such that
\begin{align*}
\underline{\bf u}\leq {\bf u}\leq \overline{\bf u}.
\end{align*}
Following the strong maximum principle, ${\bf u}> {\bf 0}$. In the set $\{x:x\cdot e\geq R \}$ for $R>0$ large enough, one has
\begin{align*}
{\bf u}\left(t,x+e\int_0^t c_{\lambda}(s)ds\right)\leq {\bf \overline{u}	}\left(t,x+e\int_0^t c_{\lambda}(s)ds\right)
&= e^{-\lambda x\cdot e}\boldsymbol{\varphi}_{\lambda}\left(t,x+e\int_0^t c_{\lambda}(s)ds\right)\leq e^{-\lambda x\cdot e}e^{C'_{\lambda}\lambda }{\bf 1}.
\end{align*}
This implies that
\begin{align*}
\lim_{x\cdot e\rightarrow \infty}{\bf u}\left(t,x+e\int_0^t c_{\lambda}(s)ds\right)={\bf 0}\text{ uniformly with respect to }t\in \R.
\end{align*}

On the other hand, it is known that ${\bf u}(t,x)\leq {\bf \overline{u}}(t,x)\leq {\bf 1},~\forall (t,x)\in \R^{N+1}.$ Thus, to obtain the limit behavior at $-\infty$, define
\begin{align*}
\rho
:=
\lim_{r\to-\infty}
\inf_{\substack{x\cdot e\leq r\\ t\in\R}}
\min_{i=1,2}
u_i\bigl(t,x+eS_\lambda(t)\bigr).
\end{align*}
The lower bound provided by $\underline{\bf u}$ implies that
$0<\omega\leq\rho\leq1$. We claim that $\rho=1$. Suppose, by contradiction, that $\rho<1$. There exist sequences
$\{x_n\}_{n\in\N}\subset\R^N$, $\{t_n\}_{n\in\N}\subset\R$,
and, after passing to a subsequence, an index $i_0\in\{1,2\}$ such that $x_n\cdot e\to-\infty,~
u_{i_0}\bigl(t_n,x_n+eS_\lambda(t_n)\bigr)\to\rho.$ Choose $k_n\in\prod_{j=1}^N L_j\Z$ such that $y_n:=x_n+eS_\lambda(t_n)-k_n
\in[{\bf 0},{\bf L})$ and define ${\bf v}_n(t,x):={\bf u}(t+t_n,x+k_n).$ Up to a subsequence, $y_n\to y\in[{\bf 0},{\bf L}]$.
By spatial periodicity, for $i=1,2$,
\begin{align*}
\partial_t v_{n,i}
-\operatorname{Tr}\bigl({\bf A}^i(t+t_n,x)D^2v_{n,i}\bigr)
+{\bf q}_i(t+t_n,x)\cdot Dv_{n,i}
=
F_i(t+t_n,x,{\bf v}_n)
\end{align*}
on $\R\times\R^N$. Since ${\bf 0}\leq{\bf v}_n\leq{\bf 1}$, the local parabolic
estimates and a diagonal extraction yield
${\bf v}_n\to{\bf v}$ locally uniformly in $\R^{N+1}$,
with weak convergence in $W^{1,2}_{p,\loc}(\R^{N+1})$
for every $1<p<\infty$.
Moreover, up to subsequences,
${\bf A}^i(t+t_n,x)$ and ${\bf q}_i(t+t_n,x)$ converge locally
uniformly to $\widetilde{\bf A}^i$ and $\widetilde{\bf q}_i$,
respectively, while $F_i(t+t_n,x,{\bf v}_n)$ converges weakly
in $L^p_{\loc}(\R^{N+1})$ to some $g_i$.
Consequently,
\begin{align*}
\partial_t v_i
-\operatorname{Tr}\bigl(\widetilde{\bf A}^i(t,x)D^2v_i\bigr)
+\widetilde{\bf q}_i(t,x)\cdot Dv_i
=
g_i(t,x) \text{ in $\R^{N+1}$.} 
\end{align*}

We next show that ${\bf v}\geq\rho{\bf 1}$.
For $(t,x)\in\R^{N+1}$, set
\begin{align*}
z_n(t,x)
:=x+k_n-eS_\lambda(t+t_n)=x+x_n-y_n
+e\bigl(S_\lambda(t_n)-S_\lambda(t+t_n)\bigr).
\end{align*}
Since $c_\lambda\in L^\infty(\R)$, $z_n(t,x)\cdot e
\leq
x_n\cdot e+|x-y_n|+\|c_\lambda\|_\infty|t|
\to-\infty$ locally uniformly in $(t,x)$.
By the definition of $\rho$, we have
\begin{align*}
{\bf v}(t,x)
=
\lim_{n\to\infty}
{\bf u}\bigl(t+t_n,z_n(t,x)+eS_\lambda(t+t_n)\bigr)
\geq\rho{\bf 1}.
\end{align*}
On the other hand, $v_{n,i_0}(0,y_n)
=
u_{i_0}\bigl(t_n,x_n+eS_\lambda(t_n)\bigr)\to\rho$ and hence $v_{i_0}(0,y)=\rho$. By the diagonal positivity assumption \eqref{KPP}, $m_\rho
:=
\min\limits_{i=1,2}
\inf\limits_{(t,x)\in\R^{N+1}}
F_i(t,x,\rho,\rho)>0.$ Let $K>0$ be a uniform Lipschitz constant
for $F_i$ with respect to $u_i$ on $[{\bf 0},{\bf 1}]$.
Cooperativity implies that, for every
${\bf a}=(a_1,a_2)\in[\rho,1]^2$,
\begin{align*}
F_i(t,x,{\bf a})
\geq
F_i(t,x,\rho,\rho)-K(a_i-\rho)
\geq
m_\rho-K(a_i-\rho),~i=1,2.
\end{align*}
Since ${\bf v}_n\to{\bf v}$ locally uniformly and the reactions
are uniformly Lipschitz in ${\bf u}$, one has
\begin{align*}
F_i(t+t_n,x,{\bf v}_n)
-
F_i(t+t_n,x,{\bf v})
\to0 \text{ locally uniformly. }
\end{align*}
Passing to the weak limit in the previous inequality, with ${\bf a}={\bf v}(t,x)$, yields $g_i(t,x)\geq m_\rho-K\bigl(v_i(t,x)-\rho\bigr) \text{ in $\R^{N+1}$.}$ Thus, $w:=v_{i_0}-\rho$ is nonnegative, satisfies $w(0,y)=0$ and
\begin{align*}
\partial_t w
-\operatorname{Tr}\bigl(\widetilde{\bf A}^{i_0}(t,x)D^2w\bigr)
+\widetilde{\bf q}_{i_0}(t,x)\cdot Dw
+Kw
\geq m_\rho>0 \text{ in $\R^{N+1}$.}
\end{align*}
The strong maximum principle implies that
$w\equiv0$ on $(-\infty,0]\times\R^N$, contradicting
the preceding strict inequality.
Therefore, $\rho=1$. Finally, since ${\bf u}\leq{\bf 1}$, the definition of $\rho$
implies that
\begin{align*}
\lim_{x\cdot e\to-\infty}
{\bf u}\left(t,x+e\int_0^t c_\lambda(s)\,ds\right)
=
{\bf 1} \text{ uniformly with respect to $t\in\R$.}
\end{align*}
 We obtain the desired result.
\end{proof}

\section{Nonexistence Criteria}\label{sec4}

\subsection{Lower spectral threshold}
In this section, we construct a subsolution to system \eqref{eq:main1} to derive a lower bound for the wave speed function $c=c(t)$ of transition fronts. We recall that
\begin{align*}
\begin{cases}
\partial_t u_1  - \text{Tr}({\bf A}^1(t,x)D^2 u_1) + {\bf q}_1(t,x)\cdot D u_1 = \mu_{11}(t,x)u_1+\mu_{12}(t,x)u_2, & t\in \R,\ x\in \R^N, \\[1ex]
\partial_t u_2  - \text{Tr} ({\bf A}^2(t,x)D^2 u_2) + {\bf q}_2(t,x)\cdot D u_2 = \mu_{21}(t,x)u_1+\mu_{22}(t,x)u_2, & t\in \R,\ x\in \R^N,
\end{cases}
\end{align*}
Consider the form ${\bf u}(t,x):=e^{-\lambda (x\cdot e -ct)}\boldsymbol{\phi}(t,x)$ for some constants $\lambda,~c$. Direct calculation yields that
\begin{align*}
({\bf P}_{\lambda}+c\lambda )\boldsymbol{\phi}=0
\end{align*}
where ${\bf P}_{\lambda}=(P_{\lambda,1},P_{\lambda,2})$ is defined, for $(t,x)\in \R\times \R^N$ and ${\bf w}=(w_1,w_2)$, by
\begin{align}\label{operatorP}
\begin{cases}
P_{\lambda,1}{\bf w}:=\partial_t w_1
- \text{Tr} \big( {\bf A}^1(t,x) D^2 w_1 \big)
+\big( {\bf q}_1(t,x) + 2\lambda {\bf A}^1(t,x)e \big) \cdot D w_1- K_1^\lambda w_1-\mu_{12}(t,x)w_2,\\
P_{\lambda,2}{\bf w}:=\partial_t w_2
- \text{Tr} \big( {\bf A}^2(t,x) D^2 w_2 \big)
+ \big( {\bf q}_2(t,x) + 2\lambda {\bf A}^2(t,x)e \big) \cdot D w_2-K_2^\lambda w_2-\mu_{21}(t,x)w_1,\\
\end{cases}
\end{align}
where $K_{i}^\lambda:=\mu_{ii}(t,x) + \lambda^2\, e \cdot {\bf A}^i(t,x)e + \lambda\, {\bf q}_i(t,x)\cdot e,~i=1,2.$ Recall the generalized principal eigenvalue.
\begin{align}\label{GPE}
\kappa(\lambda):=\inf\left\{ k\in\R:
\exists {\bf w}\in\mathcal X,\quad
\inf_{\R\times\R^N}{\bf w}>{\bf0},~
{\bf P}_\lambda{\bf w}\leq k{\bf w}
\quad\text{in }\R\times\R^N\right\}.
\end{align}
where $\mathcal X:=\Bigl[
W^{1,2,N+1}_{\mathrm{loc}}(\R\times\R^N)
\cap W^{0,1,\infty}(\R\times\R^N)
\Bigr]^2.$ This means that ${\bf w}$ and $D{\bf w}$ are bounded on $\R\times \R^N$. Taking ${\bf w}\equiv {\bf 1}$ in the definition, we obtain
\begin{align}\label{boundedabove}
\kappa(\lambda)\leq -C^{-1} \lambda^2 +|\lambda|\max_{i=1,2}\|{\bf q}_i\|_{\infty}-2\min_{i,j=1,2}\inf\limits_{\R\times\R^N} \mu_{ij},
\end{align}
where $C>0$ is in \eqref{Aassump}. Next, we prove that $\kappa(\lambda)>-\infty$. Suppose otherwise, $\kappa(\lambda)=-\infty$. Then, there exist ${\bf w}$ and $k$ as in the definition of $\kappa(\lambda)$ such that $k<-C\lambda^2-|\lambda|\max\limits_{i=1,2}\|{\bf q}_i\|_{\infty}-2\max\limits_{i,j=1,2}\sup\limits_{\R\times\R^N} \mu_{ij}.$ Now, for $\beta>0$, define ${\bf v}(t,x):=e^{-\beta t}(1,1)$. Direct calculation yields
\begin{align*}
P_{\lambda,i}{\bf v}\geq \left(-\beta -C\lambda^2-|\lambda|\max_{i=1,2}\|{\bf q}_i\|_{\infty}-2\max_{i,j=1,2}\sup\limits_{\R\times\R^N} \mu_{ij}\right)v_i,~\forall i=1,2.
\end{align*}
For sufficiently small $\beta$, the right-hand side is greater than $k$, which implies that ${\bf P}_{\lambda}{\bf v}\geq k {\bf v}$. Since ${\bf w}$ is bounded, we deduce that ${\bf v}\geq {\bf w}$ for any $t\leq t_0$ for some $t_0\in \R$. By the comparison principle, ${\bf w}(t,x)\leq e^{-\beta t},~\forall (t,x)\in \R\times\R^N.$ This implies that $0<\inf\limits_{\R\times \R^N} {\bf w}\leq e^{-\beta t}$. Contradiction. Consequently, $\kappa(\lambda)>-\infty$. Define
\begin{align}\label{criticalbelow}
c_*:=-\max_{\lambda>0}\dfrac{\kappa(\lambda)}{\lambda}.
\end{align}
This definition is well-posed if $\kappa(0)<0$ since $\kappa(\lambda)/\lambda \rightarrow -\infty $ as $\lambda \rightarrow \infty$ by \eqref{boundedabove}. Furthermore, by calculations similar to those in \cite[Proposition 4.3]{BerestyckiNadinMemoirs2022}, it is known that $\lambda \mapsto \kappa(\lambda)$ is locally Lipschitz-continuous on $\R$. Now, we show that the two-sided estimates \eqref{muestiamte} imply $\kappa(0) < 0$. For $\sigma \in W^{1,\infty}(\R)$, consider $\boldsymbol{\varphi}(t,x)=e^{-\sigma(t)}{\bf 1}$. Then, one has
\begin{align*}
{\bf P}_0\boldsymbol{\varphi}=-e^{-\sigma(t)}(\sigma'(t)+\mu_{11}+\mu_{12},\sigma'(t)+\mu_{21}+\mu_{22})\leq -\boldsymbol{\varphi}\left(\sigma'(t) +\min_{i=1,2}\inf\limits_{\R\times\R^N} \displaystyle\sum_{j=1}^2 \mu_{ij}\right).
\end{align*}
Thanks to \eqref{keylemma}, for any $\epsilon>0$, there exists $\sigma \in W^{1,\infty}(\R)$ such that
\begin{align*}
{\bf P}_0\boldsymbol{\varphi} \leq -\boldsymbol{\varphi}\left(\sigma'(t) +\min_{i=1,2}\inf\limits_{\R\times\R^N} \displaystyle\sum_{j=1}^2 \mu_{ij}\right) \leq -\boldsymbol{\varphi}\left( \left\lfloor \min_{i=1,2}\inf\limits_{\R\times\R^N} \displaystyle\sum_{j=1}^2 \mu_{ij}\right\rfloor-\epsilon\right).
\end{align*}
The two-sided estimates \eqref{muestiamte} imply that $\left\lfloor \min\limits_{i=1,2}\inf\limits_{\R\times\R^N} \displaystyle\sum_{j=1}^2 \mu_{ij}\right\rfloor>0$. We then choose $\epsilon>0$ small enough to derive that $\kappa(0)<0$. Here, $\boldsymbol{\varphi}$ satisfies the definition of $\kappa(\lambda)$. Hence, \eqref{criticalbelow} is well-defined in the sense that there exists $\lambda_*>0$ such that
\begin{align*}
c_*=-\max_{\lambda>0}\dfrac{\kappa(\lambda)}{\lambda}=-\dfrac{\kappa(\lambda_*)}{\lambda_*}.
\end{align*}

Next, following the idea of \cite[Section 4A]{nadin_2015}, we construct a positive subsolution with prescribed exponential behavior at both ends for \eqref{eq:main1}.
\begin{lemma}\label{lemma7}
Let $c,\lambda\in\mathbb R$ satisfy $\kappa(\lambda)+c\lambda<0$. Set $\rho:=x\cdot e-ct$ and $R_z:=1+|z|$ for $z\in \R$. Then, there exist $\varepsilon>0$ and $M>1$ such that, for every $z\in\mathbb R$, \eqref{eq:main1} admits a
subsolution $\underline{\mathbf v}$ with the following properties.
\begin{align*}
\begin{cases}
M^{-1}e^{-(\lambda+\varepsilon)\rho}\mathbf 1
\leq \underline{\mathbf v}(t,x)
\leq Me^{-(\lambda+\varepsilon)\rho}\mathbf 1,
&\text{if }\rho\geq z,\\
M^{-1}e^{-(\lambda-\varepsilon)\rho}\mathbf 1
\leq \underline{\mathbf v}(t,x)
\leq Me^{-(\lambda-\varepsilon)\rho}\mathbf 1,
&\text{if }\rho\leq z-R_z,\\
\inf\limits_{\substack{(t,x)\in\mathbb R^{N+1}\\z-R_z<\rho<z}}
\min\limits_{i=1,2}\underline v_i(t,x)>0.
\end{cases}
\end{align*}
\end{lemma}

\begin{proof}
By the definition of $\kappa(\lambda)$, there exist $\mathbf w$ and
$k<-c\lambda$ such that
\[
\inf_{\mathbb R^{N+1}}\min_{i=1,2}w_i>0,~
\mathbf P_\lambda\mathbf w\leq k\mathbf w
\quad\text{in }\mathbb R^{N+1}.
\]
Set
\(
m:=\inf\limits_{\mathbb R^{N+1}}\min\limits_{i=1,2}w_i>0,~
W:=\sup\limits_{\mathbb R^{N+1}}\max\limits_{i=1,2}w_i<\infty.
\) Choose a smooth nondecreasing function
$\theta:\mathbb R\to[0,1]$ such that $\theta(s)=0$ for $s\leq0$
and $\theta(s)=1$ for $s\geq1$. For each $z\in\mathbb R$, define
\[
R_z:=1+|z|,~
\zeta_z(\rho)
:=
\lambda-\varepsilon+
2\varepsilon\,
\theta\!\left(\frac{\rho-z+R_z}{R_z}\right).
\]
Then $\zeta_z=\lambda-\varepsilon$ for $\rho\leq z-R_z$ and
$\zeta_z=\lambda+\varepsilon$ for $\rho\geq z$. On the transition
interval $I_z=[z-R_z,z]$,
\[
|\rho|\leq 2R_z,~
|\zeta_z'|\leq C_\theta\frac{\varepsilon}{R_z},~
|\zeta_z''|\leq C_\theta\frac{\varepsilon}{R_z^2},
\]
where $C_\theta$ depends only on $\theta$. Recall $\rho:=x\cdot e-ct$ and define $\underline{\mathbf v}(t,x)
:=
e^{-\rho\zeta_z(\rho)}\mathbf w(t,x).$
For $i=1,2$, write
\begin{align*}
\mathcal R_i(\mathbf v)&:=\partial_t v_i
-\operatorname{Tr}({\bf A}^iD^2v_i)+{\bf q}_i\cdot Dv_i
-\sum_{j=1}^2\mu_{ij}v_j.
\end{align*}
Direct calculation yields, componentwise,
\begin{equation*}
\big[\mathcal R_i(\underline{\mathbf v})\big]
 e^{\rho\zeta_z(\rho)}
\leq\big(P_{\zeta_z,i}+c\zeta_z\big){\bf w}+CK_z(\rho),
\qquad i=1,2.
\end{equation*}
where $K_z(\rho)
:=
\bigl[(1+|\rho|+|\rho\zeta_z|+\rho^2|\zeta_z'|)
|\zeta_z'|+|\rho\zeta_z''|\bigr].$ Here $\zeta_z,\zeta_z',\zeta_z''$ are evaluated at $\rho$, and $C$
depends on the coefficient bounds, $c$, $N$, and the bounds for
$\mathbf w$ and $D\mathbf w$, but not on $z$. Since $|\zeta_z|\leq|\lambda|+\varepsilon$ and $|\rho|\leq2R_z$ on
$I_z$, the derivative bounds imply
\(
K_z(\rho)\leq C_{\theta,\lambda}(\varepsilon+\varepsilon^2),~ \rho\in I_z,
\)
with a constant independent of $z$. Outside $I_z$, the derivatives of
$\zeta_z$ vanish, so $K_z=0$. Moreover, $|\zeta_z-\lambda|\leq\varepsilon$ and
$|\zeta_z^2-\lambda^2|
\leq 2|\lambda|\varepsilon+\varepsilon^2$. Hence, for $i=1,2$,
\[
(P_{\zeta_z,i}+c\zeta_z)\mathbf w
\leq
(P_{\lambda,i}+c\lambda)\mathbf w
+C\bigl(\varepsilon+2|\lambda|\varepsilon+\varepsilon^2\bigr)
\leq
(k+c\lambda)w_i
+C\bigl(\varepsilon+2|\lambda|\varepsilon+\varepsilon^2\bigr).
\]
Thus, increasing $C$ if necessary, we deduce that
Writing $S_z=e^{-\rho\zeta_z(\rho)}$ and $R_i=S_z[(k+c\lambda)w_i+C(\varepsilon+\varepsilon^2)]$, the preceding two bounds become
\begin{equation*}
\begin{cases}
\partial_t\underline v_1-\operatorname{Tr}({\bf A}^1D^2\underline v_1)+{\bf q}_1\cdot D\underline v_1-\mu_{11}\underline v_1-\mu_{12}\underline v_2\leq R_1,\\
\partial_t\underline v_2-\operatorname{Tr}({\bf A}^2D^2\underline v_2)+{\bf q}_2\cdot D\underline v_2-\mu_{21}\underline v_1-\mu_{22}\underline v_2\leq R_2.
\end{cases}
\end{equation*}
The coefficients in these inequalities are evaluated at $(t,x)$.
Since $k+c\lambda<0$ and $w_i\geq m$, choose $\varepsilon\in(0,1)$ small enough that \(C(\varepsilon+\varepsilon^2)<-(k+c\lambda)m.\) This choice is independent of $z$. Thus, $\underline{\mathbf v}$ is a subsolution of \eqref{eq:main1}. Outside $I_z$, $\zeta_z$ is constant. Therefore, the two tail estimates hold with \(M:=\max\{1,W,m^{-1}\},\) which is independent of $z$. For each fixed $z$, the interval $I_z$ is bounded and $\zeta_z$ is bounded there. Since $w_i\geq m>0$, the last assertion follows as well. Note that no global boundedness is claimed at this stage. If $\lambda-\varepsilon>0$, the lower estimate on the left implies that $\underline{\mathbf v}$ is unbounded as $\rho\to-\infty$. Boundedness is obtained only after the gluing construction in Lemma~\ref{lemma8}.
\end{proof}

\begin{lemma}\label{lemma8}
Let $c,\underline{\lambda},\overline{\lambda}\in \R$ satisfy $\underline{\lambda}\leq 0\leq \overline{\lambda}$, $\underline{\lambda}\neq \overline{\lambda}$ and
\begin{align*}
\max_{\lambda \in [\underline{\lambda},\overline{\lambda}]}(\kappa(\lambda)+c\lambda)<0.
\end{align*}
Then, \eqref{eq:main1} admits a generalized, bounded subsolution $\underline{{\bf v}}$ such that
\begin{align*}
\lim_{r\rightarrow -\infty} \sup_{x\cdot e-ct <r} \underline{{\bf v}}(t,x)e^{\underline{\lambda}(x\cdot e-ct)}={\bf 0},\quad \lim_{r\rightarrow \infty} \sup_{x\cdot e -ct>r} \underline{{\bf v}}(t,x)e^{\overline{\lambda}(x\cdot e-ct)}={\bf 0},~\forall r_1<r_2,~ \inf_{r_1\leq x\cdot e-ct \leq r_2}\underline{{\bf v}}(t,x)>{\bf 0}.
\end{align*}
\end{lemma}

\begin{proof}

For every $\lambda\in[\underline\lambda,\overline\lambda]$, let
$\varepsilon_\lambda>0$ and $M_\lambda>1$ be as in
Lemma~\ref{lemma7}, and put
$I_\lambda=(\lambda-\varepsilon_\lambda,
\lambda+\varepsilon_\lambda)$. These intervals cover the compact
interval $[\underline\lambda,\overline\lambda]$. Choose a finite
subcover and remove intervals that are contained in other selected
intervals or do not meet $[\underline\lambda,\overline\lambda]$.
Retain a subcover of minimal cardinality, and number its intervals
from right to left. Define
\(
a_i=\lambda_i-\varepsilon_i,~
b_i=\lambda_i+\varepsilon_i,~
\varepsilon_i=\varepsilon_{\lambda_i},~
M_i=M_{\lambda_i},
\)
we have $b_1>\overline\lambda$, $a_n<\underline\lambda$, and,
if $n\ge2$,
\[
a_{i+1}<a_i<b_{i+1}<b_i,~ i=1,\ldots,n-1.
\]
Indeed, the left and right endpoints are strictly ordered because
no selected interval is contained in another. Two consecutive
intervals must overlap, since otherwise the selected intervals
would leave a point of $[\underline\lambda,\overline\lambda]$
uncovered. In particular, $d_i:=b_{i+1}-a_i>0$.

For each $i$, let $\mathbf v_i$ denote the positive subsolution
constructed in Lemma~\ref{lemma7}, with $\lambda=\lambda_i$
and a translation parameter $z_i$ to be chosen. Set
$z_1=0$, $k_1=1$, and $\ell_i=z_i-R_{z_i}$.
We choose $z_2,\ldots,z_n$ successively so that
\[
z_{i+1}<\ell_i,~
\exp\bigl(d_i(\ell_i-z_{i+1})\bigr)
>4(M_iM_{i+1})^2,~ i=1,\ldots,n-1,
\]
and define
\[
k_{i+1}:=
\frac{k_i\exp(d_i\ell_i)}{2M_iM_{i+1}}.
\]
The choices are possible because $d_i>0$, and the constants
$\varepsilon_i,M_i$ are independent of $z_i$. For every $(t,x)$ with $\rho=x\cdot e-ct=\ell_i$, the left
exponential estimate for $\mathbf v_i$ and the right estimate
for $\mathbf v_{i+1}$ imply, componentwise,
\[
\frac{k_{i+1}v_{i+1,j}(t,x)}{k_iv_{i,j}(t,x)}
\le \frac{k_{i+1}}{k_i}M_iM_{i+1}
\exp(-d_i\ell_i)=\frac12,
\qquad j=1,2.
\]
At $\rho=z_{i+1}$, the same estimates imply
\[
\frac{k_{i+1}v_{i+1,j}(t,x)}{k_iv_{i,j}(t,x)}
\ge \frac{k_{i+1}}{k_i}(M_iM_{i+1})^{-1}
\exp(-d_iz_{i+1})>2,
\qquad j=1,2.
\]
Thus $k_i\mathbf v_i$ dominates at the right endpoint
$\ell_i$, whereas $k_{i+1}\mathbf v_{i+1}$ dominates at the
left endpoint $z_{i+1}$. These inequalities are uniform in
$(t,x)$ on the corresponding moving hyperplanes.

For $n=1$, set $\underline{\mathbf v}=\mathbf v_1$.
For $n\ge2$, define $\underline{\mathbf v}$ by the following
nonoverlapping regions, with $\rho=x\cdot e-ct$:
\[
\underline{\mathbf v}(t,x)=
\begin{cases}
k_1\mathbf v_1(t,x),&\rho\ge\ell_1,\\
\max\{k_i\mathbf v_i(t,x),k_{i+1}\mathbf v_{i+1}(t,x)\},
&z_{i+1}\le\rho\le\ell_i,\quad 1\le i<n,\\
k_i\mathbf v_i(t,x),&\ell_i<\rho<z_i,
\quad 2\le i<n,\\
k_n\mathbf v_n(t,x),&\rho<z_n.
\end{cases}
\]
The maximum is taken componentwise. The inequalities at
$\rho=\ell_i$ and $\rho=z_{i+1}$ show that adjacent
expressions coincide at every interface. In fact, their strict
uniform separation there and the continuity of the constructed
subsolutions ensure that, locally near each point of an interface, the glued
function agrees with the corresponding single branch.
On every joining strip it is the componentwise maximum of two
positive subsolutions of a cooperative linear system, hence it
is a generalized subsolution. To verify the last assertion,
consider a component on which one branch is active. Its own
differential inequality remains valid because every other
component of the maximum is no smaller and $\mu_{ij}>0$
for $i\ne j$. At the switching set the same inequality follows
in the weak sense from the scalar parabolic Kato inequality,
or equivalently by smoothing the maximum and passing to the
limit. Consequently, the assembled function is a generalized
subsolution on $\mathbb R^{N+1}$.

Its right tail is $k_1\mathbf v_1$, and $b_1>\overline\lambda$;
its left tail is $k_n\mathbf v_n$, and
$a_n<\underline\lambda$. By the estimates in Lemma~\ref{lemma7},
\[
\sup_{t\in\mathbb R}\sup_{x\cdot e-ct>r}
\bigl|\underline{\mathbf v}(t,x)\bigr|
\exp\bigl(\overline\lambda(x\cdot e-ct)\bigr)
\rightarrow0\quad(r\to+\infty),
\]
and the analogous weighted limit with $\underline\lambda$
holds as $r\to-\infty$. Each of the finitely many joining
strips has bounded $\rho$-width; its constituent functions
are uniformly bounded above and below by positive constants
there. Hence $\underline{\mathbf v}$ is bounded on
$\mathbb R^{N+1}$, and for every $r_1<r_2$,
\[
\inf_{\substack{(t,x)\in\mathbb R^{N+1}\\
r_1\le x\cdot e-ct\le r_2}}
\min_{j=1,2}\underline v_j(t,x)>0.
\]
This completes the proof.

\end{proof}

\begin{theorem}\label{theo:9}
Assume \eqref{assumeq}-\eqref{KPP} hold. Suppose further that $\kappa(0)<0$. Then, for any $({\bf u},c)$ such that ${\bf u}$ is a nonnegative ${[0,1]^2}$-valued supersolution of \eqref{eq:main} and $c\in L^{\infty}(\R)$ is a wave speed function satisfying \eqref{limits}, one has
\begin{align*}
\liminf_{t\rightarrow \infty}\dfrac{1}{t}\int_0^tc(s)ds \geq c_*=-\max_{\lambda>0}\dfrac{\kappa(\lambda)}{\lambda},
\end{align*}
where $\kappa(\lambda)$ is in \eqref{GPE}.

\end{theorem}
\begin{proof}

\begingroup
Put
\[
S(t)=\int_0^t c(s)\,ds,~
\boldsymbol\phi(t,x)=\mathbf u(t,x+eS(t)),~
\Omega_R=\{x\in\mathbb R^N:x\cdot e>-R-1\}.
\]
By \eqref{limits}, the convergence of \(\boldsymbol\phi(t,x)\) to
\(\mathbf1\) as \(x\cdot e\to-\infty\) is uniform in \(t\).
Choose \(R>0\) so large that the two components of
\(\boldsymbol\phi\) are bounded away from zero on
\(\{(t,x):x\cdot e\le-R\}\), and define the scalar
\[
\beta=\min\left\{\frac12,
\inf_{\substack{t\in\mathbb R\\-R-1\le x\cdot e\le-R}}
\min_{i=1,2}\phi_i(t,x)\right\}>0.
\]
For each \(n\in\mathbb N\), the following auxiliary functions are
\textup{defined only on} \((-n,\infty)\times\overline{\Omega_R}\):
\begin{align*}
\mathbf h_n(t,x)&=w_n(t,x)\mathbf1,\\
w_n(t,x)&=
\begin{cases}
\displaystyle \beta\left(1-\frac{x\cdot e+R+1}{t+n+1}\right)^{\lambda(t+n+1)},
&0\le x\cdot e+R+1<t+n+1,\\
0,&x\cdot e+R+1\ge t+n+1.
\end{cases}
\end{align*}
Here \(\lambda>2\) will be chosen independently of \(n\).
Writing \(s=x\cdot e+R+1\) and \(T=t+n+1\), we have
\(s\ge0\), \(T>1\) and
\(0\le w_n\le\beta<1\) throughout this half-cylinder.
In particular, each reaction term is evaluated in its prescribed
domain \([0,1]^2\). In the region \(0<s<T\), put
\(\rho_n=1-s/T\in(0,1)\).
\endgroup
For $i=1,2$, direct calculation yields that
\begin{align}\label{sub1}
\begin{split}
&\dfrac{\rho_n^2}{w_n}
\left(
\begin{aligned}
&\partial_t h_{n,i}(t,x)
-\operatorname{Tr}\!\left({\bf A}^i(t,x+eS(t))D^2h_{n,i}\right)\\
&+\left({\bf q}_i(t,x+eS(t))-c(t)e\right)\cdot Dh_{n,i}
-F_i\!\left(t,x+eS(t),h_{n,1},h_{n,2}\right)
\end{aligned}
\right)\\
&=
\lambda\rho_n^2\ln\rho_n
+\dfrac{x\cdot e+R+1}{t+n+1}\lambda\rho_n
-e\cdot{\bf A}^i(t,x+eS(t))e
\left(
\lambda^2-\dfrac{\lambda}{t+n+1}
\right)\\
&\quad
-\lambda
\left({\bf q}_i(t,x+eS(t))-c(t)e\right)\cdot e\,\rho_n
-\dfrac{\rho_n^2}{w_n}
F_i\!\left(t,x+eS(t),h_{n,1},h_{n,2}\right)\\
&\leq
\lambda
\left(
1-C^{-1}\lambda+C
+\|{\bf q}_i\|_{\infty}
+\|c\|_{\infty}
\right).
\end{split}
\end{align}
\begingroup
Here \(C\) is the ellipticity constant in \eqref{Aassump}. Since
\(0<w_n\le\beta e^{-\lambda s}<1\) for \(0<s<T\),
\eqref{KPP} implies
\(F_i(t,x+eS(t),w_n,w_n)>0\).
Also \(\rho_n^2\log\rho_n\le0\), \(s/T\le1\), and
\(e\cdot\mathbf A^i e\ge C^{-1}\). Hence the right-hand side
of \eqref{sub1} is negative whenever
\[
\lambda=\lambda_0>\max\left\{2,
C\left(1+C+\max_{i=1,2}\|\mathbf q_i\|_\infty
+\|c\|_\infty\right)\right\}.
\]
The choice is independent of \(n\). On the open set \(s>T\),
\(\mathbf h_n=\mathbf0\) and \(\mathbf F(t,x,\mathbf0)=\mathbf0\).
At \(s=T\), the exponent \(\lambda_0 T>2\) ensures that
\(w_n\), its time derivative and its first and second spatial
derivatives tend to zero from the region \(s<T\). Thus the
cut-off introduces no interface measure, and \(\mathbf h_n\)
is a bounded generalized subsolution on
\((-n,\infty)\times\Omega_R\).

We check the parabolic boundary data explicitly. On the lateral
boundary \(x\cdot e=-R-1\), one has \(s=0\) and
\(
\mathbf h_n(t,x)=\beta\mathbf1\le\boldsymbol\phi(t,x)
~(t\ge-n).
\) At \(t=-n\), the function extends continuously by
\[
\mathbf h_n(-n,x)=
\begin{cases}
\beta(1-s)^{\lambda_0}\mathbf1,&0\le s<1,\\
\mathbf0,&s\ge1.
\end{cases}
\]
Whenever this function is nonzero, \(-R-1\le x\cdot e<-R\),
so that \(\mathbf h_n(-n,x)\le\beta\mathbf1
\le\boldsymbol\phi(-n,x)\). The inequality is immediate
where \(s\ge1\). The profile \(\boldsymbol\phi\) is a bounded supersolution of
\[
\partial_t\phi_i-\operatorname{Tr}\bigl(
\mathbf A^i(t,x+eS(t))D^2\phi_i\bigr)
+\bigl(\mathbf q_i(t,x+eS(t))-c(t)e\bigr)\cdot D\phi_i
=F_i(t,x+eS(t),\boldsymbol\phi).
\]
The system is cooperative and the reaction is uniformly
Lipschitz on \([0,1]^2\). For each fixed \(T_1>-n\),
the bounded comparison principle on
\([-n,T_1]\times\Omega_R\) applies to \(\mathbf h_n\) and
\(\boldsymbol\phi\). For completeness, one may exhaust the
half-space by bounded domains and use
\(\varepsilon e^{K(t+n)}(1+|x|^2)\) to control the artificial
outer boundary; the ellipticity, drift and Lipschitz bounds
allow a fixed sufficiently large \(K\) on every finite time
interval. Letting \(\varepsilon\downarrow0\) and then
\(T_1\uparrow\infty\), we obtain
\[
\mathbf h_n(t,x)\leq\boldsymbol\phi(t,x),~ t\ge-n,~ x\in\Omega_R.
\]
For fixed \((t,x)\in\mathbb R\times\Omega_R\), the inequality
holds for all sufficiently large \(n\). Since \(s>0\) is fixed,
\[
\lim_{n\to\infty}\left(1-\frac{s}{t+n+1}\right)^{
\lambda_0(t+n+1)}=e^{-\lambda_0s},
\]
and consequently
\(
\boldsymbol\phi(t,x)\ge
\beta e^{-\lambda_0(x\cdot e+R+1)}\mathbf1,~ t\in\mathbb R,~ x\in\Omega_R.
\)
In particular,
\[
\inf_{\substack{t\in\mathbb R\\x\cdot e>-R-1}}
\min_{i=1,2}\left(\phi_i(t,x)e^{\lambda_0x\cdot e}\right)
\ge\beta e^{-\lambda_0(R+1)}>0.
\]
\endgroup
By the definition of $c_*$, choosing $\underline{\lambda}=0$, $\overline{\lambda}=\lambda_0$ and $c = c_* - \varepsilon$ for any given $\varepsilon > 0$, we can verify that these choices satisfy the hypotheses of Lemma~\ref{lemma8}. On the other hand, if we consider the modified nonlinear term $\mathbf{F}_\delta(t,x,u_1,u_2):= \mathbf{F}(t,x,u_1,u_2) - \delta (u_1,u_2)= \big(F_1(t,x,u_1,u_2) - \delta u_1,\; F_2(t,x,u_1,u_2) - \delta u_2\big)$ for some sufficiently small $\delta > 0$, then the problem \eqref{eq:main1} associated with $\mathbf{F}_\delta$ also admits a generalized subsolution $\underline{\mathbf{v}}$ as in Lemma~\ref{lemma8}, with $\underline{\lambda} = 0$, $\overline{\lambda} = \lambda_0$ and $c = c_* - \varepsilon$ for any given $\varepsilon > 0$, since the process only increases the principal eigenvalue $\kappa(\lambda)$ by $\delta$. Therefore, with the help of \eqref{near_zerp}, one can prove that $h\underline{\mathbf{v}}$ is a subsolution of \eqref{eq:main} with $h>0$ small enough. Note that $\underline{\bf v}$ is bounded and decays faster than $e^{-\lambda_0 x\cdot e}$ as  $x\cdot e \rightarrow \infty$. Thus, with smaller $h$, we obtain that $h\underline{\mathbf{v}}(0,x)<{\bf u}(0,x)$ for $x\in \R^N$. By the parabolic comparison principle, one has $h\underline{\mathbf{v}}(t,x)<{\bf u}(t,x),~\forall t\geq 0,~x\in \R^N.$ This implies that, with $c=c_*-\varepsilon$ and $t>0$, $\forall r_1<r_2,~ \inf\limits_{r_1\leq x\cdot e-ct \leq r_2}{\bf u}(t,x)>{\bf 0}.$ With $x=(c_*-\varepsilon) te$, one has
\begin{align*}
{\bf 0}<\inf_{t>0}{\bf u}(t,(c_*-\varepsilon) te)=\inf_{t>0}{\bf u}\left(t,\left[(c_*-\varepsilon) t-\int_0^tc(s)ds\right] e+e\int_0^tc(s)ds \right).
\end{align*}
It follows from the limiting behavior \eqref{limits} that $\limsup\limits_{t\to\infty} \left[ (c_*-\varepsilon)t - \displaystyle\int_{0}^{t} c(s)\,ds \right]<\infty.$ This concludes the proof due to the arbitrariness of $\varepsilon$.
\end{proof}
\begin{remark}\label{remark2}
Consider the case, for $m\in (0,1]$, 
\begin{align*}
F_i(t,x,u_1,u_2)=(g_{i1}(t,x)u_1+g_{i2}(t,x)u_2)(m-u_i),~(t,x)\in \R^{N+1},~i=1,2,
\end{align*}
for some $g_{ij}\in L^{\infty}(\R^{N+1})$ such that $\inf\limits_{\R^{N+1}}g_{ij}>0$, $i,j=1,2$. In this case, one can prove that ${\bf h}_n$ is a generalized subsolution similar to \eqref{sub1} and Theorem \ref{theo:9} still holds for this type of ${\bf F}=(F_1,F_2)$.

\end{remark}

Next, we recall an important result in \cite[Proposition 4.4]{nadin_2015} as follows
\begin{align}\label{pro4.4}
\forall g\in L^{\infty}(\R),~\lfloor g \rfloor=\min_{\widetilde{g}\in \omega_g}\left(\lim_{t\rightarrow \infty}\dfrac{1}{t}\int_0^t\widetilde{g}(s)ds\right),
\end{align}
where $\omega_g$ is the $\omega$-limit set of $g$ in the weak-$\star$ topology of $L^\infty$.

\subsection{Nonexistence of transition fronts}
We now prove the lower bound for existence of transition fronts, which leads to the nonexistence criteria.

\begin{theorem}\label{theo:10}
Assume that \eqref{assumeq}-\eqref{KPP} hold. Then, for any $e\in \mathbb{S}^{N-1}$, there exists $c_*\in \R$ such that, for any generalized transition wave $({\bf u},c)$ of \eqref{eq:main} in the sense of Definition \ref{wavedef}, one has  $\left\lfloor c\right\rfloor\geq c_*$.
\end{theorem}

\begin{proof}
Let ${\bf u}$ be a generalized transition wave with the speed $c$. By \eqref{pro4.4}, there exists $\widetilde{c}\in \omega_c$ such that
\begin{align*}
\lfloor c \rfloor =\lim_{T\rightarrow \infty}\dfrac{1}{T}\int_0^T\widetilde{c}(s)ds.
\end{align*}
By the definition of $\omega_c$, there exists a sequence $\{t_n\}_{n\in \N}$ in $\R$ such that $c(\cdot + t_n) \xrightharpoonup{\star}\widetilde{c} $ in the sense of $L^\infty$ weak-$\star$ topology. For $n\in \N,~i,j=1,2$, consider
\begin{align}\label{trans}
\begin{split}
&{\bf A}^i_n (t,x):={\bf A}^i \left(t+t_n,x+e\int^{t_n}_0c(s)ds\right),~{\bf q}^i_n (t,x):={\bf q}_i \left(t+t_n,x+e\int^{t_n}_0c(s)ds\right),\\
&\mu_{ij}^n(t,x):=\mu_{ij}\left(t+t_n,x+e\int^{t_n}_0c(s)ds\right),~{\bf u}_n (t,x):={\bf u} \left(t+t_n,x+e\int^{t_n}_0c(s)ds\right).
\end{split}
\end{align}
Furthermore, fix
$0<\varepsilon<\min\{1,(2M)^{-1}\}$, where $M$ is in
\eqref{muestiamte}.
Since $F_1,F_2\in C^1([0,\delta_0)^2)$ uniformly with respect
to $(t,x)\in\R^{N+1}$ by \eqref{Fnonlinear}, there exists
$m\in(0,\min\{1,\delta_0\})$ such that
\begin{align*}
F_i(t,x,u_1,u_2)
\geq
\sum_{j=1}^2(\mu_{ij}(t,x)-\varepsilon)u_j,
~
(t,x)\in\R^{N+1},~
(u_1,u_2)\in[0,m]^2,~ i=1,2.
\end{align*}
Since $\mu_{ij}-\varepsilon\geq(2M)^{-1}>0$, it follows that
\begin{align*}
F_i(t,x,u_1,u_2)
\geq
\left(\sum_{j=1}^2(\mu_{ij}(t,x)-\varepsilon)u_j\right)
\left(1-\dfrac{u_i}{m}\right),
~ (u_1,u_2)\in[0,m]^2,~ i=1,2.
\end{align*}

For $n\in\N$, define $z_{n,i}:=\min\left\{\dfrac{u_{n,i}}{m},1\right\},
~ i=1,2,~
{\bf z}_n:=(z_{n,1},z_{n,2}).$ Denote by $F_i^n$ the translated reaction term in the
equation for $u_{n,i}$.
On the set $\{u_{n,i}<m\}$, cooperativity and the preceding
estimate imply that
\begin{align*}
F_i^n(t,x,{\bf u}_n)
&\geq F_i^n(t,x,m{\bf z}_n)\geq
m\sum_{j=1}^2(\mu_{ij}^n(t,x)-\varepsilon)z_{n,j}\geq
m\left(\sum_{j=1}^2
(\mu_{ij}^n(t,x)-\varepsilon)z_{n,j}\right)(1-z_{n,i}).
\end{align*}
On the set $\{u_{n,i}>m\}$, one has $z_{n,i}=1$ and the
last reaction term vanishes.
The minimum construction therefore yields, in the
generalized supersolution sense,
\begin{align*}
\begin{cases}
\partial_t z_{n,1}
-\text{Tr}({\bf A}^1_n(t,x)D^2z_{n,1})
+{\bf q}^1_n(t,x)\cdot Dz_{n,1}
\geq
\big[(\mu_{11}^n(t,x)-\varepsilon)z_{n,1}
+(\mu_{12}^n(t,x)-\varepsilon)z_{n,2}\big](1-z_{n,1}),
\\[1ex]
\partial_t z_{n,2}
-\text{Tr}({\bf A}^2_n(t,x)D^2z_{n,2})
+{\bf q}^2_n(t,x)\cdot Dz_{n,2}
\geq
\big[(\mu_{21}^n(t,x)-\varepsilon)z_{n,1}
+(\mu_{22}^n(t,x)-\varepsilon)z_{n,2}\big](1-z_{n,2}).
\end{cases}
\end{align*}
On the other hand, since ${\bf 0}\leq{\bf u}_n\leq{\bf 1}$
and the reaction terms are uniformly bounded, the $L^p$
parabolic interior estimates ensure that the sequences
$\{\partial_t{\bf u}_n\}_{n\in\N}$,
$\{D{\bf u}_n\}_{n\in\N}$ and
$\{D^2{\bf u}_n\}_{n\in\N}$ are bounded in $L^p(Q)$
for every $p\in(1,\infty)$ and every bounded parabolic cylinder $Q\subset\R^{N+1}$.
Hence, by the parabolic Sobolev embedding and a diagonal
extraction, ${\bf u}_n$ converges locally uniformly in
$\R^{N+1}$ to some $\widetilde{\bf u}$, with weak convergence
of the time derivatives and the first and second spatial
derivatives in $L^p_{\loc}(\R^{N+1})$.
Consequently,
\begin{align*}
{\bf z}_n\rightarrow\widetilde{\bf z}
:=
\left(
\min\left\{\dfrac{\widetilde u_1}{m},1\right\},
\min\left\{\dfrac{\widetilde u_2}{m},1\right\}
\right) \text{ locally uniformly in $\R^{N+1}$.}
\end{align*}
 Similar to Lemma \ref{lm:Haranack-type}, up to subsequences,
${\bf A}^i_n$ and ${\bf q}^i_n$ converge locally uniformly
to $\widetilde{\bf A}^i$ and $\widetilde{\bf q}_i$,
respectively, while $\mu_{ij}^n
\overset{\star}{\rightharpoonup}
\widetilde{\mu}_{ij}
~\text{in }L^\infty(\R^{N+1}),~ i,j=1,2.$ To justify the passage to the limit, we first use the equations for ${\bf u}_n$. The uniformly bounded functions
$F_i^n(t,x,{\bf u}_n)$ converge, up to subsequences,
weakly-$\star$ in $L^\infty(\R^{N+1})$ to some $g_i$.
Thus,
\begin{align*}
\partial_t\widetilde u_i
-\text{Tr}(\widetilde{\bf A}^iD^2\widetilde u_i)
+\widetilde{\bf q}_i\cdot D\widetilde u_i
=g_i
~\text{ in }\R^{N+1}.
\end{align*}
On the open set $\{\widetilde u_i<m\}$, the local uniform
convergence and the preceding cooperativity estimate imply
\begin{align*}
g_i
\geq
m\sum_{j=1}^2
(\widetilde{\mu}_{ij}-\varepsilon)\widetilde z_j
\geq
m\left(\sum_{j=1}^2
(\widetilde{\mu}_{ij}-\varepsilon)\widetilde z_j\right)
(1-\widetilde z_i).
\end{align*}
Applying the same minimum construction, we obtain
\begin{align}\label{limitu}
\begin{cases}
\partial_t\widetilde z_1
-\text{Tr}(\widetilde{\bf A}^1(t,x)D^2\widetilde z_1)
+\widetilde{\bf q}_1(t,x)\cdot D\widetilde z_1
\geq
\big[(\widetilde{\mu}_{11}(t,x)-\varepsilon)\widetilde z_1
+(\widetilde{\mu}_{12}(t,x)-\varepsilon)\widetilde z_2\big]
(1-\widetilde z_1),
\\[1ex]
\partial_t\widetilde z_2
-\text{Tr}(\widetilde{\bf A}^2(t,x)D^2\widetilde z_2)
+\widetilde{\bf q}_2(t,x)\cdot D\widetilde z_2
\geq
\big[(\widetilde{\mu}_{21}(t,x)-\varepsilon)\widetilde z_1
+(\widetilde{\mu}_{22}(t,x)-\varepsilon)\widetilde z_2\big]
(1-\widetilde z_2),
\end{cases}
\end{align}
in the generalized supersolution sense.
Here, the derivative estimates above concern ${\bf u}_n$;
no $W^{1,2}_p$ regularity is asserted for the capped
functions ${\bf z}_n$. Furthermore,
$\widetilde{\bf A}^i\in
BC^{\alpha/2,\alpha}(\R\times\R^N)$
thanks to \eqref{Aassump}.
By \eqref{muestiamte} and weak-$\star$ convergence, $M^{-1}\leq\widetilde{\mu}_{ij}(t,x)\leq M,
~(t,x)\in\R^{N+1},~ i,j=1,2.$ Moreover, the uniform limits defining the transition wave
are preserved by the truncation and the translation limit.
Define
\begin{align*}
H_i^\varepsilon(t,x,{\bf z})
:=
\left(\sum_{j=1}^2
(\widetilde{\mu}_{ij}(t,x)-\varepsilon)z_j\right)(1-z_i),~ i=1,2.
\end{align*}
These reactions are cooperative on $[0,1]^2$, vanish at
${\bf 0}$ and ${\bf 1}$, and satisfy $\partial_{z_j}H_i^\varepsilon(t,x,{\bf 0})
=
\widetilde{\mu}_{ij}(t,x)-\varepsilon.$ Furthermore, for $0<r<1$ and ${\bf z}\in[0,r]^2$,
\begin{align*}
H_i^\varepsilon(t,x,{\bf z})
\geq
\sum_{j=1}^2
(\widetilde{\mu}_{ij}(t,x)-\varepsilon)z_j
-2Mrz_i.
\end{align*}
Consequently, for every $\eta>0$, choosing
$r\leq\eta/(2M)$ implies the lower comparison with the linearized system penalized by $-\eta I$. Together with the positivity of $H_i^\varepsilon(t,x,s,s)$ for $s\in(0,1)$, this allows the small-subsolution argument in Theorem \ref{theo:9} to be applied to the generalized supersolution $\widetilde{\bf z}$. For this purpose, consider the operator
$\widetilde{\bf P}_\lambda
=(\widetilde P_{\lambda,1},\widetilde P_{\lambda,2})$
defined by
\begin{align}\label{operatorP~}
\begin{cases}
\widetilde P_{\lambda,1}{\bf w}
:=
\partial_t w_1
-\text{Tr}\big(\widetilde{\bf A}^1(t,x)D^2w_1\big)
+\big(\widetilde{\bf q}_1(t,x)
+2\lambda\widetilde{\bf A}^1(t,x)e\big)\cdot Dw_1
-\widetilde K_1^\lambda w_1
-(\widetilde{\mu}_{12}(t,x)-\varepsilon)w_2,
\\[1ex]
\widetilde P_{\lambda,2}{\bf w}
:=
\partial_t w_2
-\text{Tr}\big(\widetilde{\bf A}^2(t,x)D^2w_2\big)
+\big(\widetilde{\bf q}_2(t,x)
+2\lambda\widetilde{\bf A}^2(t,x)e\big)\cdot Dw_2
-\widetilde K_2^\lambda w_2
-(\widetilde{\mu}_{21}(t,x)-\varepsilon)w_1,
\end{cases}
\end{align}
for $x\in\R^N$, $t\in\R$, where $\widetilde K_i^\lambda
:=
\widetilde{\mu}_{ii}(t,x)-\varepsilon
+\lambda^2 e\cdot\widetilde{\bf A}^i(t,x)e
+\lambda\widetilde{\bf q}_i(t,x)\cdot e,~i=1,2.$ Denote by $\widetilde{\kappa}_\varepsilon(\lambda)$
its generalized principal eigenvalue, defined as in
\eqref{GPE}. We claim that, for every $\lambda_0>0$, there exist constants
$C_{\lambda_0}>0$ and $\varepsilon_0>0$ such that
\begin{align}\label{estimatekappa}
\widetilde{\kappa}_\varepsilon(\lambda)
\leq \kappa(\lambda)+C_{\lambda_0}\varepsilon,~
\lambda\in[0,{\lambda_0}],~
0<\varepsilon<\varepsilon_0.
\end{align}
Let us postpone the proof of this claim for the moment and derive the final conclusion. Since $\kappa(0)<0$, \eqref{estimatekappa} implies that
$\widetilde{\kappa}_\varepsilon(0)<0$ for $\varepsilon>0$
small enough. Thanks to Remark \ref{remark2}, Theorem
\ref{theo:9} yields
\begin{align*}
\lfloor c\rfloor
=
\liminf_{T\rightarrow\infty}
\dfrac{1}{T}\int_0^T\widetilde{c}(s)\,ds
\geq
\inf_{\lambda>0}
\dfrac{-\widetilde{\kappa}_\varepsilon(\lambda)}{\lambda}.
\end{align*}

By uniform ellipticity and boundedness of the coefficients,
testing the definition of
$\widetilde{\kappa}_\varepsilon(\lambda)$ with ${\bf 1}$
implies constants $a>0$ and $B,C>0$, independent of
$\varepsilon>0$ small enough, such that $\widetilde{\kappa}_\varepsilon(\lambda)
\leq -a\lambda^2+B\lambda+C,~ \lambda\geq0.$ Moreover, since $\kappa(0)<0$ and $\kappa$ is continuous
at $0$, \eqref{estimatekappa} implies the existence of constants $\eta,\delta^0>0$ such that
$\widetilde{\kappa}_\varepsilon(\lambda)\leq-\eta$
for $\lambda\in[0,\delta^0]$ and $\varepsilon>0$ small enough.
Consequently, there exist $0<\lambda_-<\lambda_+$,
independent of $\varepsilon>0$ small enough, such that
\begin{align*}
\dfrac{-\widetilde{\kappa}_\varepsilon(\lambda)}{\lambda}
\geq c_*+1,~
\lambda\in(0,\lambda_-)\cup(\lambda_+,\infty).
\end{align*}
Applying \eqref{estimatekappa} with $\lambda_0=\lambda_+$
yields
\begin{align*}
\dfrac{-\widetilde{\kappa}_\varepsilon(\lambda)}{\lambda}
&\geq
\dfrac{-\kappa(\lambda)}{\lambda}
-\dfrac{C_{\lambda_+}\varepsilon}{\lambda}\geq
c_*-\dfrac{C_{\lambda_+}}{\lambda_-}\varepsilon,~\lambda\in[\lambda_-,\lambda_+],
\end{align*}
where $c_*=\inf\limits_{\lambda>0}(-\kappa(\lambda)/\lambda)$.
It follows that
\begin{align*}
\lfloor c\rfloor
\geq c_*-\dfrac{C_{\lambda_+}}{\lambda_-}\varepsilon.
\end{align*}
Letting $\varepsilon\rightarrow0^+$, we obtain $\lfloor c\rfloor\geq c_*$. It remains to establish \eqref{estimatekappa}.\par\noindent\textbf{Proof of \eqref{estimatekappa}.} 
Fix ${\lambda_0}>0$ and $\lambda\in[0,{\lambda_0}]$.
Let $\widetilde{\bf P}_{\lambda,0}$ denote the operator
with the limiting coefficients
$\widetilde{\bf A}^i$, $\widetilde{\bf q}_i$ and
$\widetilde{\mu}_{ij}$, without the $\varepsilon$-perturbation.
Then
\begin{align*}
\bigl(\widetilde{\bf P}_{\lambda,\varepsilon}{\bf z}\bigr)_i
=
\bigl(\widetilde{\bf P}_{\lambda,0}{\bf z}\bigr)_i
+\varepsilon(z_1+z_2),~i=1,2.
\end{align*}
Choose $\varepsilon_0>0$ sufficiently small that
$\widetilde{\mu}_{ij}-\varepsilon>0$ for $i\neq j$
and $0<\varepsilon<\varepsilon_0$.

Let $\kappa(\lambda)<k<\kappa(\lambda)+1$.
By the definition of $\kappa(\lambda)$, there exists ${\bf w}\in
W^{1,2,N+1}_{\loc}(\R\times\R^N)
\cap W^{0,1,\infty}(\R\times\R^N)$ such that $\inf\limits_{\R\times\R^N}{\bf w}>{\bf 0},~{\bf P}_\lambda{\bf w}\leq k{\bf w}.$ By multiplying ${\bf w}$ by a positive constant, we may assume
that ${\bf w}\leq{\bf 1}$. We first replace ${\bf w}$ by a spatially periodic solution
of a semilinear system.
Choose $\chi\in C^\infty(\R,[0,1])$ such that
$\chi(s)=0$ for $s\leq1$ and $\chi(s)=1$ for $s\geq2$,
and set $g_k(s):=\bigl(k-B_{\lambda_0}\chi(s)\bigr)s,~
{\bf g}_k({\bf z}):=(g_k(z_1),g_k(z_2)).$ The spectral bounds for $\kappa(\lambda)$ and the coefficient
bounds allow $B_{\lambda_0}>0$ to be chosen independently of
$\lambda\in[0,{\lambda_0}]$ and
$k\in(\kappa(\lambda),\kappa(\lambda)+1)$ so that $B_{\lambda_0}
\geq
k+K_i^\lambda(t,x)+\mu_{ij}(t,x),~ i\neq j.$ Consequently, ${\bf w}$ and ${\bf 2}$ are respectively
a subsolution and a supersolution of ${\bf P}_\lambda{\bf v}={\bf g}_k({\bf v}).$ Indeed, ${\bf g}_k({\bf w})=k{\bf w}$, whereas
\begin{align*}
\bigl({\bf P}_\lambda({\bf 2})\bigr)_i
=
-2\bigl(K_i^\lambda+\mu_{ij}\bigr)
\geq 2(k-B_{\lambda_0})=g_k(2),~ i\neq j.
\end{align*}

For each $n\in\N$, solve this semilinear system for $t>-n$
with initial datum ${\bf 2}$. By comparison, the resulting solution lies between
${\bf w}$ and ${\bf 2}$.
Since the coefficients and the initial datum are spatially
periodic, this solution is also spatially periodic.
Letting $n\to\infty$ and using local parabolic compactness,
we obtain a spatially periodic entire solution ${\bf v}$ satisfying
\begin{align*}
{\bf P}_\lambda{\bf v}={\bf g}_k({\bf v}),~
{\bf 0}<{\bf w}\leq{\bf v}\leq{\bf 2}.
\end{align*}

For $i\neq j$, the equation for $v_i$ can be written as
\begin{align*}
\partial_t v_i
-\operatorname{Tr}({\bf A}^iD^2v_i)
+\bigl({\bf q}_i+2\lambda{\bf A}^ie\bigr)\cdot Dv_i
=
\bigl(K_i^\lambda+k-B_{\lambda_0}\chi(v_i)\bigr)v_i
+\mu_{ij}v_j.
\end{align*}
For each fixed ${\bf v}$, define
\(
a_i(t,x):=K_i^\lambda(t,x)+k-B_{\lambda_0}\chi(v_i(t,x)),~ i=1,2.
\)
Then, ${\bf v}$ satisfies a linear cooperative system.
The coefficients $a_i$ are uniformly bounded, independently of $k$ and
${\bf w}$ and the off-diagonal coefficients $\mu_{ij}$ are uniformly
bounded below by a positive constant. Repeating the Harnack and
maximum-principle argument of Lemma \ref{lm:Haranack-type} for this
linear system, we obtain a constant $H_{\lambda_0}\geq1$, independent
of $k$ and ${\bf w}$, such that
\begin{align*}
\max_{i=1,2}\max_{x\in[0,{\bf L}]}v_i(t,x)
\leq
H_{\lambda_0}
\min_{i=1,2}\min_{x\in[0,{\bf L}]}v_i(t,x),~ t\in\R.
\end{align*}
In particular, $v_j(t,x)\leq H_{\lambda_0} v_i(t,x),~ i,j=1,2~(t,x)\in\R^{N+1}.$ We next establish a global gradient bound.
Fix $p>N+2$ and write $\partial_t v_i
-\operatorname{Tr}({\bf A}^iD^2v_i)
+\bigl({\bf q}_i+2\lambda{\bf A}^ie\bigr)\cdot Dv_i
=h_i,$ where $h_i:=K_i^\lambda v_i+\mu_{ij}v_j+g_k(v_i)$
for $i\neq j$.
The construction implies
$\|v_i\|_\infty\leq2$ and $\|h_i\|_\infty\leq M_{\lambda_0}$,
uniformly in $k$ and ${\bf w}$. For $z=(t_0,x_0)$, set
$Q_r(z):=(t_0-r^2,t_0+r^2)\times B_r(x_0)$.
By uniform ellipticity and the global H\"older bounds
on the diffusion coefficients, the interior parabolic
estimates yield
\begin{align*}
\|v_i\|_{W^{1,2}_p(Q_1(z))}
&\leq C_{{\lambda_0},p}
\bigl(\|v_i\|_{L^p(Q_2(z))}
+\|h_i\|_{L^p(Q_2(z))}\bigr)\leq C_{{\lambda_0},p}|Q_2|^{1/p}(2+M_{\lambda_0}),
\end{align*}
with constants independent of $z$.
Since $p>N+2$, the parabolic Sobolev embedding implies
\begin{align*}
\|Dv_i\|_{L^\infty(Q_{1/2}(z))}
\leq C_p\|v_i\|_{W^{1,2}_p(Q_1(z))}
\leq C_{\lambda_0}'.
\end{align*}
Taking the supremum over $z\in\R^{N+1}$ yields
$\|D{\bf v}\|_{L^\infty(\R^{N+1})}\leq C_{\lambda_0}'$. Apply to ${\bf v}$ the same space-time translations as in
\eqref{trans}.
The local parabolic estimates yield, up to a subsequence,
a locally uniform limit $\widetilde{\bf v}$, with weak
convergence of the time derivatives and second spatial
derivatives in $L^p_{\loc}$.
Passing to the limit, we obtain $\widetilde{\bf P}_{\lambda,0}\widetilde{\bf v}
=
{\bf g}_k(\widetilde{\bf v}).$ Moreover, $\inf\limits_{\R^{N+1}}\widetilde{\bf v}>{\bf 0},~
\widetilde{\bf v}\leq{\bf 2},~\|D\widetilde{\bf v}\|_{L^\infty(\R^{N+1})}\leq C_{\lambda_0}'$ and
$\widetilde v_j\leq H_{\lambda_0}\widetilde v_i$
for $i,j=1,2$.
Here, the positive lower bound is inherited from $\inf\limits_{\R^{N+1}}{\bf w}>{\bf 0}$. Since $g_k(s)\leq ks$ for $s\geq0$, it follows that
\begin{align}\label{semilinear}
\begin{split}
\bigl(\widetilde{\bf P}_{\lambda,\varepsilon}
\widetilde{\bf v}\bigr)_i
&=
g_k(\widetilde v_i)
+\varepsilon(\widetilde v_1+\widetilde v_2)\leq
\bigl(k+(1+H_{\lambda_0})\varepsilon\bigr)\widetilde v_i,~ i=1,2.
\end{split}
\end{align}
Thus, $\widetilde{\bf v}$ is admissible in the definition of
$\widetilde\kappa_\varepsilon(\lambda)$ and
\begin{align*}
\widetilde\kappa_\varepsilon(\lambda)
\leq k+(1+H_{\lambda_0})\varepsilon.
\end{align*}
Since $H_{\lambda_0}$ is independent of $k$, letting
$k\downarrow\kappa(\lambda)$ proves \eqref{estimatekappa}
with $C_{\lambda_0}=1+H_{\lambda_0}$.
\end{proof}

\begin{remark}\label{remark51}
In the $T_0$-periodic case, the results of Berestycki and Nadin~\cite{BerestyckiNadinMemoirs2022} imply that $k(\lambda)=-\kappa(\lambda)$ and
\begin{align*}
c^*=\min_{\lambda>0}\dfrac{k(\lambda)}{\lambda}=\min_{\lambda>0}\dfrac{-\kappa(\lambda)}{\lambda}=-\max _{\lambda>0}\dfrac{\kappa(\lambda)}{\lambda}=c_*,
\end{align*}  
where $k(\lambda)$ is the principal eigenvalue in Remark \ref{periodicity}.

\end{remark}

\section{Minimal Speed $c^*(e)=c_*(e)$ in Almost-Periodic Media}\label{section 5}

In this section, we consider a framework with almost periodicity in time and periodicity in space. 
Fixing $e\in\mathbb{S}^{N-1}$, we prove that $c^*(e)=c_*(e)$, where $c^*(e)=c^*$ and $c_*(e)=c_*$ are defined in Sections~\ref{section:ex2} and \ref{sec4}, respectively. This equality provides a sharp threshold for the existence of generalized transition waves in the sense of Definition~\ref{wavedef}. We recall the  generalized principal eigenvalue.
\begin{align}\label{GPE1}
&\underline{k}(\lambda):=\inf\left\{k\in \R:\exists {\bf w} \in \big[W^{1,2,N+1}_{loc}\big(\R\times \R^N\big) \big]^2,~\inf\limits_{\R\times \R^N} {\bf w}>{\bf 0},~\|{\bf w}\|_\infty<\infty,~{\bf P}_{\lambda}{\bf w}\leq k{\bf w} \text{ in }\R\times \R^N\right\}.
\end{align}

We begin the section with the following result

\begin{lemma}
Assume that \eqref{assumeq} to \eqref{KPP} and \eqref{Fregularity} hold. Consider $\kappa(\lambda)$ in Section~\ref{sec4}. Then, $\kappa(\lambda)= \underline{k}(\lambda)$.

\end{lemma}

\begin{proof}

Under conditions \eqref{assumeq} to \eqref{KPP} and \eqref{Fregularity}, $\kappa(\lambda)\geq \underline{k}(\lambda)$. This is due to the fact that no global boundedness assumption on $D\mathbf{w}$ (with respect to $x$) is imposed in the definition of $\underline{k}(\lambda)$. We prove the reverse inequality without differentiating an arbitrary admissible subsolution.
Fix $\lambda\in\R$ and choose $k>\underline{k}(\lambda)$ so that
there exists $\mathbf w\in[W^{1,2,N+1}_{\mathrm{loc}}(\R^{N+1})]^2$
with $\inf_{\R^{N+1}}\mathbf w>\mathbf0$, $\|\mathbf w\|_\infty<\infty$, and
$\mathbf P_\lambda\mathbf w\le k\mathbf w$.
By scaling, we may assume $\mathbf0<\mathbf w\le\mathbf1$.
Choose $B>0$ satisfying
\[
B\ge k+\max_{i=1,2}\sup_{(t,x)\in\R^{N+1}}
\left(K_i^\lambda(t,x)+\sum_{j\ne i}\mu_{ij}(t,x)\right).
\]
Let $\chi\in C^\infty(\R,[0,1])$ satisfy $\chi=0$ on $(-\infty,1]$
and $\chi=1$ on $[2,\infty)$, and set
$g_k(s)=(k-B\chi(s))s$ for $s\ge0$.
Since $g_k(w_i)=kw_i$ and $g_k(2)=2(k-B)$, we have
\[
\mathbf P_\lambda\mathbf w\le\mathbf g_k(\mathbf w),\qquad
\mathbf P_\lambda(2\mathbf1)\ge\mathbf g_k(2\mathbf1),
\quad \mathbf g_k(\mathbf z)=(g_k(z_1),g_k(z_2)).
\]
For each integer $n>0$, let $\mathbf v_n$ solve the cooperative
semilinear problem $\mathbf P_\lambda\mathbf v_n=\mathbf g_k(\mathbf v_n)$
on $(-n,\infty)\times\R^N$, with $\mathbf v_n(-n,\cdot)=2\mathbf1$.
Comparison yields $\mathbf w\le\mathbf v_n\le2\mathbf1$ for $t\ge-n$.
Since the coefficients and the initial datum are $\mathbf L$-periodic
in space, $\mathbf v_n$ is spatially periodic. Local parabolic
compactness on each time slab and a diagonal subsequence yield an entire
spatially periodic solution $\mathbf v$ with
\[
\mathbf0<\inf_{\R^{N+1}}\mathbf w\le\mathbf v\le2\mathbf1,
\qquad \mathbf P_\lambda\mathbf v=\mathbf g_k(\mathbf v)\le k\mathbf v.
\]
It remains to verify the global spatial-gradient condition in the
definition of $\kappa(\lambda)$. Write the $i$-th equation as
\[
\partial_tv_i-\operatorname{Tr}(\mathbf A^iD^2v_i)
+(\mathbf q_i+2\lambda\mathbf A^ie)\cdot Dv_i=h_i,
\quad h_i=K_i^\lambda v_i+\sum_{j\ne i}\mu_{ij}v_j+g_k(v_i).
\]
The uniform boundedness of $\mathbf v$, the coefficients, and $g_k$ on
$[0,2]$ implies $\|h_i\|_{L^\infty(\R^{N+1})}<\infty$.
Fix $p>N+2$. Uniform interior parabolic $W^{1,2}_p$ estimates on cylinders
$Q_1(z)\Subset Q_2(z)$, valid because $\mathbf A^i$ are uniformly elliptic
and globally H\"older continuous, show that
\[
\sup_{z\in\R^{N+1}}\|v_i\|_{W^{1,2}_p(Q_1(z))}<\infty.
\]
The parabolic Sobolev embedding, applied on $Q_{1/2}(z)$, now implies
$\|Dv_i\|_{L^\infty(\R^{N+1})}<\infty$. Consequently $\mathbf v$ is an
admissible test function for $\kappa(\lambda)$, so
$\kappa(\lambda)\le k$. Letting $k\downarrow\underline{k}(\lambda)$
and recalling the opposite inequality proves
$\kappa(\lambda)=\underline{k}(\lambda)$.
\end{proof}

Recall that $\boldsymbol{\eta}_{\lambda}$ is the solution of \eqref{div}. We first revisit the results in \cite[Proposition 3.12]{BerestyckiNadinMemoirs2022}.

\begin{lemma}\label{lem:least-mean-growth}
Assume \eqref{assumeq} to \eqref{KPP} and \eqref{Fregularity}. Then, 
\begin{align*}
\underline{k}(\lambda)=-\lim_{T \to +\infty} \inf_{t \in \mathbb{R}} 
\frac{1}{T} \ln\dfrac{ \|\boldsymbol{\eta}_{\lambda}(t+T,\cdot)\|_{\infty}}{ \|\boldsymbol{\eta}_{\lambda}(t,\cdot)\|_{\infty}}=-\lambda \lfloor c_{\lambda} \rfloor,
\end{align*}
where $\lfloor c_{\lambda} \rfloor$ is in \eqref{clambdaest}.
\end{lemma}

\begin{proof}
We show that $\underline{k}(\lambda)=-\lambda \lfloor c_{\lambda} \rfloor=:-b$. Recall $\boldsymbol{\varphi}_{\lambda}(t,x)
=e^{-\lambda S_{\lambda}(t)}\boldsymbol{\eta}_{\lambda}(t,x).$ The  estimate \eqref{estimatevarphi} implies that $\boldsymbol{\varphi}_{\lambda}$ is bounded above and bounded away from zero. Moreover, ${\bf P}_{\lambda}\boldsymbol{\varphi}_{\lambda}=-\lambda c_{\lambda}(t)\boldsymbol{\varphi}_{\lambda}\text{ in }\mathbb{R}^{N+1}.$ For every $\varepsilon>0$, \eqref{keylemma} provides
$\sigma_{\varepsilon}\in W^{1,\infty}(\mathbb{R})$ such that $\sigma_{\varepsilon}'(t)+\lambda c_{\lambda}(t)
\geq b-\varepsilon
 \text{ for }t\in\mathbb{R}.$ Thus, $\boldsymbol{w}_{\varepsilon}
:=e^{-\sigma_{\varepsilon}(t)}\boldsymbol{\varphi}_{\lambda}$
is admissible in the definition of $\underline{k}(\lambda)$ and satisfies
\begin{align*}
{\bf P}_{\lambda}\boldsymbol{w}_{\varepsilon}
=-\big(\sigma_{\varepsilon}'+\lambda c_{\lambda}\big)
\boldsymbol{w}_{\varepsilon}
\leq(-b+\varepsilon)\boldsymbol{w}_{\varepsilon}.
\end{align*}
It follows that $\underline{k}(\lambda)\leq-b+\varepsilon$.
Letting $\varepsilon\to0^+$ yields $\underline{k}(\lambda)\leq-b$. Conversely, let $k\in\mathbb{R}$ and $\boldsymbol{w}$ be any admissible pair satisfying
\begin{align*}
{\bf P}_{\lambda}\boldsymbol{w}\leq k\boldsymbol{w},~
m_w:=\inf_{(t,x)\in\mathbb{R}^{N+1}}\min_{i=1,2}w_i(t,x)>0,~
M_w:=\|\boldsymbol{w}\|_{\infty}<\infty.
\end{align*}
Set $M(t):=\|\boldsymbol{\eta}_{\lambda}(t,\cdot)\|_{\infty}$. \eqref{Harnack-zero}  implies that
$\eta_{\lambda,i}(t,x)\geq C_0M(t)$ for $i=1,2$. Hence, for each fixed $t\in\mathbb{R}$,
\begin{align*}
\boldsymbol{w}(t,x)
\leq\frac{M_w}{C_0M(t)}\boldsymbol{\eta}_{\lambda}(t,x).
\end{align*}
Since ${\bf P}_{\lambda}\boldsymbol{\eta}_{\lambda}=0$ and ${\bf P}_{\lambda}\big(e^{-k(s-t)}\boldsymbol{w}(s,x)\big)
=e^{-k(s-t)}
\big({\bf P}_{\lambda}\boldsymbol{w}-k\boldsymbol{w}\big)(s,x)
\leq0,$ the comparison principle implies that
\begin{align*}
e^{-kT}\boldsymbol{w}(t+T,x)
\leq\frac{M_w}{C_0M(t)}
\boldsymbol{\eta}_{\lambda}(t+T,x),~\forall T>0.
\end{align*}
Consequently, $\dfrac{M(t+T)}{M(t)}
\geq\dfrac{C_0m_w}{M_w}e^{-kT}$ and therefore $\inf\limits_{t\in\mathbb{R}}\dfrac1T\ln\dfrac{M(t+T)}{M(t)}
\geq-k+\dfrac1T\ln\dfrac{C_0m_w}{M_w}.$ Letting $T\to+\infty$, we obtain $b\geq-k$. Taking the infimum over all admissible $k$ yields
$\underline{k}(\lambda)\geq-b$. Thus, $\underline{k}(\lambda)=-b=-\lambda \lfloor c_{\lambda} \rfloor$. The second equality follows from \eqref{clambdaest}. This completes the proof.
\end{proof}

It is known that Lemmas \ref{lemma continuous} and \ref{monotonec} imply that $\lfloor c_{\lambda} \rfloor\geq \lfloor c_{\lambda^*} \rfloor=c^*(e),~\forall \lambda \in (0,\lambda^*]$. 
It then follows that
\begin{align*}
c^*(e)=-\sup_{\lambda\in (0,\lambda^*	)}\dfrac{\underline{k}(\lambda)}{\lambda}\geq -\max_{\lambda>0}\dfrac{\kappa(\lambda)}{\lambda} =c_*(e).
\end{align*}

Now, in the almost-periodic (a.p.) setting, we shall prove that $\lambda \mapsto -\underline{k}(\lambda)$ is convex, which yields that $\lambda \mapsto -\underline{k}(\lambda)/\lambda$ is decreasing on $(0,\lambda^*)$ and either increasing on $(\lambda^*,\infty)$ or constant on $[\lambda^*,\lambda_3]$ for some $\lambda_3$ and then is increasing on $(\lambda_3,\infty)$, showing that the minimum is attained at $\lambda^*$ and hence that $c^*(e)=c_*(e)$. For this purpose, define 
\begin{align*}
a(t,T):=\ln\dfrac{ \|\boldsymbol{\eta}_{\lambda}(t+T,\cdot)\|_{\infty}}{ \|\boldsymbol{\eta}_{\lambda}(t,\cdot)\|_{\infty}}.
\end{align*}
It is known that $(t,T)\mapsto a(t,T)$ is continuous due to the parabolic estimates for $\boldsymbol{\eta}_\lambda $. 

Next, for $\lambda>0$, recall $\mathcal A_\lambda(t)$ denote the spatial operator
in the right-hand side of \eqref{div}; equivalently,
$\mathbf P_\lambda=\partial_t-\mathcal A_\lambda(t)$, $\mathcal C=[\mathbf0,\mathbf L]$ and set
\begin{align*}
\mathcal J(\mathbf z):=\frac1{|\mathcal C|}
\sum_{i=1}^2\int_{\mathcal C}z_i(x)\,dx.
\end{align*}
The normalized principal section is denoted by $\mathbf p_\lambda$;
its instantaneous growth rate is denoted by $b_\lambda$.

Before proving the key result, we make the following observation.

\begin{remark}\label{mu_ap1}
By \eqref{Fnonlinear} and \eqref{APassume}, the functions
$\mu_{ij}(t,x)=\partial_{u_j}F_i(t,x,0,0)$,
$i,j=1,2$, are almost periodic in $t$ uniformly with
respect to $x\in\mathbb R^N$. Indeed, this follows
from the uniform convergence of the corresponding
difference quotients and the uniform-limit property
of almost periodic functions.
\end{remark}

\begin{lemma}\label{lem:AP-normalized-bundle}
Assume \eqref{assumeq}, \eqref{Aassump}, \eqref{Fnonlinear}
and \eqref{Fregularity}. Suppose further that the almost periodicity \eqref{APassume}.
No KPP inequality or near-zero nonlinear
estimate is needed for the following assertions about the
linearized evolution.
For each $\lambda>0$, there exists a unique positive spatially
periodic entire solution $\mathbf p_\lambda$ of
\begin{align}\label{eq:v2-normalized}
\partial_t\mathbf p_\lambda
=\mathcal A_\lambda(t)\mathbf p_\lambda
-b_\lambda(t)\mathbf p_\lambda,~ \mathcal J(\mathbf p_\lambda(t,\cdot))=1,
\end{align}
where
\begin{align}\label{eq:v2-growth}
b_\lambda(t)=\mathcal J\big(
\mathcal A_\lambda(t)\mathbf p_\lambda(t,\cdot)\big).
\end{align}
Both $\mathbf p_\lambda$ and $b_\lambda$ are almost periodic
in $t$, uniformly in $x$ for $\mathbf p_\lambda$. The mappings $\lambda\mapsto\mathbf p_\lambda$ and
$\lambda\mapsto b_\lambda$ are locally $C^1$ in the corresponding uniform norms. For each compact interval $I\Subset(0,\infty)$, there are
$0<m_I<M_I<\infty$ such that
\begin{align}\label{eq:v2-harnack}
m_I\le p_{\lambda,i}(t,x)\le M_I,~ |\partial_\lambda p_{\lambda,i}(t,x)|
+|\partial_\lambda b_\lambda(t)|\le M_I
\end{align}
for $\lambda\in I$, $t\in\R$, $x\in\R^N$ and $i=1,2$.
%
\end{lemma}

\begin{proof}
Write $U_\lambda(t,s)$ for the strongly positive evolution
operator generated by $\mathcal A_\lambda(t)$ on the periodic cell.
We first prove a uniform two-sided bound for the matrix kernel of
$U_\lambda$. Fix $I\Subset(0,\infty)$ and let
$S_{\lambda,i}(t,s)$ be the scalar evolution obtained by retaining
only the $i$-th diagonal equation, including its zero-order term
$\mu_{ii}$. Uniform ellipticity, the uniform H\"older bounds of the
principal coefficients, and boundedness of the lower-order terms
imply the following fixed-positive-time scalar kernel bounds: for
any $0<\tau_0<\tau_1<\infty$ there are constants
$c_I,C_I>0$ such that, on the periodic cell,
\[
c_I\le G_{\lambda,i}(t,s;x,y)\le C_I
\quad\text{if }\tau_0\le t-s\le\tau_1,
\quad i=1,2,\ \lambda\in I.
\]
All constants are independent of $s,t$ and of the time translate of
the coefficients. Indeed, the scalar parabolic Gaussian upper and
lower estimates on the covering space, summed over the periodic
lattice, give these bounds on the compact torus; the scalar
zero-order potentials are uniformly bounded and are incorporated by
comparison with their essential bounds.

Set $\mu_0=\min_{i\ne j}\inf_{t,x}\mu_{ij}>0$ and fix
$\tau>0$. Positivity and variation of constants give, for $i\ne j$
and nonnegative $f$,
\[
(U_\lambda(t+3\tau,t))_{ij}f
\ge\mu_0\int_{t+\tau}^{t+2\tau}
S_{\lambda,i}(t+3\tau,r)
S_{\lambda,j}(r,t)f\,dr.
\]
Both scalar factors in this integral have durations between $\tau$
and $2\tau$. Their positive kernels are bounded below uniformly,
so each off-diagonal entry of the full matrix kernel is bounded
below by a strictly positive constant. The diagonal entries
dominate $S_{\lambda,i}(t+3\tau,t)$ and enjoy the same property.
For the upper bound, let $M$ dominate the absolute values of the
zero-order coefficients. The scalar upper bounds and comparison with
the positive constant coupling matrix having every entry $M$ give a
uniform upper bound for all four kernels on this fixed time block.
Integration against $\mathbf z\ge\mathbf0$ therefore yields, with
$T=3\tau$ and a rescaling by the measure of the periodic cell,
constants $0<a_I<A_I$ satisfying
\begin{align}\label{eq:v2-doeblin}
a_I\mathcal J(\mathbf z)\mathbf1
\le U_\lambda(t+T,t)\mathbf z
\le A_I\mathcal J(\mathbf z)\mathbf1
\end{align}
for every $\mathbf z\ge\mathbf0$, $\mathbf z\not\equiv\mathbf0$,
$t\in\R$ and $\lambda\in I$. Here, all inequalities are componentwise.
In particular, the normalized time-$T$ maps have uniformly finite
projective diameter. Birkhoff's contraction inequality implies a
uniform contraction factor $\vartheta_I<1$ in the Hilbert projective
metric. Iterating these maps from time $t-nT$ and normalizing by
$\mathcal J$ determines a unique projective entire section.
The inequalities \eqref{eq:v2-doeblin} imply \eqref{eq:v2-harnack}
without its derivative terms. Normalizing at each time yields
\eqref{eq:v2-normalized}--\eqref{eq:v2-growth}.

We next prove the differentiability of the normalized positive solutions with respect to $\lambda$. Fix $I\Subset(0,\infty)$ and let
$X=C_{\rm per}(\mathcal C)^2$. For $T>0$ as in
\eqref{eq:v2-doeblin}, put
\begin{align*}
\mathcal N_{\lambda,t}(z)
=\frac{U_\lambda(t+T,t)z}
{\mathcal J(U_\lambda(t+T,t)z)},~
\mathcal J(z)=1.
\end{align*}
After one time block, every normalized positive vector lies in the
order interval $a_I/A_I\le z_i\le A_I/a_I$ (with a harmless
rescaling of the constants). 
On this interval, set $X_0:=\{h\in X:\mathcal J(h)=0\}$.
For a normalized function $z$ and $h\in X_0$ the infinitesimal
Hilbert seminorm is
\begin{align*}
\|h\|_{H,z}=\max_{i,x}\frac{h_i(x)}{z_i(x)}
-\min_{i,x}\frac{h_i(x)}{z_i(x)}.
\end{align*}
Writing $m_I^{\mathrm{p}}=a_I/A_I$ and $M_I^{\mathrm{p}}=A_I/a_I$, the normalization
$\mathcal J(h)=0$ implies that $h_i/z_i$ takes values of both
signs unless $h=0$. Consequently,
\begin{align*}
\|h\|_\infty\le M_I^{\mathrm{p}}\|h\|_{H,z},\qquad
\|h\|_{H,z}\le\frac{2}{m_I^{\mathrm{p}}}\|h\|_\infty.
\end{align*}
The Hilbert contraction factor for one time block is bounded by
$\vartheta_I<1$, uniformly in the coefficient hull and
$\lambda\in I$. Thus, for $m$ consecutive blocks,
\begin{align*}
\big\|D_z\mathcal N^{(m)}_{\lambda,t}h\big\|_\infty
\le\frac{2M_I^{\mathrm{p}}}{m_I^{\mathrm{p}}}\vartheta_I^m\|h\|_\infty,
\qquad h\in X_0.
\end{align*}
We choose $m$ so large that $\frac{2M_I^{\mathrm{p}}}{m_I^{\mathrm{p}}}\vartheta_I^m<1$.
The projective contraction estimate therefore implies

\begin{align}\label{eq:v10-contraction}
\sup_{t\in\R,\lambda\in I}
\big\|D_z\mathcal N^{(m)}_{\lambda,t}\big\|_{X_0\to X_0}
\le\theta<1,
\end{align}
where $\mathcal N^{(m)}_{\lambda,t}$ denotes the normalized
evolution over $mT$; the derivative is evaluated along normalized
positive states in that interval. The bound follows also by
differentiating Birkhoff's metric inequality at two nearby
normalized states.

The mapping $\lambda\mapsto U_\lambda(t+mT,t)$ are $C^1$ in
operator norm on $X$, uniformly in $t$ and locally uniformly in
$\lambda\in I$. 
To see the parameter bound, write
$\mathcal A_\lambda(t)=\mathcal A_0(t)+\lambda B_1(t)+\lambda^2B_2(t)$,
where $B_1$ has order at most one and $B_2$ is a multiplication
operator. The $B_j$ have uniformly bounded coefficients on every
compact parameter interval. For any fixed time $T>0$,
variation of constants yields, initially on smooth functions,
\begin{align*}
\partial_\lambda U_\lambda(t+T,t)
=\int_t^{t+T}U_\lambda(t+T,r)
\big(B_1(r)+2\lambda B_2(r)\big)U_\lambda(r,t)\,dr.
\end{align*}
The short-time parabolic gradient estimate
$\|DU_\lambda(r,t)z\|_\infty
\le C_I(1+(r-t)^{-1/2})\|z\|_\infty$
for $0<r-t\le T$, together with the evolution bounds on a
finite time block, makes the integral uniformly convergent in
operator norm. The same estimate applied to parameter differences
proves continuity of the derivative, uniformly in the initial
time and on compact $\lambda$-intervals. This also justifies
Duhamel's identity on $X$ by density.
To specify the regularity used below, fix $0<\beta<\alpha$ and
$0<\tau<T/4$. The parabolic interior estimates on the periodic cell,
applied after a positive time interval, imply
\begin{align}\label{eq:v15-uniform-smoothing}
\sup_{t\in\mathbb R,\lambda\in I}
\left(\|U_\lambda(t+T,t)\|_{X\to C^{2+\beta}_{\rm per}}
+\|\partial_\lambda U_\lambda(t+T,t)\|_{X\to C^{2+\beta}_{\rm per}}\right)
\le C_{I,T,\beta}.
\end{align}
We justify the estimate on the left of
\eqref{eq:v15-uniform-smoothing}, fix $0<\tau<T/4$ and write
$z(s)=U_\lambda(s,t)z_0$ and
$y(s)=(\partial_\lambda U_\lambda(s,t))z_0$.
On $[t,t+2\tau]$, the Duhamel identity and the short-time gradient
estimate imply
\begin{align*}
\|y(t+2\tau)\|_\infty
&\le C_{I,T}\|z_0\|_\infty
\int_0^{2\tau}(1+s^{-1/2})\,ds
\le C_{I,T,\tau}\|z_0\|_\infty.
\end{align*}
For $s>t$, $y$ satisfies
\begin{align*}
\partial_s y=\mathcal A_\lambda(s)y+
(B_1(s)+2\lambda B_2(s))z,~ y(t)=0.
\end{align*}
The positive-time regularization of $z$ on
$[t+\tau,t+T]$ controls its parabolic
$C^{1+\beta/2,2+\beta}$ norm on
$[t+2\tau,t+T]$, after reducing $\beta\in(0,\alpha)$ if
necessary. In particular, the inhomogeneous term is uniformly
$C^{\beta/2,\beta}$ on $[t+2\tau,t+T]$.
The interior parabolic Schauder estimate for the equation of $y$
on $[t+2\tau,t+T]$, with the $C^0$ bound already established,
now controls $\|y(t+T)\|_{C^{2+\beta}}$.
The same procedure applied to $z$ proves
\eqref{eq:v15-uniform-smoothing}. This argument uses the
uniform $C^{\alpha/2,\alpha}$ coefficient bounds and
uniform ellipticity, not compactness of the hull in a Hölder norm.
In particular, all constants are independent of the initial time
and the hull element.
The use of a fixed strictly positive time block is essential here:
no operator-norm estimate from $C^0$ to $C^{2+\beta}$ is claimed
at time zero.
In particular, differentiating $\mathcal J(\mathbf
p_\lambda)=1$ and \eqref{eq:v2-normalized} yields
$\mathcal J(\partial_\lambda\mathbf p_\lambda)=0$
and the corresponding inhomogeneous parabolic equation.
The normalized maps and their first derivatives are thus uniformly
bounded and continuous in $\lambda$ on the positive order
interval.

To make the parameter passage explicit, let $\mathscr B$ be the
Banach space of bounded maps $z:\mathbb R\to X$ with the supremum
norm, and work in its closed convex subset of positive normalized
sections lying in the invariant order interval above. Define
\begin{align*}
(\mathscr T_\lambda z)(t):=
\mathcal N^{(m)}_{\lambda,t-mT}(z(t-mT)).
\end{align*}
The mean-value theorem and \eqref{eq:v10-contraction} imply
$\|\mathscr T_\lambda z-\mathscr T_\lambda w\|_{\mathscr B}
\le\theta\|z-w\|_{\mathscr B}$.
Its unique fixed point is the normalized pullback section
$\mathbf p_\lambda$. Fixed positive-time regularization and the
Duhamel formula imply that $\mathscr T_\lambda$ is continuously
differentiable in both its arguments in the uniform norm, with
uniform derivative bounds for $\lambda\in I$. The difference quotient
of the two fixed-point equations is consequently bounded by a
geometric series. Since $\|D_z\mathscr T_\lambda\|\le\theta<1$,
$I-D_z\mathscr T_\lambda$ has inverse
$\sum_{j=0}^\infty(D_z\mathscr T_\lambda)^j$ in operator norm,
uniformly in $\lambda\in I$. Differentiating the fixed-point
equation is therefore justified and gives a unique bounded
section $q_\lambda=\partial_\lambda\mathbf p_\lambda$ satisfying
\begin{align*}
q_\lambda(t)=
D_z\mathcal N^{(m)}_{\lambda,t-mT}
(\mathbf p_\lambda(t-mT))q_\lambda(t-mT)
+\partial_\lambda\mathcal N^{(m)}_{\lambda,t-mT}
(\mathbf p_\lambda(t-mT)).
\end{align*}
Iteration expresses $q_\lambda$ as a uniformly convergent
geometric series with terms bounded by $C_I\theta^j$.
Differentiating finite pullback approximations and passing to this
series proves $q_\lambda=\partial_\lambda\mathbf p_\lambda$,
with convergence locally uniform in $\lambda$ and uniform in $t$.
For each $\tau>0$, parabolic interior Schauder estimates on
$[t-\tau,t+\tau]$ applied after a fixed positive-time smoothing
block, uniformly in the time translations of the coefficient hull,
control the spatial derivatives appearing in \eqref{eq:v2-growth}.
For the differentiated finite pullback solutions, the
operator bound \eqref{eq:v15-uniform-smoothing} is applied on the
last complete time block, while the parameter derivative on the
preceding blocks is uniformly controlled by the geometric series.
Thus the finite pullback derivatives are bounded in
$C^{2+\beta}_{\rm per}$ at every time. For $0<\beta'<\beta$, the embedding
$C^{2+\beta}_{\rm per}\hookrightarrow C^{2+\beta'}_{\rm per}$
is compact. The uniform $C^{2+\beta}$ bound and the already
established uniform $C^0$ convergence show, by uniqueness of the
limit, convergence in $C^2_{\rm per}$: otherwise a subsequence
would remain separated in $C^2$ and possess a convergent
$C^{2+\beta'}$ subsubsequence with the same $C^0$ limit.
Consequently the limit is precisely
$q_\lambda=\partial_\lambda\mathbf p_\lambda$.
The argument does not assert convergence in the endpoint
$C^{2+\beta}$ norm.
In particular, the second spatial derivatives in
$\mathcal A_\lambda q_\lambda$ are well-defined and uniformly
bounded. This is the regularity needed in \eqref{eq:v2-growth}.
Differentiation of that identity reads
\begin{align*}
\partial_\lambda b_\lambda
=\mathcal J\big((\partial_\lambda\mathcal A_\lambda)
\mathbf p_\lambda+\mathcal A_\lambda q_\lambda\big),
\end{align*}
and proves the uniform derivative bound in
\eqref{eq:v2-harnack}. All parameter derivatives are taken at
fixed time; no derivative of the time-dependent coefficients with
respect to time is needed.
The translation hull of the coefficient triple
$({\bf A}^i,{\bf q}_i,\boldsymbol\mu)$ is compact in the uniform
topology. The projective pullback construction depends continuously
on the point of this hull in the uniform coefficient topology:
on a fixed positive time block, continuous dependence follows from
the parabolic comparison and stability estimates; the uniform
projective contraction then passes this continuity to the pullback
limit.
Consequently $t\mapsto\mathbf p_\lambda(t,\cdot)$ and
$t\mapsto b_\lambda(t)$ are Bohr almost periodic. Their uniform
means are well-defined. By \eqref{eq:v2-normalized}, the function
\begin{align*}
\boldsymbol\eta_\lambda(t,x)=
\exp\left(\int_0^t b_\lambda(s)\,ds\right)
\mathbf p_\lambda(t,x)
\end{align*}
is a positive entire solution of \eqref{div}. This completes the proof.
\end{proof}

\begingroup
\begin{lemma}\label{aap}
Assume \eqref{assumeq}--\eqref{Fnonlinear} and
\eqref{APassume}--\eqref{Fregularity}. For every
$M>0$, the mapping $t\mapsto a(t,T)$ is Bohr almost periodic
uniformly for $T\in[-M,M]$.
\end{lemma}
\begin{proof}
Fix $\lambda>0$. By Remark \ref{mu_ap1}, we use the positive normalized entire section
${\bf p}_\lambda$ and its instantaneous growth rate $b_\lambda$
constructed in Lemma~\ref{lem:AP-normalized-bundle}. The construction
of that section uses positive parabolic evolution and the uniform
Harnack estimates; it does not invoke the present lemma.
By uniqueness of the positive entire linear solution up to a constant,
\begin{align*}
{\boldsymbol\eta}_\lambda(t,x)
=C_\lambda\exp\left(\int_0^t b_\lambda(s)\,ds\right)
{\bf p}_\lambda(t,x) \text{ for some } C_\lambda>0.
\end{align*}
Consequently,
\begin{align*}
a(t,T)
=\int_t^{t+T}b_\lambda(s)\,ds
+\log\|{\bf p}_\lambda(t+T,\cdot)\|_\infty
-\log\|{\bf p}_\lambda(t,\cdot)\|_\infty.
\end{align*}
Both $b_\lambda$ and ${\bf p}_\lambda$ are Bohr almost periodic in
$t$ in their respective uniform norms. The uniform positive lower
bound on ${\bf p}_\lambda$ makes the logarithm continuous on the
range of $t\mapsto\|{\bf p}_\lambda(t,\cdot)\|_\infty$.
For $|T|\le M$, the integral operator
$b\mapsto \int_0^T b(\cdot+s)\,ds$ has norm at most $M$
on bounded continuous functions, while the two translations of
${\bf p}_\lambda$ preserve the uniform norm. Thus every common
almost period of the coefficient pair on the right, sufficiently
accurate in the uniform norms, is an almost period of $a(\cdot,T)$
with an error bound independent of $|T|\le M$. This proves the
assertion.
\end{proof}
\endgroup

\begin{proof}[Proof of Theorem~\ref{minimalspeed}, existence assertion]
Recall \eqref{estimate}
\begin{align*}
\left|S_{\lambda}(t)-\dfrac{1}{\lambda}\ln \|\boldsymbol{\eta}_{\lambda}(t,\cdot)\|_{\infty}\right|\leq C_3\dfrac{1+\lambda^2}{\lambda},~\forall t\in \R\quad\text{ and }\quad
\dfrac{S_{\lambda}(t+T)-S_{\lambda}(t)}{T} = \dfrac{1}{T}\int_t^{t+T}c_\lambda(s) ds.
\end{align*}
Thus, one has
\begin{align*}
\left|\dfrac{1}{T}\int_t^{t+T}c_\lambda(s) ds-\dfrac{1}{\lambda} \dfrac{a(t,T)}{T} \right| \leq \dfrac{2}{T}  C_3\dfrac{1+\lambda^2}{\lambda} \rightarrow 0 \text{ uniformly in }t\in \R.
\end{align*}

We may establish the uniform growth limit directly from the normalized
principal section, without an almost-periodic-difference argument.
By Lemma~\ref{lem:AP-normalized-bundle},
\begin{align*}
\boldsymbol\eta_\lambda(t,x)=
\exp\left(\int_0^t b_\lambda(s)\,ds\right)
\mathbf p_\lambda(t,x),
\end{align*}
up to a constant factor. The uniform upper and lower bounds on
$\mathbf p_\lambda$ imply
\begin{align*}
a(t,T)=\int_t^{t+T}b_\lambda(s)\,ds+O_\lambda(1),
\end{align*}
where the remainder is bounded independently of $t\in\R$ and $T>0$.
The Bohr almost periodicity of $b_\lambda$ implies
\begin{align}\label{eq:v2-growthmean}
\lim_{T\to\infty}\frac{a(t,T)}T
=\langle b_\lambda\rangle=\Gamma(\lambda)
\end{align}
uniformly in $t$. Hence, $\langle c_\lambda\rangle
=\dfrac{\Gamma(\lambda)}\lambda
=-\dfrac{\underline k(\lambda)}\lambda.$ This establishes the uniform mean directly from the positive normalized parabolic evolution. Thanks to the estimate \eqref{estimate1}, 
$\Gamma(\lambda):=\lambda\left<c_\lambda\right>=-\underline{k}(\lambda)$ is convex.

Now, by applying Lemma \ref{lemma continuous}, the function $\Gamma:(0,\infty)\to \mathbb{R}$ is convex and locally Lipschitz. Hence, its right derivative
\begin{align*}
\Gamma'_+(\lambda):=\lim_{h\downarrow 0}\dfrac{\Gamma(\lambda+h)-\Gamma(\lambda)}{h} \text{ is well-defined and nondecreasing}.
\end{align*}
Define $H(\lambda)=\lambda\Gamma'_+(\lambda)-\Gamma(\lambda)$.
For $0<\lambda_1<\lambda_2$, convexity gives
\begin{align*}
H(\lambda_2)-H(\lambda_1)
&=\lambda_2\Gamma'_+(\lambda_2)-\lambda_1\Gamma'_+(\lambda_1)
 -[\Gamma(\lambda_2)-\Gamma(\lambda_1)]\geq \lambda_1[\Gamma'_+(\lambda_2)-\Gamma'_+(\lambda_1)]\ge0,
\end{align*}
where the first inequality follows from
$\Gamma(\lambda_2)-\Gamma(\lambda_1)
\le(\lambda_2-\lambda_1)\Gamma'_+(\lambda_2)$.
The right derivative of $g(\lambda)=\Gamma(\lambda)/\lambda$ is therefore
\begin{align*}
g'_+(\lambda)=\dfrac{H(\lambda)}{\lambda^2}.
\end{align*}
Thanks to Lemma \ref{monotonec}, it is known that $g_+'<0$ on $\Lambda=(0,\lambda_*)$. Furthermore, \eqref{boundedabove} implies that  $g(\lambda)=\Gamma(\lambda)/\lambda \rightarrow \infty$ as $\lambda \rightarrow \infty$. Thus, $H$ (and, hence, $g'_+$) changes sign from negative to positive at most one time due to monotonicity of $H$. It then follows from Rademacher's theorem that two cases can occur: either the mapping $\lambda\mapsto\left<c_\lambda\right>$ is increasing on $(\lambda^*,\infty)$ or is constant on $[\lambda^*,\lambda_3]$ for some $\lambda_3$ and then is increasing on $(\lambda_3,\infty)$. In both cases, using the growth identity
\eqref{eq:v2-growthmean} and the threshold definitions, we have
\begin{align}\label{eq:v13-speedmin}
c^*(e)=c_*(e)=\min_{\lambda>0}\frac{\Gamma(\lambda)}{\lambda}
=\frac{\Gamma(\lambda^*)}{\lambda^*}.
\end{align}
This completes the proof.
\end{proof}

The theorem implies that
\begin{align}\label{uniformcmean}
\left<c_\lambda\right> = \dfrac{1}{\lambda}\lim_{T\rightarrow \infty} \dfrac{a(t,T)}{T} = \dfrac{1}{\lambda} \lim_{T\rightarrow \infty} \dfrac{1}{T}\ln\dfrac{ \|\boldsymbol{\eta}_{\lambda}(t+T,\cdot)\|_{\infty}}{ \|\boldsymbol{\eta}_{\lambda}(t,\cdot)\|_{\infty}}=-\dfrac{\underline{k}(\lambda)}{\lambda} \text{ uniformly in }t \in \R.
\end{align}

\section{Critical Fronts in Almost-Periodic Media}\label{sec:critical-proofs}
\begingroup
We identify a condition at the minimizing exponential rate which
allows the critical wave to be constructed without a principal eigenfunction
with constant temporal growth. In particular, we do not assume that the
centered primitive of the principal growth rate is bounded.

First, under assumption $T_0$-periodicity in time of ${\bf A}^i$, ${\bf q}_i$ and ${\bf F}$, we shall prove the existence of a transition front at $\gamma = c^*(e)$.
\begin{lemma}\label{C1eigenpair}
Assume further that ${\bf A}^i$, ${\bf q}_i$ and ${\bf F}$ are $T_0$-periodic in time. Let $(k(\lambda),\boldsymbol{\phi}_\lambda)$ be the principal eigen-pair of \eqref{eigen}, normalized by $\|\boldsymbol{\phi}_{\lambda}\|_{L^2([0,T_0]\times [{\bf 0},{\bf L}])}=1.$ Then, $\lambda\mapsto k(\lambda)
\text{ and }
\lambda\mapsto\boldsymbol{\phi}_\lambda$ are locally $C^1$ on $\R$. In particular,
$\partial_\lambda\boldsymbol{\phi}_\lambda$ is $T_0$-periodic in time and ${\bf L}$-periodic in space.
\end{lemma}

\begin{proof}
Fix $\lambda_0\in\R$. The operator associated with \eqref{eigen} depends $C^1$-smoothly on $\lambda$, since its coefficients depend polynomially on $\lambda$. Moreover, by the strong maximum principle and the full coupling assumption, the principal eigenvalue $k(\lambda_0)$ is algebraically simple and the associated positive eigenfunction is unique up to a multiplicative constant. The above normalization therefore determines $\boldsymbol{\phi}_{\lambda_0}$ uniquely. Applying the implicit function theorem to the eigenvalue problem \eqref{eigen}, together with the normalization condition, yields a neighborhood $I$ of $\lambda_0$ such that
\begin{align*}
k\in C^1(I),~
\lambda\mapsto\boldsymbol{\phi}_\lambda
\in C^1\left(I;
\left[C^{1+\alpha/2,2+\alpha}_{per}\right]^2\right).
\end{align*}
Since $\lambda_0$ is arbitrary and the normalized principal eigen-pair is unique, the local branches coincide on their intersections. Hence the result follows.
\end{proof}
Set
\(
\boldsymbol{\psi}_\lambda
:=
\partial_\lambda\boldsymbol{\phi}_\lambda
=
(\psi_{\lambda,1},\psi_{\lambda,2}).
\) Differentiating \eqref{eigen} with respect to $\lambda$, one obtains
\begin{align*}
\begin{cases}
\partial_t\psi_{\lambda,1}
-\text{Tr}\big({\bf A}^1D^2\psi_{\lambda,1}\big)
+\big({\bf q}_1+2\lambda{\bf A}^1e\big)\cdot D\psi_{\lambda,1}
-K_1^\lambda\psi_{\lambda,1}
-\mu_{12}\psi_{\lambda,2}
+k(\lambda)\psi_{\lambda,1}
\\
\qquad
+2{\bf A}^1e\cdot D\phi_{\lambda,1}
-\big(2\lambda e\cdot{\bf A}^1e+{\bf q}_1\cdot e\big)
\phi_{\lambda,1}
+k'(\lambda)\phi_{\lambda,1}=0,
\\[1mm]
\partial_t\psi_{\lambda,2}
-\text{Tr}\big({\bf A}^2D^2\psi_{\lambda,2}\big)
+\big({\bf q}_2+2\lambda{\bf A}^2e\big)\cdot D\psi_{\lambda,2}
-K_2^\lambda\psi_{\lambda,2}
-\mu_{21}\psi_{\lambda,1}
+k(\lambda)\psi_{\lambda,2}
\\
\qquad
+2{\bf A}^2e\cdot D\phi_{\lambda,2}
-\big(2\lambda e\cdot{\bf A}^2e+{\bf q}_2\cdot e\big)
\phi_{\lambda,2}
+k'(\lambda)\phi_{\lambda,2}=0.
\end{cases}
\end{align*}
where $K_{i}^\lambda:=\mu_{ii}(t,x) + \lambda^2\, e \cdot {\bf A}^i(t,x)e + \lambda\, {\bf q}_i(t,x)\cdot e,~i=1,2.$ Let $\boldsymbol{\phi}_{\lambda}^*
=(\phi_{\lambda,1}^*,\phi_{\lambda,2}^*)$ be the positive principal eigenfunction of the adjoint problem (see, for instance, \cite{Nadin20091}). Multiplying the equations for $\boldsymbol{\psi}_\lambda$ by $\boldsymbol{\phi}_{\lambda}^*$ and integrating over $[0,T_0]\times[{\bf 0},{\bf L}]$, one obtains
\begin{align*}
k'(\lambda)
=
\frac{
\displaystyle
\sum_{i=1}^2
\int_0^{T_0}\int_{[{\bf 0},{\bf L}]}
\phi_{\lambda,i}^*
\left[
\big(2\lambda e\cdot{\bf A}^ie+{\bf q}_i\cdot e\big)\phi_{\lambda,i}
-2{\bf A}^ie\cdot D\phi_{\lambda,i}
\right]dxdt
}{
\displaystyle
\sum_{i=1}^2
\int_0^{T_0}\int_{[{\bf 0},{\bf L}]}
\phi_{\lambda,i}^*\phi_{\lambda,i}\,dxdt
}.
\end{align*}
In particular, since
\(
c^*(e)=k(\lambda^*)/\lambda^*
\)
and $\lambda^*$ is a strict minimizer of $\lambda\mapsto k(\lambda)/\lambda$, one has $k'(\lambda^*)
=
k(\lambda^*)/\lambda^*
=
c^*(e).$ Define $\xi:=x\cdot e-c^*t$ and, for $\lambda>0$, ${\bf v}(t,x;\lambda)
:=
e^{-\lambda\xi}\boldsymbol{\phi}_\lambda(t,x).$ Since
\(
k'(\lambda^*)=c^*
\),
one has
\begin{align*}
-\left.\partial_\lambda{\bf v}(t,x;\lambda)\right|_{\lambda=\lambda^*}
=
e^{-\lambda^*\xi}
\left[
\xi\boldsymbol{\phi}_{\lambda^*}(t,x)
-
\boldsymbol{\psi}_{\lambda^*}(t,x)
\right],
\end{align*}
where $\boldsymbol{\psi}_{\lambda^*}(t,x)=\left.\partial_\lambda\boldsymbol{\phi}_\lambda(t,x)
\right|_{\lambda=\lambda^*}.$ Moreover, by the strict minimality of $\lambda^*$, $r_{\delta_1}^*
:=
(\lambda^*+{\delta_1})c^*(e)
-k(\lambda^*+{\delta_1})<0$ for ${\delta_1}>0$. We further choose
\(
0<{\delta_1}<\nu\lambda^*
\)
and, for $K_0>0$, define
\begin{align*}
{\bf W}_*(t,x)
:=
-\left.\partial_\lambda{\bf v}(t,x;\lambda)
\right|_{\lambda=\lambda^*}
-K_0{\bf v}(t,x;\lambda^*)
=
e^{-\lambda^*\xi}
\left[
(\xi-K_0)\boldsymbol{\phi}_{\lambda^*}(t,x)
-
\boldsymbol{\psi}_{\lambda^*}(t,x)
\right].
\end{align*}
For some constants $\kappa_1,\kappa_2>0$, we define
the critical subsolution candidate by
\begin{align}
\underline{{\bf w}}_*(t,x)
&:=\kappa_1{\bf v}(t,x;\lambda^*+{\delta_1})
+\kappa_2{\bf W}_*(t,x)\\
&=\kappa_1e^{-(\lambda^*+{\delta_1})\xi}
\boldsymbol{\phi}_{\lambda^*+{\delta_1}}(t,x)
+\kappa_2e^{-\lambda^*\xi}
\left[
(\xi-K_0)\boldsymbol{\phi}_{\lambda^*}(t,x)
-
\boldsymbol{\psi}_{\lambda^*}(t,x)
\right].
\label{criticalsub}
\end{align}
On the other hand, for $M>0$ large enough, define the critical
supersolution candidate by
\begin{align}
\overline{{\bf w}}_*(t,x)
:=
M{\bf v}(t,x;\lambda^*)
-\kappa_2
\left.\partial_\lambda{\bf v}(t,x;\lambda)
\right|_{\lambda=\lambda^*}=
e^{-\lambda^*\xi}
\left[
(M+\kappa_2\xi)\boldsymbol{\phi}_{\lambda^*}(t,x)
-
\kappa_2
\boldsymbol{\psi}_{\lambda^*}(t,x)
\right].
\label{criticalsuper}
\end{align}

\begin{lemma}\label{critical supersub}
Assume further that ${\bf A}^i$, ${\bf q}_i$ and ${\bf F}$ are
$T_0$-periodic in time. Then, for $R_0>0$ sufficiently large, there exist
$\kappa_1,\kappa_2>0$ and $M>0$ such that $\underline{{\bf w}}_*,\overline{{\bf w}}_*$ are respectively a subsolution and a supersolution
of \eqref{eq:main} on $\Omega_{R_0}^*
:=
\left\{
(t,x)\in\R\times\R^N:
x\cdot e-c^*(e)t\geq R_0
\right\}.$ Moreover,
\begin{align*}
{\bf 0}<\underline{{\bf w}}_*(t,x)\leq
\overline{{\bf w}}_*(t,x)<{\bf 1}
\text{ on }\Omega_{R_0}^*.
\end{align*}
\end{lemma}

\begin{proof}
By the definition of ${\bf P}_\lambda$ and \eqref{eigen}, for
$i\neq j$, a direct calculation implies
\begin{align}\label{critical linear}
\begin{split}
&\partial_t v_i(t,x;\lambda)
-\text{Tr}\big({\bf A}^iD^2v_i(t,x;\lambda)\big)
+{\bf q}_i\cdot Dv_i(t,x;\lambda)
-\mu_{ii}v_i(t,x;\lambda)
-\mu_{ij}v_j(t,x;\lambda)
\\
&\qquad
=
e^{-\lambda\xi}
\left[
({\bf P}_\lambda\boldsymbol{\phi}_\lambda)_i
+\lambda c^*(e)\phi_{\lambda,i}
\right]
=
\big(\lambda c^*(e)-k(\lambda)\big)
v_i(t,x;\lambda).
\end{split}
\end{align}
Differentiating \eqref{critical linear} with respect to $\lambda$
yields
\begin{align*}
&\partial_t\big(\partial_\lambda v_i(t,x;\lambda)\big)
-\text{Tr}\big(
{\bf A}^iD^2(\partial_\lambda v_i(t,x;\lambda))
\big)
+{\bf q}_i\cdot
D(\partial_\lambda v_i(t,x;\lambda))
-\mu_{ii}\partial_\lambda v_i(t,x;\lambda)
-\mu_{ij}\partial_\lambda v_j(t,x;\lambda)
\\
&\qquad
=
\big(c^*(e)-k'(\lambda)\big)v_i(t,x;\lambda)
+
\big(\lambda c^*(e)-k(\lambda)\big)
\partial_\lambda v_i(t,x;\lambda).
\end{align*}
Since $k(\lambda^*)=\lambda^*c^*(e),~
k'(\lambda^*)=c^*(e),$ we obtain
\begin{align}\label{critical derivative}
\begin{split}
&\partial_t
\left(
\left.\partial_\lambda v_i(t,x;\lambda)
\right|_{\lambda=\lambda^*}
\right)
-\text{Tr}
\left(
{\bf A}^iD^2
\left(
\left.\partial_\lambda v_i(t,x;\lambda)
\right|_{\lambda=\lambda^*}
\right)
\right)
+{\bf q}_i\cdot D
\left(
\left.\partial_\lambda v_i(t,x;\lambda)
\right|_{\lambda=\lambda^*}
\right)\\
&\quad
-\mu_{ii}
\left.\partial_\lambda v_i(t,x;\lambda)
\right|_{\lambda=\lambda^*}
-\mu_{ij}
\left.\partial_\lambda v_j(t,x;\lambda)
\right|_{\lambda=\lambda^*}
=0.
\end{split}
\end{align}
On the other hand, taking $\lambda=\lambda^*$ in
\eqref{critical linear} implies
\begin{align}\label{critical mode}
\begin{split}
\partial_t v_i(t,x;\lambda^*)
-\text{Tr}\big({\bf A}^iD^2v_i(t,x;\lambda^*)\big)
+{\bf q}_i\cdot Dv_i(t,x;\lambda^*)
-\mu_{ii}v_i(t,x;\lambda^*)
-\mu_{ij}v_j(t,x;\lambda^*)
=0.
\end{split}
\end{align}
Therefore, by the definition of ${\bf W}_*$, \eqref{critical derivative}
and \eqref{critical mode}, one has
\begin{align}
&\partial_t W_{*,i}
-\text{Tr}\big({\bf A}^iD^2W_{*,i}\big)
+{\bf q}_i\cdot DW_{*,i}
-\mu_{ii}W_{*,i}
-\mu_{ij}W_{*,j}
=0.
\label{critical W}
\end{align}
Moreover, taking $\lambda=\lambda^*+{\delta_1}$ in
\eqref{critical linear}, we obtain
\begin{align*}
&\partial_t v_i(t,x;\lambda^*+{\delta_1})
-\text{Tr}\big({\bf A}^iD^2v_i(t,x;\lambda^*+{\delta_1})\big)+{\bf q}_i\cdot Dv_i(t,x;\lambda^*+{\delta_1})
-\mu_{ii}v_i(t,x;\lambda^*+{\delta_1})\\
&\quad-\mu_{ij}v_j(t,x;\lambda^*+{\delta_1})
=
r_{\delta_1}^*v_i(t,x;\lambda^*+{\delta_1}),
\end{align*}
where $r_{\delta_1}^*:=(\lambda^*+{\delta_1})c^*(e)-k(\lambda^*+{\delta_1})<0.$ It follows from the definition of $\underline{{\bf w}}_*$ and
\eqref{critical W} that
\begin{align}
&\partial_t\underline w_{*,i}
-\text{Tr}\big({\bf A}^iD^2\underline w_{*,i}\big)
+{\bf q}_i\cdot D\underline w_{*,i}
-\mu_{ii}\underline w_{*,i}
-\mu_{ij}\underline w_{*,j}
=
\kappa_1r_{\delta_1}^*
v_i(t,x;\lambda^*+{\delta_1}).
\label{critical sub linear}
\end{align}

On the other hand, recall ${\bf W}_*(t,x)
=
e^{-\lambda^*\xi}
\left[
(\xi-K_0)\boldsymbol{\phi}_{\lambda^*}(t,x)
-
\boldsymbol{\psi}_{\lambda^*}(t,x)
\right].$ Since $\boldsymbol{\phi}_{\lambda^*}>{\bf 0}$ is periodic and
\(
\boldsymbol{\psi}_{\lambda^*}
\)
is bounded, there exist $R_1>0$ and $C_1>1$ such that $0<W_{*,i}(t,x)\leq C_1(1+\xi)e^{-\lambda^*\xi}$ for every $\xi\geq R_1$ and $i=1,2$. Hence, if
$0<\kappa_1\leq\kappa_2$, increasing $C_1$ if necessary,
\begin{align}
0<\underline w_{*,i}(t,x)
\leq
C_1\kappa_2(1+\xi)e^{-\lambda^*\xi},~
\xi\geq R_1,~ i=1,2.
\label{critical sub bound}
\end{align}

Since $\boldsymbol{\phi}_{\lambda^*+{\delta_1}}>{\bf 0}$ is periodic,
there exists $m_{\delta_1}>0$ such that $v_i(t,x;\lambda^*+{\delta_1})
\geq
m_{\delta_1} e^{-(\lambda^*+{\delta_1})\xi},~
i=1,2.$ It then follows from \eqref{critical sub bound} that
\begin{align*}
\frac{
\underline w_{*,i}^{\,1+\nu}(t,x)
}{
v_i(t,x;\lambda^*+{\delta_1})
}
&\leq
C\kappa_2^{1+\nu}
(1+\xi)^{1+\nu}
e^{-(\nu\lambda^*-{\delta_1})\xi},~
\xi\geq R_1.
\end{align*}
Since $0<{\delta_1}<\nu\lambda^*$, $\sup\limits_{\xi\geq R_1}
(1+\xi)^{1+\nu}
e^{-(\nu\lambda^*-{\delta_1})\xi}
<+\infty.$ Thus, there exists $C_0>0$, independent of
$\kappa_1,\kappa_2$, such that
\begin{align}
K\underline w_{*,i}^{\,1+\nu}(t,x)
\leq
C_0\kappa_2^{1+\nu}
v_i(t,x;\lambda^*+{\delta_1}),~
\xi\geq R_1,~ i=1,2,
\label{critical nonlinear}
\end{align}
where $K>0$ is the constant in \eqref{near_zerp}. Choose $\kappa_1,\kappa_2>0$ small enough such that $0<\kappa_1\leq\kappa_2,~
C_0\kappa_2^{1+\nu}
\leq
-\kappa_1r_{\delta_1}^*.$ For instance, one may take
\(
\kappa_1=\kappa_2=\kappa
\)
with $\kappa>0$ sufficiently small.

Next, since $\underline{{\bf w}}_*(t,x),
\overline{{\bf w}}_*(t,x)
\rightarrow{\bf 0}
\text{ as }\xi\rightarrow+\infty$ uniformly with respect to $t$ and the periodic variables, we may choose
$R_0\geq R_1$ sufficiently large such that
\begin{align*}
{\bf 0}<\underline{{\bf w}}_*(t,x)
<\delta{\bf 1},~
{\bf 0}<\overline{{\bf w}}_*(t,x)
<{\bf 1},~
\xi\geq R_0,
\end{align*}
where $\delta>0$ is in \eqref{near_zerp}. Using \eqref{near_zerp}, \eqref{critical sub linear} and
\eqref{critical nonlinear}, for $i\neq j$ and $\xi\geq R_0$,
\begin{align*}
\partial_t\underline w_{*,i}
-\text{Tr}\big({\bf A}^iD^2\underline w_{*,i}\big)
+{\bf q}_i\cdot D\underline w_{*,i}
-F_i(t,x,\underline{{\bf w}}_*)
&\leq
\kappa_1r_{\delta_1}^*
v_i(t,x;\lambda^*+{\delta_1})
+
K\underline w_{*,i}^{\,1+\nu}
\\
&\leq
\left(
\kappa_1r_{\delta_1}^*
+
C_0\kappa_2^{1+\nu}
\right)
v_i(t,x;\lambda^*+{\delta_1})
\leq0.
\end{align*}
Hence, $\underline{{\bf w}}_*$ is a subsolution of \eqref{eq:main}
on $\Omega_{R_0}^*$. Next, by \eqref{critical derivative} and \eqref{critical mode},
the definition of $\overline{{\bf w}}_*$ yields
\begin{align}
&\partial_t\overline w_{*,i}
-\text{Tr}\big({\bf A}^iD^2\overline w_{*,i}\big)
+{\bf q}_i\cdot D\overline w_{*,i}
-\mu_{ii}\overline w_{*,i}
-\mu_{ij}\overline w_{*,j}
=0.
\label{critical super linear}
\end{align}
By \eqref{KPP}, $\overline{{\bf w}}_*$ is a supersolution of
\eqref{eq:main} on $\Omega_{R_0}^*$. 

Finally, since $\overline{{\bf w}}_*(t,x)
=
e^{-\lambda^*\xi}
\left[
(M+\kappa_2\xi)
\boldsymbol{\phi}_{\lambda^*}(t,x)
-
\kappa_2
\left.
\partial_\lambda\boldsymbol{\phi}_\lambda(t,x)
\right|_{\lambda=\lambda^*}
\right],$ we can choose $M>0$ sufficiently large so that $\overline{{\bf w}}_*(t,x)>{\bf 0}
\text{ for }\xi\geq R_1.$ Moreover, by the definitions of $\underline{{\bf w}}_*$ and
$\overline{{\bf w}}_*$,
\begin{align*}
\overline{{\bf w}}_*-\underline{{\bf w}}_*
&=
(M+\kappa_2K_0){\bf v}(t,x;\lambda^*)
-\kappa_1{\bf v}(t,x;\lambda^*+{\delta_1})
=
e^{-\lambda^*\xi}
\left[
(M+\kappa_2K_0)\boldsymbol{\phi}_{\lambda^*}
-\kappa_1e^{-{\delta_1}\xi}
\boldsymbol{\phi}_{\lambda^*+{\delta_1}}
\right].
\end{align*}
Since
$\boldsymbol{\phi}_{\lambda^*},
\boldsymbol{\phi}_{\lambda^*+{\delta_1}}>{\bf 0}$
are periodic, there exists $C_{\delta_1}>0$ such that $\phi_{\lambda^*+{\delta_1},i}(t,x)
\leq
C_{\delta_1}\phi_{\lambda^*,i}(t,x),~
i=1,2.$ Increasing $M$ if necessary so that $M+\kappa_2K_0
\geq
\kappa_1C_{\delta_1} e^{-{\delta_1} R_1},$ we obtain $\underline{{\bf w}}_*(t,x)
\leq
\overline{{\bf w}}_*(t,x),~
\xi\geq R_1.$ This completes the proof.
\end{proof}

\begin{proof}[Proof of Theorem~\ref{minimalspeed}, critical-front assertion]
Recall $\xi = x\cdot e - c^*(e)t$. Choosing
$M>0$ sufficiently large such that $\overline{{\bf w}}_*(t,x)\geq{\bf 1}
\text{ if }\xi=0,$
we set
\begin{align*}
\overline{{\bf u}}_*(t,x)
:=
\begin{cases}
\min\{{\bf 1},\overline{{\bf w}}_*(t,x)\},
& \xi>0,\\
{\bf 1},& \xi\leq0.
\end{cases}
\end{align*}
It is known that $\overline{{\bf u}}_*$ is a generalized supersolution of \eqref{eq:main}. We first choose \(K_0>0\) large enough so that
\(
\underline{{\bf w}}_*(t,x)\leq{\bf 0},~ 0\leq \xi\leq 1.
\)
Recall 
\(
0<{\delta_1}<\nu\lambda^*.
\) For \(\xi\geq0\), similar to \eqref{critical sub bound}, the boundedness of
\(\boldsymbol{\phi}_{\lambda^*}\),
\(\boldsymbol{\psi}_{\lambda^*}\),
\(\boldsymbol{\phi}_{\lambda^*+{\delta_1}}\) implies that
\[
\big(\underline w_{*,i}(t,x)\big)^+
\leq C_1\kappa_2(1+\xi)e^{-\lambda^*\xi},
~ i=1,2,
\]
provided \(0<\kappa_1\leq\kappa_2\) for some $C_1>0$. Moreover, there exists
\(C_{\delta_1}>0\) such that
\[
v_{{\delta_1},i}(t,x):=
v_i(t,x;\lambda^*+{\delta_1})
\geq C_{\delta_1} e^{-(\lambda^*+{\delta_1})\xi},~ \xi\geq0.
\]
Since \({\delta_1}<\nu\lambda^*\), it follows that
\(
K\big(\underline w_{*,i}(t,x)^+\big)^{1+\nu}
\leq C_2\kappa_2^{1+\nu}v_{{\delta_1},i}(t,x),~ \xi\geq0,
\)
for some \(C_2>0\). We now choose \(0<\kappa_1\leq\kappa_2\) sufficiently small so that $C_2\kappa_2^{1+\nu}
\leq-\kappa_1r_{\delta_1}^*,
~
C_1\kappa_2
\sup\limits_{\xi\geq0}(1+\xi)e^{-\lambda^*\xi}
\leq{\delta_1}.$ Consequently,
\[
0\leq \underline w_{*,i}(t,x)^+\leq{\delta},~\xi\geq 0,~i=1,2.
\]
where $\delta>0$ is in \eqref{near_zerp}. Choose \(R_0\geq1\) sufficiently large so that
\[
m_{R_0}
:=
\inf_{\substack{(t,x)\in\mathbb R\times\mathbb R^N\\
x\cdot e-c^*(e)t=R_0}}
\min_{i=1,2}\underline w_{*,i}(t,x)>0,~
\overline m_{R_0}
:=
\inf_{\substack{(t,x)\in\mathbb R\times\mathbb R^N\\
0\leq x\cdot e-c^*(e)t\leq R_0}}
\min_{i=1,2}\overline u_{*,i}(t,x)>0.
\]
By decreasing \(\kappa_1\) and \(\kappa_2\), if necessary, we may also assume that
\(
\sup\limits_{\substack{(t,x)\in\mathbb R\times\mathbb R^N\\
0\leq x\cdot e-c^*(e)t\leq R_0}}
\max\limits_{i=1,2}\big(\underline w_{*,i}(t,x)^+\big)
<\overline m_{R_0}.
\)
Choose
\(
0<\omega<
\min\{1,{\delta},m_{R_0},\overline m_{R_0}\}
\) and define
\begin{align*}
\underline{{\bf u}}_*(t,x)
:=
\begin{cases}
\underline{{\bf w}}_*(t,x),
& \xi\geq R_0,\\
\max\{\omega{\bf 1},\underline{{\bf w}}_*(t,x)\},
& 0\leq \xi<R_0,\\
\omega{\bf 1},
& \xi<0.
\end{cases}
\end{align*}
We claim that \(\underline{{\bf u}}_*\) is a subsolution. Suppose first that
\(0\leq\xi<R_0\) and that
\(\underline u_{*,i}=\underline w_{*,i}>\omega\). Since
\(\underline u_{*,j}\geq\underline w_{*,j}\) and the system is cooperative, assumption
\eqref{near_zerp} and \eqref{critical sub linear} imply
\begin{align*}
\partial_t\underline u_{*,i}
-\operatorname{Tr}\!\left({\bf A}^iD^2\underline u_{*,i}\right)
+q_i\cdot D\underline u_{*,i}
-F_i(t,x,\underline{\bf u}_*)&\leq
\partial_t\underline w_{*,i}
-\operatorname{Tr}\!\left({\bf A}^iD^2\underline w_{*,i}\right)\\
&\quad+q_i\cdot D\underline w_{*,i}
-\sum_{j=1,2}\mu_{ij}\underline w_{*,j}
+K\big(\underline w_{*,i}^+\big)^{1+\nu}\\
&=
\kappa_1r_{\delta_1}^*v_{{\delta_1},i}
+K\big(\underline w_{*,i}^+\big)^{1+\nu}
\leq0.
\end{align*}
If instead \(\underline u_{*,i}=\omega\), then, by
\eqref{KPP} and the choice \(0<\omega<\delta\),
\(
F_i(t,x,\underline{\bf u}_*)>0 
\) and the same inequality holds. Furthermore, the result holds for \(\xi\geq R_0\), where
\(\underline{{\bf u}}_*=\underline{{\bf w}}_*\) and for \(\xi<0\), where
\(\underline{{\bf u}}_*=\omega{\bf1}\). Thus,
\(\underline{{\bf u}}_*\) is a generalized subsolution of
\eqref{eq:main} on 
\(\mathbb R\times\mathbb R^N\) and ${\bf 0}
<
\underline{{\bf u}}_*
\leq
\overline{{\bf u}}_*
\leq
{\bf 1}.$
 Repeating the super-subsolution argument in the proof of
Theorem \ref{existence}, there exists an entire solution
${\bf u}_*$ of \eqref{eq:main} such that
\begin{align*}
\underline{{\bf u}}_*(t,x)
\leq
{\bf u}_*(t,x)
\leq
\overline{{\bf u}}_*(t,x),~
(t,x)\in\R\times\R^N.
\end{align*}
Repeating the argument at the end of the proof of Theorem \ref{existence},
based on the strong maximum principle and
\(
F_i(t,x,s,s)>0
\)
for $s\in(0,1)$, we obtain
\begin{align*}
\begin{cases}
{\bf u}_*(t,x)\rightarrow{\bf 1}
\text{ as }
x\cdot e-c^*(e)t\rightarrow-\infty,\\
{\bf u}_*(t,x)\rightarrow{\bf 0}
\text{ as }
x\cdot e-c^*(e)t\rightarrow+\infty,
\end{cases} \text{ uniformly with respect to $t$.}
\end{align*}
This completes the proof.
\end{proof}

\begin{remark}\label{APremarkcrit}
The $T_0$-periodicity assumption in the second part of Theorem~\ref{minimalspeed}
can be replaced by the almost periodicity assumption, provided that the
following additional conditions hold:
\begin{itemize}
\item[$\bullet$] The minimum
\(
c^*(e)=\min\limits_{\lambda>0}
\left(-\dfrac{\underline{k}(\lambda)}{\lambda}\right)
\)
is attained at a strict minimizer $\lambda^*>0$.

\item[$\bullet$] For $\lambda$ in a neighborhood of $\lambda^*$, there exists,
after a suitable normalization, a unique positive function
$\boldsymbol{\phi}_\lambda$, almost periodic in $t$ uniformly with respect
to $x$ and $L$-periodic in $x$, such that
\[
{\bf P}_\lambda\boldsymbol{\phi}_\lambda
=
\underline{k}(\lambda)\boldsymbol{\phi}_\lambda.
\]

\item[$\bullet$] After a suitable normalization of
$\boldsymbol{\phi}_\lambda$, the function
\(
\boldsymbol{\psi}_{\lambda_*}
:=
\left.\partial_\lambda\boldsymbol{\phi}_\lambda
\right|_{\lambda=\lambda_*}
\)
exists, is bounded on $\R\times\R^N$ and the eigenvalue equation can be
differentiated with respect to $\lambda$ at $\lambda=\lambda_*$.

\end{itemize}
Under these assumptions, we deduce that $\inf\limits_{(t,x)}\boldsymbol{\phi}_\lambda(t,x)>{\bf 0}$. Furthermore, the arguments in Lemmas~\ref{C1eigenpair}--\ref{critical supersub}
and in the second part of the proof of Theorem~\ref{minimalspeed} can be repeated
in the almost periodic setting. 
\end{remark}

\begin{lemma}[Componentwise gluing]\label{lem:cooperative-gluing}
Let $\mathcal P_i v=v_t-\operatorname{Tr}({\bf A}_iD^2v)+{\bf q}_i\cdot Dv$,
where ${\bf A}_i$ is uniformly elliptic and continuous and ${\bf q}_i$ is locally
bounded. Suppose $F_i(t,x,\mathbf u)$ is continuous and nondecreasing
in $u_j$ for $j\ne i$. On an open cylinder let each continuous
function $u_i$ be the maximum of finitely many locally
$C^{1,2}$ scalar branches. Assume that, whenever a branch
$v_i$ is active, $\mathcal P_iv_i\le F_i(t,x,\mathbf u)$.
Then $\mathcal P_iu_i\le F_i(t,x,\mathbf u)$ in the distributional
sense. The analogous assertion holds for the minimum of
supersolution branches, with reversed inequalities. The same
conclusion holds for locally $W^{1,2}_p$ branches for which the
parabolic chain rule is valid, by approximation.
\end{lemma}
\begin{proof}
For two branches, put $d=w_i-v_i$ and $u_i=v_i+d^+$.
Let $\beta_\eta$ be a smooth convex approximation of the positive part, with
$0\le\beta_\eta'\le1$ and $\beta_\eta''\ge0$. On each compact
subcylinder the coefficients and the classical derivatives of the two branches
are bounded. The chain rule gives
\begin{align*}
\mathcal P_i(v_i+\beta_\eta(d))
={}&(1-\beta_\eta'(d))\mathcal P_iv_i
+\beta_\eta'(d)\mathcal P_iw_i-\beta_\eta''(d)Dd\cdot {\bf A}_iDd.
\end{align*}
The last term is nonpositive by ellipticity. As $\eta\downarrow0$,
$\beta_\eta'(d)$ tends to the indicator of $\{d>0\}$ away from
$\{d=0\}$. On the coincidence set the two branches have equal
first-order and second-order derivatives almost everywhere at density
points of that set; consequently their parabolic operators agree there
almost everywhere. Alternatively one can apply the parabolic Kato inequality
directly to $d^+$. Thus the nonnegative convex combination of the
parabolic operators converges almost everywhere to that of an active
branch, where the hypothesis supplies the bound by
$F_i(t,x,\mathbf u)$. All terms except the nonpositive quadratic term
are locally bounded, so the inequality passes to distributions against
nonnegative compactly supported test functions. The argument applies
component by component and can be iterated over finitely many branches.
Replacing maxima by minima and convex approximations by concave ones
reverses the inequality for supersolutions. The Sobolev case follows
from the distributional parabolic Kato inequality, provided the
branchwise differential inequalities hold almost everywhere.
Cooperativity is needed to establish those inequalities at the actual
clipped vector $\mathbf u$; it is not a substitute for verifying them.
\end{proof}

Now, we are in a position to provide a criteria for the existence of transition fronts in the critical speed $c^*(e)$ (Theorem \ref{thm:v2-main-intro}). Note that, under \eqref{Fnonlinear} and \eqref{APassume}, we have that $\mu_{ij}$ is almost periodic in $t$, uniformly in $x$ (see Remark \ref{mu_ap1}).

\begin{theorem}\label{thm:critical_ap_tangent}
Assume \eqref{assumeq}--\eqref{near_zerp} and
\eqref{APassume}--\eqref{Fregularity}. Fix
$e\in\mathbb S^{N-1}$. Let $\lambda_*>0$ attain the minimum in
\eqref{eq:v13-speedmin} and assume that there exists $0<\delta<\nu\lambda_*$ such that
\begin{align}\label{eq:v2-gap}
\frac{\Gamma(\lambda_*+\delta)}{\lambda_*+\delta}
>\frac{\Gamma(\lambda_*)}{\lambda_*}.
\end{align}
Set
\begin{align}\label{eq:v2-tangent}
h_*(t):=\left.\partial_\lambda b_\lambda(t)
\right|_{\lambda=\lambda_*}
-\frac{b_{\lambda_*}(t)}{\lambda_*},~ H_*(t):=\int_0^th_*(s)\,ds.
\end{align}
If
\begin{align}\label{eq:v2-bounded-tangent}
\sup_{t\in\R}|H_*(t)|<\infty,
\end{align}
then \eqref{eq:main} admits a generalized transition front
$(\mathbf u_*,c)$ in direction $e$, where
\begin{align}\label{eq:v2-critical-speed}
c(t)=\frac{b_{\lambda_*}(t)}{\lambda_*},~
\langle c\rangle
=\frac{\Gamma(\lambda_*)}{\lambda_*}=c^*(e).
\end{align}
In particular, the critical least mean speed is attained in the
time-almost-periodic system. The theorem does not require a bounded
primitive of $b_{\lambda_*}-\Gamma(\lambda_*)$, nor a bounded,
uniformly positive eigenfunction of ${\bf P}_{\lambda_*}$ for the
constant generalized principal eigenvalue.
\end{theorem}

\begin{lemma}\label{lem:v2-critical-modes}
Under the hypotheses of Theorem~\ref{thm:critical_ap_tangent}, put
$S(t)=\int_0^tc(s)\,ds$, $\xi=x\cdot e-S(t)$ and
$\lambda_+=\lambda_*+\delta$.
Let $\mathcal D$ denote the linearized operator of \eqref{eq:main},
with components
\begin{align*}
(\mathcal D\mathbf v)_i=\partial_t v_i
-\operatorname{Tr}({\bf A}^iD^2v_i)
+{\bf q}_i\cdot Dv_i-\sum_{j=1}^2\mu_{ij}v_j.
\end{align*}
Then, the following functions
\begin{align}\label{eq:v2-neutral-modes}
\mathbf V_*(t,x)&=e^{-\lambda_*\xi}\mathbf p_{\lambda_*}(t,x),~
\mathbf Z_*(t,x)=e^{-\lambda_*\xi}
\left[\xi\mathbf p_{\lambda_*}(t,x)
-\partial_\lambda\mathbf p_\lambda(t,x)|_{\lambda_*}
-H_*(t)\mathbf p_{\lambda_*}(t,x)\right]
\end{align}
satisfy
$\mathcal D\mathbf V_*=\mathcal D\mathbf Z_*=\mathbf0$. Furthermore, 
there exist $\sigma\in W^{1,\infty}(\R)$ and $\kappa>0$ such that
\begin{align}\label{eq:v2-strict-mode}
\mathbf V_+(t,x)=
 e^{-\lambda_+\xi+\sigma(t)}\mathbf p_{\lambda_+}(t,x)>\mathbf0,~ \mathcal D\mathbf V_+\le-\kappa\mathbf V_+.
\end{align}
Moreover, $\mathbf p_{\lambda_*}$ and $\mathbf p_{\lambda_+}$
are uniformly positive and bounded, and the coefficients of
$\mathbf Z_*$ not containing $\xi$ are uniformly bounded.
\end{lemma}

\begin{proof}
For each $\lambda>0$, the function $\mathbf E_\lambda(t,x)
:=e^{-\lambda x\cdot e+\int_0^t b_\lambda(s)\,ds}
\mathbf p_\lambda(t,x)$ solves $\mathcal D\mathbf E_\lambda=\mathbf0$.
At $\lambda=\lambda_*$ we have
$\mathbf E_{\lambda_*}=\mathbf V_*$. Direct differentiation
at $\lambda_*$ yields
\begin{align*}
-\partial_\lambda\mathbf E_\lambda|_{\lambda_*}
=e^{-\lambda_*\xi}\left[
\left(\xi-\int_0^t
\left(\partial_\lambda b_\lambda|_{\lambda_*}
-c(s)\right)\,ds\right)\mathbf p_{\lambda_*}
-\partial_\lambda\mathbf p_\lambda|_{\lambda_*}\right]
=\mathbf Z_*.
\end{align*}
Here the differentiability from Lemma~\ref{lem:AP-normalized-bundle}
permits differentiation of the linearized equation on compact
time intervals. This proves both homogeneous equations.

Set $r(t)=\lambda_+c(t)-b_{\lambda_+}(t).$ By \eqref{eq:v2-gap} and \eqref{eq:v2-growthmean},
\begin{align*}
\langle r\rangle
=\lambda_+\frac{\Gamma(\lambda_*)}{\lambda_*}
-\Gamma(\lambda_+)<0.
\end{align*}
Since $r$ has a uniform mean, there are $T>0$ and $\kappa>0$
such that
$T^{-1}\int_{t-T}^tr(s)\,ds\le-\kappa$ for every $t$.
Define
\begin{align*}
\sigma(t):=-\int_0^T\left(1-\frac{s}{T}\right)r(t-s)\,ds.
\end{align*}
The function $\sigma$ is bounded and Lipschitz and
\begin{align*}
\sigma'(t)+r(t)=\frac1T\int_{t-T}^tr(s)\,ds\le-\kappa.
\end{align*}
Substitution in the linearized equation yields
$\mathcal D\mathbf V_+\le-\kappa\mathbf V_+$.
The remaining statements follow from
\eqref{eq:v2-harnack} and \eqref{eq:v2-bounded-tangent}.
\end{proof}

\begin{proof}[Proof of Theorem~\ref{thm:critical_ap_tangent}]
The linear construction is contained in
Lemma~\ref{lem:AP-normalized-bundle}. The passage to a nonlinear
front uses separately the positivity and cooperativity of
\eqref{Fnonlinear}, the KPP upper bound and positivity along the
diagonal in \eqref{KPP}, and the near-zero estimate
\eqref{near_zerp}. We follow the sub- and supersolution arrangement of Lemma~\ref{critical supersub},
with the constant-speed Floquet modes replaced by the three functions
in Lemma~\ref{lem:v2-critical-modes}. We include the estimates that
are needed when the coefficients are not periodic in time.

Fix $K_0>0$ sufficiently large and let $\varepsilon>0$. Put $\mathbf w_-=\varepsilon\mathbf V_+
+\varepsilon(\mathbf Z_*-K_0\mathbf V_*),~\mathbf w_+=M\mathbf V_*+\varepsilon\mathbf Z_*.$ The uniform bounds in \eqref{eq:v2-harnack},
\eqref{eq:v2-bounded-tangent} and \eqref{eq:v2-strict-mode}
show that there are constants $C_1,C_2>0$ and $R_0>1$,
independent of $\varepsilon$ for $0<\varepsilon\le1$, such that
\begin{align}\label{eq:v2-tail}
0< w_{-,i}(t,x)\le C_1\varepsilon(1+\xi)e^{-\lambda_*\xi},~
V_{+,i}(t,x)\ge C_2^{-1}e^{-\lambda_+\xi}
\end{align}
for $\xi\ge R_0$. Increasing $K_0$ and choosing $M$ sufficiently
large, we arrange
\begin{align*}
\mathbf w_-\le\mathbf0\ \text{for }0\le\xi\le1,
\quad
\mathbf w_-\le\mathbf w_+\ \text{for }\xi\ge0,
\quad
\mathbf w_+(t,x)\ge\mathbf1\ \text{if }\xi=0.
\end{align*}
The constants are chosen in this order. Write $\mathbf Z_*=e^{-\lambda_*\xi}(\xi\mathbf p_{\lambda_*}+\mathbf r_*)$ with $\|\mathbf r_*\|_\infty\leq C_*$ and $\mathbf p_{\lambda_*}\geq m_*\mathbf1$. Taking $R_0>2(C_*+K_0M_*)/m_*$ makes $\mathbf Z_*-K_0\mathbf V_*\geq\frac12m_*\xi e^{-\lambda_*\xi}\mathbf1$ on $\xi\geq R_0$. Moreover, the cancellation\[\mathbf w_+-\mathbf w_-=(M+\varepsilon K_0)\mathbf V_*-\varepsilon\mathbf V_+\]shows that the ordering follows from the bounded ratio $\mathbf V_+/\mathbf V_*$ on $\xi\geq0$. A sufficiently large $M$ also ensures $\mathbf w_+(\xi=0)\geq\mathbf1$. After fixing $K_0,R_0,M$, we reduce $\varepsilon$ without changing them.
Using \eqref{eq:v2-neutral-modes}--\eqref{eq:v2-strict-mode},
we obtain $\mathcal D\mathbf w_-\le-\varepsilon\kappa\mathbf V_+,~ \mathcal D\mathbf w_+=\mathbf0.$ Since $\delta<\nu\lambda_*$, the estimates in
\eqref{eq:v2-tail} imply, for $\xi\ge R_0$,
\begin{align}\label{eq:v2-nonlinear-bound}
\frac{K(w_{-,i}^+)^{1+\nu}}{\varepsilon V_{+,i}}
&\le C\varepsilon^\nu(1+\xi)^{1+\nu}
 e^{-(\nu\lambda_*-\delta)\xi}
\le C'\varepsilon^\nu.
\end{align}
On $0\le\xi\le R_0$ the same estimate follows
from the uniform boundedness of the normalized sections and their
parameter derivatives, together with the positive lower bound of
$e^{-\lambda_+\xi+\sigma(t)}p_{\lambda_+,i}$.
Thus, with a larger constant the estimate holds for every
$\xi\ge0$ at which $w_{-,i}>0$. More precisely, the coefficient of
$e^{-\lambda_*\xi}$ in $\mathbf Z_*-K_0\mathbf V_*$ has
at most linear growth in $\xi$, uniformly in $(t,x)$.
Consequently, for some $C_3>0$ independent of $\varepsilon$,
\begin{align}\label{eq:v15-global-smallness}
\sup_{\substack{t\in\mathbb R,\ x\in\mathbb R^N\\\xi\ge0}}
\max_{i=1,2}w_{-,i}^+(t,x)\le C_3\varepsilon.
\end{align}
Write $\delta_{\rm nl}>0$ for the neighborhood radius in
\eqref{near_zerp}, to distinguish it from the spectral increment
$\delta$ in \eqref{eq:v2-gap}. We decrease $\varepsilon$ so that
$C'\varepsilon^\nu\le\kappa$ and
$C_3\varepsilon<\delta_{\rm nl}/2$. Whenever
$w_{-,i}>0$, the near-zero inequality and
\eqref{eq:v2-nonlinear-bound} imply, provided both components
of the chosen barrier remain below $\delta_{\rm nl}$,
\begin{align}\label{eq:v15-active-subsolution}
(\mathcal D\mathbf w_-)_i
+K(w_{-,i})^{1+\nu}\le
-\varepsilon\kappa V_{+,i}
+K(w_{-,i})^{1+\nu}\le0.
\end{align}
The comparison of the other component with $w_{-,j}$ is
performed using the nonnegative off-diagonal coefficient
$\mu_{ij}$. In particular, if the other component is replaced
by a larger nonnegative barrier value, the reaction inequality
becomes no less favorable. Finally, the KPP upper bound in
\eqref{KPP} proves the supersolution inequality for
$\mathbf w_+$ at points where the selected component lies in
$(0,1)$; no value of $F_i$ outside $[0,1]^2$ is needed.

Set $\overline{\mathbf u}=\mathbf1$ for $\xi\le0$ and
$\overline{\mathbf u}=\min\{\mathbf1,\mathbf w_+\}$ for $\xi>0$,
with componentwise minimum. By cooperativity, this is a generalized
supersolution. Choose $R\ge R_0$ so that
$\mathbf w_->\mathbf0$ on $\xi\ge R$, and choose
$0<\omega<\min\{\delta,\delta_{\rm nl}/2\}$ smaller than
every component of $\mathbf w_-$ on $\xi=R$ and of
$\overline{\mathbf u}$ on $0\le\xi\le R$.
By \eqref{eq:v15-global-smallness}, we may also require
$\max\{\omega,C_3\varepsilon\}<\delta_{\rm nl}$.
The lower bounds used to select $\omega$ do not depend on $t$:
on $\xi=R$ this follows from the uniform positivity of
$\mathbf p_{\lambda_*}$ and the leading $R$-term of
$\mathbf Z_*$; on the compact strip $0\le\xi\le R$ it follows
from the bounds for $\mathbf w_+$ and the value $\mathbf1$.
Define
\begin{align*}
\underline{\mathbf u}(t,x)=
\begin{cases}
\omega\mathbf1,&\xi<0,\\
\max\{\omega\mathbf1,\mathbf w_-(t,x)\},&0\le\xi<R,\\
\mathbf w_-(t,x),&\xi\ge R.
\end{cases}
\end{align*}
All maxima are componentwise. We prove the inequalities
at each selected component and at the interfaces. If
$\underline u_i=\omega$, the other component is at least
$\omega$; hence cooperativity and \eqref{KPP} imply
\begin{align*}
F_i(t,x,\underline{\mathbf u})
\ge F_i(t,x,\omega,\omega)>0
=\partial_t\omega-\operatorname{Tr}({\bf A}^iD^2\omega)
+{\bf q}_i\cdot D\omega.
\end{align*}
If $\underline u_i=w_{-,i}>0$, its other component lies in
$[0,\delta_{\rm nl})$ by \eqref{eq:v15-global-smallness}
and the choice of $\omega$, and is at least $w_{-,j}$.
Using \eqref{near_zerp} at the actual nonnegative barrier
vector, and then $\mu_{ij}\ge0$, we obtain
\begin{align*}
F_i(t,x,\underline{\mathbf u})
&\ge\sum_{j=1}^2\mu_{ij}\underline u_j
-K(w_{-,i})^{1+\nu}\ge\sum_{j=1}^2\mu_{ij}w_{-,j}
-K(w_{-,i})^{1+\nu}.
\end{align*}
The last expression is at least the scalar diffusion--transport
operator applied to $w_{-,i}$ by
\eqref{eq:v15-active-subsolution}. Thus the differential
subsolution inequality holds for whichever smooth branch is
active in component $i$.
The interfaces at $\xi=0$ and $\xi=R$ produce no surface
measure: $\mathbf w_-<\omega\mathbf1$ on a uniform strip
$0\le\xi\le1$ by the initial choice of $K_0$ and the
subsequent choice of $\omega$, whereas
$\mathbf w_->\omega\mathbf1$ on a uniform strip about
$\xi=R$ by uniform continuity and the strictly positive margin
in the choice of $\omega$. Thus the same smooth branch is active
on both sides of either prescribed interface. At an interior
switching set, Lemma~\ref{lem:cooperative-gluing}
applies to the two scalar branches. In particular, the convexity
term in the chain rule has the favorable sign and no positive
interface measure appears. The preceding componentwise estimates, evaluated at
the actual vector $\underline{\mathbf u}$, then preserve its
nonlinear cooperative right-hand side. Since the branches are
continuous and coincide at the switching interface, no jump
term is introduced in the solution itself. If the branch
$\overline u_i=w_{+,i}<1$ is active, all components of
$\overline{\mathbf u}$ lie in $[0,1]$ and
$\overline u_j\le w_{+,j}$. The KPP inequality and
$\mu_{ij}>0$ imply
\begin{align*}
F_i(t,x,\overline{\mathbf u})
\le\sum_{j=1}^2\mu_{ij}\overline u_j
\le\sum_{j=1}^2\mu_{ij}w_{+,j}
=\partial_t w_{+,i}
-\operatorname{Tr}({\bf A}^iD^2w_{+,i})
+{\bf q}_i\cdot Dw_{+,i}.
\end{align*}
If instead $\overline u_i=1$, the other component is at most
$1$ and cooperativity implies
$F_i(t,x,1,\overline u_j)\le F_i(t,x,1,1)=0$.
Thus each smooth active branch satisfies the supersolution
inequality without evaluating $\mathbf F$ outside $[0,1]^2$.
Lemma~\ref{lem:cooperative-gluing}, applied with reversed
inequalities, then proves the supersolution inequality across
the switching set. The preceding estimates verify its
active-branch hypotheses against the actual clipped vector. Accordingly,
$\underline{\mathbf u}$ and $\overline{\mathbf u}$ form an ordered
subsolution--supersolution pair and
\begin{align*}
\mathbf0<\underline{\mathbf u}\le
\overline{\mathbf u}\le\mathbf1.
\end{align*}

Solving the Cauchy problem at times $-n$ with initial values
$\underline{\mathbf u}(-n,\cdot)$ and using comparison,
parabolic interior estimates and a diagonal argument, we obtain an
entire solution $\mathbf u_*$ between the two barriers.
The upper barrier converges uniformly to zero when $\xi\to+\infty$.
On the other hand, $\mathbf u_*\ge\omega\mathbf1$ for
$\xi\le R$. Suppose that convergence to $\mathbf1$ as
$\xi\to-\infty$ is not uniform. Then there are $\eta>0$,
$t_n\in\mathbb R$ and $x_n\in\mathbb R^N$ such that
$\xi_n=x_n\cdot e-S(t_n)\to-\infty$ and
$\min_i u_{*,i}(t_n,x_n)\le1-\eta$.
Translate the solutions by $(t_n,x_n)$; spatial periodicity,
compactness of the temporal translation hull and interior parabolic
estimates provide a subsequence converging locally to an entire
solution $\mathbf U$. The bound $c\in L^\infty(\mathbb R)$ ensures
$\xi(t_n+t,x_n+x)\to-\infty$ uniformly on each compact
set of $(t,x)$. Hence
$\omega\mathbf1\le\mathbf U\le\mathbf1$ everywhere and
$\min_i U_i(0,0)\le1-\eta$.
Let $s=\inf_{t,x,i}U_i(t,x)\in[\omega,1-\eta]$.
Taking a further sequence of translates approaching this infimum,
and using the same hull compactness, produces a positive entire
solution $\mathbf W\ge s\mathbf1$ and an index $i$ such that
$W_i(0,0)=s$. At this space--time minimum,
$\partial_tW_i(0,0)=0$, $DW_i(0,0)=0$, and
$D^2W_i(0,0)$ is nonnegative definite. The equation therefore
gives
\[
0\ge F_i(0,0,\mathbf W(0,0))
\ge F_i(0,0,s,s)>0,
\]
where the second inequality uses cooperativity and the last one
uses \eqref{KPP} for $0<s<1$. This contradiction proves the
uniform convergence to $\mathbf1$ and hence both limits of
Definition~\ref{wavedef}.
Finally, $c_*\in L^\infty(\R)$ is almost periodic and
\begin{align*}
\langle c_*\rangle
=\frac{\langle b_{\lambda_*}\rangle}{\lambda_*}
=\frac{\Gamma(\lambda_*)}{\lambda_*}=c^*(e).
\end{align*}
\end{proof}

\begin{corollary}\label{cor:v2-quasiperiodic}
Assume the hypotheses of Theorem~\ref{thm:critical_ap_tangent},
except for \eqref{eq:v2-bounded-tangent}, and retain
\eqref{eq:v2-gap}. Suppose that, for some frequency vector
$\omega\in\R^d$ and a function $\mathfrak h\in C^s(\mathbb T^d)$,
\begin{align*}
h_*(t)=\mathfrak h(\omega t),\qquad
\int_{\mathbb T^d}\mathfrak h(\theta)\,d\theta=0,
\end{align*}
where $s>d+\tau$ and
\begin{align}\label{eq:v2-diophantine}
|k\cdot\omega|\ge\gamma |k|^{-\tau}
\qquad(k\in\mathbb Z^d\setminus\{0\})
\end{align}
for some $\gamma>0$ and $\tau\ge0$.
Then \eqref{eq:main} admits a critical generalized transition
front with speed \eqref{eq:v2-critical-speed}. In particular,
no time-periodicity assumption is needed.
\end{corollary}

\begin{proof}
With $\mathbb T^d=\mathbb R^d/\mathbb Z^d$, we use the Fourier
basis $e^{2\pi\mathrm i k\cdot\theta}$. The assumption $s>d+\tau$
and the Fourier coefficient bound
$|\widehat{\mathfrak h}(k)|\le C(1+|k|)^{-s}$ imply
\begin{align*}
\sum_{k\in\mathbb Z^d\setminus\{0\}}
\frac{|\widehat{\mathfrak h}(k)|}{|k\cdot\omega|}
<\infty.
\end{align*}
Consequently the series
\begin{align*}
\mathfrak H(\theta)=
\sum_{k\in\mathbb Z^d\setminus\{0\}}
\frac{\widehat{\mathfrak h}(k)}{2\pi\mathrm i k\cdot\omega}
 e^{2\pi\mathrm i k\cdot\theta}
\end{align*}
converges absolutely and uniformly. Integrating the Fourier
series of $\mathfrak h$ along the orbit $t\mapsto\omega t$,
first for trigonometric polynomials and then by uniform
convergence, shows that
\begin{align*}
\int_0^t h_*(r)\,dr
=\mathfrak H(\omega t)-\mathfrak H(0).
\end{align*}
Thus \eqref{eq:v2-bounded-tangent} holds and
Theorem~\ref{thm:critical_ap_tangent} applies.
\end{proof}

\begin{remark}
The condition on $\mathfrak h$ is a sufficient regularity
condition on $h_*$, not a claim that all smooth
quasiperiodic coefficients determine a function $h_*$ of that
regularity. Such a deduction requires separate estimates for
parameter and torus-phase derivatives of the normalized principal
section.
\end{remark}

\begin{proof}[Proof of Remark~\ref{rem:v2-examples-intro}]
At an interior minimizing rate, differentiability of $\Gamma$ implies
$\Gamma'(\lambda_*)=\Gamma(\lambda_*)/\lambda_*$. By the uniform
$C^1$ dependence established in Lemma~\ref{lem:AP-normalized-bundle},
we have $\langle h_*\rangle=0$. If the coefficients are periodic,
then $b_{\lambda_*}$ and its $\lambda$-derivative are periodic.
Thus the primitive of $h_*$ is periodic and bounded.

Suppose next that $b_\lambda(t)=a(t)\gamma(\lambda)+\lambda\beta(t)$
near $\lambda_*$. Then
\begin{align*}
h_*(t)=a(t)\left(\gamma'(\lambda_*)-
\frac{\gamma(\lambda_*)}{\lambda_*}\right)=0,
\end{align*}
because $\lambda_*$ is an interior minimizer of
$\gamma(\lambda)/\lambda$. For
$\gamma(\lambda)=d\lambda^2+r$ with $d,r>0$,
$\gamma(\lambda)/\lambda=d\lambda+r/\lambda$ has its
unique minimum at $\sqrt{r/d}$, with strict separation to the
right. Since $\Gamma(\lambda)=\langle a\rangle\gamma(\lambda)
+\lambda\langle\beta\rangle$, its minimizing rate is the same,
and the claimed mean speed follows.

If $h_*$ is a finite trigonometric sum with nonzero frequencies,
termwise integration yields
\begin{align*}
H_*(t)=\sum_{j=1}^m\left[
\frac{A_j}{\omega_j}\sin(\omega_jt)
+\frac{B_j}{\omega_j}\big(1-\cos(\omega_jt)\big)\right],
\end{align*}
which is bounded regardless of rational relations between the
frequencies. If the Fourier coefficients of an almost-periodic
$h_*$ satisfy the summability condition stated after the theorem,
then
\begin{align*}
H_*(t)=\sum_{k\ne0}\frac{\widehat h_*(k)}
{2\pi\mathrm i k\cdot\omega}\big(e^{2\pi\mathrm i k\cdot\omega t}-1\big)
\end{align*}
converges uniformly in $t$. Absolute summability of the original
Fourier coefficients also permits termwise integration from
finite partial sums, and hence identifies this series with the
primitive of $h_*$. This proves the asserted sufficient condition.

Finally, consider $h(t)=\displaystyle\sum_{n=1}^\infty 2^{-n}\cos(4^{-n}t).$ The defining series and all its differentiated series converge
uniformly on $\R$. The function $h$ is therefore smooth and
Bohr almost periodic. It has zero uniform mean, because each
summand has zero mean and the series converges uniformly. For
fixed $t$ its primitive is
\begin{align*}
H(t)=\sum_{n=1}^\infty2^n\sin(4^{-n}t),
\end{align*}
where the series converges absolutely for each $t$.
At $t_j=(\pi/2)4^j$, all terms with $n<j$ vanish, the term $n=j$
equals $2^j$, and every term with $n>j$ is nonnegative. Thus
$H(t_j)\ge 2^j\to\infty$.
\end{proof}

\endgroup

\begingroup
For coefficients obtained from an autonomous system by time-dependent changes of variables, the existence statement in Theorem~\ref{minimalspeed} extends beyond temporal
periodicity without requiring a bounded almost-periodic principal eigenfunction in Remark \ref{APremarkcrit}.
We consider a class for which the minimal speed and a critical front can both
be identified. The distinction between the existence of a critical front and
the existence of an almost-periodic eigenfunction is essential.

Next, let ${\bf A}^i_0(y)$, ${\bf q}_{i,0}(y)$ and
${\bf F}_0(y,{\bf u})$ be independent of time and ${\bf L}$-periodic in $y$.
Assume that they satisfy \eqref{assumeq}--\eqref{near_zerp} and the
regularity conditions used for the time-independent case of
Theorem~\ref{minimalspeed}. In the statements below, recall $\mathcal C=[{\bf0},{\bf L}]$.
We write $c_0^*(e)$ for the minimal speed of this autonomous system.

\begin{theorem}
\label{thm:critical_ap_conjugacy}
Suppose that $a\in C_b^1(\R)$ and
$B\in C_b^2(\R;\R^N)$ satisfy
\begin{align}\label{eq:critical_ap_assumptions}
0<a_-\le a(t)\le a_+<\infty,
\qquad a,\ B,\ B'\text{ are Bohr almost periodic on }\R.
\end{align}
For $i=1,2$, define
\begin{align}\label{eq:critical_ap_coefficients}
\begin{split}
{\bf A}^i(t,x)&=a(t){\bf A}^i_0(x-B(t)),\\
{\bf q}_i(t,x)&=a(t){\bf q}_{i,0}(x-B(t))+B'(t),\\
F_i(t,x,{\bf u})&=a(t)F_{i,0}(x-B(t),{\bf u}).
\end{split}
\end{align}
Then these coefficients are spatially periodic and almost periodic in time,
uniformly in their respective spatial and state variables. For every
$e\in\mathbb S^{N-1}$, system \eqref{eq:main} admits a generalized
transition front $({\bf u}_e,c_e)$ with
\begin{align}\label{eq:critical_ap_speed}
c_e(t)=a(t)c_0^*(e)+B'(t)\cdot e,
\qquad
\langle c_e\rangle=\langle a\rangle c_0^*(e)=c^*(e).
\end{align}
In particular, this front attains the minimal least mean speed of the
almost-periodic system. No boundedness condition is imposed on
$\int_0^t(a(s)-\langle a\rangle)\,ds$.
\end{theorem}

\begin{proof}
Set $\tau(t)=\int_0^t a(r)\,dr,~ y=x-B(t).$ Since $a\ge a_->0$, $\tau$ is a strictly increasing bijection of $\R$.
If ${\bf U}$ is a solution of the autonomous system, define
${\bf u}(t,x)={\bf U}(\tau(t),x-B(t))$. Direct differentiation implies
\begin{align*}
\partial_t u_i
 &=a(t)\partial_s U_i(\tau(t),y)-B'(t)\cdot D_yU_i(\tau(t),y),~
D_xu_i=D_yU_i,~ D_x^2u_i=D_y^2U_i.
\end{align*}
Consequently,
\begin{align*}
&\partial_tu_i-\operatorname{Tr}({\bf A}^i(t,x)D_x^2u_i)
       +{\bf q}_i(t,x)\cdot D_xu_i=a(t)\big[\partial_sU_i-
\operatorname{Tr}({\bf A}^i_0(y)D_y^2U_i)
+{\bf q}_{i,0}(y)\cdot D_yU_i\big]
=a(t)F_{i,0}(y,{\bf U}).
\end{align*}
The coefficients in \eqref{eq:critical_ap_coefficients} inherit their
spatial periodicity and their structural hypotheses from the autonomous
coefficients. Their temporal almost periodicity follows from the uniform
continuity of the periodic spatial coefficients and the relative compactness
of the translation families of $a,B,B'$.

By the time-independent case of Theorem~\ref{minimalspeed}, there is a
critical front ${\bf U}_e$ of the autonomous system whose interface is
$s\mapsto c_0^*(e)s$. Consider
\begin{align*}
{\bf u}_e(t,x)={\bf U}_e(\tau(t),x-B(t)),\qquad
S_e(t)=c_0^*(e)\tau(t)+(B(t)-B(0))\cdot e.
\end{align*}
The function $S_e$ is locally absolutely continuous and
$S_e'=c_e\in L^\infty(\R)$. Furthermore, for
$x\cdot e\to\pm\infty$ the points
$x+S_e(t)e-B(t)$ have signed longitudinal distance from the autonomous
interface equal to
\begin{align*}
\big(x+S_e(t)e-B(t)\big)\cdot e-c_0^*(e)\tau(t)
=x\cdot e-B(0)\cdot e.
\end{align*}
Thus both front limits are uniform in $t$, exactly as for ${\bf U}_e$.
Since $B$ is bounded and $a$ has a uniform mean, the speed in
\eqref{eq:critical_ap_speed} has the indicated uniform mean.

It remains to verify minimality. Let $\gamma_0(\lambda)$ and
$\boldsymbol\phi^0_\lambda>0$ be the principal eigenvalue and an
${\bf L}$-periodic principal eigenfunction of the stationary tilted
linearized operator in direction $e$. Its positive entire linearized
solution is
$e^{-\lambda y\cdot e+\gamma_0(\lambda)s}
\boldsymbol\phi^0_\lambda(y)$.
The conjugacy above transforms this solution into
\begin{align}\label{eq:critical_ap_tilted_solution}
e^{-\lambda x\cdot e}\boldsymbol\eta_\lambda(t,x),\text{ where }
\boldsymbol\eta_\lambda(t,x)=
 e^{\lambda B(t)\cdot e+\gamma_0(\lambda)\tau(t)}
 \boldsymbol\phi^0_\lambda(x-B(t)).
\end{align}
Spatial translation leaves the supremum norm of the periodic
$\boldsymbol\phi^0_\lambda$ unchanged. Hence, for every $T>0$,
\begin{align*}
\frac1T\log\frac{\|\boldsymbol\eta_\lambda(t+T,\cdot)\|_\infty}
 {\|\boldsymbol\eta_\lambda(t,\cdot)\|_\infty}
&=\gamma_0(\lambda)\frac1T\int_t^{t+T}a(r)\,dr
+\lambda e\cdot\frac{B(t+T)-B(t)}T.
\end{align*}
The right-hand side tends to
$\langle a\rangle\gamma_0(\lambda)$ uniformly in $t$.
The principal growth characterization in
the spectral characterization of growth in Section~\ref{section 5} therefore yields
\begin{align*}
c^*(e)=\inf_{\lambda>0}
\frac{\langle a\rangle\gamma_0(\lambda)}\lambda
=\langle a\rangle c_0^*(e).
\end{align*}
This proves that the constructed front has the minimal least mean speed.
\end{proof}

\begin{corollary}\label{cor:v2-conjugate-tangent}
For the coefficients obtained from an autonomous periodic system by
\eqref{eq:critical_ap_coefficients}, assume that the stationary
minimal speed is attained at an interior minimizing rate $\lambda_*$
and that \eqref{eq:v2-gap} holds. Then the tangent function in
\eqref{eq:v2-tangent} vanishes identically. In particular,
Theorem~\ref{thm:critical_ap_tangent} applies, even when the
centered primitive of the temporal multiplier $a$ is unbounded.
\end{corollary}
\begin{proof}
Normalize the stationary principal eigenfunction
$\boldsymbol\phi^0_\lambda$ by $\mathcal J(\boldsymbol\phi^0_\lambda)=1$.
The transformed positive linearized solution in
\eqref{eq:critical_ap_tilted_solution} has normalized section
$\mathbf p_\lambda(t,x)=\boldsymbol\phi^0_\lambda(x-B(t))$.
Consequently
\begin{align*}
b_\lambda(t)=a(t)\gamma_0(\lambda)+\lambda B'(t)\cdot e.
\end{align*}
Since $\gamma_0'(\lambda_*)=\gamma_0(\lambda_*)/\lambda_*$
at an interior differentiable minimizer,
\begin{align*}
\left.\partial_\lambda b_\lambda(t)\right|_{\lambda_*}
-\frac{b_{\lambda_*}(t)}{\lambda_*}=0
\end{align*}
for every $t$. Thus $H_*\equiv0$, as asserted.
\end{proof}
\begin{proposition}
\label{prop:critical_ap_coboundary}
Under the hypotheses of Theorem~\ref{thm:critical_ap_conjugacy}, set
\begin{align*}
G(t)=\int_0^t\big(a(s)-\langle a\rangle\big)\,ds.
\end{align*}
Fix $\lambda>0$ with $\gamma_0(\lambda)\ne0$. Assume in addition that
the stationary tilted operator and its adjoint have strictly positive,
periodic principal eigenfunctions, with sufficient spatial regularity
for periodic integration by parts. Then the following statements are equivalent:
\begin{enumerate}[label=(\roman*)]
\item $G$ is bounded on $\R$;
\item the almost-periodic tilted operator has a positive,
${\bf L}$-periodic eigenfunction $\boldsymbol\Phi_\lambda(t,x)$,
bounded above and bounded away from zero on $\R\times\R^N$, satisfying
\begin{align}\label{eq:critical_ap_eigenpair}
{\bf P}_\lambda\boldsymbol\Phi_\lambda
=-\langle a\rangle\gamma_0(\lambda)
\boldsymbol\Phi_\lambda;
\end{align}
\item such an eigenfunction can be chosen Bohr almost periodic in $t$,
uniformly in $x$.
\end{enumerate}
Thus, the existence of a critical front in
Theorem~\ref{thm:critical_ap_conjugacy} does not require the
almost-periodic eigenfunction postulated in Remark~\ref{APremarkcrit}.
\end{proposition}

\begin{proof}
Suppose first that $G$ is bounded. The Bohl--Bohr theorem for
real-valued Bohr almost-periodic functions (see \cite{Fink1974}) implies
that $G$ is almost periodic. Define
\begin{align}\label{eq:critical_ap_eigenfunction}
\boldsymbol\Phi_\lambda(t,x)=
\exp\big(\gamma_0(\lambda)G(t)+\lambda B(t)\cdot e\big)
\boldsymbol\phi^0_\lambda(x-B(t)).
\end{align}
This function is almost periodic in $t$, uniformly in $x$, is positive,
and has positive uniform lower and finite upper bounds. Differentiating
\eqref{eq:critical_ap_eigenfunction}, or using
\eqref{eq:critical_ap_tilted_solution}, proves
\eqref{eq:critical_ap_eigenpair}. This establishes (i)$\Rightarrow$(iii),
while (iii)$\Rightarrow$(ii) is immediate.

Conversely, suppose (ii) holds. Let
$\boldsymbol\phi^{0,*}_\lambda>0$ be a periodic principal
eigenfunction of the adjoint stationary tilted operator and define
\begin{align*}
Z(t)=e^{-\lambda B(t)\cdot e}
\sum_{i=1}^2\int_{\mathcal C}
\phi^{0,*}_{\lambda,i}(y)
\Phi_{\lambda,i}(t,y+B(t))\,dy.
\end{align*}
The bounds on $B$ and $\boldsymbol\Phi_\lambda$ imply
$0<\inf_t Z(t)\le\sup_tZ(t)<\infty$.
Upon differentiating $Z$, translating the equation
\eqref{eq:critical_ap_eigenpair} back to the stationary coordinates,
and integrating by parts on $\mathcal C$, the adjoint eigenvalue
equation implies
\begin{align*}
Z'(t)=\gamma_0(\lambda)
\big(a(t)-\langle a\rangle\big)Z(t).
\end{align*}
Therefore
$Z(t)=Z(0)\exp(\gamma_0(\lambda)G(t))$.
As $\gamma_0(\lambda)\ne0$, the uniform upper and lower bounds on
$Z$ imply that $G$ is bounded. This proves (ii)$\Rightarrow$(i).
\end{proof}

\begin{corollary}
\label{cor:critical_ap_phase}
Assume the hypotheses of Theorem~\ref{thm:critical_ap_conjugacy} and
$c_0^*(e)\ne0$. For the front constructed in that theorem, let
$S_e(t)=\int_0^tc_e(r)\,dr$. Then
\begin{align*}
\sup_{t\in\R}|S_e(t)-c^*(e)t|<\infty
\text{ if and only if }
\sup_{t\in\R}|G(t)|<\infty.
\end{align*}
If in addition the hypotheses of
Proposition~\ref{prop:critical_ap_coboundary} hold at a minimizing
$\lambda_* >0$, then these conditions are equivalent to the existence
of a uniformly positive, bounded almost-periodic eigenfunction at
$\lambda_*$. If $G$ is unbounded, the constructed front still has
minimal mean speed, but its interface has an unbounded displacement
from the linear function $t\mapsto c^*(e)t$.
\end{corollary}

\begin{proof}
By \eqref{eq:critical_ap_speed}, $S_e(t)-c^*(e)t
=c_0^*(e)G(t)+(B(t)-B(0))\cdot e.$ The second term is bounded and $c_0^*(e)\ne0$, proving the first
equivalence. At a minimizing $\lambda_*$,
$\gamma_0(\lambda_*)=\lambda_*c_0^*(e)\ne0$; the second equivalence
follows from Proposition~\ref{prop:critical_ap_coboundary}.
\end{proof}

\begin{corollary}
\label{cor:critical_ap_example}
For $N\ge1$, set
\begin{align*}
a(t):=2+\sum_{n=1}^\infty2^{-n}\cos(4^{-n}t)
\end{align*}
and consider the cooperative KPP system
\begin{align}\label{eq:critical_ap_explicit}
\begin{cases}
\partial_tu_1-a(t)\Delta u_1
=a(t)(2u_1+u_2-3u_1^2),\\
\partial_tu_2-a(t)\Delta u_2
=a(t)(u_1+2u_2-3u_2^2).
\end{cases}
\end{align}
For every $e\in\mathbb S^{N-1}$, this system possesses a critical
transition front with speed
\begin{align*}
c(t)=2\sqrt3\,a(t),\qquad
\langle c\rangle=c^*(e)=4\sqrt3.
\end{align*}
Nevertheless, for every $\lambda>0$ there is no positive spatially
periodic eigenfunction of ${\bf P}_\lambda$ at the generalized
principal eigenvalue which is simultaneously bounded above and
bounded away from zero on $\R\times\R^N$.
\end{corollary}

\begin{proof}
The series defining $a$ and each of its differentiated series converge
uniformly. Thus $a$ is smooth and Bohr almost periodic,
$1\le a(t)\le3$, and $\langle a\rangle=2$.
The diagonal subspace $u_1=u_2$ is invariant. Let $W$ be the classical
critical Fisher--KPP profile, satisfying
\begin{align*}
W''+2\sqrt3\,W'+3W(1-W)=0,
\quad W(-\infty)=1,\quad W(+\infty)=0.
\end{align*}
For
$\tau(t)=\int_0^t a(s)\,ds$, the function
\begin{align*}
{\bf u}(t,x)=
\big(W(x\cdot e-2\sqrt3\,\tau(t)),
W(x\cdot e-2\sqrt3\,\tau(t))\big)
\end{align*}
is a critical entire front of \eqref{eq:critical_ap_explicit}.
Indeed, the stationary linearized matrix is
$\left(\begin{matrix}2&1\\1&2\end{matrix}\right)$,
whose Perron eigenvalue is $3$; hence the stationary minimal speed is
$2\sqrt3$, and Theorem~\ref{thm:critical_ap_conjugacy} applies.

On the other hand,
\begin{align*}
G(t)=\int_0^t(a(s)-2)\,ds
=\sum_{n=1}^\infty2^n\sin(4^{-n}t).
\end{align*}
The series on the right converges for each fixed $t$, since its terms
are bounded in absolute value by $2^{-n}|t|$ for large $n$.
For $t_k=(\pi/2)4^k$, all the terms with $n<k$ vanish, the term
with $n=k$ equals $2^k$, and every term with $n>k$ is nonnegative.
Therefore $G(t_k)\ge2^k\to\infty$.
For this system $\gamma_0(\lambda)=\lambda^2+3>0$, and the
nonexistence assertion follows from
Proposition~\ref{prop:critical_ap_coboundary}.
\end{proof}

\begin{corollary}
\label{cor:critical_ap_periodic_base}
Let the reference coefficients ${\bf A}^i_0(s,y)$,
${\bf q}_{i,0}(s,y)$ and ${\bf F}_0(s,y,{\bf u})$ be
$T_0$-periodic in $s$, spatially periodic in $y$, and satisfy the
hypotheses of the time-periodic critical-front result in
Theorem~\ref{minimalspeed}. Suppose \eqref{eq:critical_ap_assumptions}
holds and that $G(t)=\int_0^t(a(r)-\langle a\rangle)\,dr$
is bounded. Define the new coefficients by
\begin{align*}
{\bf A}^i(t,x)&=a(t){\bf A}^i_0(\tau(t),x-B(t)),\\
{\bf q}_i(t,x)&=a(t){\bf q}_{i,0}(\tau(t),x-B(t))+B'(t),\\
F_i(t,x,{\bf u})&=a(t)F_{i,0}(\tau(t),x-B(t),{\bf u}).
\end{align*}
Then these coefficients are Bohr almost periodic in $t$ and the new
system admits a critical generalized transition front with mean speed
$\langle a\rangle c_0^*(e)$, where $c_0^*(e)$ is the minimal speed of
the reference time-periodic system.
\end{corollary}

\begin{proof}
By the Bohl--Bohr theorem, $G$ is almost periodic and
$\tau(t)=\langle a\rangle t+G(t)$. Composition with continuous
periodic coefficients preserves almost periodicity, as follows by
uniform approximation with trigonometric polynomials; the same holds
for spatial translation by $B(t)$. The change of variables in
Theorem~\ref{thm:critical_ap_conjugacy} maps the reference front to
one with speed $a(t)c_0^*(e)+B'(t)\cdot e$ and uniform mean
$\langle a\rangle c_0^*(e)$. The associated Floquet exponent is
multiplied by $\langle a\rangle$, and Theorem~\ref{minimalspeed}
identifies this speed as critical.
\end{proof}
\endgroup

\section{Small- and Large-Diffusion Asymptotics}\label{sec:memoirs-diffusion-regimes}
The large-advection results below vary the amplitude of the advective field while the diffusion operators remain fixed. Variational spectral descriptions of minimal KPP front speeds in periodic shears are classical; see, for example, \cite{BerestyckiHamelNadirashvili2005,NolenXin2005}. We keep the advective field fixed and vary a common diffusion rate in both components. This distinction is essential: the two limits are not obtained by simply replacing the large-advection parameter in Proposition~\ref{prop:memoirs-second-asympt}. The setting below is a subclass of \eqref{eq:main} in which the two diffusion matrices coincide and are constant in space and time.

Let $d>0$ be the diffusion rate and assume ${\bf A}^1_d={\bf A}^2_d=d{\bf A}_0$, where ${\bf A}_0$ is a fixed symmetric positive definite matrix. We first consider shears for which, in the longitudinal and transverse coordinates $(s,z)$, the constant-coefficient differential operator has the diagonal form $a\partial_{ss}+D\partial_{zz}$, with $a,D>0$. This alignment condition holds in both explicit examples of this section when the physical diffusion matrices are $dI_2$: the transformed coefficients are $(a,D)=(1,1)$ in Case~1 and $(a,D)=(1,2)$ in Case~2. No conclusion below is asserted for a general non-diagonal constant matrix in these coordinates without additional analysis.

Fix $r,\gamma,h>0$ and $k=2\pi$. The associated cooperative system is
\begin{equation}\label{eq:memoirs-diffusion-model}
\begin{cases}
(u_1)_t=d(a(u_1)_{ss}+D(u_1)_{zz})-h\sin(kz)(u_1)_s
 +(r-\gamma)u_1+\gamma u_2-ru_1^2,\\
(u_2)_t=d(a(u_2)_{ss}+D(u_2)_{zz})-h\sin(kz)(u_2)_s
 +\gamma u_1+(r-\gamma)u_2-ru_2^2,
\end{cases}
\end{equation}
where $s\in\mathbb R$ and $z\in\mathbb T=\mathbb R/\mathbb Z$. Thus the two components have the \textup{same constant diffusion matrix}. The reaction matrix is diagonalized by the vectors $(1,1)$ and $(1,-1)$, with eigenvalues $r$ and $r-2\gamma$, respectively. The diffusion and shear terms act in the same way on each component. Since $\gamma>0$, the positive principal eigenfunction lies in the symmetric mode $(\phi,\phi)$, for which the tilted eigenvalue equation is scalar. Denote by $\nu_d(q)$ the largest periodic eigenvalue of
\begin{equation}\label{eq:memoirs-diffusion-eigen}
 dD\phi''+(r+q\sin(kz))\phi=\nu_d(q)\phi,
 \qquad \phi(z+1)=\phi(z)>0.
\end{equation}
The minimal spectral spreading speed in the positive $s$-direction has the exact characterization
\begin{equation}\label{eq:memoirs-diffusion-speed}
 c^*(d)=\inf_{\lambda>0}\left\{ad\lambda+\frac{\nu_d(h\lambda)}{\lambda}\right\}.
\end{equation}
In this formula $h$ is fixed as $d\downarrow0$ or $d\to\infty$. For the examples of this section one has $h=Bb_j$ and $B=5$; in particular, $h$ is not the varying diffusion rate.

\begin{lemma}[Two periodic eigenvalue regimes]\label{lem:memoirs-diffusion-eigen}
Let $D,r>0$, $k=2\pi$, and let $\nu_d(q)$ be defined by \eqref{eq:memoirs-diffusion-eigen}. Then $\nu_d(-q)=\nu_d(q)$ and $\nu_d(q)\ge r$. The following asymptotic estimates hold.
\begin{enumerate}[label=\textup{(\roman*)},leftmargin=2em]
\item For $|q|/(dD)$ sufficiently small,
\begin{equation}\label{eq:memoirs-diffusion-perturb}
\nu_d(q)=r+\frac{q^2}{2dDk^2}
 +O\left(\frac{q^4}{(dD)^3}\right).
\end{equation}
\item For $q/(dD)\to+\infty$,
\begin{equation}\label{eq:memoirs-diffusion-localized}
\nu_d(q)=r+q-\frac{k}{\sqrt2}\sqrt{dDq}+O(dD),
\end{equation}
with a constant independent of $q$ and $d$ whenever $q/(dD)$ is sufficiently large.
\end{enumerate}
\end{lemma}
\begin{proof}
Translation by half a period interchanges $q$ with $-q$ without changing the spectrum, and the constant test function in the Rayleigh quotient shows $\nu_d(q)\ge r$. For \textup{(i)}, divide the operator by $dD$ and regard the potential $q\sin(kz)/(dD)$ as a perturbation of $\partial_{zz}$. The constant normalized eigenfunction has eigenvalue zero, isolated from the remaining eigenvalues $-k^2n^2$, $n\in\mathbb Z\setminus\{0\}$. The first-order coefficient vanishes because $\int_0^1\sin(kz)\,dz=0$. Its second-order coefficient is
\begin{equation*}
\sum_{n\ne0}\frac{|\langle\sin(kz),e^{inkz}\rangle|^2}{k^2n^2}
=\frac{1}{2k^2},
\end{equation*}
since only the modes $n=\pm1$ contribute. Analytic perturbation of this isolated simple eigenvalue is valid for small $|q|/(dD)$, and the half-period symmetry makes its Taylor series even. Rescaling the eigenvalue gives \eqref{eq:memoirs-diffusion-perturb}.

For \textup{(ii)}, let $z_0=1/4$, so that $\sin(kz_0)=1$, and introduce
\begin{equation*}
\varepsilon=\sqrt{\frac{dD}{q}},\qquad
\mathcal H_\varepsilon=-\varepsilon^2\frac{d^2}{dz^2}
 +1-\sin(kz)
\end{equation*}
on the periodic unit circle. If $E_0(\varepsilon)$ is its lowest eigenvalue, the Rayleigh principle gives $\nu_d(q)=r+q-qE_0(\varepsilon)$. Near its unique minimum modulo the period, the potential satisfies
\begin{equation*}
1-\sin(k(z_0+y))=\frac{k^2}{2}y^2+O(y^4).
\end{equation*}
The harmonic-oscillator approximation yields
\begin{equation}\label{eq:memoirs-diffusion-oscillator}
E_0(\varepsilon)=\frac{k}{\sqrt2}\varepsilon+O(\varepsilon^2).
\end{equation}
For completeness, we justify the order of the remainder. A Gaussian proportional to
$\exp(-k y^2/(2\sqrt2\,\varepsilon))$, cut off within a fixed neighborhood of $y=0$, gives the upper bound in \eqref{eq:memoirs-diffusion-oscillator}: its second and fourth moments are $O(\varepsilon)$ and $O(\varepsilon^2)$, respectively. For the lower bound, write $V(y)=1-\cos(ky)$ on the circle, centered at $y=0$, and let $\psi_\varepsilon$ be its normalized positive ground state. The upper bound already proved yields $E_0(\varepsilon)\leq C\varepsilon$. Choose a smooth periodic function $S\geq0$ such that $S(y)=y^2$ near zero, while $S$ is positive and constant outside a fixed neighborhood of zero. Since $V$ has a unique nondegenerate zero on the circle, there are constants $c,C>0$ with $V\geq cS$ and $|S'|^2\leq CS$. For a fixed sufficiently small $\eta>0$, put $W_\varepsilon=\eta S/\varepsilon$; then $\varepsilon^2|W_\varepsilon'|^2\leq V/4$. Multiplying the eigenvalue equation by $e^{2W_\varepsilon}\psi_\varepsilon$ and integrating by parts gives the exact identity
\begin{equation*}
\varepsilon^2\int_{\mathbb T}|(e^{W_\varepsilon}\psi_\varepsilon)'|^2+\int_{\mathbb T}(V-E_0(\varepsilon)-\varepsilon^2|W_\varepsilon'|^2)e^{2W_\varepsilon}\psi_\varepsilon^2=0.
\end{equation*}
On the set $V\geq4C\varepsilon$ the second integrand is bounded below by $(V/2)e^{2W_\varepsilon}\psi_\varepsilon^2$, after choosing $C$ so that $E_0(\varepsilon)\le C\varepsilon$. On its complement, $S\leq C_1\varepsilon$ and $W_\varepsilon$ is bounded independently of $\varepsilon$. The negative part of the integral is consequently bounded by $C_2\varepsilon$, since $\|\psi_\varepsilon\|_2=1$. Hence $\int_{\{V\geq4C\varepsilon\}}Ve^{2W_\varepsilon}\psi_\varepsilon^2\leq C_4\varepsilon$. The inequality $V\geq4C\varepsilon$ on this set, together with the boundedness of $W_\varepsilon$ on its complement, now yields $\int_{\mathbb T}e^{2W_\varepsilon}\psi_\varepsilon^2\leq C_3$, uniformly as $\varepsilon\downarrow0$. Since $S(y)=y^2$ near zero and $S$ has a positive lower bound away from a fixed neighborhood of zero, the inequalities $y^4e^{-2\eta y^2/\varepsilon}\leq C\varepsilon^2$ near zero and $e^{-2\eta S/\varepsilon}\leq C_\delta e^{-c_\delta/\varepsilon}$ away from zero imply
\begin{equation*}
\int_{\mathbb T}\operatorname{dist}(y,\mathbb Z)^4\psi_\varepsilon^2\,dy=O(\varepsilon^2),\qquad
\int_{|y|\geq\delta}\psi_\varepsilon^2\,dy=O(e^{-c_\delta/\varepsilon})
\end{equation*}
for every fixed sufficiently small $\delta>0$. Let $\chi$ be a smooth cutoff equal to one near $y=0$ and supported where $V(y)=k^2y^2/2+O(y^4)$. The IMS identity, the moment estimate and the exponential tail show that the quadratic form at $\psi_\varepsilon$ is bounded below by that of the full-line harmonic oscillator $-\varepsilon^2\partial_{yy}+k^2y^2/2$ at $\chi\psi_\varepsilon$, minus $C\varepsilon^2$. Its ground-state energy is $k\varepsilon/\sqrt2$, and $\|\chi\psi_\varepsilon\|_2^2=1+O(e^{-c/\varepsilon})$. Consequently $E_0(\varepsilon)\geq k\varepsilon/\sqrt2-C\varepsilon^2$. Together with the Gaussian upper bound this proves the two-sided estimate. This proves \eqref{eq:memoirs-diffusion-oscillator}, and multiplication by $q$ gives \eqref{eq:memoirs-diffusion-localized}.
\end{proof}

We now prove the diffusion-rate theorem stated in the introduction. The spectral reduction above applies because the coupling matrix is symmetric and the two diffusion matrices coincide. In particular, the positive spectral mode has eigenvalue $r$ for every $\gamma>0$, even when $r-\gamma\leq0$; no positivity of the individual diagonal reaction coefficients is used here.

\begin{proof}[Proof of Theorem~\ref{thm:memoirs-diffusion-two-regimes}]
We begin with the large-diffusion limit. At every $d>0$, the infimum in \eqref{eq:memoirs-diffusion-speed} is attained: the lower bound $\nu_d(q)\ge r$ implies divergence of the quotient as $\lambda\downarrow0$, while $ad\lambda$ tends to infinity as $\lambda\to\infty$. The same lower bound and Lemma~\ref{lem:memoirs-diffusion-eigen}\textup{(i)}, tested at $\lambda=\sqrt{r/(ad)}$, yield
\begin{equation*}
c^*(d)\le2\sqrt{ra\,d}+O(d^{-3/2}).
\end{equation*}
On the other hand,
\begin{equation*}
c^*(d)\ge ad\lambda_d+\frac r{\lambda_d}
\ge2\sqrt{ra\,d}.
\end{equation*}
Writing $\lambda_d=x_d/\sqrt d$, these inequalities imply $ax_d+r/x_d\le2\sqrt{ra}+O(d^{-2})$. Since $x\mapsto ax+r/x$ has a unique strict minimum at $x_0=\sqrt{r/a}$ and tends to infinity at both endpoints, $x_d\to x_0$; indeed $|x_d-x_0|=O(d^{-1})$ by its positive second derivative near $x_0$. Thus $h\lambda_d/(dD)=O(d^{-3/2})$, and the eigenvalue expansion is uniform on a fixed compact interval containing all sufficiently large-$d$ minimizers. There,
\begin{equation}\label{eq:memoirs-diffusion-effective}
\begin{aligned}
ad\lambda+\frac{\nu_d(h\lambda)}\lambda
={}&\left(ad+\frac{h^2}{2dDk^2}\right)\lambda
 +\frac r\lambda
 +O\left(\frac{h^4\lambda^3}{(dD)^3}\right).
\end{aligned}
\end{equation}
Since $\lambda\asymp d^{-1/2}$, the remainder is $O(d^{-9/2})$ uniformly near both the original and the reduced minimizing exponents. The explicit minimum of the leading terms is
\begin{equation*}
2\sqrt{r\left(ad+\frac{h^2}{2dDk^2}\right)}.
\end{equation*}
Taylor's formula for this square root proves \eqref{eq:memoirs-diffusion-large}.

For $d\downarrow0$, set $\lambda=x/\sqrt d$. The upper bound $\nu_d(h\lambda)\le r+h\lambda$ gives
\begin{equation}\label{eq:memoirs-diffusion-small-upper}
c^*(d)\le h+2\sqrt{ra\,d}.
\end{equation}
A localized Gaussian test function, as in Lemma~\ref{lem:memoirs-diffusion-eigen}\textup{(ii)}, yields $\nu_d(q)\ge r+q-C\sqrt{dDq}$ for $q\ge dD$, with $C$ independent of $q$ and $d$. Choose a minimizing exponent $\lambda_d$. If $h\lambda_d<dD$, then $\nu_d(h\lambda_d)\ge r$ and the quotient is at least $r/\lambda_d>rh/(dD)$, contradicting \eqref{eq:memoirs-diffusion-small-upper} for small $d$. Hence $h\lambda_d\ge dD$. Put $x_d=\lambda_d\sqrt d$. For a constant $C_1>0$ independent of $d$, the Gaussian estimate now implies
\begin{equation*}
 h+\sqrt d\left(ax_d+\frac r{x_d}-C_1d^{1/4}x_d^{-1/2}\right)
 \le c^*(d)\le h+2\sqrt{rad}.
\end{equation*}
For $0<x\le1$ and small $d$, $r/x-C_1d^{1/4}/\sqrt x\ge r/(2x)$; for $x\ge1$, $ax-C_1d^{1/4}/\sqrt x\ge ax/2$. Therefore the displayed inequalities exclude $x_d\downarrow0$ and $x_d\to\infty$, and there exist $0<m<M<\infty$ such that $m\le x_d\le M$. On this interval $h\lambda_d/(dD)\asymp d^{-3/2}$, and Lemma~\ref{lem:memoirs-diffusion-eigen}\textup{(ii)} applies with a uniform remainder. More precisely,
$\varepsilon=\sqrt{dD/(h\lambda)}\le C d^{3/4}$ for
$\lambda=x/\sqrt d$ and $x\in[m,M]$. The proof of the spectral
lemma bounds $|E_0(\varepsilon)-k\varepsilon/\sqrt2|$ by
$C\varepsilon^2$ with a constant depending only on the fixed
periodic potential. After multiplication by $q=h\lambda$ and
division by $\lambda$, the remainder in the speed functional
is at most $C dD/\lambda\le C' d^{3/2}$ uniformly on
$[m,M]$. Thus the minimization below takes place inside the
interval on which the error bound is uniform.

Put $K=k\sqrt{Dh}/\sqrt2$. The eigenvalue estimate and \eqref{eq:memoirs-diffusion-speed} now imply
\begin{equation}\label{eq:memoirs-diffusion-rescaled}
c^*(d)=h+\sqrt d\min_{m\le x\le M}
\left\{ax+\frac r x-d^{1/4}Kx^{-1/2}\right\}
+O(d^{3/2}).
\end{equation}
The unperturbed function $F(x)=ax+r/x$ has its unique nondegenerate minimum at $x_0=\sqrt{r/a}$, with
\begin{equation*}
F(x_0)=2\sqrt{ra},\qquad
F''(x_0)=\frac{2a^{3/2}}{\sqrt r}.
\end{equation*}
Set $G(x)=-Kx^{-1/2}$ and $\varepsilon=d^{1/4}$. Taylor expansion at $x_0$, or equivalently the implicit-function theorem for $F'(x)+\varepsilon G'(x)=0$, gives
\begin{equation}\label{eq:memoirs-diffusion-min-perturb}
\min_x\{F(x)+\varepsilon G(x)\}
=F(x_0)+\varepsilon G(x_0)
-\frac{\varepsilon^2[G'(x_0)]^2}{2F''(x_0)}
+O(\varepsilon^3).
\end{equation}
Direct calculation yields
\begin{equation*}
G(x_0)=-K\left(\frac ar\right)^{1/4},~
\frac{[G'(x_0)]^2}{2F''(x_0)}
=\frac{K^2}{16r}=\frac{k^2Dh}{32r}.
\end{equation*}
Substitution into \eqref{eq:memoirs-diffusion-rescaled} proves \eqref{eq:memoirs-diffusion-small}, since the perturbative remainder after multiplication by $\sqrt d$ is $O(d^{5/4})$. The limits in \eqref{eq:memoirs-diffusion-limits} follow.
\end{proof}

\section{Large-Advection Asymptotics}\label{sec:large-advection-proofs}
In the almost-periodic setting, consider the advection terms ${\bf q}_i=A{\bf v}_i$, $i=1,2$, where $A>0$ and ${\bf v}_i={\bf v}_i(t,x)$ are fixed. In this case, $c^*=c^*(e,A)$. Recall \eqref{div}.
\begin{align}\label{dirfcase1}
\begin{cases}
\partial_t \eta_{\lambda,1}
&= \text{Tr} \big( {\bf A}^1 D^2 \eta_{\lambda,1} \big)
- \big( A {\bf v}_1 + 2\lambda {\bf A}^1e \big) \cdot D \eta_{\lambda,1}+ \lambda A \left({\bf v}_1\cdot e\right)\eta_{\lambda,1} + H_1^\lambda\eta_{\lambda,1}+\mu_{12}\eta_{\lambda,2},\\
\partial_t \eta_{\lambda,2}
&= \text{Tr} \big( {\bf A}^2  D^2 \eta_{\lambda,2} \big)
- \big( A {\bf v}_2  + 2\lambda {\bf A}^2 e \big) \cdot D \eta_{\lambda,2}+\lambda A \left({\bf v}_2 \cdot e\right) \eta_{\lambda,2}+H_2^\lambda\eta_{\lambda,2}+\mu_{21} \eta_{\lambda,1},
\end{cases}
\end{align}
 where $H_{i}^\lambda:=\mu_{ii}  + \lambda^2\, e \cdot {\bf A}^i e,~i=1,2.$ By the additional regularity assumptions \eqref{Fregularity} and standard parabolic Schauder estimates,
\(
\boldsymbol{\eta}_\lambda
\in
\left[C_{\mathrm{loc}}^{1+\alpha/2,2+\alpha}
(\R\times\R^N)\right]^2.
\) 

We first investigate the case where the drift amplitude $A\rightarrow \infty$ for $\boldsymbol{\eta}_{\lambda}$. Consider
\begin{align}\label{sub2}
\begin{cases}
\partial_t \psi_1
= \lambda A v_1 \psi_1,\\
\partial_t \psi_2
= \lambda A v_2 \psi_2,
\end{cases}
\end{align}
where $v_i(t) = \inf\limits_{x\in \R^N}\left({\bf v}_i(t,x) \cdot e\right)$. For an initial constant vector with nonnegative components, the solution of \eqref{sub2} is spatially constant. Since
$\lambda A\bigl(\mathbf v_i(t,x)\cdot e-v_i(t)\bigr)\ge0$
and the remaining reaction and coupling terms in the tilted equation are nonnegative, this vector is a subsolution of \eqref{dirfcase1}. Thus, for $i=1,2$ and $m_i(t) := \min \limits_{x\in [{\bf 0},{\bf L}]} \eta_{\lambda,i}(t,x)$, the parabolic comparison principle implies that $m_i(t+T)\geq m_i(t) \exp\left(\lambda A\displaystyle \int_t^{t+T} v_i(s)ds\right).$ Using \eqref{below} or \eqref{Harnack-zero}, we get
\begin{align}\label{below1}
a(t,T) \geq \lambda A  \int_t^{t+T} v_i(s)ds + \ln( C_0(A)).
\end{align}

To obtain the upper estimate, we use \eqref{super1} and write
\begin{align*}
M_i(t)&=\max_{x\in[\mathbf0,\mathbf L]}\eta_{\lambda,i}(t,x),~ M(t)=\max_{i=1,2}M_i(t),~
w_i(t)=\sup_{x\in\R^N}({\bf v}_i(t,x)\cdot e).
\end{align*}
Then,
\begin{align}\label{upper1}
M_i(t+T) \leq M(t)\exp\left(\int_t^{t+T} \left(\lambda A \max_{i=1,2}w_i(s) + \mu_\infty(s)\right)ds + C\lambda^2 T\right),
\end{align}
where $C$ is in \eqref{Aassump}. Note that if $t\mapsto g(t,x)$ is almost periodic in time uniformly with respect to $x\in \R^N$, then both $t\mapsto \sup\limits_{x\in\R^N}g(t,x),\inf\limits_{x\in\R^N}g(t,x)$ are almost periodic in time. It then follows from \eqref{below1} and \eqref{upper1} that
\begin{align}\label{estimatec}
A  \max_{i=1,2}\left<  v_i \right> \leq \left<c_\lambda\right>\leq C\lambda +A \left<  \max_{i=1,2}w_i \right> + \lambda^{-1}\left< \mu_{\infty}\right>,~\forall \lambda>0.
\end{align}
The minimal-speed characterization in Theorem~\ref{minimalspeed} is applied for each fixed $A>0$.

This implies that for each $\lambda>0$,
\begin{align*}
\max_{i=1,2}\left<  v_i \right>\leq \liminf\limits_{A \rightarrow \infty} \dfrac{\left<c_\lambda\right>}{A}\leq \limsup\limits_{A \rightarrow \infty} \dfrac{\left<c_\lambda\right>}{A} \leq \left<   \max_{i=1,2} w_i \right>.
\end{align*}
From Theorem \ref{minimalspeed}, one has $\left<c_{\lambda^*}\right> \leq \left<c_\lambda\right> $, where $\lambda^*=\lambda^*(A)$ is a global minimizer over $(0,\infty)$. Thus,
\begin{align*}
\limsup_{A \rightarrow \infty}\dfrac{\left<c_{\lambda^*}\right>}{A} \leq  \limsup\limits_{A \rightarrow \infty} \dfrac{\left<c_\lambda\right>}{A} \leq \left< \max_{i=1,2} w_i \right>.
\end{align*}
By choosing $\lambda=\lambda^*(A)$ in \eqref{estimatec}, we obtain
\begin{align}\label{welldeflim}
\max_{i=1,2}\left< \inf\limits_{x\in \R^N} ({\bf v}_i\cdot e) \right>\leq \liminf\limits_{A \rightarrow \infty}\dfrac{c^*(e,A)}{A} \leq \limsup\limits_{A \rightarrow \infty} \dfrac{c^*(e,A)}{A}  \leq \left<   \max_{i=1,2} \sup\limits_{x\in \R^N}({\bf v}_i\cdot e )\right>.
\end{align}
\begingroup
The same bounds hold for ${\bf q}_i=A{\bf v}_i+{\bf r}_i$ with ${\bf r}_i$ bounded independently of $A$. Indeed, the exponentially tilted linear equation then contains the diagonal zeroth-order term $\lambda(A{\bf v}_i+{\bf r}_i)\cdot e$. The spatially constant comparison functions used above eliminate the transport derivatives, including ${\bf r}_i\cdot D\eta_{\lambda,i}$. Let $R=\max_i\|{\bf r}_i\cdot e\|_\infty$. The same subsolution and supersolution estimates, now retaining the contribution of ${\bf r}_i\cdot e$, imply
\begin{equation}\label{eq:v5-remainder}
A\max_i\left\langle\inf_x({\bf v}_i\cdot e)\right\rangle-R
\le\langle c_\lambda\rangle\le
A\left\langle\max_i\sup_x({\bf v}_i\cdot e)\right\rangle+R+C\lambda+\frac{C_1}{\lambda},\qquad\lambda>0,
\end{equation}
where $C,C_1$ depend only on the ellipticity bounds and the reaction coefficients, not on $A$. Since the lower bound is independent of $\lambda$ and the upper bound may be evaluated at any fixed $\lambda>0$, minimizing with respect to $\lambda$ and dividing by $A$ proves \eqref{eq:v5-general} in the claimed setting. 
\endgroup
Thus, the liminf and limsup are finite. In particular, if ${\bf v}_i \cdot e =0$ on $\R^{N+1}$ for $i=1,2$, one has
\begin{align*}
\lim\limits_{A \rightarrow \infty} \dfrac{c^*(e,A)}{A} =0.
\end{align*}


Next, we estimate the value of these limits. 
\subsection{Case 1: ${\bf v}_i={\bf v}_i(t)$ for $i=1,2$.}

Since the velocity fields do not depend on $x$, we write
$v_i(t)={\bf v}_i(t)\cdot e$. In direction $e$, estimate
\eqref{estimatec} reads
\begin{align}\label{estimatec12}
A\max_{i=1,2}\langle v_i\rangle
\le \langle c_\lambda\rangle
\le C\lambda+A\langle\max_{i=1,2}v_i\rangle
+\lambda^{-1}\langle\mu_\infty\rangle.
\end{align}
Suppose that ${\bf v}_1(t)={\bf v}_2(t)$ for all $t>T(e)$.
The difference ${\bf v}_1-{\bf v}_2$ is Bohr almost periodic
and vanishes on a half-line; therefore it vanishes identically.
Indeed, for any $s\in\R$ and any $\varepsilon>0$, choose an
$\varepsilon$-almost period $\tau>T(e)-s$. Then
$|{\bf v}_1(s)-{\bf v}_2(s)|\le\varepsilon$, and let
$\varepsilon\downarrow0$. Thus, the assumption is equivalent to
${\bf v}_1={\bf v}_2={\bf v}$ on the whole time axis.
The lower and upper estimates above consequently have the same
leading term and imply
\begin{align*}
\lim_{A\to\infty}\frac{c^*(e,A)}A
=\langle{\bf v}\cdot e\rangle.
\end{align*}
This is the common-velocity case, not an independent weaker
assumption.

\subsection{Case 2: ${\bf v}_1={\bf v}_2={\bf v}(t,x)$, $\nabla_x\cdot{\bf v}=0$ and ${\bf v}$ has zero spatial mean.}


 We investigate the case where the drift amplitude $A\rightarrow \infty$ under condition \eqref{driftterm}. Recall, from \eqref{div}, that
\begin{align}\label{hatsys}
\begin{cases}
\partial_t \eta_{\lambda,1}
&= \operatorname{Tr} \big(  {\bf A}^1 D^2 \eta_{\lambda,1} \big)
- \big( A  {\bf v} + 2\lambda  {\bf A}^1e \big) \cdot D \eta_{\lambda,1}+ \lambda A  ({\bf v}\cdot e)\eta_{\lambda,1} +  {H}_1^\lambda\eta_{\lambda,1}+ {\mu}_{12}\eta_{\lambda,2},\\
\partial_t \eta_{\lambda,2}
&= \operatorname{Tr} \big(  {\bf A}^2  D^2 \eta_{\lambda,2} \big)
- \big( A  {\bf v}  + 2\lambda  {\bf A}^2 e \big) \cdot D \eta_{\lambda,2}+\lambda A  ({\bf v} \cdot e) \eta_{\lambda,2}+ {H}_2^\lambda\eta_{\lambda,2}+ {\mu}_{21} \eta_{\lambda,1},
\end{cases}
\end{align}
where  $ {H}_{i}^\lambda= {\mu}_{ii}  + \lambda^2\, e \cdot  {\bf A}^i e>0,~i=1,2.$ Thanks to the uniqueness in Lemma \ref{lm:Haranack-type}, we may assume that $\|\boldsymbol{\eta}_{\lambda}(0,\cdot)\|_{\infty}=1,~\forall \lambda>0$. 
Recall the set of first integrals 
\begin{align*}
\mathfrak{I}:=\left\{ \boldsymbol{\phi}\in [H_{per}^1(\R^N)]^2:\boldsymbol{\phi}=(\phi_1,\phi_2),~ {\bf v}\cdot D \phi_i=0 \text{ on }\R\times \R^N,~\forall i=1,2 \right\},
\end{align*}
and
\begin{align*}
\mathfrak{I}_{{\bf A}}^{\mu}:=\left\{\boldsymbol{\phi}=(\phi_1,\phi_2)\in \mathfrak{I}:\sum_{i=1,2}\left<\int_{[{\bf 0},{\bf L}]} D\phi_i\cdot   {\bf A}^iD\phi_i dx\right>\leq \sum_{i,j=1,2}\left<\int_{[{\bf 0},{\bf L}]} {\mu}_{ij} \phi_i\phi_j dx\right>  \right\}
\end{align*}
Note that both sets $\mathfrak{I}$ and $\mathfrak{I}_{\mathbf{A}}^{\mu}$ are nonempty since they contain constant functions.

Let $\boldsymbol{\phi}=(\phi_1,\phi_2) \in \mathfrak{I}_{{\bf A}}^{\mu}\setminus\{{\bf 0}\}$. Define $\boldsymbol{\psi}=(\psi_1,\psi_2)=\ln{\boldsymbol{\eta}}_{\lambda} = (\ln \eta_{\lambda,1},\ln \eta_{\lambda,2})$. We obtain the following lemma.
\begin{lemma}

The following estimate holds
\begin{align}\label{cidentity}
\left<\int_{[{\bf 0},{\bf L}]}\phi^2_i\partial_t \psi_jdx\right>= \int_{[{\bf 0},{\bf L}]}\phi^2_i\left<\partial_t \psi_j\right>dx=\lambda \left< {c}_\lambda\right> \int_{[{\bf 0},{\bf L}]}\phi_i^2(x) dx,~\forall i,j=1,2.
\end{align}

\end{lemma}

\begin{proof}
For $i,j=1,2$,  consider
\begin{align*}
l(t) := \ln \| {\boldsymbol{\eta}}_{\lambda}(t,\cdot)\|_\infty,~G_{i,j}(t) := \int_{[{\bf 0},{\bf L}]}\phi_i^2(x) \psi_j(t,x) dx,~B_i := \int_{[{\bf 0},{\bf L}]}\phi_i^2(x) dx.
\end{align*}
These functions are well-defined since $\phi_i \in L^2([{\bf 0},{\bf L}])$ and $x\mapsto\psi_j(t,x)$ is bounded for each $t>0$. Note that
\begin{align*}
\left<\int_{[{\bf 0},{\bf L}]}\phi^2_i\partial_t \psi_jdx\right>=\lim_{T\rightarrow \infty}\dfrac{1}{T}\int_{[{\bf 0},{\bf L}]}\phi^2_i(x) \int_t^{t+T}\partial_t \psi_j(s,x)dsdx=\lim_{T\rightarrow \infty}\dfrac{1}{T}(G_{i,j}(t+T)-G_{i,j}(t)),
\end{align*}
uniformly in $t\in \R$. Similarly to Lemma \ref{lm:Haranack-type} and \eqref{Harnack-zero}, there exists $ {C}_A\in (0,1)$ small enough so that
\begin{align*}
\min_{i=1,2}\min\limits_{x\in \R^{N}}\eta_{\lambda,i}(t,x)\geq  {C}_A\max_{i=1,2}\max_{x\in \R^{N}}\eta_{\lambda,i}(t,x)= {C}_A\| {\boldsymbol{\eta}}_{\lambda}(t,\cdot)\|_\infty,~\forall t\in \R.
\end{align*}
This implies that $l(t)\geq \psi_j(t,x)>l(t)-K_A,~\forall (t,x)\in \R\times \R^N$ for some $K_A>0$. Direct calculation yields that
\begin{align*}
|G_{i,j}(t+T)-G_{i,j}(t)-B_i(l(t+T)-l(t))|\leq K_AB_i.
\end{align*}
From \eqref{uniformcmean}, one has
\begin{align*}
\lim_{T\rightarrow \infty}\dfrac{G_{i,j}(t+T)-G_{i,j}(t)}{T}=\lim_{T\rightarrow \infty}B_i\dfrac{l(t+T)-l(t)}{T}=B_i\lambda \left< {c}_\lambda\right> \text{ uniformly in $t\in \R$. }
\end{align*}
This completes the proof.
\end{proof}

We are now in a position to prove the large-drift limit.

\begin{proof}[Proof of Theorem~\ref{limitspeed1}\textup{(b)}, lower bound]
For $\boldsymbol{\phi}\in \mathfrak{I}\setminus\{\bf 0\}$ and $i,j\in\{1,2\}$ with $i\neq j$, we have
\begin{align}\label{phi_i^2}
&\int_{[{\bf 0},{\bf L}]}\phi^2_i\left[\partial_t \psi_i
-\text{Tr} \big(  {\bf A}^i D^2 \psi_i\big)
-D\psi_i\cdot {\bf A}^i D\psi_i + \big( A  {\bf v} + 2\lambda  {\bf A}^ie \big) \cdot D \psi_i- \lambda A  ({\bf v}\cdot e )-  {H}_i^\lambda- {\mu}_{ij}\dfrac{\eta_{\lambda,j}}{\eta_{\lambda,i}}\right]dx=0.
\end{align}
It follows from \eqref{regularityofA}-\eqref{driftterm} and the definition of $\mathfrak{I}$ that
\begin{align*}
&\int_{[{\bf 0},{\bf L}]}\phi^2_iA {\bf v} \cdot D\psi_idx  = -A\int_{[{\bf 0},{\bf L}]}\psi_i\left[ {\bf v}\cdot D (\phi^2_i)+\phi^2_i\nabla_x\cdot  {\bf v} \right]dx=0, \\
&-\int_{[{\bf 0},{\bf L}]}\phi^2_i\text{Tr} \big(  {\bf A}^i D^2 \psi_i\big)dx=2\int_{[{\bf 0},{\bf L}]}\phi_i D \phi_i\cdot  {\bf A}^i D\psi_idx,
\end{align*}
using \(\bf L\)-periodicity, \({\bf v}\cdot D\phi_i=0\), \(\nabla_x\cdot {\bf v}=0\), and \(\nabla_x\cdot {\bf A}^i=0\). Thus, one has
\begin{align}\label{eq:test}
\begin{split}
 \int_{[{\bf 0},{\bf L}]}\phi^2_i\partial_t \psi_idx &= \int_{[{\bf 0},{\bf L}]}\left(\phi^2_iD\psi_i\cdot {\bf A}^i D\psi_i-2\phi_i D \phi_i\cdot  {\bf A}^i D\psi_i
 -  2\lambda \phi^2_i  {\bf A}^ie \cdot D \psi_i\right)dx \\
 &+ \lambda A \int_{[{\bf 0},{\bf L}]}( {\bf v}\cdot e) \phi^2_i dx +   \int_{[{\bf 0},{\bf L}]}(\lambda^2 e \cdot  {\bf A}^i e)\phi^2_idx+ \int_{[{\bf 0},{\bf L}]}\left[ {\mu}_{ii}  +  {\mu}_{ij}\dfrac{\eta_{\lambda,j}}{\eta_{\lambda,i}}\right]\phi^2_idx.
 \end{split}
\end{align}
For $\theta\in(0,1)$, set $X=\phi_iD\psi_i$, $Y=D\phi_i$ and $Z=\lambda\phi_i e$. Completing the square in the scalar product induced by ${\bf A}^i$, and then applying the weighted Young inequality, gives
\[
X\cdot{\bf A}^iX-2X\cdot{\bf A}^i(Y+Z)
\geq -(Y+Z)\cdot{\bf A}^i(Y+Z)
\geq-\theta^{-1}Y\cdot{\bf A}^iY
-(1-\theta)^{-1}Z\cdot{\bf A}^iZ.
\]
Integration over the periodic cell yields
\begin{align*}
&\int_{[{\bf 0},{\bf L}]}\left(\phi^2_iD\psi_i\cdot {\bf A}^i D\psi_i-2\phi_i D \phi_i\cdot  {\bf A}^i D\psi_i
 -  2\lambda \phi^2_i  {\bf A}^i e \cdot D \psi_i\right)dx\\
&\geq -\int_{[{\bf 0},{\bf L}]} \left(
\dfrac{1}{\theta}D\phi_i\cdot {\bf A}^iD\phi_i
+\dfrac{\lambda^2}{1-\theta}\phi_i^2 e\cdot{\bf A}^ie
\right)dx. 
\end{align*}
Choose $\theta=(1+\lambda \alpha)^{-1}$ for some $\alpha>0$, we get
\begin{align*}
\int_{[{\bf 0},{\bf L}]}\phi^2_i\partial_t \psi_idx &\geq \int_{[{\bf 0},{\bf L}]} \left({\mu}_{ii} \phi^2_i- D\phi_i\cdot   {\bf A}^iD\phi_i\right) dx - \lambda\int_{[{\bf 0},{\bf L}]}\left( \dfrac{1}{\alpha}\phi_i^2 e \cdot  {\bf A}^ie + \alpha D\phi_i\cdot  {\bf A}^i D\phi_i\right) dx
\\
&+ \lambda A \int_{[{\bf 0},{\bf L}]}( {\bf v}\cdot e) \phi^2_i dx 
+ \int_{[{\bf 0},{\bf L}]}  {\mu}_{ij}\dfrac{\eta_{\lambda,j}}{\eta_{\lambda,i}}\phi^2_idx.
\end{align*}
This implies that
\begin{align*}
\sum_{i=1,2}\int_{[{\bf 0},{\bf L}]}\phi^2_i\partial_t \psi_i\,dx
&\geq \sum_{i=1,2}\int_{[{\bf 0},{\bf L}]} \left({\mu}_{ii} \phi^2_i- D\phi_i\cdot {\bf A}^iD\phi_i\right) dx\\
&\quad- \lambda \sum_{i=1,2}\int_{[{\bf 0},{\bf L}]}\left( \dfrac{1}{\alpha}\phi_i^2 e \cdot {\bf A}^ie + \alpha D\phi_i\cdot {\bf A}^iD\phi_i\right) dx
\\
&+ \lambda A \int_{[{\bf 0},{\bf L}]}( {\bf v}\cdot e) |\boldsymbol{\phi}|^2 dx 
+ \int_{[{\bf 0},{\bf L}]} 2 \sqrt{ {\mu}_{12} {\mu}_{21}}\phi_1\phi_2 dx.
\end{align*}
Here, we use $2ab\leq a^2+b^2$ for the last term. Taking the uniform mean and using $\mu_{12}=\mu_{21}$, \eqref{cidentity} implies that
\begin{align*}
\lambda \left< {c}_\lambda\right> \int_{[{\bf 0},{\bf L}]}|\boldsymbol{\phi}|^2 dx 
&\geq \sum_{i,j=1,2}\left<\int_{[{\bf 0},{\bf L}]} {\mu}_{ij} \phi_i\phi_j dx\right>- \sum_{i=1,2}\left<\int_{[{\bf 0},{\bf L}]} D\phi_i\cdot   {\bf A}^iD\phi_i dx\right> \\
&- \lambda\sum_{i=1,2}\left<\int_{[{\bf 0},{\bf L}]}\left( \dfrac{1}{\alpha}\phi_i^2 e \cdot  {\bf A}^ie + \alpha D\phi_i\cdot  {\bf A}^i D\phi_i\right) dx\right>
+ \lambda A \left<\int_{[{\bf 0},{\bf L}]}( {\bf v}\cdot e) |\boldsymbol{\phi}|^2 dx\right>.
\end{align*}

Since $\boldsymbol{\phi}=(\phi_1,\phi_2) \in \mathfrak{I}_{{\bf A}}^{\mu}\setminus\{{\bf 0}\}$, direct calculation yields that
\begin{align*}
\dfrac{\left< {c}_\lambda\right>}{A} \int_{[{\bf 0},{\bf L}]}|\boldsymbol{\phi}|^2 dx \geq -\dfrac{1}{A}\sum_{i=1,2}\left<\int_{[{\bf 0},{\bf L}]} \left( \dfrac{1}{\alpha}\phi_i^2 e \cdot  {\bf A}^ie + \alpha D\phi_i\cdot  {\bf A}^iD\phi_i\right) dx\right>+  \left<\int_{[{\bf 0},{\bf L}]}( {\bf v}\cdot e)|\boldsymbol{\phi}|^2 dx\right>.
\end{align*}
By choosing $\lambda = \lambda^*(A)$ and passing to the limit $A\rightarrow \infty$, we get
\begin{align*}
\liminf_{A\rightarrow \infty}\dfrac{ {c}^*(e,A)}{A} \geq \max_{\boldsymbol{\phi}\in \mathfrak{I}_{{\bf A}}^{\mu}\setminus\{\bf 0\}}\dfrac{\left<\displaystyle\int_{[{\bf 0},{\bf L}]}( {\bf v}\cdot e)|\boldsymbol{\phi}|^2 dx\right>}{\displaystyle\int_{[{\bf 0},{\bf L}]}|\boldsymbol{\phi}|^2 dx},
\end{align*} 
which concludes the proof.
\end{proof}

\begin{remark}\label{remark7}
If we do not assume the symmetric property, i.e., $\mu_{12}\neq \mu_{21}$, one can repeat the calculation so that
\begin{align*}
\liminf_{A\rightarrow \infty}\dfrac{ {c}^*(e,A)}{A} \geq \max_{\boldsymbol{\phi}\in \mathfrak{I}_{{\bf A},ns}^{\mu}\setminus\{\bf 0\}}\dfrac{\left<\displaystyle\int_{[{\bf 0},{\bf L}]}( {\bf v}\cdot e)|\boldsymbol{\phi}|^2 dx\right>}{\displaystyle\int_{[{\bf 0},{\bf L}]}|\boldsymbol{\phi}|^2 dx}.
\end{align*} 
where 
\begin{align*}
\mathfrak{I}_{{\bf A},ns}^{\mu}:=\left\{\boldsymbol{\phi}=(\phi_1,\phi_2)\in \mathfrak{I}:\sum_{i=1,2}\left<\int_{[{\bf 0},{\bf L}]} D\phi_i\cdot   {\bf A}^iD\phi_i dx\right>\leq \sum_{i=1,2}\left<\int_{[{\bf 0},{\bf L}]} {\mu}_{ii} \phi_i^2dx\right>  \right\}.
\end{align*}
In both cases, since ${\bf v}$ has zero mean, we deduce that
\begin{align*}
\liminf\limits_{A\rightarrow \infty}\dfrac{ {c}^*(e,A)}{A} \geq 0 \text{ implies that }\limsup\limits_{A\rightarrow \infty}\dfrac{ {c}^*(e,A)}{A} \geq 0.
\end{align*}

\end{remark}

Next, we assume \eqref{regularityofA} and $T_0$-periodicity. Recall
\begin{align*}
 \mathfrak{I}_{\text{time},p}=\Big\{\boldsymbol{\phi}\in [L^2\bigl(0,T_0; H^1_{\mathrm{per}}([{\bf 0},{\bf L}])\bigr)]^2 :
\boldsymbol{\phi}\ \text{is $T_0$- and $\mathbf{L}$-periodic},\
\mathbf{v}(t,x)\cdot D\phi_i=0,\ \forall i=1,2
\Big\},
\end{align*}
and the norm for $[L^2\bigl(0,T_0; H^1_{\mathrm{per}}([{\bf 0},{\bf L}])]^2$ 
\begin{align*}
\|{\bf u}\|^2 =\|{\bf u}\|^2_2+ \|D{\bf u}\|_{2,{\bf A}^1,{\bf A}^2}^2,~{\bf u} =(u_1,u_2) \in [L^2\bigl(0,T_0; H^1_{\mathrm{per}}([{\bf 0},{\bf L}])]^2
\end{align*}
where $\|D{\bf u}\|_{2,{\bf A}^1,{\bf A}^2}^2= \displaystyle\int_{[0,T_0]\times[{\bf 0},{\bf L}]}Du_1 \cdot {\bf A}^1Du_1 dx dt + \displaystyle\int_{[0,T_0]\times[{\bf 0},{\bf L}]}Du_2\cdot {\bf A}^2Du_2 dx dt$.

By Remarks \ref{periodicity} and \ref{remark51}, for $\lambda \in \R$, recall that $(k(\lambda,A),{\bf v}_\lambda)$ is a principal eigen-pair of \eqref{eigen}, that is,
\begin{align}\label{almosteigenvalue1}
\begin{cases}
\partial_t v_{\lambda,1}
&= \text{div} \big(  {\bf A}^1 D v_{\lambda,1} \big)
- \big( A  {\bf v} + 2\lambda  {\bf A}^1e \big) \cdot D v_{\lambda,1}+ \lambda A ( {\bf v}\cdot e)v_{\lambda,1} \\
&+ ( {\mu}_{11}  + \lambda^2\, e \cdot  {\bf A}^1 e)v_{\lambda,1}+ {\mu}_{12}v_{\lambda,2}-k(\lambda,A)v_{\lambda,1},\\
\partial_t v_{\lambda,2}
&= \text{div} \big(  {\bf A}^2  D v_{\lambda,2} \big)
- \big( A  {\bf v}  + 2\lambda  {\bf A}^2 e \big) \cdot D v_{\lambda,2}+\lambda A ( {\bf v} \cdot e) v_{\lambda,2}\\
&+( {\mu}_{22}  + \lambda^2\, e \cdot  {\bf A}^2 e)v_{\lambda,2}+ {\mu}_{21} v_{\lambda,1}-k(\lambda,A)v_{\lambda,2}.
\end{cases}
\end{align}
Without loss of generality, we assume that $\|{\bf v}_
\lambda\|_2=1$.

\begin{proof}[Proof of Theorem~\ref{limitspeed1}\textup{(b)}, upper bound]

If
$\limsup_{A\to\infty}c^*(e,A)/A=0$, the asserted upper
bound follows by choosing the constant test function
$\boldsymbol\phi=\mathbf1$ in $\mathfrak I_{\mathrm{time},p}$:
by the zero spatial mean of $\mathbf v$, its quotient is zero.
We may therefore assume $\limsup_{A\to\infty}c^*(e,A)/A>0$.
Choose \(\omega\) such that
\begin{align*}
\limsup\limits_{A\rightarrow \infty}\dfrac{ {c}^*(e,A)}{A}>\omega>0.
\end{align*}
Then, there exists $A_n\rightarrow \infty$ as $n\rightarrow \infty$ such that $ {c}^*(e,A_n)\geq \omega A_n.$ Based on Theorem \ref{minimalspeed}, one has $k(\lambda,A_n)\geq  \omega \lambda  A_n.$ For any $\varepsilon>0$, choose $\lambda_n=(\varepsilon A_n)^{-1}$ and $k_n = k(\lambda_n,A_n)$. It then follows from \eqref{uniformcmean} and \eqref{estimatec} that 
\begin{align}\label{boundedofeige}
\varepsilon^{-1}\omega\leq k_n \leq C\lambda_n^2 + \varepsilon^{-1} \left<  \sup\limits_{x\in \R^N}( {\bf v}\cdot e)\right> +\left< \mu_{\infty}\right>.
\end{align}
Without loss of generality, we assume that $\lambda_n \leq 1$. Let ${\bf v}_n(t,x)$ be the solution of \eqref{almosteigenvalue1} associated with $A=A_n$ and $\lambda=\lambda_n$ for each $n\in \N$. Thus, we get
\begin{align}\label{almosteigenvalueperi}
\begin{cases}
\partial_t v_{n,1}
&= \text{div} \big(  {\bf A}^1 D v_{n,1} \big)
- \big( A_n  {\bf v} + 2\lambda_n  {\bf A}^1e \big) \cdot D v_{n,1}+ \varepsilon^{-1}( {\bf v}\cdot e)v_{n,1}\\
& + ( {\mu}_{11}  + \lambda^2_n\, e \cdot  {\bf A}^1 e)v_{n,1}+ {\mu}_{12}v_{n,2}-k_nv_{n,1},\\
\partial_t v_{n,2}
&= \text{div} \big(  {\bf A}^2 D v_{n,2} \big)
- \big( A_n  {\bf v} + 2\lambda_n  {\bf A}^2e \big) \cdot D v_{n,2}+ \varepsilon^{-1}( {\bf v}\cdot e)v_{n,2}\\
& + ( {\mu}_{22}  + \lambda^2_n\, e \cdot  {\bf A}^2 e)v_{n,2}+ {\mu}_{21}v_{n,1}-k_nv_{n,2}.
\end{cases}
\end{align}
For $i,j=1,2$, multiplying the $i$-th equation by $v_{n,i}$ and integrating over $Q:=[0,T_0]\times [{\bf 0},{\bf L}]$, 
one has
\begin{align}\label{estimateeii}
\begin{split}
 \|D{\bf v}_n\|_{2,{\bf A}^1,{\bf A}^2}^2 + k_n &= \sum_{i=1,2}\int_{Q}\mu_{ii}v_{n,i}^2  dxdt+\sum_{\substack{i,j=1,2\\ i\neq j}}
\int_{Q}
\mu_{ij}\, v_{n,i}v_{n,j}\, dx\,dt\\
&+\lambda_n^2\sum_{i=1,2}\int_{Q}\, (e \cdot  {\bf A}^i e)  v_{n,i}^2dx dt +
\varepsilon^{-1}\int_{Q} ( {\bf v}\cdot e) \left|{\bf v}_n\right|^2dxdt,\\
\end{split}
\end{align}
which yields
\begin{align}\label{boundedabove123}
\left\| D{\bf v}_n\right\|_{2,{\bf A}^1,{\bf A}^2}^2
+ k_n
\leq   \sup_{i,j=1,2}\left\| \mu_{ii}\right\|_\infty  +   \sup_{i,j=1,2,i\neq j}\left\|\mu_{ij} \right\|_\infty 
+ \varepsilon^{-1}\left\|({\bf v}\cdot e)\right\|_\infty 
+ C,
\end{align}
where $C>0$ is in \eqref{Aassump}. Here, we use $2ab\leq (a^2+b^2)$.

We next derive a uniform bound for the spatial $L^1$-mass of
${\bf v}_n$. Set
\(
m_n(t)
:=
\displaystyle\sum_{i=1,2}
\int_{[{\bf 0},{\bf L}]}
v_{n,i}(t,x)\,dx .
\)
Since ${\bf v}_n$ is a positive principal eigenfunction,
$m_n(t)>0$ for all $t\in\mathbb R$. Integrating the equations
in \eqref{almosteigenvalue1} with respect to $x$ and using again
the spatial periodicity,
$\nabla_x\cdot{\bf v}=0$, and
$\nabla_x\cdot({\bf A}^ie)=0$, we obtain
\begin{align*}
m_n'(t)
={}&
\sum_{i=1,2}
\int_{[{\bf 0},{\bf L}]}
\left[
\varepsilon^{-1}({\bf v}\cdot e)
+\mu_{ii}
+\lambda_n^2(e\cdot{\bf A}^ie)
-k_n
\right]
v_{n,i}\,dx
+
\int_{[{\bf 0},{\bf L}]}
\left(
\mu_{12}v_{n,2}
+\mu_{21}v_{n,1}
\right)\,dx .
\end{align*}
Therefore, for each fixed $\varepsilon>0$, there exists
$C_\varepsilon>0$, independent of $n$, such that $|m_n'(t)|
\leq
C_\varepsilon m_n(t),~t\in\mathbb R.$ Hence,
\begin{align}\label{masscomparison}
m_n(t)
\leq
e^{C_\varepsilon T_0}m_n(s),~
s,t\in[0,T_0].
\end{align}

On the other hand, by the Cauchy--Schwarz inequality,
\(
m_n(t)^2
\leq
2|[{\bf 0},{\bf L}]|
\displaystyle\int_{[{\bf 0},{\bf L}]}
|{\bf v}_n(t,x)|^2\,dx.
\)
Since $\|{\bf v}_n\|_2=1$, we obtain
\(
\displaystyle\int_0^{T_0}m_n(t)^2\,dt
\leq
2|[{\bf 0},{\bf L}]|,
\) which yields
\(
m_n(s_n)
\leq
\left(
\dfrac{2|[{\bf 0},{\bf L}]|}{T_0}
\right)^{1/2}
\) for some $s_n\in[0,T_0]$. Together with \eqref{masscomparison}, there exists $C_\varepsilon>0$ such that
\begin{align}\label{uniformL1vn}
\sup_{n\in\mathbb N}
\sup_{t\in[0,T_0]}
\sum_{i=1,2}
\|v_{n,i}(t,\cdot)\|_{L^1([{\bf 0},{\bf L}])}
\leq
C_\varepsilon .
\end{align}

Let
\(
q:=2+\frac{2}{N}.
\)
By the Gagliardo--Nirenberg inequality on the periodic cell,
there exists $C>0$ such that
\[
\|f\|_{L^q([{\bf 0},{\bf L}])}^{q}
\leq
C
\|Df\|_{L^2([{\bf 0},{\bf L}])}^{2}
\|f\|_{L^1([{\bf 0},{\bf L}])}^{2/N}
+
C
\|f\|_{L^1([{\bf 0},{\bf L}])}^{q}.
\]
Applying this inequality to $v_{n,i}(t,\cdot)$, integrating with
respect to $t$, and using \eqref{boundedabove123} and
\eqref{uniformL1vn}, we get
\begin{align}\label{higherintegrabilityvn}
\sup_{n\in\mathbb N}
\|v_{n,i}\|_{L^{2+2/N}(Q)}
<\infty,~i=1,2.
\end{align}

Set
\(
\rho_{n,i}:=v_{n,i}^2,~ i=1,2.
\)
It follows from \eqref{higherintegrabilityvn} that
$\{\rho_{n,i}\}_{n\in\mathbb N}$ is bounded in
$L^{1+1/N}(Q)$. Hence, up to a subsequence, there exist
$\rho_i\in L^{1+1/N}(Q)$, $\rho_i\geq0$, such that
\begin{align}\label{densityconvergence}
\rho_{n,i}
\rightharpoonup
\rho_i\text{ in }L^{1+1/N}(Q),~ i=1,2.
\end{align}
Here and below, the functions are identified with functions on the space--time torus, and hence the limits $\rho_i$ are $T_0$-periodic in $t$ and ${\bf L}$-periodic in $x$. Moreover, since \(1\in \big(L^{1+1/N}(Q)\big)^*,\) \eqref{densityconvergence} and $\|{\bf v}_n\|_2=1$ imply that
\begin{align}\label{densitynormalization}
\sum_{i=1,2}\int_Q\rho_i\,dxdt
=
\lim_{n\to\infty}
\int_Q|{\bf v}_n|^2\,dxdt
=1.
\end{align}

We next prove that $\rho_i$ is a first integral of ${\bf v}$.
Multiplying the $i$-th equation of
\eqref{almosteigenvalueperi} by $2v_{n,i}$ implies
\begin{align}\label{rhoequation}
\begin{split}
\partial_t\rho_{n,i}
={}&
\operatorname{div}
\left(
{\bf A}^iD\rho_{n,i}
\right)
-
2Dv_{n,i}\cdot{\bf A}^iDv_{n,i}
-
A_n{\bf v}\cdot D\rho_{n,i}-
2\lambda_n{\bf A}^ie\cdot D\rho_{n,i}
\\
&+
2
\left[
\varepsilon^{-1}({\bf v}\cdot e)
+\mu_{ii}
+\lambda_n^2(e\cdot{\bf A}^ie)
-k_n
\right]
\rho_{n,i}+
2\mu_{ij}v_{n,i}v_{n,j},~ i\neq j.
\end{split}
\end{align}
Indeed,
\(
2v_{n,i}
\operatorname{div}({\bf A}^iDv_{n,i})
=
\operatorname{div}({\bf A}^iD\rho_{n,i})
-
2Dv_{n,i}\cdot{\bf A}^iDv_{n,i}.
\) Dividing \eqref{rhoequation} by $A_n$, all the terms except
${\bf v}\cdot D\rho_{n,i}$ tend to zero in the sense of
distributions. To see this, observe first that
\[
\int_Q
|D\rho_{n,i}|\,dxdt
=
2\int_Q
v_{n,i}|Dv_{n,i}|\,dxdt
\leq
2
\|v_{n,i}\|_{L^2(Q)}
\|Dv_{n,i}\|_{L^2(Q)}
\leq
C_\varepsilon,
\]
whereas, by \eqref{boundedabove123}, $\left\| D {\bf v}_{n} \right\|_{2,{\bf A}^1,{\bf A}^2}^2\leq C_\varepsilon.$ Furthermore, $\{\rho_{n,i}\}$ is bounded in $L^1(Q)$,
$\{k_n\}$ is bounded for fixed $\varepsilon$, and
\(
\dfrac{\lambda_n}{A_n}
=
\dfrac{1}{\varepsilon A_n^2}
\rightarrow0.
\)
Therefore, passing to the limit of \eqref{rhoequation}/$A_n$ in the distributional sense, 
\begin{align}\label{rhofirstintegral}
{\bf v}(t,x)\cdot D\rho_i(t,x)=0
\quad\text{in }\mathcal D'(Q),~ i=1,2.
\end{align}

It remains to prove that $\sqrt{\rho_i}$ has the required
spatial $H^1$-regularity. 
By Mazur's lemma, for $p:=1+\frac1N$, there exist convex combinations
\[
\overline\rho_{m,i}
:=
\sum_{n=m}^{N_m}\alpha_{m,n}\rho_{n,i},~
\alpha_{m,n}\geq0,~
\sum_{n=m}^{N_m}\alpha_{m,n}=1,
\]
such that
\begin{align}\label{mazur-density}
\overline\rho_{m,i}\rightarrow\rho_i
\quad\text{strongly in }L^p(Q) \text{ as }m\to \infty.
\end{align}
Since $Q$ has finite measure, this convergence also holds
strongly in $L^1(Q)$. We claim that
$\{\sqrt{\overline\rho_{m,i}}\}_{m\in\mathbb N}$
is bounded in
$L^2(0,T_0;H^1_{\rm per}([{\bf0},{\bf L}]))$.
Indeed, since $\rho_{n,i}=v_{n,i}^2$ and
\(
D\rho_{n,i}=2v_{n,i}Dv_{n,i},
\) one has
\[
D\overline\rho_{m,i}
=
2\sum_{n=m}^{N_m}
\alpha_{m,n}v_{n,i}Dv_{n,i}.
\]
By the Cauchy--Schwarz inequality,
\begin{align*}
|D\overline\rho_{m,i}|^2
&\leq
4
\left(
\sum_{n=m}^{N_m}
\alpha_{m,n}v_{n,i}^2
\right)
\left(
\sum_{n=m}^{N_m}
\alpha_{m,n}|Dv_{n,i}|^2
\right)=
4\overline\rho_{m,i}
\sum_{n=m}^{N_m}
\alpha_{m,n}|Dv_{n,i}|^2.
\end{align*}
Therefore, $|D\sqrt{\overline\rho_{m,i}}|^2
\leq
\displaystyle\sum_{n=m}^{N_m}
\alpha_{m,n}|Dv_{n,i}|^2
\text{ in }Q.$
Consequently, by \eqref{boundedabove123},
\begin{align*}
&\int_Q
|D\sqrt{\overline\rho_{m,i}}|^2\,dxdt
\leq
\sum_{n=m}^{N_m}
\alpha_{m,n}
\int_Q|Dv_{n,i}|^2\,dxdt
\leq C_\varepsilon\text{ and }\int_Q
|\sqrt{\overline\rho_{m,i}}|^2\,dxdt
=
\int_Q\overline\rho_{m,i}\,dxdt
\leq1.
\end{align*}
Thus,
\(
\{\sqrt{\overline\rho_{m,i}}\}_{m\in\mathbb N}
\text{ is bounded in }
L^2(0,T_0;H^1_{\rm per}([{\bf0},{\bf L}])).
\)
Using \eqref{mazur-density}, we therefore obtain
\begin{align*}
\left\|
\sqrt{\overline\rho_{m,i}}
-
\sqrt{\rho_i}
\right\|_{L^2(Q)}^2
&\leq
\|\overline\rho_{m,i}-\rho_i\|_{L^1(Q)}
\rightarrow0.
\end{align*}
Here, we use $|\sqrt a-\sqrt b|^2\leq|a-b|$ for all $a,b\geq0$. Hence, $\sqrt{\overline\rho_{m,i}}
\rightarrow
\sqrt{\rho_i}
\text{ strongly in }L^2(Q).$ By weak compactness in
$L^2(0,T_0;H^1_{\rm per}([{\bf0},{\bf L}]))$,
up to a subsequence,
\[
\sqrt{\overline\rho_{m,i}}
\rightharpoonup
\sqrt{\rho_i}
\quad\text{in }
L^2(0,T_0;H^1_{\rm per}([{\bf0},{\bf L}])).
\]
Consequently,
\(
\phi_i:=\sqrt{\rho_i}
\in
L^2(0,T_0;H^1_{\rm per}([{\bf 0},{\bf L}])),~\forall i=1,2.
\) Furthermore, by weak lower semicontinuity,
\[
\int_Q|D\sqrt{\rho_i}|^2\,dxdt
\leq
\liminf_{m\to\infty}
\int_Q
|D\sqrt{\overline\rho_{m,i}}|^2\,dxdt
 <\infty.
\]

Next, set
\(
\boldsymbol{\phi}:=(\phi_1,\phi_2).
\)
Since $\rho_i=\phi_i^2$, \eqref{rhofirstintegral} implies
\(
2\phi_i\,{\bf v}\cdot D\phi_i=0
~\text{ in }Q.
\)
On the other hand, $D\phi_i=0$ a.e. on the set
$\{\phi_i=0\}$. Hence
\(
{\bf v}\cdot D\phi_i=0~\text{ in }Q,~ i=1,2,
\)
which yields $\boldsymbol{\phi}
\in
\mathfrak{I}_{\mathrm{time},p}.$ Moreover, by \eqref{densitynormalization}, $\|\boldsymbol{\phi}\|_{L^2(Q)}^2
=
\displaystyle\sum_{i=1,2}
\displaystyle\int_Q\rho_i\,dxdt
=1.$ In particular,
$\boldsymbol{\phi}\not\equiv{\bf0}$. We now return to \eqref{estimateeii}. Since the diffusion term is nonnegative,
\[
k_n
\leq
C_\mu
+
C_A\lambda_n^2
+
\varepsilon^{-1}
\int_Q
({\bf v}\cdot e)|{\bf v}_n|^2\,dxdt.
\]
Together with $k_n\geq\varepsilon^{-1}\omega$, this implies
\begin{align}\label{omegaestimate}
\omega
\leq
\varepsilon C_\mu
+
\varepsilon C_A\lambda_n^2
+
\int_Q
({\bf v}\cdot e)|{\bf v}_n|^2\,dxdt.
\end{align}
Since
\(
|{\bf v}_n|^2
=
\rho_{n,1}+\rho_{n,2}
\)
and ${\bf v}\cdot e\in L^\infty(Q)$, it follows from
\eqref{densityconvergence} that
\begin{align*}
\int_Q
({\bf v}\cdot e)|{\bf v}_n|^2\,dxdt
&\rightarrow
\sum_{i=1,2}
\int_Q
({\bf v}\cdot e)\rho_i\,dxdt
=
\int_Q
({\bf v}\cdot e)
|\boldsymbol{\phi}|^2\,dxdt.
\end{align*}
Letting $n\to\infty$ in
\eqref{omegaestimate} yields $\omega
\leq
\varepsilon C_\mu
+
\displaystyle\int_Q
({\bf v}\cdot e)
|\boldsymbol{\phi}|^2\,dxdt.$ Since $\boldsymbol{\phi}
\in
\mathfrak{I}_{\mathrm{time},p}$ and $\|\boldsymbol{\phi}\|_{L^2(Q)}=1$, one has
\begin{align*}
\omega
&\leq
\varepsilon C_\mu
+
\sup_{\boldsymbol{\psi}\in
\mathfrak{I}_{\mathrm{time},p}\setminus\{{\bf0}\}}
\frac{
\displaystyle
\int_Q
({\bf v}\cdot e)|\boldsymbol{\psi}|^2\,dxdt
}{
\displaystyle
\int_Q|\boldsymbol{\psi}|^2\,dxdt
}.
\end{align*}
Since $\varepsilon>0$ is arbitrary 
and $\omega$ can be chosen arbitrarily close to
\(
\limsup\limits_{A\to\infty}\dfrac{c^*(e,A)}{A},
\)
we conclude that
\begin{align*}
\limsup_{A\to\infty}
\frac{c^*(e,A)}{A}
\leq
\sup_{\boldsymbol{\psi}\in
\mathfrak{I}_{\mathrm{time},p}\setminus\{{\bf0}\}}
\frac{
\left\langle
\displaystyle
\int_{[{\bf 0},{\bf L}]}
({\bf v}\cdot e)
|\boldsymbol{\psi}|^2\,dx
\right\rangle
}{
\left\langle
\displaystyle
\int_{[{\bf 0},{\bf L}]}
|\boldsymbol{\psi}|^2\,dx
\right\rangle
}.
\end{align*}
This completes the proof.
\end{proof}

Next, we assume that ${\bf A}^i={\bf A}^i(x),{\bf q}_i={\bf q}_i(x),{\bf F}={\bf F}(x,u_1,u_2)$. By uniqueness of the principal eigenvalue, let $(k(\lambda,A),{\bf v}_\lambda)$ be a principal eigen-pair of \eqref{eigen} as follows
\begin{equation}\label{eigenspace}
\begin{cases}
\begin{aligned}
\operatorname{div}({\bf A}^1Dv_{\lambda,1})-(A{\bf v}+2\lambda{\bf A}^1e)\cdot Dv_{\lambda,1}
+\lambda A({\bf v}\cdot e)v_{\lambda,1}\\
\hspace{2em}+(\mu_{11}+\lambda^2e\cdot{\bf A}^1e)v_{\lambda,1}+\mu_{12}v_{\lambda,2}
&=k(\lambda,A)v_{\lambda,1},
\end{aligned}\\
\begin{aligned}
\operatorname{div}({\bf A}^2Dv_{\lambda,2})-(A{\bf v}+2\lambda{\bf A}^2e)\cdot Dv_{\lambda,2}
+\lambda A({\bf v}\cdot e)v_{\lambda,2}\\
\hspace{2em}+(\mu_{22}+\lambda^2e\cdot{\bf A}^2e)v_{\lambda,2}+\mu_{21}v_{\lambda,1}
&=k(\lambda,A)v_{\lambda,2}.
\end{aligned}
\end{cases}
\end{equation}
Without loss of generality, we assume that $\|{\bf v}_
\lambda\|_2^2=\displaystyle\int_{[{\bf 0},{\bf L}]}|{\bf v}_
\lambda|^2dx=1$. Recall $V=|[{\bf 0},{\bf L}]|$ and
\begin{align}\label{paramerke}
\kappa_e(\lambda):=\sup_{{\bf w}=(w_1,w_2)\in \mathfrak{I}\setminus \{\bf 0\}}
\left\{
\dfrac{
\lambda \displaystyle\int_{[{\bf 0},{\bf L}]}
({\bf v}\cdot e)\,|{\bf w}|^2 \, dx
+\sum_{\substack{i,j=1,2}}
\int_{[{\bf 0},{\bf L}]}
\mu_{ij}\, w_i\, w_j  \, dx -  \left\| D {\bf w} \right\|_{2,{\bf A}^1,{\bf A}^2}^2
}{
\|{\bf w}\|_2^2
}
\right\}.
\end{align}
The constant vector $(1,1)$ belongs to $\mathfrak I$. Taking this vector in the Rayleigh quotient defining $\kappa_e(\lambda)$ and using $\int_{[\mathbf0,\mathbf L]}\mathbf v\cdot e\,dx=0$ yields, for every $\lambda\in\mathbb R$, \[\kappa_e(\lambda)\geq\frac1{2V}\sum_{i,j=1}^2\int_{[\mathbf0,\mathbf L]}\mu_{ij}(x)\,dx>0.\]

\begin{lemma}[Convex quotients]\label{lem:convex-infimum}
Let $f_n:[0,\infty)\to\mathbb R$ be finite convex functions which
converge locally uniformly to a convex function $f$. Assume
$f_n(0)>0$ for all sufficiently large $n$ and $f(0)>0$, and set
$r=\inf_{\lambda>0}f(\lambda)/\lambda\ge0$.
Then,
\begin{align*}
\lim_{n\to\infty}\inf_{\lambda>0}\frac{f_n(\lambda)}\lambda=r.
\end{align*}
\end{lemma}
\begin{proof}
Testing at any fixed $\lambda>0$ gives the upper bound. For the
lower bound fix $\eta>0$. Choose $\delta>0$ such that
$f(\lambda)\ge f(0)/2$ on $[0,\delta]$ and
$f(0)/(4\delta)>\max\{r-\eta,0\}$. Local uniform convergence
implies $f_n(\lambda)\ge f(0)/4$ on this interval for all large
$n$. Choose $L>\delta$ with $f(0)/L<\eta/3$.
By definition of $r$, 
\begin{align*}
\frac{f(L)-f(0)}L\ge r-\frac{f(0)}L>r-\eta/3.
\end{align*}
Local uniform convergence gives, for sufficiently large $n$,
$(f_n(L)-f_n(0))/L\ge r-2\eta/3$.
For $\lambda\ge L$, convexity and $f_n(0)>0$ yield
\begin{align*}
\frac{f_n(\lambda)}\lambda
&\ge\frac{f_n(\lambda)-f_n(0)}\lambda
\ge\frac{f_n(L)-f_n(0)}L\ge r-\eta.
\end{align*}
On $[\delta,L]$, local uniform convergence gives
$f_n(\lambda)/\lambda\ge r-\eta$ for large $n$; on
$(0,\delta]$ the preceding positive lower bound gives the
same inequality. Thus,
$\liminf_n\inf_{\lambda>0}f_n(\lambda)/\lambda\ge r-\eta$.
Letting $\eta\downarrow0$ proves the assertion, including $r=0$.
\end{proof}

\begin{lemma}\label{lemm17}
Assume that ${\bf A}^i={\bf A}^i(x),{\bf q}_i={\bf q}_i(x),{\bf F}={\bf F}(x,u_1,u_2)$ and $\mu_{12}=\mu_{21}$. Then, \eqref{paramerke} has a maximizer and
\begin{align*}
\lim_{A\rightarrow \infty}\dfrac{c^*(e,A)}{A}=\lim_{A\rightarrow \infty}\inf_{\lambda>0}\dfrac{k(\lambda/A,A)}{\lambda}=\inf_{\lambda>0}\dfrac{\kappa_e(\lambda)}{\lambda}.
\end{align*}
Furthermore, the maximizer is unique up to a multiplicative constant.

\end{lemma}

\begin{proof}

Recall from Remark \ref{remark51} that
\begin{align*}
c^*(e,A)=\min_{\lambda>0}\dfrac{k(\lambda,A)}{\lambda} \text{ implies that }\dfrac{c^*(e,A)}{A}=\min_{\lambda>0}\dfrac{k(\lambda/A,A)}{\lambda}.
\end{align*}

Let \(\mathbf v_\lambda^A\) denote the normalized positive eigenfunction associated with \(k(\lambda/A,A)\). From \eqref{eigenspace}, one has
\begin{align}\label{eigenvalueoverA}
\begin{cases}
&\text{div} \big(  {\bf A}^1 D v_{\lambda,1}^A \big)
- \big( A  {\bf v} + \dfrac{2\lambda}{A}  {\bf A}^1e \big) \cdot D v_{\lambda,1}^A+ \lambda ( {\bf v}\cdot e)v_{\lambda,1}^A\\
& + \left[ {\mu}_{11}  + \left(\dfrac{\lambda}{A}\right)^2\, e \cdot  {\bf A}^1 e \right]v_{\lambda,1}^A+ {\mu}_{12}v_{\lambda,2}^A-k(\lambda/A,A)v_{\lambda,1}^A=0,\\
&\text{div} \big(  {\bf A}^2  D v_{\lambda,2}^A \big)
- \big( A  {\bf v}  + \dfrac{2\lambda}{A}  {\bf A}^2 e \big) \cdot D v_{\lambda,2}^A+\lambda  ( {\bf v} \cdot e) v_{\lambda,2}^A\\
&+\left[ {\mu}_{22}  + \left(\dfrac{\lambda}{A}\right)^2\, e \cdot  {\bf A}^2 e\right]v_{\lambda,2}^A+ {\mu}_{21} v_{\lambda,1}^A-k(\lambda/A,A)v_{\lambda,2}^A=0.
\end{cases}
\end{align}
For $i=1,2$, multiplying the $i$-th equation by $(v_{\lambda,i}^A)^{-1}$ and integrate over $[{\bf 0},{\bf L}]$, one has
\begin{align}\label{logestimate}
0\leq \left\|D\ln v_{\lambda,i}^A\right\|_{2,{\bf A}^i}^2
\leq Vk(\lambda/A,A).
\end{align}
To verify this bound, the periodic integration-by-parts identities are
\begin{align*}
\int_{[\mathbf0,\mathbf L]}
\frac{\operatorname{div}({\bf A}^i Dv_{\lambda,i}^A)}{v_{\lambda,i}^A}
=\int_{[\mathbf0,\mathbf L]}
D\ln v_{\lambda,i}^A\cdot{\bf A}^iD\ln v_{\lambda,i}^A,~
\int_{[\mathbf0,\mathbf L]}\mathbf v\cdot D\ln v_{\lambda,i}^A
=0,~
\int_{[\mathbf0,\mathbf L]}{\bf A}^i e\cdot D\ln v_{\lambda,i}^A=0.
\end{align*}
The last two equalities use, respectively,
$\operatorname{div}\mathbf v=0$ and
$\operatorname{div}({\bf A}^ie)=0$, together with periodicity.
Also $\int_{[\mathbf0,\mathbf L]}\mathbf v\cdot e=0$.
Thus the equation for $v_{\lambda,i}^A$, divided by
$v_{\lambda,i}^A$ and integrated, yields the exact relation
\begin{align*}
\left\|D\ln v_{\lambda,i}^A\right\|_{2,{\bf A}^i}^2
={}&Vk(\lambda/A,A)-\int_{[\mathbf0,\mathbf L]}
\left(\mu_{ii}+\frac{\lambda^2}{A^2}e\cdot{\bf A}^ie
+\mu_{ij}\frac{v_{\lambda,j}^A}{v_{\lambda,i}^A}\right)dx,
~ j\ne i.
\end{align*}
The integrand on the last line is nonnegative under the
standing assumptions. This proves \eqref{logestimate}.

 Similarly, multiplying the $i$-th equation by $v_{\lambda,i}^A$ yields that
\begin{align}\label{estimateei}
\begin{split}
 \left\| D {\bf v}_{\lambda}^A \right\|_{2,{\bf A}^1,{\bf A}^2}^2 + k(\lambda/A,A) &= \sum_{i=1,2}\int_{ [{\bf 0},{\bf L}]}\mu_{ii}(v_{\lambda,i}^A)^2  dx+\sum_{\substack{i,j=1,2\\ i\neq j}}
\int_{ [{\bf 0},{\bf L}]}
\mu_{ij}\, v_{\lambda,i}^A\, v_{\lambda,j}^A \, dx\\
&+\sum_{i=1,2}\int_{ [{\bf 0},{\bf L}]}\left[\left(\dfrac{\lambda}{A}\right)^2\, e \cdot  {\bf A}^i e + \lambda ( {\bf v}\cdot e) \right](v_{\lambda,i}^A)^2dx,\\
\end{split}
\end{align}
which implies that
\begin{align}\label{boundedabove12}
\left\| D {\bf v}_{\lambda}^A \right\|_{2,{\bf A}^1,{\bf A}^2}^2
+ k(\lambda/A, A)
\leq   \sup_{i,j=1,2}\left\| \mu_{ii}\right\|_\infty  + 2  \sup_{i,j=1,2,i\neq j}\left\|\mu_{ij} \right\|_\infty 
+ |\lambda |\left\|({\bf v}\cdot e)\right\|_\infty 
+ C \left( \dfrac{\lambda}{A} \right)^2,
\end{align}
where $C>0$ is in \eqref{Aassump}. Here, we use $2ab\leq (a^2+b^2)$.
It is well-known that there exist $A_n\rightarrow \infty$ and ${\bf w}_0=(w_{1,0},w_{2,0})\in [H^1_{per}]^2$ such that ${\bf v}_{\lambda}^{A_n} \rightarrow{\bf w}_0\text{ in }[L^2_{per}]^2,~D{\bf v}_{\lambda}^{A_n} \rightharpoonup D{\bf w}_0\text{ in }[L^2_{per}]^2.$ Moreover, dividing \eqref{eigenvalueoverA} by $A_n$ and passing to the limit, in the sense of distributions, yields
\(
{\bf v}\cdot Dw_{i,0}=0, \text{ in }[{\bf 0},{\bf L}],~ i=1,2.
\)
Hence, ${\bf w}_0\in \mathfrak I$. 

We first show that the drift term vanishes when tested against the limiting first integral. Since
$\nabla_x\cdot\mathbf v=0$ and
$\mathbf v\cdot Dw_{i,0}=0$ in distributions,
periodic integration by parts, initially for smooth
approximations and then by density, shows that
\begin{align*}
\int_{[\mathbf0,\mathbf L]}
A_n(\mathbf v\cdot Dv_{\lambda,i}^{A_n})w_{i,0}\,dx
=-A_n\int_{[\mathbf0,\mathbf L]}
 v_{\lambda,i}^{A_n}(\mathbf v\cdot Dw_{i,0})\,dx=0.
\end{align*}
The drift term vanishes for every $n$.
The remaining first-order term, with coefficient
$2\lambda/A_n$, is $O(A_n^{-1})$ because the eigenfunctions are
bounded in $H^1_{\rm per}$.
Next, the weak formulation of the system \eqref{eigenvalueoverA}, tested by ${\bf w}_0$, is
\begin{align*}
k(\lambda/A_n,A_n)+o(1)
&= -\sum_{i=1,2}\int_{[{\bf0},{\bf L}]}
Dv_{\lambda,i}^{A_n}\cdot{\bf A}^i Dw_{i,0}\,dx+\lambda\int_{[{\bf0},{\bf L}]}
({\bf v}\cdot e)|{\bf w}_0|^2\,dx\\
&\quad+\sum_{i,j=1,2}\int_{[{\bf0},{\bf L}]}
\mu_{ij}w_{i,0}w_{j,0}\,dx+o(1).
\end{align*}
where $o(1)$ tends to $0$ as $A_n\rightarrow \infty$. To compare the two energy identities, put
\begin{align*}
\mathcal E(u,z)
=\sum_{i=1}^2\int_{[\mathbf0,\mathbf L]}
Du_i\cdot\mathbf A^iDz_i\,dx.
\end{align*}
The strong $L^2$ convergence of $\mathbf v_\lambda^{A_n}$,
its unit normalization and the boundedness of $k(\lambda/A_n,A_n)$
show that the zero-order, reaction and transport integrals in the
identity tested by $\mathbf v_\lambda^{A_n}$ differ by $o(1)$
from the corresponding mixed terms in the identity tested by
$\mathbf w_0$. The large-drift terms vanish exactly in both
identities, the first by $\nabla\cdot\mathbf v=0$ and the
second by the first-integral condition. Consequently,
\begin{align*}
\mathcal E(\mathbf v_\lambda^{A_n},\mathbf v_\lambda^{A_n})
-\mathcal E(\mathbf v_\lambda^{A_n},\mathbf w_0)=o(1).
\end{align*}
Weak $H^1$ convergence implies
$\mathcal E(\mathbf v_\lambda^{A_n},\mathbf w_0)
\to\mathcal E(\mathbf w_0,\mathbf w_0)$.
We therefore have convergence of the gradient-energy norms.
Uniform ellipticity then implies strong convergence
of the gradients in $L^2$; combined with the strong $L^2$
convergence already established, this proves strong
$H^1_{\rm per}$ convergence. The weak convergence implies that
 $\left\| D {\bf v}_{\lambda}^{A_n} \right\|_{2,{\bf A}^1,{\bf A}^2}^2 \rightarrow \left\| D {\bf w}_0 \right\|_{2,{\bf A}^1,{\bf A}^2}^2.$ Thus, one has 
\begin{align}\label{strongconverges}
&{\bf v}_{\lambda}^{A_n} \rightarrow{\bf w}_0\text{ in }[H^1_{per}]^2.
\end{align}
Passing to a subsequence, \eqref{estimateei} implies that
\begin{align*}
k_0:=\lim_{n\rightarrow \infty}k(\lambda/A_n,A_n) &=\sum_{\substack{i,j=1,2}}
\int_{ [{\bf 0},{\bf L}]}
\mu_{ij}\, w_{i,0}w_{j,0}\,  \, dx+\lambda\int_{ [{\bf 0},{\bf L}]} ( {\bf v}\cdot e) \left|{\bf w}_0\right|^2dx - \left\| D {\bf w}_0 \right\|_{2,{\bf A}^1,{\bf A}^2}^2.
\end{align*}

Consider ${\bf w}=(w_1,w_2)\in  \mathfrak{I}\cap [L^\infty(  [{\bf 0},{\bf L}])]^2 $. Multiplying the $i$-th equation of \eqref{eigenvalueoverA} for $A=A_n$ by $w_{i}^2/v_{\lambda,i}^{A_n} $, one has
\begin{align}\label{estimatekmax}
\begin{split}
k(\lambda/A_n,A_n)\|{\bf w}\|_2^2 &= \sum_{\substack{i=1,2}}\int_{ [{\bf 0},{\bf L}]} w_i^2 D\ln v_{\lambda,i}^{A_n} \cdot {\bf A}^i  D\ln v_{\lambda,i}^{A_n} -2w_i Dw_i {\bf A}^i D\ln v_{\lambda,i}^{A_n} dx\\
&-\dfrac{2\lambda}{A_n} \sum_{i=1,2}\int_{ [{\bf 0},{\bf L}]} w_{i}^2  {\bf A}^i e \cdot D \ln v_{\lambda,i}^{A_n}dx + \lambda \int_{ [{\bf 0},{\bf L}]} ({\bf v}\cdot e)\, |{\bf w}|^2 dx\\
&+ \sum_{i,j=1,2} \int_{ [{\bf 0},{\bf L}]} \mu_{ij}\, \frac{v_{\lambda,j}^{A_n}}{v_{\lambda,i}^{A_n}}\, w_i^2 dx
+ \left(\frac{\lambda}{A_n}\right)^2\sum_{\substack{i=1,2}} \int_{ [{\bf 0},{\bf L}]} (e \cdot {\bf A}^i e)\, w^2_i dx.
\end{split}
\end{align}
The inequality used to pass to the limit follows directly from the exact identity above. Write $z_i=v_{\lambda,i}^{A_n}>0$. For each $i$ the two gradient terms satisfy, pointwise,
\begin{align*}
w_i^2D\log z_i\cdot {\bf A}^iD\log z_i
-2w_iDw_i\cdot{\bf A}^iD\log z_i
=|Dw_i-w_iD\log z_i|_{{\bf A}^i}^2-|Dw_i|_{{\bf A}^i}^2
\ge -|Dw_i|_{{\bf A}^i}^2.
\end{align*}
Because $\mu_{12}=\mu_{21}>0$, the two off-diagonal terms obey the pointwise estimate
\begin{align*}
\mu_{12}\left(\frac{z_2}{z_1}w_1^2+\frac{z_1}{z_2}w_2^2\right)
\ge 2\mu_{12}|w_1w_2|\ge 2\mu_{12}w_1w_2.
\end{align*}
The diagonal terms are unchanged. The term involving $2\lambda/A_n$ tends to zero: since ${\bf w}$ is bounded and the coefficients ${\bf A}^i$ are bounded, its absolute value is at most
\begin{align*}
\frac{C|\lambda|}{A_n}\sum_{i=1}^2
\|w_i\|_\infty^2\|D\log z_i\|_{L^2},
\end{align*}
which is $O(A_n^{-1})$ by \eqref{logestimate} and \eqref{boundedabove12}. Inserting these inequalities in \eqref{estimatekmax}, then letting $n\to\infty$, gives
\begin{align}\label{maxvalue}
k_0\|{\bf w}\|_2^2 & \geq \lambda \displaystyle\int_{ [{\bf 0},{\bf L}]}
({\bf v}\cdot e)\,|{\bf w}|^2 \, dx\
+\sum_{\substack{i,j=1,2}}
\int_{ [{\bf 0},{\bf L}]}
\mu_{ij}\, w_i\, w_j  \, dx -  \left\| D {\bf w} \right\|_{2,{\bf A}^1,{\bf A}^2}^2,
\end{align}

For $\mathbf w=(w_1,w_2)\in\mathfrak I$, set
$T_N(s)=\operatorname{sgn}(s)\min\{|s|,N\}$ and
$w_{i,N}=T_N(w_i)$. The Sobolev chain rule gives
\begin{align*}
Dw_{i,N}=\mathbf1_{\{|w_i|<N\}}Dw_i\quad\text{a.e.},
\qquad \mathbf v\cdot Dw_{i,N}=0.
\end{align*}
Thus $\mathbf w_N\in\mathfrak I\cap[L^\infty]^2$, and
$\mathbf w_N\to\mathbf w$ in $[H^1_{\rm per}]^2$ by dominated convergence.
The bounded-coefficient terms in \eqref{maxvalue} pass to the limit, so \eqref{maxvalue} holds for every
${\bf w}\in  \mathfrak{I}$, which implies that $k_0 = \kappa_e(\lambda)$ and ${\bf w}_0$ is a maximizer of \eqref{paramerke}. On the other hand, if we take another $B_n\rightarrow \infty$ so that \(k(\lambda/B_n,B_n) \rightarrow k_1\) for some $k_1$, then we can repeat the process so that
$k_1=k_0=\kappa_e(\lambda)$ and
\begin{align}\label{convergofeige}
\kappa_e(\lambda) = \lim_{A\rightarrow \infty}k(\lambda/A,A) = \sum_{\substack{i,j=1,2}}
\int_{ [{\bf 0},{\bf L}]}
\mu_{ij}\, w_{i,0}w_{j,0}\,  \, dx+\lambda\int_{ [{\bf 0},{\bf L}]} ( {\bf v}\cdot e) \left|{\bf w}_0\right|^2dx - \left\| D {\bf w}_0 \right\|_{2,{\bf A}^1,{\bf A}^2}^2.
\end{align}
Thus, $\kappa_e(\lambda)$ is convex as the supremum of affine
functions of $\lambda$.  Moreover, $\kappa_e(\lambda)>0$ for all
$\lambda$, by testing against the constant vector $(1,1)$ and
using $\int_{[{\bf0},{\bf L}]}\mathbf v\,dx=\mathbf0$. 
This further implies that \eqref{convergofeige} is uniform on any bounded interval. Hence, for $\delta>0$ sufficiently small and $M>0$ sufficiently large, there exists $K>0$ such that
\begin{align*}
k(\lambda/A,A)\geq \dfrac{1}{2V}\sum_{i,j=1}^2\int_{[{\bf0},{\bf L}]}\mu_{ij}(x)\,dx>0,~\forall \lambda \in [-{\delta}, M],~A>K.
\end{align*}
Thus, there exists $\delta_2>0$ so that
\begin{align*}
\inf_{\lambda>0}\dfrac{k(\lambda/A,A)}{\lambda}=\inf_{\lambda\in [\delta_2,\infty)}\dfrac{k(\lambda/A,A)}{\lambda},~\inf_{\lambda>0}\dfrac{\kappa_e(\lambda)}{\lambda}=\inf_{\lambda\in [\delta_2,\infty)}\dfrac{\kappa_e(\lambda)}{\lambda}.
\end{align*}

Next, we show that
\begin{align*}
\lim_{A\rightarrow \infty}\inf_{\lambda>0}\dfrac{k(\lambda/A,A)}{\lambda}=\inf_{\lambda>0}\dfrac{\kappa_e(\lambda)}{\lambda}.
\end{align*}
The convergence of the scalar infima now follows from
Lemma~\ref{lem:convex-infimum}. Indeed, define
$f_A(\lambda)=k(\lambda/A,A)$ and $f(\lambda)=\kappa_e(\lambda)$.
The preceding argument proves local uniform convergence.
Both functions are convex in $\lambda$. Moreover,
$f_A(0)>0$ uniformly in $A$ by the strict cooperativity and
the constant positive subsolution for the untilted periodic
operator; and $f(0)=\kappa_e(0)>0$ by testing the constant
vector $(1,1)$. The nonnegativity of
$\inf_{\lambda>0}\kappa_e(\lambda)/\lambda$ follows from
$\kappa_e(\lambda)>0$ for every $\lambda>0$. Applying
Lemma~\ref{lem:convex-infimum} along an arbitrary sequence
$A_n\to\infty$ proves
\begin{align*}
\lim_{A\to\infty}\inf_{\lambda>0}
\frac{k(\lambda/A,A)}\lambda
=\inf_{\lambda>0}\frac{\kappa_e(\lambda)}\lambda.
\end{align*}
In particular this argument also covers a limiting infimum
equal to zero; it does not require the minimizing exponents to
remain in a fixed compact interval.

To conclude the lemma, we show that any maximizer of
\eqref{paramerke} is a multiple of ${\bf w}_0$.
Denote ${\bf v}_n:={\bf v}_{\lambda}^{A_n}$.
From \eqref{strongconverges}, one has
${\bf v}_n\rightarrow{\bf w}_0$ strongly in $[H^1_{per}]^2$. Up to a subsequence, ${\bf v}_n\rightarrow{\bf w}_0$ and $D{\bf v}_n\rightarrow D{\bf w}_0$ a.e. Let $\mathcal{R}_\lambda({\bf w})$ denote the functional in the definition of $\kappa_e(\lambda)$. First, $w_{i,0}\not\equiv0$ for $i=1,2$. Indeed, ${\bf w}_0\geq{\bf 0}$ is a maximizer of \eqref{paramerke} with $\|{\bf w}_0\|_2=1$. If $w_{1,0}\equiv0$, with $\tau>0$, define ${\bf w}^{\tau}:= \dfrac{(\tau,w_{2,0})}{\sqrt{1+\tau^2V}}$. Then, ${\bf w}^{\tau}$ is admissible. Since both components are first integrals, ${\bf w}^{\tau}\in\mathfrak I$. Since $\mu_{12}=\mu_{21}>0$ and $w_{2,0}\geq0$ is nonzero, the two coupling terms contribute $2\tau B$. The gradient of the first component is zero. Therefore, direct calculation yields that
\begin{align*}
R_\lambda({\bf w}^{\tau})-\kappa_e(\lambda)
=
\frac{
2\tau B+\tau^2\bigl(C_\lambda-V\kappa_e(\lambda)\bigr)
}{
1+\tau^2V
}
>0
\end{align*}
for all sufficiently small $\tau>0$, where $B
:=
\displaystyle\int_{ [{\bf 0},{\bf L}]} \mu_{12}w_{2,0}\,dx>0,~
C_\lambda
:=
\displaystyle\int_{ [{\bf 0},{\bf L}]} \bigl(\mu_{11}+\lambda{\bf v}\cdot e\bigr)\,dx.$
This contradicts maximality. The case $w_{2,0}\equiv0$
is excluded similarly. 

Next, let ${\bf c}_n$ be the average of $\ln{\bf v}_n$.
From \eqref{logestimate} and Poincar\'e's inequality,
up to a subsequence,
$\ln{\bf v}_n-{\bf c}_n\rightarrow\boldsymbol{\omega}$
strongly in $[L^2_{per}]^2$, weakly in $[H^1_{per}]^2$,
and a.e. For each $i=1,2$, the set $\{w_{i,0}>0\}$ has positive
measure. Choose a point $x_i$ in this set at which the
above pointwise convergences hold. Then
\begin{align*}
c_{n,i}
&=\ln v_{n,i}(x_i)
  - \left(\ln v_{n,i}(x_i)-c_{n,i}\right)\rightarrow
\ln w_{i,0}(x_i)-\omega_i(x_i)=:c_i\in\R.
\end{align*}
Thus, ${\bf c}_n\rightarrow{\bf c}\in\R^2$ and
$\ln{\bf v}_n\rightarrow{\bf g}:=
\boldsymbol{\omega}+{\bf c}$ a.e., where
${\bf g}\in[H^1_{per}]^2$.
Since ${\bf v}_n\rightarrow{\bf w}_0$ almost everywhere
along this subsequence, exponentiation gives
${\bf w}_0=e^{\bf g}>{\bf0}$ almost everywhere, and hence
$\ln{\bf w}_0={\bf g}\in[H^1_{\rm per}]^2$. The strong
$H^1$ convergence in \eqref{strongconverges} permits, after
extracting the same subsequence further, almost-everywhere
convergence of both $\mathbf v_n$ and $D\mathbf v_n$.
At these points, positivity of the limit implies
\[
D\ln{\bf v}_n
=\frac{D{\bf v}_n}{{\bf v}_n}
\longrightarrow
\frac{D{\bf w}_0}{{\bf w}_0}=D\ln{\bf w}_0
\quad\text{almost everywhere}.
\]
The uniform $L^2$ bound from \eqref{logestimate} also implies
weak convergence of these logarithmic gradients in $L^2$.
No strong $L^2$ convergence of the quotients is used. In particular,
the bounds obtained so far do not imply convergence of their norms;
weak convergence, the almost-everywhere identification and the subsequent
rigidity argument are kept separate.

Let ${\bf w}$ be another nonzero maximizer of
\eqref{paramerke}. Replacing each component by its absolute value
preserves the gradient energies and diagonal terms and cannot
decrease the coupling contribution, because $\mu_{12}=\mu_{21}>0$.
Thus its componentwise absolute value is also a maximizer, and we
may first consider ${\bf w}\geq{\bf0}$. Define $W_i=(w_{0,i}^2+w_i^2)^{1/2},~ i=1,2.$ The chain rule for Sobolev functions and the definition of
$\mathfrak I$ imply that ${\bf W}\in\mathfrak I$ and
$\|{\bf W}\|_2^2=\|{\bf w}_0\|_2^2+\|{\bf w}\|_2^2$.
For each $i$, the Cauchy--Schwarz inequality in the scalar
product induced by ${\bf A}^i$ implies $DW_i\cdot {\bf A}^iDW_i
\le Dw_{0,i}\cdot {\bf A}^iDw_{0,i}
+Dw_i\cdot {\bf A}^iDw_i.$ Moreover, pointwise, $W_1W_2\ge w_{0,1}w_{0,2}+w_1w_2.$ The diagonal terms of the quadratic form are additive under
${\bf W}$ and the off-diagonal terms have the indicated
sign since $\mu_{12}=\mu_{21}>0$. Consequently,
\begin{align*}
\kappa_e(\lambda)\|{\bf W}\|_2^2
&\ge \int_{[\bf0,\bf L]}
\left[\lambda({\bf v}\cdot e)|{\bf W}|^2
+\sum_{i,j=1}^2\mu_{ij}W_iW_j\right]dx
-\|D{\bf W}\|_{2,{\bf A}^1,{\bf A}^2}^2\geq\kappa_e(\lambda)
\big(\|{\bf w}_0\|_2^2+\|{\bf w}\|_2^2\big).
\end{align*}
Both inequalities are equalities. The nonnegative difference of
the gradient energies therefore vanishes for each component.
Since $w_{0,i}>0$, the equality case of the
Cauchy--Schwarz inequality yields
\begin{align}\label{eq:v14-rigidity}
Dw_i=w_iD\log w_{0,i}\quad\text{ on }[\mathbf0,\mathbf L],~ i=1,2.
\end{align}
The identity holds on the set $\{w_{0,i}^2+w_i^2>0\}$,
which has full measure.

We now verify the positivity needed to use logarithms. Fix an index
$i$. Since $w_i\geq0$ and $w_i\not\equiv0$ (the preceding
coupling perturbation also applies to $\mathbf w$), there are
$\delta_i>0$ and a measurable set $E_i\subset[\mathbf0,\mathbf L]$
with $|E_i|>0$ and $w_i\geq\delta_i$ on $E_i$. Put
$u_{i,\varepsilon}=\log(w_i+\varepsilon)$,
$0<\varepsilon\le1$. By \eqref{eq:v14-rigidity},
\begin{align*}
D u_{i,\varepsilon}
=\frac{w_i}{w_i+\varepsilon}D\ln w_{0,i},
\qquad
\|D u_{i,\varepsilon}\|_2\le\|D\ln w_{0,i}\|_2.
\end{align*}
Let $m_{i,\varepsilon}$ be the spatial mean of
$u_{i,\varepsilon}$. The periodic Poincar\'e inequality implies
\begin{align*}
\|u_{i,\varepsilon}-m_{i,\varepsilon}\|_2\le C.
\end{align*}
Since $u_{i,\varepsilon}\ge\log\delta_i$ on $E_i$, we get $m_{i,\varepsilon}
\ge\log\delta_i-C|E_i|^{-1/2}.$ The inequality $u_{i,\varepsilon}\le w_i+1$ also bounds
$m_{i,\varepsilon}$ from above. Thus
$\{u_{i,\varepsilon}\}_{0<\varepsilon\le1}$ is bounded in
$H^1_{\rm per}$. Along a sequence $\varepsilon\downarrow0$ it
converges weakly in $H^1_{\rm per}$ and strongly in $L^2$;
pointwise convergence identifies its limit with $\log w_i$.
Consequently, $w_i>0$ almost everywhere and
$\log w_i\in H^1_{\rm per}$. Letting $\varepsilon\downarrow0$
in the gradient identity above, we obtain
\begin{align*}
D\log w_i=D\log w_{0,i}\quad\text{a.e.}
\end{align*}
Connectedness of the periodic cell then implies
$w_i=d_iw_{0,i}$ for some $d_i>0$, $i=1,2$.

It remains to prove $d_1=d_2$. Write
\begin{align*}
&a_i:=\int_{[\mathbf0,\mathbf L]}w_{0,i}^2\,dx>0,~B:=\int_{[\mathbf0,\mathbf L]}\mu_{12}w_{0,1}w_{0,2}\,dx>0,\\
&Q_i:=\int_{[\mathbf0,\mathbf L]}
\left((\lambda\mathbf v\cdot e+\mu_{ii})w_{0,i}^2
-Dw_{0,i}\cdot\mathbf A^iDw_{0,i}\right)dx.
\end{align*}
For $r>0$, define $f(r)=\mathcal R_\lambda(rw_{0,1},w_{0,2})=\dfrac{Q_1r^2+2Br+Q_2}{a_1r^2+a_2}.$ Since $\mathbf w_0$ is a maximizer, $f'(1)=0$, and hence
$Q_1a_2-Q_2a_1=B(a_1-a_2)$. Direct differentiation yields
\begin{align}\label{eq:v13-positive-critical}
f'(r)=\frac{2B(1-r)(a_1r+a_2)}{(a_1r^2+a_2)^2},
\qquad r>0.
\end{align}
Thus, $r=1$ is the unique positive critical point and the unique maximum on $(0,\infty)$. Since $\mathbf w$ is a maximizer and $\mathcal R_\lambda(\mathbf w)=\mathcal R_\lambda((d_1/d_2)w_{0,1},w_{0,2}),$ we conclude that $d_1=d_2$.

Finally, let $\mathbf w$ be an arbitrary real maximizer and
$\mathbf w_\pm=\max\{\pm\mathbf w,\mathbf0\}$ componentwise.
Both $\mathbf w_+$ and $\mathbf w_-$ belong to $\mathfrak I$.
The gradient energies split, while the cross-component reaction
terms with $\mu_{12}=\mu_{21}>0$ satisfy the corresponding
one-sided inequality. Thus, omitting a term whenever its
numerator and denominator both vanish,
\begin{align*}
\mathcal R_\lambda(\mathbf w)
&\le\frac{\|\mathbf w_+\|_2^2}{\|\mathbf w\|_2^2}
 \mathcal R_\lambda(\mathbf w_+)
+\frac{\|\mathbf w_-\|_2^2}{\|\mathbf w\|_2^2}
 \mathcal R_\lambda(\mathbf w_-)\le\kappa_e(\lambda)=\mathcal R_\lambda(\mathbf w).
\end{align*}
Every nonzero part is therefore a nonnegative maximizer and hence
strictly positive almost everywhere in both components by the
preceding argument. If both parts were nonzero, then for each
$i=1,2$ we would have $w_{+,i}>0$ and $w_{-,i}>0$ almost
everywhere, contradicting $w_{+,i}w_{-,i}=0$. Hence exactly one
part is nonzero. Every real maximizer therefore has a constant sign
and is a scalar multiple of $\mathbf w_0$.
\end{proof}

We are now in a position to prove the main results.
 
\begin{proof}[Proof of Theorem~\ref{limitspeed1}\textup{(b)}, stationary case]

Inequality $\geq$ follows immediately from Lemma \ref{lemm17} and the definition of $\kappa_e(\lambda)$. Next, let ${\bf w}_0(\lambda) = (w_{0,1}(\lambda), w_{0,2}(\lambda))$ be a normalized nonnegative maximizer of $\kappa_e(\lambda)$. To obtain the reverse inequality, it is sufficient to find $\lambda>0$ such that $\gamma(\lambda)= 0$,
where 
\begin{align*}
\gamma(\lambda) := \sum_{i=1,2}\int_{[{\bf 0},{\bf L}]}Dw_{0,i}(\lambda)\cdot   {\bf A}^iDw_{0,i}(\lambda)\, dx-  \sum_{i,j=1,2}\int_{[{\bf 0},{\bf L}]}{\mu}_{ij} w_{0,i}(\lambda)w_{0,j}(\lambda)\, dx.
\end{align*}

For this purpose, $\gamma(0)=-\kappa_e(0)<0$.  We first
verify continuity of $\gamma$. If $\lambda_n\to\lambda_0\ge0$,
the normalized maximizers $\mathbf w_n=\mathbf w_0(\lambda_n)$
are bounded in $[H^1_{per}]^2$, since their Rayleigh quotients equal
$\kappa_e(\lambda_n)$ and the latter is locally bounded.
After extraction, $\mathbf w_n\rightharpoonup\boldsymbol\omega$
in $[H^1_{per}]^2$ and $\mathbf w_n\to\boldsymbol\omega$
in $[L^2_{per}]^2$, where $\boldsymbol\omega\in\mathfrak I$
and $\|\boldsymbol\omega\|_2=1$. Define
\begin{align*}
\mathcal V_e(\mathbf w):=\int_{[\mathbf0,\mathbf L]}
(\mathbf v\cdot e)|\mathbf w|^2\,dx,~
\mathcal B(\mathbf w):=\sum_{i,j=1}^2
\int_{[\mathbf0,\mathbf L]}\mu_{ij}w_iw_j\,dx,~
\mathcal E(\mathbf w):=\|D\mathbf w\|_{2,\mathbf A^1,\mathbf A^2}^2.
\end{align*}
The $L^2$ convergence, boundedness of the coefficients and
continuity of $\kappa_e$ imply
\begin{align*}
\lim_{n\to\infty}\mathcal E(\mathbf w_n)
&=\lambda_0\mathcal V_e(\boldsymbol\omega)
 +\mathcal B(\boldsymbol\omega)-\kappa_e(\lambda_0).
\end{align*}
On the other hand, the variational inequality for
$\boldsymbol\omega$ and lower semicontinuity imply, respectively,
\begin{align*}
\mathcal E(\boldsymbol\omega)&\ge
\lambda_0\mathcal V_e(\boldsymbol\omega)
 +\mathcal B(\boldsymbol\omega)-\kappa_e(\lambda_0),~
\mathcal E(\boldsymbol\omega)\le
\liminf_{n\to\infty}\mathcal E(\mathbf w_n).
\end{align*}
Thus both inequalities are equalities. Uniform ellipticity now
implies strong $H^1$ convergence. By the simplicity of the positive
principal maximizer established above,
$\boldsymbol\omega=\mathbf w_0(\lambda_0)$ after the chosen
normalization. Therefore $\gamma$ is continuous.

Set
\begin{align*}
\mathcal M_e=\sup_{\boldsymbol\phi\in
\mathfrak I_{\bf A}^{\mu}\setminus\{\mathbf0\}}
\frac{\mathcal V_e(\boldsymbol\phi)}{\|\boldsymbol\phi\|_2^2},
\qquad
\mathcal V_e^*=\sup_{\boldsymbol\phi\in
\mathfrak I\setminus\{\mathbf0\}}
\frac{\mathcal V_e(\boldsymbol\phi)}{\|\boldsymbol\phi\|_2^2}.
\end{align*}
For every $\lambda>0$, the variational characterization implies
$\kappa_e(\lambda)/\lambda\ge\mathcal M_e$.  If
$\gamma(\lambda_0)=0$ for some $\lambda_0>0$, then
$\mathbf w_0(\lambda_0)\in\mathfrak I_{\bf A}^{\mu}$ and
\begin{align*}
\inf_{\lambda>0}\frac{\kappa_e(\lambda)}\lambda
\le\frac{\kappa_e(\lambda_0)}{\lambda_0}
=\mathcal V_e(\mathbf w_0(\lambda_0))\le\mathcal M_e.
\end{align*}

It remains to consider the possibility that
$\gamma(\lambda)<0$ for every $\lambda>0$. By the uniform
bound for the quadratic reaction form on the $L^2$ unit sphere,
\begin{align}\label{eq:v9-recession}
\lim_{\lambda\to\infty}\frac{\kappa_e(\lambda)}\lambda
=\mathcal V_e^*.
\end{align}
Indeed, the upper estimate follows from
$\mathcal B(\mathbf w)-\mathcal E(\mathbf w)\le C$
for $\|\mathbf w\|_2=1$; testing against any fixed
$\mathbf w\in\mathfrak I$ establishes the reverse estimate.
Since the positive principal eigenvalue is simple,
$\kappa_e$ is differentiable and the Hellmann--Feynman identity reads
\begin{align*}
\kappa_e'(\lambda)=\mathcal V_e(\mathbf w_0(\lambda)),~ \gamma(\lambda)=\lambda\kappa_e'(\lambda)-\kappa_e(\lambda).
\end{align*}
Consequently $\lambda\mapsto\kappa_e(\lambda)/\lambda$
is nonincreasing, and its infimum equals $\mathcal V_e^*$
by \eqref{eq:v9-recession}. Moreover, the normalized maximizer
belongs to $\mathfrak I_{\bf A}^{\mu}$, while
\begin{align*}
\frac{\kappa_e(\lambda)}\lambda
=\mathcal V_e(\mathbf w_0(\lambda))-
\frac{\gamma(\lambda)}\lambda
\le \mathcal V_e(\mathbf w_0(\lambda))+\frac C\lambda,
\end{align*}
because $-\gamma(\lambda)=\mathcal B(\mathbf w_0(\lambda))
-\mathcal E(\mathbf w_0(\lambda))\le C$.
Thus $\mathcal M_e\ge\mathcal V_e^*$ and equality follows.
If $\gamma$ becomes nonnegative at any point, continuity and
$\gamma(0)<0$ reduce the argument to the preceding zero case.
This completes the stationary variational identity.

Now, for $\alpha>0$, we consider the reaction term $\alpha{\bf F}$
and $c^*(e,A,\alpha{\bf F})$. Recall
$V:=|[{\bf 0},{\bf L}]|$ and $M_0=\left(\int_{[{\bf0},{\bf L}]}\mu_{ij}(x)\,dx\right)_{i,j=1,2}.$ We first use the stationary variational identity, applied to the
reaction $\alpha\mathbf F$, in the form
\begin{align*}
\mathcal L_e(\alpha)
:=\lim_{A\to\infty}\frac{c^*(e,A,\alpha\mathbf F)}{A}
=\sup_{\substack{\boldsymbol\phi\in\mathfrak I\setminus\{\mathbf0\},\\
\boldsymbol\phi\ge\mathbf0,\ \mathcal E(\boldsymbol\phi)
\le\alpha\mathcal B(\boldsymbol\phi)}}
\frac{\mathcal V_e(\boldsymbol\phi)}{\|\boldsymbol\phi\|_2^2}.
\end{align*}
The restriction to nonnegative functions does not change the
supremum. Indeed, componentwise absolute values preserve the
energy and the transport term, while they do not decrease the
reaction form, since the off-diagonal coefficients are positive.

We next prove the upper bound. Normalize an admissible
$\boldsymbol\phi$ by $\|\boldsymbol\phi\|_2=1$ and define
\begin{align*}
\mathbf a:=V^{-1}\int_{[\mathbf0,\mathbf L]}
\boldsymbol\phi\,dx,\qquad
\boldsymbol\omega:=\boldsymbol\phi-\mathbf a,\qquad
T:=\mathcal E(\boldsymbol\omega)^{1/2}.
\end{align*}
All constants below are independent of the admissible function
and of sufficiently small $\alpha>0$. Uniform boundedness of
$\mu_{ij}$ and the energy constraint imply
$T^2\le C\alpha$. The periodic Poincar\'e inequality implies
\begin{align*}
\|\boldsymbol\omega\|_2\le CT,\qquad
1=V|\mathbf a|^2+\|\boldsymbol\omega\|_2^2,\qquad
|\mathbf a|^2=V^{-1}+O(\alpha).
\end{align*}
Furthermore, the reaction form satisfies 
$\mathcal B(\boldsymbol\phi)
=\mathbf a\cdot M_0\mathbf a+O(T).$ In particular,
\begin{align}\label{eq:v21-energy-quadratic}
T^2\le\alpha\bigl(\mathbf a\cdot M_0\mathbf a+CT\bigr).
\end{align}
The normalization and nonnegativity of $\boldsymbol\phi$ imply
$\mathbf a\ge\mathbf0$. Since all entries of $M_0$ are
strictly positive, there is $m_0>0$, independent of $\alpha$,
such that $\mathbf a\cdot M_0\mathbf a\ge m_0$ for all
sufficiently small $\alpha$. Solving the quadratic inequality
\eqref{eq:v21-energy-quadratic} for $T\ge0$, we obtain
\begin{align}\label{eq:v21-energy-root}
\frac{T}{\sqrt\alpha}
\le\sqrt{\mathbf a\cdot M_0\mathbf a}+C\sqrt\alpha.
\end{align}
This estimate is valid also when $T=0$; no division by an
energy quantity tending to zero is required.

For $\xi\in\mathbb S^1$, write
\begin{align*}
Q_e(\xi):=\sup_{\substack{\boldsymbol\psi\in\mathfrak I,\\
\int_{[\mathbf0,\mathbf L]}\boldsymbol\psi\,dx=\mathbf0,\\
\mathcal E(\boldsymbol\psi)>0}}
\frac{\displaystyle\int_{[\mathbf0,\mathbf L]}
(\mathbf v\cdot e)(\xi\cdot\boldsymbol\psi)\,dx}
{\mathcal E(\boldsymbol\psi)^{1/2}}.
\end{align*}
The sign of $\boldsymbol\psi$ may be reversed, so $Q_e(\xi)\ge0$.
The Poincar\'e inequality and the boundedness of $\mathbf v$
show that $Q_e$ is finite and continuous on $\mathbb S^1$.
If $T>0$, put $\xi=\mathbf a/|\mathbf a|$. Since
$\int_{[\mathbf0,\mathbf L]}(\mathbf v\cdot e)\,dx=0$,
\begin{align*}
\mathcal V_e(\boldsymbol\phi)
&=2\int_{[\mathbf0,\mathbf L]}
(\mathbf v\cdot e)(\mathbf a\cdot\boldsymbol\omega)\,dx
+\int_{[\mathbf0,\mathbf L]}
(\mathbf v\cdot e)|\boldsymbol\omega|^2\,dx\le 2|\mathbf a|Q_e(\xi)T+CT^2.
\end{align*}
Combining this estimate with \eqref{eq:v21-energy-root} and
$|\mathbf a|^2=V^{-1}+O(\alpha)$, we obtain
\begin{align}\label{eq:v21-small-upper}
\frac{\mathcal V_e(\boldsymbol\phi)}{2\sqrt\alpha}
\le\frac1V\sqrt{\xi\cdot M_0\xi}\,Q_e(\xi)
+C\sqrt\alpha.
\end{align}
When $T=0$, the left-hand side vanishes, so the same bound
holds. The constant in \eqref{eq:v21-small-upper} is uniform
over all normalized admissible functions. Taking their supremum
yields
\begin{align*}
\limsup_{\alpha\downarrow0}
\frac{\mathcal L_e(\alpha)}{2\sqrt\alpha}
\le\frac1V
\sup_{\substack{\xi\in\mathbb S^1,\\\xi_i\ge0}}
\sqrt{\xi\cdot M_0\xi}\,Q_e(\xi).
\end{align*}

We now establish the reverse inequality. Fix $\xi\in\mathbb S^1$
with $\xi_i>0$, and choose a bounded first integral
$\boldsymbol\psi\in\mathfrak I\setminus\mathbb R^2$ such that
\begin{align*}
\int_{[\mathbf0,\mathbf L]}\boldsymbol\psi\,dx=\mathbf0,~ \mathcal E(\boldsymbol\psi)=1.
\end{align*}
Since $m(\xi):=\xi\cdot M_0\xi>0$, the scalar equation $\delta^2
=\alpha\mathcal B(\xi+\delta\boldsymbol\psi)$ has a positive solution for sufficiently small $\alpha$, with $\delta=\sqrt{\alpha m(\xi)}+O(\alpha).$ The boundedness of $\boldsymbol\psi$ and the strict positivity
of the components of $\xi$ ensure that
$\mathbf W_\alpha:=\xi+\delta\boldsymbol\psi\ge\mathbf0$
when $\alpha$ is small enough. The vector $\mathbf W_\alpha$
is therefore admissible for the variational formula, and
\begin{align*}
\mathcal V_e(\mathbf W_\alpha)
=2\delta\int_{[\mathbf0,\mathbf L]}
(\mathbf v\cdot e)(\xi\cdot\boldsymbol\psi)\,dx
+O(\delta^2),~
\|\mathbf W_\alpha\|_2^2
=V+\delta^2\|\boldsymbol\psi\|_2^2=V+O(\alpha).
\end{align*}
Consequently,
\begin{align*}
\liminf_{\alpha\downarrow0}
\frac{\mathcal L_e(\alpha)}{2\sqrt\alpha}
\ge\frac{\sqrt{\xi\cdot M_0\xi}}{V}
\int_{[\mathbf0,\mathbf L]}
(\mathbf v\cdot e)(\xi\cdot\boldsymbol\psi)\,dx.
\end{align*}
This inequality extends to arbitrary energy-normalized
$\boldsymbol\psi\in\mathfrak I\setminus\mathbb R^2$:
componentwise truncation followed by subtraction of spatial
means preserves the first-integral condition, and the
truncated functions converge in $[H^1_{\rm per}]^2$.
We may then take the supremum over $\boldsymbol\psi$.
The strict positivity restriction on $\xi$ may be removed
by approximating nonnegative unit vectors by strictly
positive ones and using the continuity of $Q_e$.

Finally, for an arbitrary $\xi\in\mathbb S^1$, componentwise
sign changes of $\boldsymbol\psi$ show that
$Q_e(\xi)=Q_e(|\xi|)$. As $M_0$ is symmetric with positive
off-diagonal entries, we have
$\xi\cdot M_0\xi\le |\xi|\cdot M_0|\xi|$.
Thus, the supremum over nonnegative unit vectors equals the
supremum stated in \eqref{eq:v5-small}. The lower and upper
bounds now coincide and establish the asserted iterated limit. This completes the proof.
\end{proof}

\subsection{A time-periodic counterexample}\label{counter}

We next show that the equalities associated with
$\mathfrak{I}_{{\bf A}}^\mu$ and
$\mathfrak{I}_{{\bf A},time}^\mu$ do not hold in general,
even if $\mu_{12}=\mu_{21}$.

Now, we construct a counterexample. Consider $T_0=1$ and, in $\R^2$, choose the periodic cell
$[0,1]^2$ and
\begin{align*}
{\bf v}(x,y)=(\sin(2\pi y),0),\quad
e=(1,0),\quad
{\bf A}^1={\bf A}^2=I_2,\quad
{\bf q}_1={\bf q}_2=A{\bf v}.
\end{align*}
For $\delta\in(0,1)$, set
$a_\delta(t):=\delta\cos(2\pi t)$ and consider
\begin{align*}
\begin{cases}
F_1(t,x,u_1,u_2)=(2+a_\delta(t))u_1+u_2-(3+a_\delta(t))u_1^2,\\
F_2(t,x,u_1,u_2)=u_1+(2-a_\delta(t))u_2-(3-a_\delta(t))u_2^2.
\end{cases}
\end{align*}
Then
\begin{align*}
M_\delta(t):=(\mu_{ij}(t))
=
\begin{pmatrix}
2+a_\delta(t)&1\\
1&2-a_\delta(t)
\end{pmatrix},
\end{align*}
and, in particular, $\mu_{12}=\mu_{21}=1$. Furthermore,
\begin{align*}
F_1(t,x,s,s)=(3+a_\delta(t))s(1-s)>0,~
F_2(t,x,s,s)=(3-a_\delta(t))s(1-s)>0,~\forall s\in (0,1).
\end{align*}

Let $\nu(\alpha)$ denote the principal eigenvalue of
\begin{align}\label{counter1d}
\begin{cases}
\phi''+\alpha\sin(2\pi y)\phi=\nu(\alpha)\phi,\quad y\in\R,\\
\phi(y+1)=\phi(y),\quad \phi(y)>0.
\end{cases}
\end{align}
For $r>0$, consider the two-dimensional periodic eigenvalue problem associated with the diffusion matrix $I_2$ and the advection term $A{\bf v}$. Since ${\bf v}(x,y)=(\sin(2\pi y),0)$, its positive principal eigenfunction is independent of $x$ and the corresponding eigenvalue problem reduces to \eqref{counter1d}. Therefore, by the variational characterization in \cite{Zlatos2010}, one has
\begin{align}\label{Srvariational}
\mathcal{S}(r)
:=
\inf_{\alpha>0}
\dfrac{\nu(\alpha)+r}{\alpha}
=
\max_{\substack{
\phi\in H^1_{per}(0,1)\setminus\{0\}\\
\int_0^1|\phi'|^2\,dy\leq r\int_0^1\phi^2\,dy}}
\dfrac{\displaystyle\int_0^1\sin(2\pi y)\phi^2\,dy}
{\displaystyle\int_0^1\phi^2\,dy}.
\end{align}

We first calculate the large-drift limit. Let $\rho_\delta$ be the
principal Floquet exponent associated with
\begin{align}\label{counterfloquet}
{\bf p}'_\delta(t)
=
M_\delta(t){\bf p}_\delta(t)
-
\rho_\delta{\bf p}_\delta(t),
\end{align}
where ${\bf p}_\delta>0$ is $1$-periodic. Let
$\phi_{\lambda A}$ be the positive eigenfunction of
\eqref{counter1d} associated with $\alpha=\lambda A$.
Since the coefficients are independent of $x$, the product $\mathbf p_\delta(t)\phi_{\lambda A}(y)$ is positive and satisfies the periodic tilted eigenvalue problem. Uniqueness of the positive principal eigenvalue then implies $k(\lambda,A)
=
\lambda^2+\nu(\lambda A)+\rho_\delta$,
\begin{align*}
\dfrac{c^*(e,A)}{A}
=
\min_{\lambda>0}\dfrac{k(\lambda,A)}{\lambda A}
=
\min_{\alpha>0}
\left(
\dfrac{\alpha}{A^2}
+
\dfrac{\nu(\alpha)+\rho_\delta}{\alpha}
\right),
\end{align*}
and
\begin{align}\label{counterlimit}
\lim_{A\rightarrow\infty}
\dfrac{c^*(e,A)}{A}
=
\mathcal{S}(\rho_\delta).
\end{align}

We next show that $3<\rho_\delta<
2+\sqrt{1+\delta^2}.$ Set $q_+:=\dfrac{1}{\sqrt{2}}(1,1),~
q_-:=\dfrac{1}{\sqrt{2}}(1,-1),~
{\bf p}_\delta(t)=X(t)q_++Y(t)q_-.$ Then \eqref{counterfloquet} becomes $X'=(3-\rho_\delta)X+\delta\cos(2\pi t)Y,~
Y'=\delta\cos(2\pi t)X+(1-\rho_\delta)Y.$ Since ${\bf p}_\delta>0$, one has
$X>0$ and $|Y|<X$. Set $z:=Y/X$. Then $|z|<1$ and
\begin{align}\label{counterz}
z'
=
\delta\cos(2\pi t)(1-z^2)-2z.
\end{align}
Furthermore, $\dfrac{X'}{X}
=
3-\rho_\delta+\delta\cos(2\pi t)z.$ Taking the temporal mean, using the periodicity of $X$ and $z$ and multiplying \eqref{counterz} by $z/(1-z^2)$, we obtain
\begin{align}\label{counterrho1}
\rho_\delta
&=
3+\delta\int_0^1\cos(2\pi t)z(t)\,dt =
3+2\int_0^1\frac{z^2}{1-z^2}\,dt
>3,
\end{align}
where the last inequality follows from $z\not\equiv0$. Moreover, from \eqref{counterfloquet}, $\dfrac{1}{2}\dfrac{d}{dt}\ln|{\bf p}_\delta|^2
=
\dfrac{{\bf p}_\delta^TM_\delta(t){\bf p}_\delta}
{|{\bf p}_\delta|^2}
-\rho_\delta.$ Taking the temporal mean, we obtain
\begin{align*}
\rho_\delta
=
\int_0^1
\dfrac{{\bf p}_\delta^TM_\delta(t){\bf p}_\delta}
{|{\bf p}_\delta|^2}dt
\leq
\int_0^1\lambda_{\max}(M_\delta(t))dt.
\end{align*}
Since $\lambda_{\max}(M_\delta(t))
=
2+\sqrt{1+\delta^2\cos^2(2\pi t)}
<
2+\sqrt{1+\delta^2}$ for almost every $t\in(0,1)$, it follows that $\rho_\delta
<
2+\sqrt{1+\delta^2}.$ 

We first calculate the part associated with the
time-independent class $\mathfrak{I}_{{\bf A}}^\mu$.
Since $\boldsymbol{\phi}\in\mathfrak{I}_{{\bf A}}^\mu$
is independent of $t$ and $\left<M_\delta(t)\right>
=
\begin{pmatrix}
2&1\\
1&2
\end{pmatrix}
=:M_0,$ one has
\begin{align*}
\sum_{i=1,2}
\int_{[0,1]^2}|D\phi_i|^2dxdy
\leq
\int_{[0,1]^2}
\boldsymbol{\phi}^TM_0\boldsymbol{\phi}\,dxdy
\leq
3\int_{[0,1]^2}
|\boldsymbol{\phi}|^2dxdy.
\end{align*}
By the Rayleigh characterization of $\nu(\alpha)$, for every
$\alpha>0$,
\begin{align*}
\alpha
\int_{[0,1]^2}
\sin(2\pi y)|\boldsymbol{\phi}|^2dxdy
-
\sum_{i=1,2}
\int_{[0,1]^2}|D\phi_i|^2dxdy
\leq
\nu(\alpha)
\int_{[0,1]^2}
|\boldsymbol{\phi}|^2dxdy.
\end{align*}
It follows that $\dfrac{
\left<
\displaystyle\int_{[0,1]^2}
({\bf v}\cdot e)|\boldsymbol{\phi}|^2dxdy
\right>
}{
\displaystyle\int_{[0,1]^2}
|\boldsymbol{\phi}|^2dxdy
}
\leq
\dfrac{\nu(\alpha)+3}{\alpha}.$ Taking the infimum with respect to $\alpha>0$, we derive
\begin{align}\label{counterIupper}
\max_{\boldsymbol{\phi}\in
\mathfrak{I}_{{\bf A}}^\mu\setminus\{{\bf 0}\}}
\dfrac{
\left<
\displaystyle\int_{[0,1]^2}
({\bf v}\cdot e)|\boldsymbol{\phi}|^2dxdy
\right>
}{
\displaystyle\int_{[0,1]^2}
|\boldsymbol{\phi}|^2dxdy
}
\leq
\mathcal{S}(3).
\end{align}

Conversely, choose
$\psi\in H^1_{per}(\R)\setminus\{0\}$ such that $\displaystyle\int_0^1|\psi'|^2dy
\leq
3\displaystyle\int_0^1\psi^2dy$ and set
$\phi_1(x,y)=\phi_2(x,y):=\psi(y)$.
Since ${\bf v}\cdot D\phi_i
=
\sin(2\pi y)\partial_x\phi_i=0$ and
\begin{align*}
\sum_{i=1,2}
\int_{[0,1]^2}|D\phi_i|^2dxdy
=
2\int_0^1|\psi'|^2dy
\leq
6\int_0^1\psi^2dy
=
\sum_{i,j=1,2}
\left<
\int_{[0,1]^2}
\mu_{ij}\phi_i\phi_jdxdy
\right>,
\end{align*}
one has
$\boldsymbol{\phi}\in\mathfrak{I}_{{\bf A}}^\mu$.
Furthermore,
\begin{align*}
\dfrac{
\left<
\displaystyle\int_{[0,1]^2}
({\bf v}\cdot e)|\boldsymbol{\phi}|^2dxdy
\right>
}{
\displaystyle\int_{[0,1]^2}
|\boldsymbol{\phi}|^2dxdy
}
=
\dfrac{
\displaystyle\int_0^1
\sin(2\pi y)\psi^2dy
}{
\displaystyle\int_0^1\psi^2dy
}.
\end{align*}
Taking the maximum with respect to $\psi$ and using
\eqref{Srvariational}, we obtain
\begin{align}\label{counterI}
\max_{\boldsymbol{\phi}\in
\mathfrak{I}_{{\bf A}}^\mu\setminus\{{\bf 0}\}}
\dfrac{
\left<
\displaystyle\int_{[0,1]^2}
({\bf v}\cdot e)|\boldsymbol{\phi}|^2dxdy
\right>
}{
\displaystyle\int_{[0,1]^2}
|\boldsymbol{\phi}|^2dxdy
}
\geq
\mathcal{S}(3), \text{ which yields }\max_{\boldsymbol{\phi}\in
\mathfrak{I}_{{\bf A}}^\mu\setminus\{{\bf 0}\}}
\dfrac{
\left<
\displaystyle\int_{[0,1]^2}
({\bf v}\cdot e)|\boldsymbol{\phi}|^2dxdy
\right>
}{
\displaystyle\int_{[0,1]^2}
|\boldsymbol{\phi}|^2dxdy
}
=
\mathcal{S}(3).
\end{align}

We now calculate the part associated with
$\mathfrak{I}_{{\bf A},time}^\mu$. Set $R_\delta
:=
\max\limits_{t\in\R}\lambda_{\max}(M_\delta(t))
=
2+\sqrt{1+\delta^2}.$ For
$\boldsymbol{\phi}\in\mathfrak{I}_{{\bf A},time}^\mu$,
one has $\displaystyle\sum_{i,j=1,2}\mu_{ij}\phi_i\phi_j
=
\boldsymbol{\phi}^TM_\delta(t)\boldsymbol{\phi}
\leq
R_\delta|\boldsymbol{\phi}|^2.$ Consequently,
\begin{align}\label{counterJenergy}
\sum_{i=1,2}
\left<
\int_{[0,1]^2}|D\phi_i|^2dxdy
\right>
\leq
R_\delta
\left<
\int_{[0,1]^2}
|\boldsymbol{\phi}|^2dxdy
\right>.
\end{align}
On the other hand, the Rayleigh characterization of
$\nu(\alpha)$ implies, for every $\alpha>0$,
\begin{align*}
\alpha
\left<
\int_{[0,1]^2}
\sin(2\pi y)|\boldsymbol{\phi}|^2dxdy
\right>
-
\sum_{i=1,2}
\left<
\int_{[0,1]^2}
|D\phi_i|^2dxdy
\right>
\leq
\nu(\alpha)
\left<
\int_{[0,1]^2}
|\boldsymbol{\phi}|^2dxdy
\right>.
\end{align*}
Together with \eqref{counterJenergy}, it follows that $\dfrac{
\left<
\displaystyle\int_{[0,1]^2}
({\bf v}\cdot e)|\boldsymbol{\phi}|^2dxdy
\right>
}{
\left<
\displaystyle\int_{[0,1]^2}
|\boldsymbol{\phi}|^2dxdy
\right>
}
\leq
\dfrac{\nu(\alpha)+R_\delta}{\alpha}.$ Taking the infimum with respect to $\alpha>0$, we derive
\begin{align}\label{counterJupper}
\sup_{\boldsymbol{\phi}\in
\mathfrak{I}_{{\bf A},time}^\mu\setminus\{{\bf 0}\}}
\dfrac{
\left<
\displaystyle\int_{[0,1]^2}
({\bf v}\cdot e)|\boldsymbol{\phi}|^2dxdy
\right>
}{
\left<
\displaystyle\int_{[0,1]^2}
|\boldsymbol{\phi}|^2dxdy
\right>
}
\leq
\mathcal{S}(R_\delta).
\end{align}

To prove the opposite inequality, fix $r\in(0,R_\delta)$.
Choose $t_0\in[0,1]$ such that
$\lambda_{\max}(M_\delta(t_0))=R_\delta$ and let
$\xi_0\in\R^2$, $|\xi_0|=1$, be an eigenvector satisfying $M_\delta(t_0)\xi_0=R_\delta\xi_0.$ Since
$t\mapsto\xi_0^TM_\delta(t)\xi_0$ is continuous and equals
$R_\delta$ at $t=t_0$, there exists a nonzero nonnegative
$1$-periodic function $\chi=\chi(t)$ such that
\begin{align}\label{counterchi}
\int_0^1
\xi_0^TM_\delta(t)\xi_0\chi^2(t)dt
\geq
r\int_0^1\chi^2(t)dt.
\end{align}
Choose
$\psi\in H^1_{per}(\R)\setminus\{0\}$ such that $\displaystyle\int_0^1|\psi'|^2dy
\leq
r\displaystyle\int_0^1\psi^2dy$ and set $\boldsymbol{\phi}(t,x,y)
:=
\chi(t)\psi(y)\xi_0.$ Since
${\bf v}\cdot D\phi_i
=
\sin(2\pi y)\partial_x\phi_i=0$,
one has
$\boldsymbol{\phi}\in\mathfrak{I}_{\text{time},p}$.
Furthermore, by \eqref{counterchi},
\begin{align*}
\sum_{i=1,2}
\left<
\int_{[0,1]^2}|D\phi_i|^2dxdy
\right>
\leq
\int_0^1
\xi_0^TM_\delta(t)\xi_0\chi^2(t)dt
\int_0^1\psi^2dy=
\sum_{i,j=1,2}
\left<
\int_{[0,1]^2}
\mu_{ij}\phi_i\phi_jdxdy
\right>.
\end{align*}
Thus,
$\boldsymbol{\phi}\in\mathfrak{I}_{{\bf A},time}^\mu$.
Moreover,
\begin{align*}
\dfrac{
\left<
\displaystyle\int_{[0,1]^2}
({\bf v}\cdot e)|\boldsymbol{\phi}|^2dxdy
\right>
}{
\left<
\displaystyle\int_{[0,1]^2}
|\boldsymbol{\phi}|^2dxdy
\right>
}
=
\dfrac{
\displaystyle\int_0^1
\sin(2\pi y)\psi^2dy
}{
\displaystyle\int_0^1\psi^2dy
}.
\end{align*}
Taking the supremum with respect to $\psi$ and then letting
$r\rightarrow R_\delta^-$ gives the inequality below. The reverse inequality established above yields equality. By
\eqref{Srvariational} that
\begin{align}\label{counterJ}
\sup_{\boldsymbol{\phi}\in
\mathfrak{I}_{{\bf A},time}^\mu\setminus\{{\bf 0}\}}
\dfrac{
\left<
\displaystyle\int_{[0,1]^2}
({\bf v}\cdot e)|\boldsymbol{\phi}|^2dxdy
\right>
}{
\left<
\displaystyle\int_{[0,1]^2}
|\boldsymbol{\phi}|^2dxdy
\right>
}
\geq
\mathcal{S}(R_\delta),~
\sup_{\boldsymbol{\phi}\in
\mathfrak{I}_{{\bf A},time}^\mu\setminus\{{\bf 0}\}}
\dfrac{
\left<
\displaystyle\int_{[0,1]^2}
({\bf v}\cdot e)|\boldsymbol{\phi}|^2dxdy
\right>
}{
\left<
\displaystyle\int_{[0,1]^2}
|\boldsymbol{\phi}|^2dxdy
\right>
}
=
\mathcal{S}(R_\delta).
\end{align}

Finally, $\mathcal{S}(r)$ is strictly increasing with respect
to $r>0$. Indeed, for $r_2>r_1>0$, let $\alpha_{r_2}>0$
be a minimizer associated with $\mathcal{S}(r_2)$. Then
\begin{align*}
\mathcal{S}(r_2)
=
\dfrac{\nu(\alpha_{r_2})+r_2}{\alpha_{r_2}}
\geq
\mathcal{S}(r_1)
+
\dfrac{r_2-r_1}{\alpha_{r_2}}
>
\mathcal{S}(r_1).
\end{align*}
Therefore, by $3<\rho_\delta<
2+\sqrt{1+\delta^2}$, $\mathcal{S}(3)
<
\mathcal{S}(\rho_\delta)
<
\mathcal{S}(R_\delta).$ Together with \eqref{counterlimit}, \eqref{counterI} and
\eqref{counterJ}, we finally obtain
\begin{align*}
&\max_{\boldsymbol{\phi}\in
\mathfrak{I}_{{\bf A}}^\mu\setminus\{{\bf 0}\}}
\dfrac{
\left<
\displaystyle\int_{[0,1]^2}
({\bf v}\cdot e)|\boldsymbol{\phi}|^2dxdy
\right>
}{
\displaystyle\int_{[0,1]^2}
|\boldsymbol{\phi}|^2dxdy
}<
\lim_{A\rightarrow\infty}
\dfrac{c^*(e,A)}{A}
<
\sup_{\boldsymbol{\phi}\in
\mathfrak{I}_{{\bf A},time}^\mu\setminus\{{\bf 0}\}}
\dfrac{
\left<
\displaystyle\int_{[0,1]^2}
({\bf v}\cdot e)|\boldsymbol{\phi}|^2dxdy
\right>
}{
\left<
\displaystyle\int_{[0,1]^2}
|\boldsymbol{\phi}|^2dxdy
\right>
}.
\end{align*}
Thus, neither of the two expected equalities holds in general,
even under the assumption $\mu_{12}=\mu_{21}$.

\section{Explicit Models and Illustrations}\label{section 6}
\begingroup
\setlength{\textfloatsep}{8pt plus 2pt minus 2pt}
\setlength{\floatsep}{7pt plus 2pt minus 2pt}
\setlength{\intextsep}{8pt plus 2pt minus 2pt}

We consider two space-periodic cooperative systems for which the spectral characterization of Theorem~\ref{minimalspeed} reduces to a one-dimensional periodic eigenvalue problem. We first examine the minimal speed and its large-advection limit. We then solve the full nonlinear system, without imposing equality of the two components, using a mass-lumped bilinear finite-element discretization and a second-order time integrator. All the numerical quantities and the surfaces below are computed from the equations stated in this section. The finite-time solutions are not identified with the entire critical fronts constructed in Section~\ref{sec:critical-proofs}.

\medskip
\noindent\textbf{Case 1. A horizontal periodic shear.}
On $[0,1]^2$, set
\begin{align*}
{\bf v}(x,y)=(\sin(2\pi y),0),\quad e=(1,0),\quad
{\bf A}^1={\bf A}^2=I_2,\quad {\bf q}_1={\bf q}_2=A{\bf v},
\end{align*}
with
\begin{align*}
\begin{cases}
F_1(u_1,u_2)=2u_1+u_2-3u_1^2,\\
F_2(u_1,u_2)=u_1+2u_2-3u_2^2.
\end{cases}
\end{align*}
The linearization at ${\bf0}$ has principal eigenvalue $3$, associated with the eigenvector $(1,1)$. The field ${\bf v}$ is periodic, divergence free and has zero mean. Writing $z=y$, the tilted principal eigenvalue satisfies
\begin{align*}
k(\lambda,A)=\lambda^2+\lambda_{1,1}(\lambda A),
\end{align*}
where
\begin{equation}\label{eq:app-case1-1d}
\phi''+(3+\alpha\sin(2\pi z))\phi=\lambda_{1,1}(\alpha)\phi,
\quad\phi(z+1)=\phi(z)>0.
\end{equation}
Consequently,
\begin{equation}\label{eq:app-case1-speed}
\frac{c^*(e,A)}{A}=\inf_{\alpha>0}
\left\{\frac{\alpha}{A^2}+\frac{\lambda_{1,1}(\alpha)}{\alpha}\right\}.
\end{equation}

\medskip
\noindent\textbf{Case 2. A diagonal periodic shear.}
We take
\begin{align*}
{\bf v}(x,y)=\bigl(\sin(2\pi(x-y)),\sin(2\pi(x-y))\bigr),\quad
 e=\frac1{\sqrt2}(1,1),\quad
{\bf A}^1={\bf A}^2=I_2,\quad {\bf q}_1={\bf q}_2=A{\bf v},
\end{align*}
and
\begin{align*}
\begin{cases}
F_1(u_1,u_2)=3u_1+2u_2-5u_1^2,\\
F_2(u_1,u_2)=2u_1+3u_2-5u_2^2.
\end{cases}
\end{align*}
The principal eigenvalue of the linearized reaction matrix is $5$. With $s=(x+y)/\sqrt2$ and $z=x-y$, there holds
$\Delta=\partial_{ss}+2\partial_{zz}$ and
${\bf v}\cdot e=\sqrt2\sin(2\pi z)$. Thus,
\begin{align*}
k(\lambda,A)=\lambda^2+\lambda_{1,2}(\lambda A),
\end{align*}
where
\begin{equation}\label{eq:app-case2-1d}
2\phi''+(5+\sqrt2\alpha\sin(2\pi z))\phi
=\lambda_{1,2}(\alpha)\phi,\quad \phi(z+1)=\phi(z)>0,
\end{equation}
and
\begin{equation}\label{eq:app-case2-speed}
\frac{c^*(e,A)}{A}=\inf_{\alpha>0}
\left\{\frac{\alpha}{A^2}+\frac{\lambda_{1,2}(\alpha)}{\alpha}\right\}.
\end{equation}
Both shear fields are shown in Figure~\ref{fig:app-flows}.

\begin{figure}[!htbp]
\centering
\includegraphics[width=.92\textwidth]{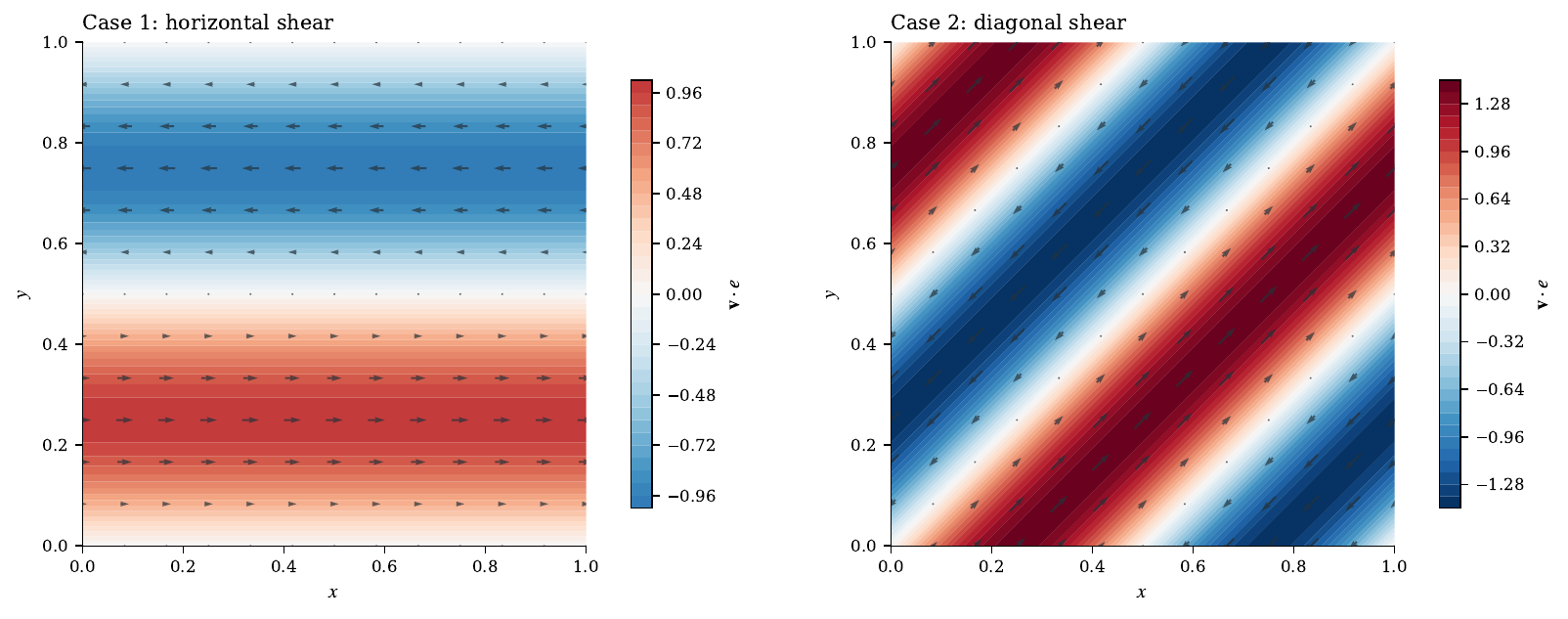}
\caption{Periodic velocity fields in the two examples. The arrows represent ${\bf v}$ and the background represents ${\bf v}\cdot e$. The direction of propagation is horizontal in Case~1 and diagonal in Case~2. Both fields are divergence free.}
\label{fig:app-flows}
\end{figure}

\medskip
\noindent\textbf{The principal eigenvalue and the minimal speed.}
For $j=1,2$, write
\begin{equation}\label{eq:app-unified}
(D_j,r_j,b_j)=\begin{cases}
(1,3,1),&j=1,\\
(2,5,\sqrt2),&j=2.
\end{cases}
\end{equation}

\subsection{Spectral reduction}
The two reaction matrices are symmetric and have the same Perron
eigendirection. This is the structural reason that the two-component
eigenvalue computation reduces to a scalar periodic problem.

\begin{lemma}\label{lem:memoirs-spectral-check}
Set $\gamma_1=1$ and $\gamma_2=2$. For $j=1,2$ and
$\alpha\geq0$, consider the periodic operator
\begin{equation}\label{eq:memoirs-tilted-system}
\mathcal B_{j,\alpha}=D_j\partial_{zz}I_2+
 b_j\alpha\sin(2\pi z)I_2+
\begin{pmatrix}r_j-\gamma_j&\gamma_j\\
\gamma_j&r_j-\gamma_j\end{pmatrix}.
\end{equation}
Its principal periodic eigenvalue equals the largest eigenvalue of
$D_j\partial_{zz}+r_j+b_j\alpha\sin(2\pi z)$.
The positive principal eigenfunction of $\mathcal B_{j,\alpha}$
may be chosen in the form $(\phi,\phi)$, where $\phi$ is the
positive periodic principal eigenfunction of the scalar operator.
The complementary antisymmetric mode has its largest eigenvalue
smaller by $2\gamma_j$.
\end{lemma}
\begin{proof}
The constant vectors $(1,1)$ and $(1,-1)$ diagonalize the
reaction matrix, with eigenvalues $r_j$ and $r_j-2\gamma_j$.
The diffusion and advection operators act identically on both
components and therefore preserve these two invariant subspaces.
On the symmetric subspace the eigenvalue equation is precisely
\begin{equation*}
D_j\phi''+(r_j+b_j\alpha\sin(2\pi z))\phi=\Lambda\phi,
\end{equation*}
whereas on the antisymmetric subspace the potential $r_j$ is
replaced by $r_j-2\gamma_j$. The two scalar operators differ by
the strictly negative constant $-2\gamma_j$. Their largest
periodic eigenvalues consequently differ by the same constant.
The symmetric eigenfunction is strictly positive in both
components and hence represents the principal eigenpair of the
cooperative system.
\end{proof}

The longitudinal exponential weight contributes $\lambda^2$
from the unit longitudinal diffusion and changes the shear
transport into the multiplication potential
$\lambda A b_j\sin(2\pi z)$. Lemma~\ref{lem:memoirs-spectral-check}
therefore verifies the identities
\begin{equation}\label{eq:memoirs-spectral-identity}
 k_j(\lambda,A)=\lambda^2+\lambda_{1,j}(A\lambda),
 \qquad j=1,2,
\end{equation}
used in the speed calculations below. Here the eigenvalues refer
to the equations for the two distinct coefficient families;
there is no assumption that their speeds coincide.

The periodic principal eigenvalue of
$D_j\partial_{zz}+r_j+b_j\alpha\sin(2\pi z)$
is approximated on $[0,1]$ by continuous piecewise-linear functions with periodic endpoint identification. Let $M_h$ and $K_h$ be the mass and stiffness matrices, and $S_h$ the weighted mass matrix associated with $\sin(2\pi z)$. The generalized symmetric eigenvalue problem reads
\begin{equation}\label{eq:app-Fourier}
(-D_j K_h+r_j M_h+b_j\alpha S_h)\Phi_h
=\lambda_{1,j,h}(\alpha)M_h\Phi_h.
\end{equation}
We use the largest generalized eigenvalue. Each potential integral is evaluated by three-point Gaussian quadrature. For each case, the eigenvalue is computed at $160$ geometrically spaced parameters in $[0.16,75]$. A cubic interpolation of these values is then minimized over the same interval. The computed minima lie in its interior. The figures use $192$ periodic elements. Doubling this number to $384$ changes the large-advection normalized-speed limits by at most $1.6\times10^{-5}$ and the normalized speeds at $A=5$ by at most $8.1\times10^{-6}$. These are spatial sensitivity comparisons, not certified finite-element error bounds. The cubic interpolation is used only for scalar minimization after the generalized eigenvalues have been computed.

\begin{table}[!htbp]
\centering
\caption{Normalized minimal speed $c^*(e,A)/A$ for the horizontal periodic shear.}
\label{table1}
\renewcommand{\arraystretch}{1.09}
\begin{tabular}{@{}c|rrrrrrr@{}}\hline
$A$ & $0.5$ & $1$ & $5$ & $10$ & $15$ & $20$ & $25$\\\hline
$c^*(e,A)/A$ & $6.939160$ & $3.485949$ & $0.793473$ & $0.516263$ & $0.445038$ & $0.416828$ & $0.402970$\\\hline
\end{tabular}
\end{table}
\begin{table}[!htbp]
\centering
\caption{Normalized minimal speed $c^*(e,A)/A$ for the diagonal periodic shear.}
\label{table2}
\renewcommand{\arraystretch}{1.09}
\begin{tabular}{@{}c|rrrrrrr@{}}\hline
$A$ & $0.5$ & $1$ & $5$ & $10$ & $15$ & $20$ & $25$\\\hline
$c^*(e,A)/A$ & $8.958417$ & $4.500345$ & $1.024684$ & $0.667590$ & $0.576241$ & $0.540192$ & $0.522528$\\\hline
\end{tabular}
\end{table}

Let
\begin{equation}\label{eq:app-limit}
\ell_j=\inf_{\alpha>0}\frac{\lambda_{1,j}(\alpha)}{\alpha},\quad
\alpha_j^*\in\operatorname*{argmin}_{\alpha>0}
\frac{\lambda_{1,j}(\alpha)}{\alpha}.
\end{equation}
By the variational characterization,
\begin{equation}\label{eq:app-bound}
0\leq\frac{c^*(e,A)}{A}-\ell_j\leq\frac{\alpha_j^*}{A^2},
\quad j=1,2.
\end{equation}
\subsection{Higher-order large-advection asymptotics}
The scalar reductions admit a sharper asymptotic expansion than the
first-order bound in \eqref{eq:app-bound}.  The result concerns the
explicit stationary shear models of this section; it does not assert
such an expansion for the general nonautonomous coefficients studied above.

\begin{proposition}\label{prop:memoirs-second-asympt}
Let $D>0$, $r>0$ and $b>0$.  Denote by $\nu(\alpha)$ the
largest periodic eigenvalue of
\begin{equation}\label{eq:memoirs-speed-schrodinger}
D\phi''+(r+b\alpha\sin(2\pi z))\phi=\nu(\alpha)\phi,
\qquad \phi(z+1)=\phi(z)>0.
\end{equation}
Then $g(\alpha):=\nu(\alpha)/\alpha$ has a unique minimizer
$\alpha_*>0$, and $g''(\alpha_*)>0$.  If
\begin{equation*}
\mathcal C(A):=\inf_{\alpha>0}\bigg\{g(\alpha)+\frac{\alpha}{A^2}\bigg\},
\qquad \ell=g(\alpha_*),
\end{equation*}
then, as $A\to\infty$,
\begin{equation}\label{eq:memoirs-higher-speed}
\mathcal C(A)=\ell+\frac{\alpha_*}{A^2}
-\frac{1}{2g''(\alpha_*)A^4}+O(A^{-6}).
\end{equation}
More precisely, the minimizing parameter satisfies
\begin{equation}\label{eq:memoirs-minimizer-expansion}
\alpha(A)=\alpha_*-\frac{1}{g''(\alpha_*)A^2}+O(A^{-4}).
\end{equation}
These statements apply to both Cases~1 and~2, with
$(D,r,b)=(D_j,r_j,b_j)$ and $\mathcal C(A)=c^*(e,A)/A$.
\end{proposition}
\begin{proof}
The operator in \eqref{eq:memoirs-speed-schrodinger} is
self-adjoint on periodic $H^2(0,1)$.  Its largest eigenvalue is
simple, has a strictly positive eigenfunction, and is analytic
with respect to $\alpha$.  Normalize this eigenfunction by
$\|\phi_\alpha\|_{L^2(0,1)}=1$.  For an orthonormal basis of
real eigenfunctions $\{\phi_{\alpha,n}\}_{n\geq0}$, with
$\phi_{\alpha,0}=\phi_\alpha$ and ordered eigenvalues
$\nu(\alpha)>\nu_1(\alpha)\geq\nu_2(\alpha)\geq\cdots$,
second-order self-adjoint perturbation theory yields
\begin{equation}\label{eq:memoirs-strict-convex}
\nu''(\alpha)=2b^2\sum_{n\geq1}
\frac{\big|\int_0^1\sin(2\pi z)
\phi_\alpha(z)\phi_{\alpha,n}(z)\,dz\big|^2}
{\nu(\alpha)-\nu_n(\alpha)}>0.
\end{equation}
The inequality is strict: otherwise
$\sin(2\pi z)\phi_\alpha$ would be proportional to
$\phi_\alpha$, which is impossible because
$\phi_\alpha>0$ and $\sin(2\pi z)$ is nonconstant.

Since $\nu(0)=r$, one has $g(\alpha)\to+\infty$ as
$\alpha\downarrow0$.  The Rayleigh characterization implies
$\nu(\alpha)\leq r+b\alpha$.  Conversely, for every
$\eta>0$ choose a nonzero smooth periodic test function
supported where $\sin(2\pi z)>1-\eta$.  Its Rayleigh quotient
shows that
\begin{equation*}
\lim_{\alpha\to\infty}\frac{\nu(\alpha)}\alpha=b.
\end{equation*}
In fact, $g(\alpha)<b$ for all sufficiently large $\alpha$.
To verify this strict inequality, observe that
\begin{equation*}
\nu(\alpha)=r+b\alpha-E(\alpha),\qquad
E(\alpha)=\inf_{\|f\|_2=1}
\bigg\{D\int_0^1|f'|^2\,dz+
 b\alpha\int_0^1(1-\sin(2\pi z))|f|^2\,dz\bigg\}.
\end{equation*}
If $E(\alpha)$ remained bounded along some
$\alpha_n\to\infty$, normalized approximate minimizers
would be bounded in $H^1_{\rm per}(0,1)$.
Compactness in $L^2(0,1)$ would produce a limit of norm one
supported on the zero set of $1-\sin(2\pi z)$, a set of
measure zero, which is impossible.  Thus
$E(\alpha)\to\infty$, and $\nu(\alpha)<b\alpha$
for all sufficiently large $\alpha$.
It follows that $g$ attains its global minimum in $(0,\infty)$.

Set $H(\alpha)=\alpha\nu'(\alpha)-\nu(\alpha)$.  By
\eqref{eq:memoirs-strict-convex},
$H'(\alpha)=\alpha\nu''(\alpha)>0$ for $\alpha>0$.
Moreover, $g'(\alpha)=H(\alpha)/\alpha^2$.
Consequently the minimizer is unique, and
\begin{equation}\label{eq:memoirs-nondegen}
 g''(\alpha_*)=\frac{\nu''(\alpha_*)}{\alpha_*}>0.
\end{equation}

Write $\varepsilon=A^{-2}$.  Since
$g(\alpha)\to+\infty$ at zero,
$g(\alpha)\to b>\ell$ at infinity, and $\alpha_*$ is
its unique minimizer, every minimizer of
$g(\alpha)+\varepsilon\alpha$ converges to $\alpha_*$ as
$\varepsilon\downarrow0$.  The analytic implicit-function theorem
applied to $g'(\alpha)+\varepsilon=0$ therefore gives
\begin{equation*}
\alpha(\varepsilon)=\alpha_*
 -\frac{\varepsilon}{g''(\alpha_*)}+O(\varepsilon^2).
\end{equation*}
Taylor expansion of
$g(\alpha(\varepsilon))+\varepsilon\alpha(\varepsilon)$
then yields
\begin{equation*}
\ell+\alpha_*\varepsilon
 -\frac{\varepsilon^2}{2g''(\alpha_*)}+O(\varepsilon^3).
\end{equation*}
This proves \eqref{eq:memoirs-higher-speed} and
\eqref{eq:memoirs-minimizer-expansion}.
\end{proof}

The expansion refines the bounds in \eqref{eq:app-bound}
for the stationary shears.  The coefficients of the correction
are characterized by the continuous periodic eigenvalue
problem; the finite-element data below are numerical
approximations and are not used in the proof.

For the mesh used here,
\begin{equation}\label{eq:app-values}
\begin{aligned}
\ell_1&\approx0.376752389,&\alpha_1^*&\approx17.124602,\\
\ell_2&\approx0.489219494,&\alpha_2^*&\approx21.695988.
\end{aligned}
\end{equation}
Figures~\ref{fig:intro-speed-selection}--\ref{fig:app-spectral} display the spectral minimization, the principal eigenfunctions and the large-advection limit. These plots concern the eigenvalue problems \eqref{eq:app-case1-1d} and \eqref{eq:app-case2-1d}, not fitted propagation curves.

\begin{figure}[!htbp]
\centering
\includegraphics[width=.93\textwidth]{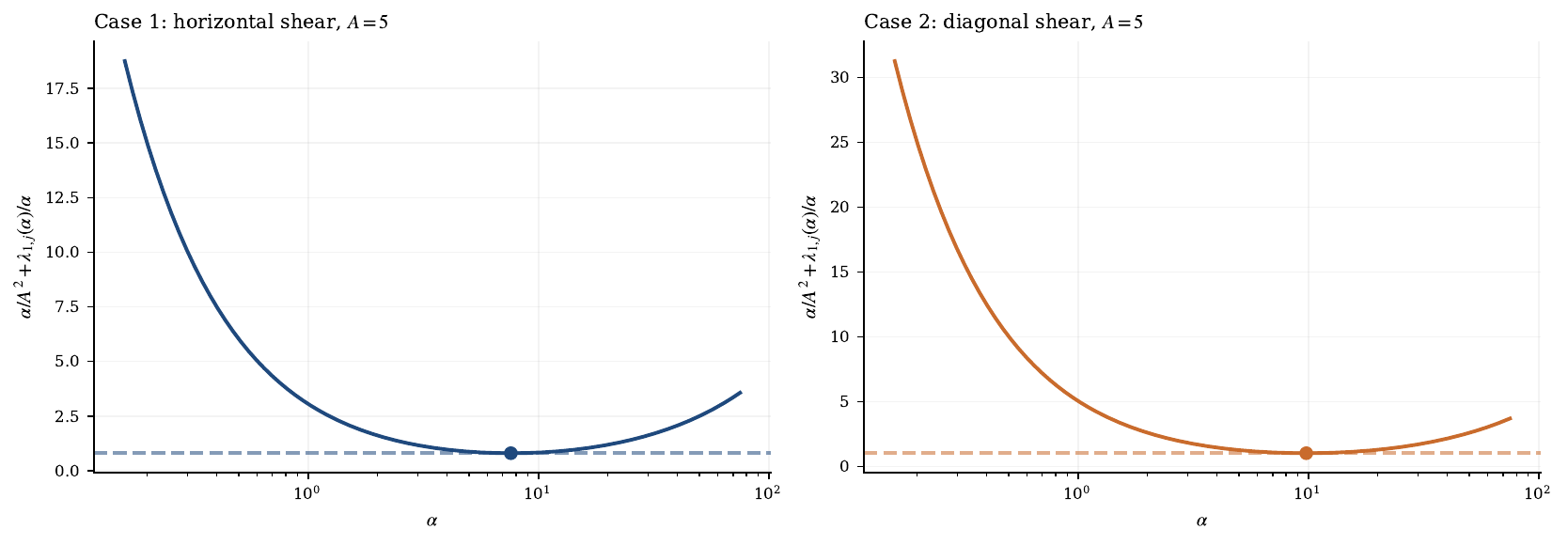}
\caption{Variational characterization of the minimal speed in Theorem~\ref{minimalspeed}. For $A=5$, the curves represent $\alpha/A^2+\lambda_{1,j}(\alpha)/\alpha$ in the two periodic media. Their minima equal $c^*(e,A)/A$.}
\label{fig:intro-speed-selection}
\end{figure}
\begin{figure}[!htbp]
\centering
\includegraphics[width=.93\textwidth]{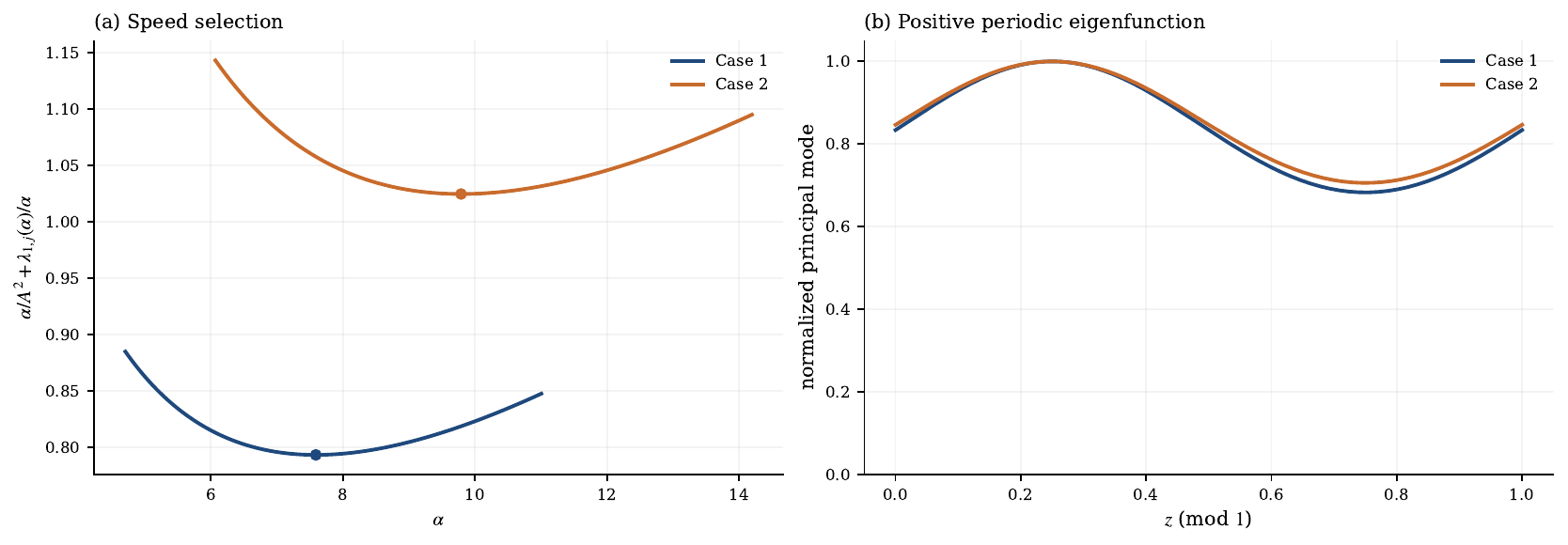}
\caption{The minimizing spectral modes for $A=5$. Panel (a) shows the variational quotients near their minima; panel (b) shows the corresponding positive periodic eigenfunctions, normalized to have maximum one.}
\label{fig:intro-critical-modes}
\end{figure}
\begin{figure}[!htbp]
\centering
\includegraphics[width=.94\textwidth]{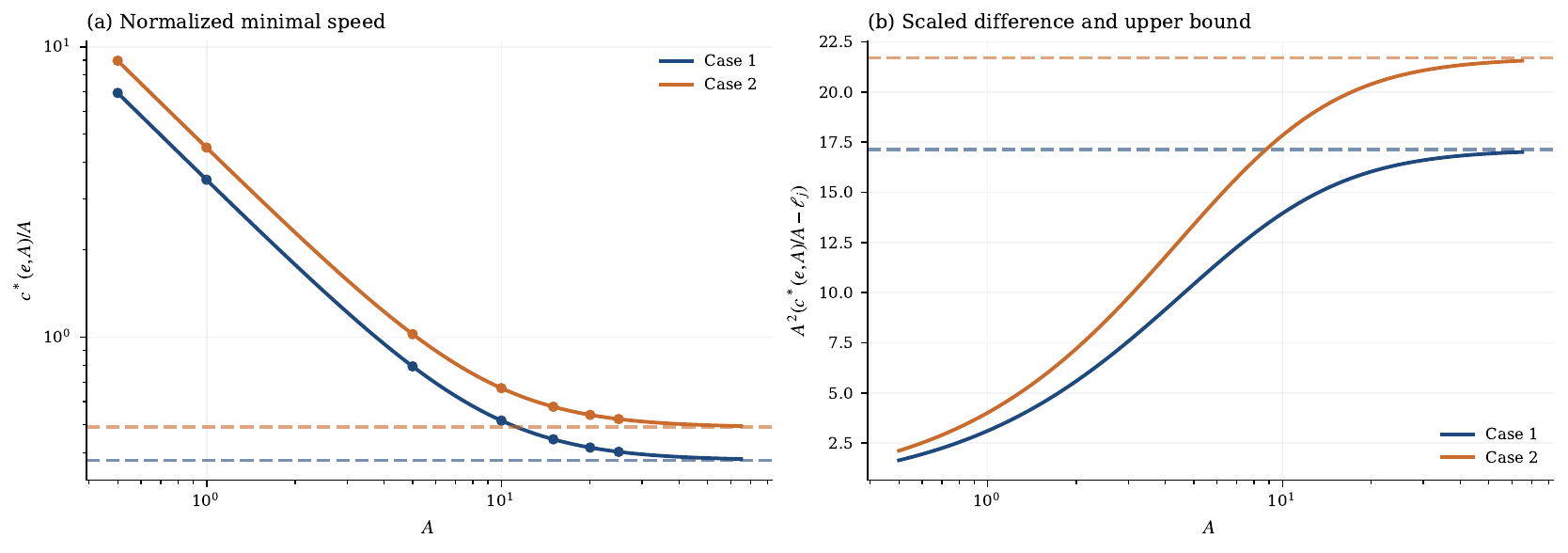}
\caption{Large-advection asymptotics in Theorem~\ref{limitspeed1}\textup{(b)}. Panel (a) shows $c^*(e,A)/A$ and its limit $\ell_j$; panel (b) shows $A^2(c^*(e,A)/A-\ell_j)$ and the upper bound $\alpha_j^*$ in \eqref{eq:app-bound}. The symbols in panel (a) correspond to Tables~\ref{table1} and~\ref{table2}.}
\label{fig:app-speed}
\end{figure}
\begin{figure}[!htbp]
\centering
\includegraphics[width=.93\textwidth]{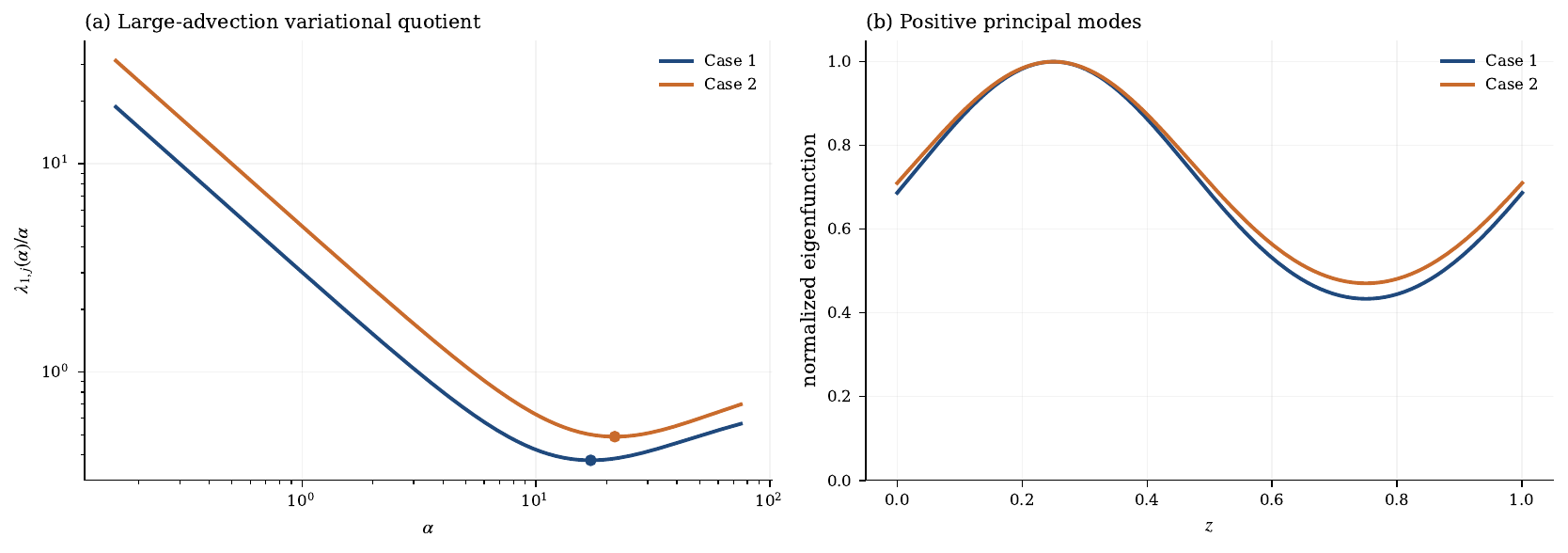}
\caption{One-dimensional spectral reduction of Theorem~\ref{limitspeed1}\textup{(b)}. Panel (a) contains the quotients $\lambda_{1,j}(\alpha)/\alpha$ and their minimizers; panel (b) contains the associated normalized positive principal eigenfunctions over one period.}
\label{fig:app-spectral}
\end{figure}

\medskip
\noindent\textbf{Exact supercritical fronts in the absence of advection.}
When $A=0$ and $u_1=u_2=u$ is independent of $z$, both systems reduce to
$u_t=u_{ss}+r_j u(1-u)$. The function
\begin{equation}\label{eq:app-exact-wave}
U_j(s)=\bigl(1+\exp(\sqrt{r_j/6}\,s)\bigr)^{-2},\quad
c_j=5\sqrt{r_j/6}
\end{equation}
satisfies $U_j''+c_j U_j'+r_j U_j(1-U_j)=0$, with limits $U_j(-\infty)=1$ and $U_j(+\infty)=0$. Hence $(U_j,U_j)$ is an exact front of the cooperative system. As $c_j>2\sqrt{r_j}$, these fronts are supercritical. They do not represent the critical fronts considered in Section~\ref{sec:critical-proofs}.
\begin{figure}[!htbp]
\centering
\includegraphics[width=.90\textwidth]{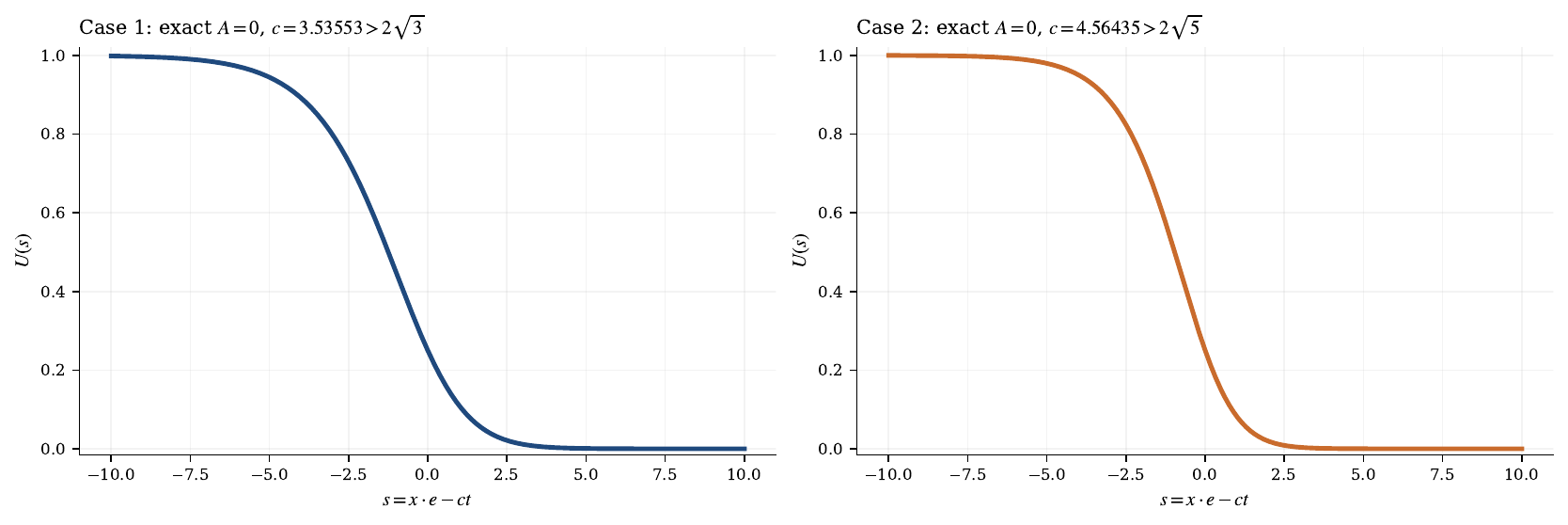}
\caption{Exact nonlinear traveling fronts \eqref{eq:app-exact-wave} when $A=0$. Panels (a) and (b) correspond to $r_1=3$ and $r_2=5$, respectively. Both components coincide, and the propagation speed is strictly above the minimal speed.}
\label{fig:app-exact}
\end{figure}

\medskip
\noindent\textbf{Nonlinear cooperative propagation.}
We now solve the full two-component system. Set $s=x\cdot e$ and take $z=y$ in Case~1, while $s=(x+y)/\sqrt2$ and $z=x-y$ in Case~2. With $\gamma_1=1$ and $\gamma_2=2$, the exact reduction for solutions depending on $(s,z)$ is
\begin{equation}\label{eq:app-full-pde}
\begin{cases}
(u_1)_t=(u_1)_{ss}+D_j(u_1)_{zz}-Ab_j\sin(2\pi z)(u_1)_s
+(r_j-\gamma_j)u_1+\gamma_j u_2-r_j u_1^2,\\
(u_2)_t=(u_2)_{ss}+D_j(u_2)_{zz}-Ab_j\sin(2\pi z)(u_2)_s
+\gamma_j u_1+(r_j-\gamma_j)u_2-r_j u_2^2.
\end{cases}
\end{equation}
For completeness, one may separate the symmetric and antisymmetric
parts of the nonlinear system without assuming that its components
coincide.
\begin{lemma}\label{lem:memoirs-mode-splitting}
Let $(u_1,u_2)$ be a classical solution of
\eqref{eq:app-full-pde}, and put $v=(u_1+u_2)/2$ and
$d=(u_1-u_2)/2$. If
\begin{equation*}
\mathscr D_j=\partial_{ss}+D_j\partial_{zz}
-Ab_j\sin(2\pi z)\partial_s,
\end{equation*}
then
\begin{equation}\label{eq:memoirs-split}
\begin{cases}
v_t=\mathscr D_j v+r_jv(1-v)-r_jd^2,\\
d_t=\mathscr D_j d+(r_j-2\gamma_j-2r_jv)d.
\end{cases}
\end{equation}
In particular, equality of the components is an invariant
constraint, but need not hold for nonsymmetric initial data.
\end{lemma}
\begin{proof}
Insert $u_1=v+d$ and $u_2=v-d$ into
\eqref{eq:app-full-pde}. On adding the two equations and dividing
by two, the reaction sum is
\begin{equation*}
 r_jv-\frac{r_j}{2}\big((v+d)^2+(v-d)^2\big)
 =r_jv(1-v)-r_jd^2.
\end{equation*}
Their difference divided by two has reaction term
\begin{equation*}
(r_j-2\gamma_j)d-
\frac{r_j}{2}\big((v+d)^2-(v-d)^2\big)
 =(r_j-2\gamma_j-2r_jv)d.
\end{equation*}
The spatial operators commute with these linear combinations.
This proves \eqref{eq:memoirs-split}. If $d$ is initially zero
and both components have identical boundary values, uniqueness
for the resulting linear equation gives $d\equiv0$.
\end{proof}
The mode decomposition by itself does not imply decay of $d$,
because its zero-order coefficient need not be negative everywhere.
The weighted estimate established below provides an additional
conclusion under the precise spectral normalization and initial
conditions of Theorem~\ref{thm:memoirs-mode-decay}.

\subsection{Consequences for the two stationary shears}\label{sec:diffusion-consequences}
The preceding diffusion-rate theorem applies to both stationary shears. We record the corresponding parameter values and the constant-coefficient limiting case.
\begin{corollary}[The two periodic shears]\label{cor:memoirs-diffusion-cases}
For the horizontal and diagonal shear examples, replace the common physical diffusion matrices $I_2$ by $dI_2$ in both components, without changing the amplitude $B>0$ of the shear. Theorem~\ref{thm:memoirs-diffusion-two-regimes} applies with
\begin{equation}\label{eq:memoirs-diffusion-case-parameters}
(a,D,r,\gamma,h)=
\begin{cases}
(1,1,3,1,B),&\text{Case 1},\\
(1,2,5,2,\sqrt2 B),&\text{Case 2}.
\end{cases}
\end{equation}
In particular the small-diffusion limiting speeds are $B$ and $\sqrt2 B$, respectively, and the large-diffusion leading terms are $2\sqrt{3d}$ and $2\sqrt{5d}$. All higher-order coefficients are obtained by substituting \eqref{eq:memoirs-diffusion-case-parameters} in \eqref{eq:memoirs-diffusion-large}--\eqref{eq:memoirs-diffusion-small}.
\end{corollary}
\begin{proof}
For Case~1, $s=x$ and $z=y$, so $\Delta=\partial_{ss}+\partial_{zz}$ and the shear projection is $B\sin(kz)$. For Case~2, $s=(x+y)/\sqrt2$ and $z=x-y$, so $\Delta=\partial_{ss}+2\partial_{zz}$ and the projection is $\sqrt2 B\sin(kz)$. Both systems retain their symmetric coupling, and the positive tilted eigenmode therefore remains the scalar symmetric mode. All the asserted parameters follow directly.
\end{proof}

\begin{proposition}[Constant-coefficient consistency check]\label{prop:memoirs-diffusion-homogeneous}
Let ${\bf A}_0$ be any constant symmetric positive definite $N\times N$ matrix, let $q_0\in\mathbb R^N$ be a common constant transport vector, and let $M$ be a constant irreducible cooperative $2\times2$ matrix with spectral bound $s(M):=\max\{\operatorname{Re}z:z\in\sigma(M)\}>0$. For the constant-coefficient linearization with diffusion $d{\bf A}_0$ in both components and transport $q_0$, its spectral linear-determinacy speed in direction $e\in\mathbb S^{N-1}$ is exactly
\begin{equation}\label{eq:memoirs-diffusion-homogeneous}
c_{\rm hom}(d,e)=q_0\cdot e
 +2\sqrt{d\,s(M)\,e\cdot{\bf A}_0 e}.
\end{equation}
This identity applies to the nonlinear KPP system whenever its minimal speed is linearly determinate.
\end{proposition}
\begin{proof}
Insertion of the exponential ansatz $e^{-\lambda(x\cdot e-ct)}p$ with the positive eigenvector $p$ corresponding to $s(M)$ gives the tilted principal eigenvalue $d\lambda^2e\cdot{\bf A}_0e+\lambda q_0\cdot e+s(M)$. Dividing by $\lambda>0$ and minimizing yields \eqref{eq:memoirs-diffusion-homogeneous}. In particular, the flow-free version of \eqref{eq:memoirs-diffusion-model} has $q_0=0$ and $s(M)=r$, yielding the exact relation $c^*(d)=2\sqrt{ra\,d}$. This agrees with the vanishing-shear specialization of the preceding formulas, although the fixed-positive-shear small-diffusion expansion was not asserted uniformly in $h$.
\end{proof}

\begin{figure}[!htbp]
\centering
\includegraphics[width=.94\textwidth]{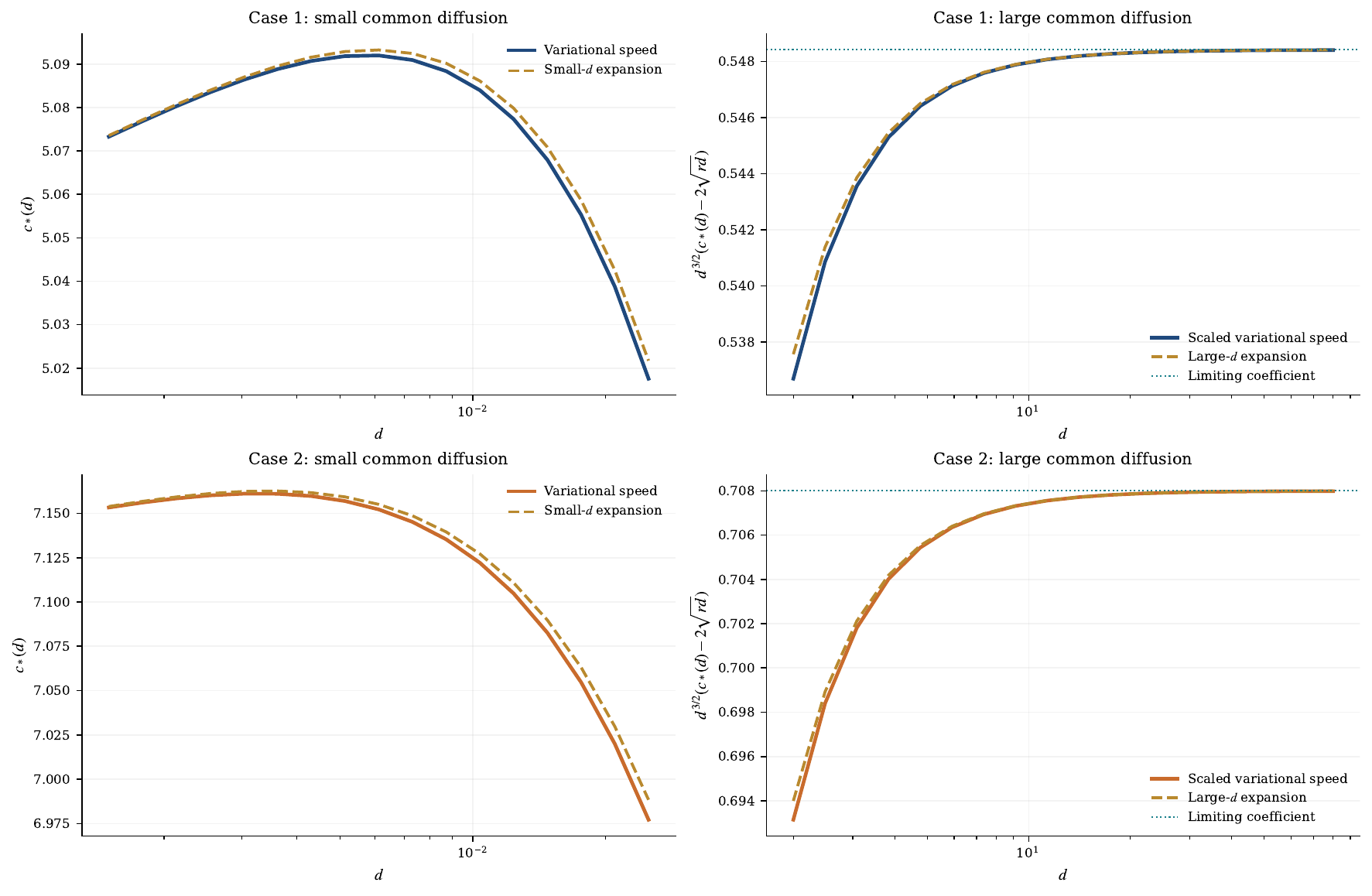}
\caption{Dependence of the spectral minimal speed on the common diffusion rate for the two stationary shears at fixed amplitude $B=5$. The left panels compare $c^*(d)$ with its small-diffusion expansion; the right panels display the scaled large-diffusion correction and its limiting coefficient. Both cases use the parameter choices of \eqref{eq:memoirs-diffusion-case-parameters}; the curves are computed from the periodic eigenvalue problem \eqref{eq:memoirs-diffusion-speed}.}
\label{fig:memoirs-diffusion-regimes}
\end{figure}

The two diffusion limits correspond to different spectral regimes. As $d\to\infty$, averaging over the transverse period leads to the effective longitudinal diffusion coefficient $ad+h^2/(2dDk^2)$, up to higher-order terms at the minimizing exponential weight. As $d\downarrow0$, the principal eigenfunction concentrates near the streamline $z=1/4$, where the projected velocity is maximal. The term of order $d^{3/4}$ arises from the ground-state energy of the limiting harmonic oscillator. The expansions are established for the stationary shear flows considered above, through their self-adjoint scalar spectral reduction.

The transverse variable $z$ has period one. Equality $u_1=u_2$ is preserved by the system, but it is not imposed in the computations. We fix $A=5$. The spectral computation on the $192$-element mesh gives
\begin{equation}\label{eq:app-cstar-A5}
c_1^*\approx3.96736688,\quad c_2^*\approx5.12341981.
\end{equation}
With $\xi=s-c_j^*t$ and $w_i(t,\xi,z)=u_i(t,s,z)$, equation \eqref{eq:app-full-pde} becomes
\begin{equation}\label{eq:app-moving-pde}
\begin{cases}
(w_1)_t=(w_1)_{\xi\xi}+D_j(w_1)_{zz}
+(c_j^*-Ab_j\sin(2\pi z))(w_1)_\xi
+(r_j-\gamma_j)w_1+\gamma_j w_2-r_jw_1^2,\\
(w_2)_t=(w_2)_{\xi\xi}+D_j(w_2)_{zz}
+(c_j^*-Ab_j\sin(2\pi z))(w_2)_\xi
+\gamma_j w_1+(r_j-\gamma_j)w_2-r_jw_2^2.
\end{cases}
\end{equation}
On $\xi\in[-24,24]$ we prescribe $w_i(t,-24,z)=1$ and $w_i(t,24,z)=0$, together with periodicity in $z$. The initial data are deliberately nonsymmetric:
\begin{equation}\label{eq:app-initial}
\begin{aligned}
w_1(0,\xi,z)&=\frac1{1+\exp(2.4[\xi+0.24-0.07\sin(2\pi z)])},\\
w_2(0,\xi,z)&=\frac1{1+\exp(2.4[\xi-0.24-0.07\sin(2\pi z)])}.
\end{aligned}
\end{equation}
The two endpoint rows are replaced by the prescribed boundary values. This initial separation of the components is part of the mathematical data and is not introduced when drawing the curves.

\subsection{Decay of the antisymmetric mode}
The exact mode splitting yields a quantitative estimate along the
spectrally selected critical frame.  The estimate does not require
the two initial components to coincide.

\begin{theorem}\label{thm:memoirs-mode-decay}
Fix $j\in\{1,2\}$ and $A>0$.  Let $\lambda_{A,j}>0$
minimize $\dfrac{\lambda^2+\lambda_{1,j}(A\lambda)}\lambda,$ so that the corresponding spectral speed is
\begin{equation}\label{eq:memoirs-critical-weight}
c_j^*(A)=\frac{\lambda_{A,j}^2+
\lambda_{1,j}(A\lambda_{A,j})}{\lambda_{A,j}}.
\end{equation}
Let $\phi_{A,j}>0$ be the periodic principal eigenfunction of
\begin{equation}\label{eq:memoirs-mode-weight-eig}
D_j\phi''+(r_j+Ab_j\lambda_{A,j}\sin(2\pi z))\phi
=\lambda_{1,j}(A\lambda_{A,j})\phi.
\end{equation}
Suppose $w_1,w_2\in[0,1]$ are classical solutions of
\eqref{eq:app-moving-pde}, with $c_j^*$ there interpreted as
$c_j^*(A)$, on either $\mathbb R\times\mathbb T$
or $[a,b]\times\mathbb T$, with periodicity in $z$.  In the
bounded-interval case, suppose their longitudinal boundary
values agree.  In the whole-line case, assume that the weighted difference
$e^{\lambda_{A,j}\xi}(w_1-w_2)/\phi_{A,j}(z)$ is bounded
on every finite time cylinder $[0,T]\times\mathbb R\times\mathbb T$.  Set
\begin{equation*}
\mathcal N_{A,j}(t)=\sup_{\xi,z}
\frac{e^{\lambda_{A,j}\xi}|w_1(t,\xi,z)-w_2(t,\xi,z)|}
{\phi_{A,j}(z)}.
\end{equation*}
Then, at every subsequent time for which the solution exists,
\begin{equation}\label{eq:memoirs-mode-decay}
\mathcal N_{A,j}(t)\leq e^{-2\gamma_j t}
\mathcal N_{A,j}(0),~
\gamma_1=1,~\gamma_2=2.
\end{equation}
In particular, the decay rates in the corresponding
exponentially weighted norms are at least $2$ and $4$.
\end{theorem}
\begin{proof}
Put $h=w_1-w_2$.  Subtracting the two equations in
\eqref{eq:app-moving-pde}, without linearizing the reaction,
we obtain
\begin{equation}\label{eq:memoirs-h-eq}
\begin{split}
h_t={}&h_{\xi\xi}+D_jh_{zz}
+(c_j^*(A)-Ab_j\sin(2\pi z))h_\xi+(r_j-2\gamma_j-r_j(w_1+w_2))h.
\end{split}
\end{equation}
Define
$\Psi(\xi,z)=e^{-\lambda_{A,j}\xi}\phi_{A,j}(z)$.
By \eqref{eq:memoirs-critical-weight} and
\eqref{eq:memoirs-mode-weight-eig}, direct differentiation gives
\begin{equation}\label{eq:memoirs-weight-exact}
\begin{split}
&\big[\partial_{\xi\xi}+D_j\partial_{zz}
+(c_j^*(A)-Ab_j\sin(2\pi z))\partial_\xi
+r_j-2\gamma_j\big]\Psi=-2\gamma_j\Psi.
\end{split}
\end{equation}
Let $C=\mathcal N_{A,j}(0)$ and
$B(t,\xi,z)=C e^{-2\gamma_jt}\Psi(\xi,z)$.
Since $r_j>0$ and $w_1+w_2\geq0$,
\eqref{eq:memoirs-weight-exact} yields
\begin{equation*}
\begin{split}
B_t-\big[\partial_{\xi\xi}+D_j\partial_{zz}
+(c_j^*(A)-Ab_j\sin(2\pi z))\partial_\xi
+r_j-2\gamma_j-r_j(w_1+w_2)\big]B=r_j(w_1+w_2)B\geq0.
\end{split}
\end{equation*}
Both $B-h$ and $B+h$ are therefore supersolutions
of the scalar linear equation associated with
\eqref{eq:memoirs-h-eq}.  They are nonnegative initially;
in the bounded interval they are also nonnegative at the
longitudinal boundaries because $h=0$ there.
The parabolic comparison principle yields $|h|\leq B$.
On the whole line, division by $\Psi$ transforms the scalar
equation into one with bounded drift coefficients and zero-order
coefficient $-2\gamma_j-r_j(w_1+w_2)\leq-2\gamma_j$.
The bounded-solution maximum principle applies to the weighted
difference on each finite time cylinder.  Dividing by $\Psi$ proves
\eqref{eq:memoirs-mode-decay}.
\end{proof}

The weighted estimate is compatible with the computed
nonsymmetric initial data.  It does not assert a global
unweighted bound on $\mathbb R$ or convergence of the entire
cooperative solution to a critical front.  On every bounded
longitudinal interval it yields a local uniform bound with
a constant depending on that interval and the initial data.

\medskip
\noindent\textbf{Spatial and temporal discretization.}
Let $V_h$ be the space of continuous functions that are bilinear on the rectangular elements of $[-24,24]\times[0,1]$, periodic in $z$ and with the prescribed endpoint values in $\xi$. We use a uniform grid with $321$ longitudinal nodes and $32$ periodic transverse nodes, so that $\Delta\xi=0.15$ and $\Delta z=1/32$. The spatial diffusion and transport forms, for a test function $v_h$ vanishing at the longitudinal endpoints, are
\begin{align}\label{eq:app-weak}
&\int\nolimits_{[-24,24]\times[0,1]} (w_i)_t v_h
+\int\nolimits_{[-24,24]\times[0,1]}
\bigl((w_i)_\xi(v_h)_\xi+D_j(w_i)_z(v_h)_z\bigr)\notag\\
&\hspace{1cm}-\int\nolimits_{[-24,24]\times[0,1]}
(c_j^*-Ab_j\sin(2\pi z))(w_i)_\xi v_h
=\int\nolimits_{[-24,24]\times[0,1]}F_i(w_1,w_2)v_h.
\end{align}
We use the bilinear tensor-product finite-element basis with vertex-based nodal quadrature for the mass, reaction, diffusion and transport forms. The mass-lumped semidiscrete scheme, including central Galerkin transport, consequently has the following five-point representation at every interior node $(\xi_k,z_l)$:
\begin{align}\label{eq:app-fe-stencil}
(\mathcal L_{j,h} W)_{k,l}={}&\frac{W_{k-1,l}-2W_{k,l}+W_{k+1,l}}{h_\xi^2}
+D_j\frac{W_{k,l-1}-2W_{k,l}+W_{k,l+1}}{h_z^2}\notag\\
&+\bigl(c_j^*-Ab_j\sin(2\pi z_l)\bigr)
\frac{W_{k+1,l}-W_{k-1,l}}{2h_\xi},
\end{align}
with $l$ interpreted modulo $N_z$; the values at $k=0$ and $k=N_\xi-1$ are fixed by the longitudinal Dirichlet conditions. In particular, \eqref{eq:app-fe-stencil} is a nodally integrated, mass-lumped $Q_1$ finite-element operator, not an independently chosen finite-difference approximation. The corresponding element P\'eclet numbers
\begin{equation}\label{eq:app-pe}
\mathrm{Pe}_j=\frac{h_\xi}{2}\max_{z\in[0,1]}
|c_j^*-Ab_j\sin(2\pi z)|
\end{equation}
are approximately $0.673$ and $0.915$, both below one; hence the longitudinal off-diagonal coefficients of $\mathcal L_{j,h}$ are nonnegative.

Let $Y$ denote the vector of free interior nodal values of one component, let $L_{j,h}$ be the matrix in \eqref{eq:app-fe-stencil} restricted to those degrees of freedom and let $g_{j,h}$ collect the prescribed endpoint contributions. The linear diffusion--transport substep is computed by Crank--Nicolson:
\begin{equation}\label{eq:app-cn}
\bigl(I-\tfrac12\Delta t L_{j,h}\bigr)Y^{\mathrm{new}}
=\bigl(I+\tfrac12\Delta t L_{j,h}\bigr)Y^{\mathrm{old}}+\Delta t\,g_{j,h}.
\end{equation}
For the nonlinear system we apply a Heun half-step to the coupled reaction terms, the Crank--Nicolson substep \eqref{eq:app-cn} separately to each component, and a second Heun half-step. For sufficiently regular semidiscrete solutions, the composition is second-order accurate in time; this order statement does not provide an unconditional positivity theorem for Crank--Nicolson. A sparse LU factorization of the left-hand side of \eqref{eq:app-cn} is reused throughout each run. We take $\Delta t=0.01$, integrate to $T=10$, and extract the computed nodal solutions at $t=0,2,4,6,8,10$. No prescribed traveling-wave profiles or post-processing translations are used in the simulation plots.

\begin{figure}[!htbp]
\centering
\includegraphics[width=.98\textwidth]{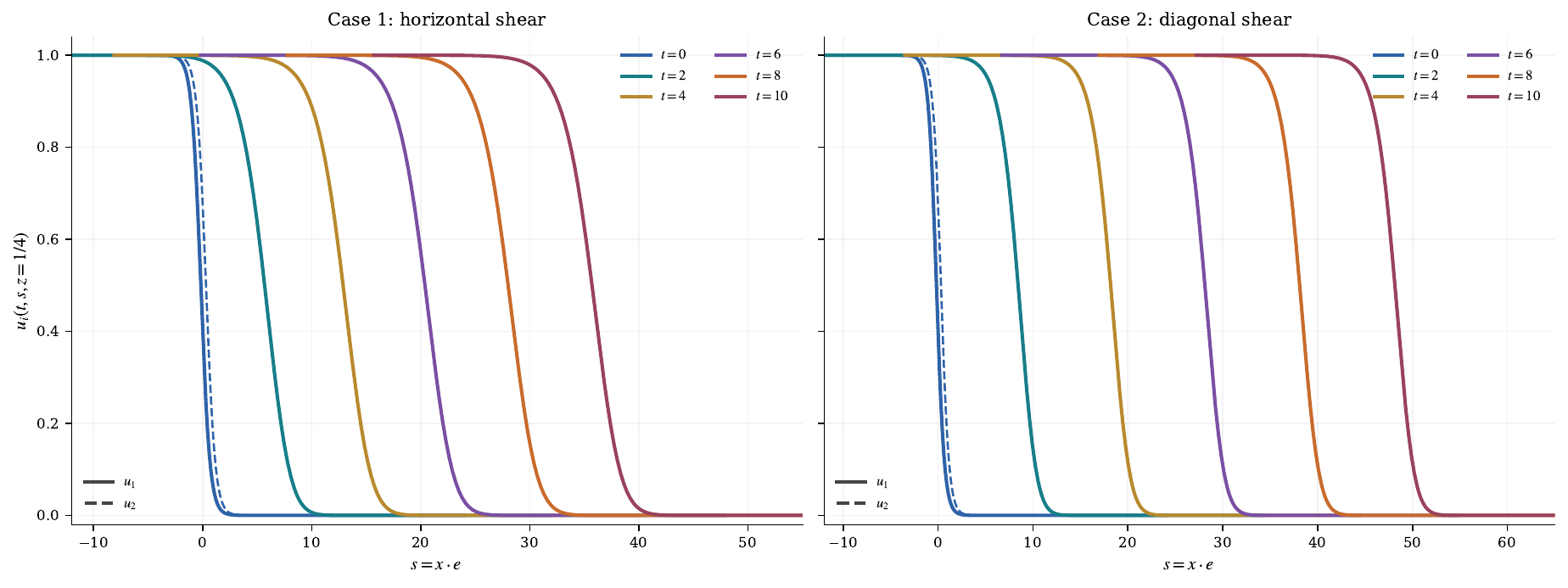}
\caption{Propagation of the cooperative system \eqref{eq:app-full-pde} at $z=1/4$. Panels (a) and (b) represent Cases~1 and~2, respectively. Each color identifies one observation time; solid and dashed curves are the two components. The physical longitudinal coordinate is $s=x\cdot e$. The observed shifts are consistent with the minimal-speed description of Theorem~\ref{minimalspeed}, but are finite-time computations.}
\label{fig:app-profiles}
\end{figure}
\begin{figure}[!htbp]
\centering
\includegraphics[width=.98\textwidth]{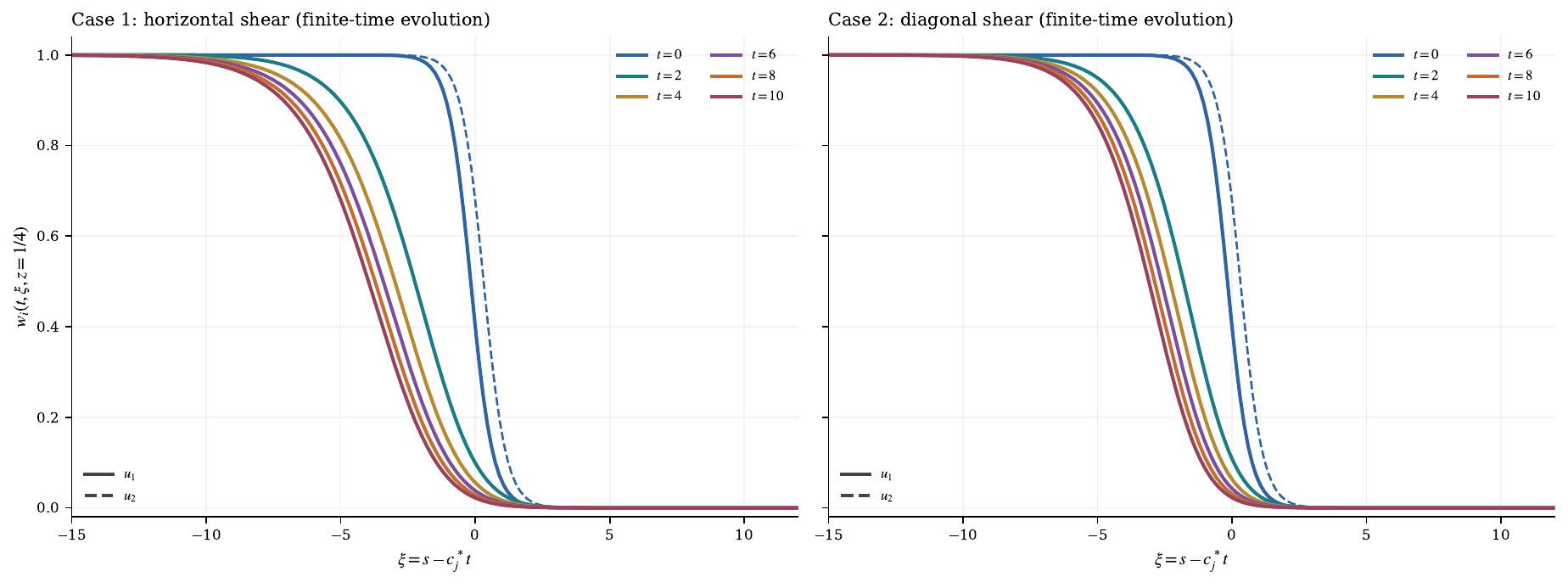}
\caption{The same solutions in the coordinates $\xi=s-c_j^*t$, with $c_j^*$ from \eqref{eq:app-cstar-A5}. Each panel contains several observation times on the section $z=1/4$; solid curves represent $w_1$ and dashed curves $w_2$. These finite-time profiles are observed in a frame moving at the minimal spectral speed; they are not asserted to be entire critical fronts.}
\label{fig:app-critical-snapshots}
\end{figure}
\begin{figure}[!htbp]
\centering
\includegraphics[width=.78\textwidth]{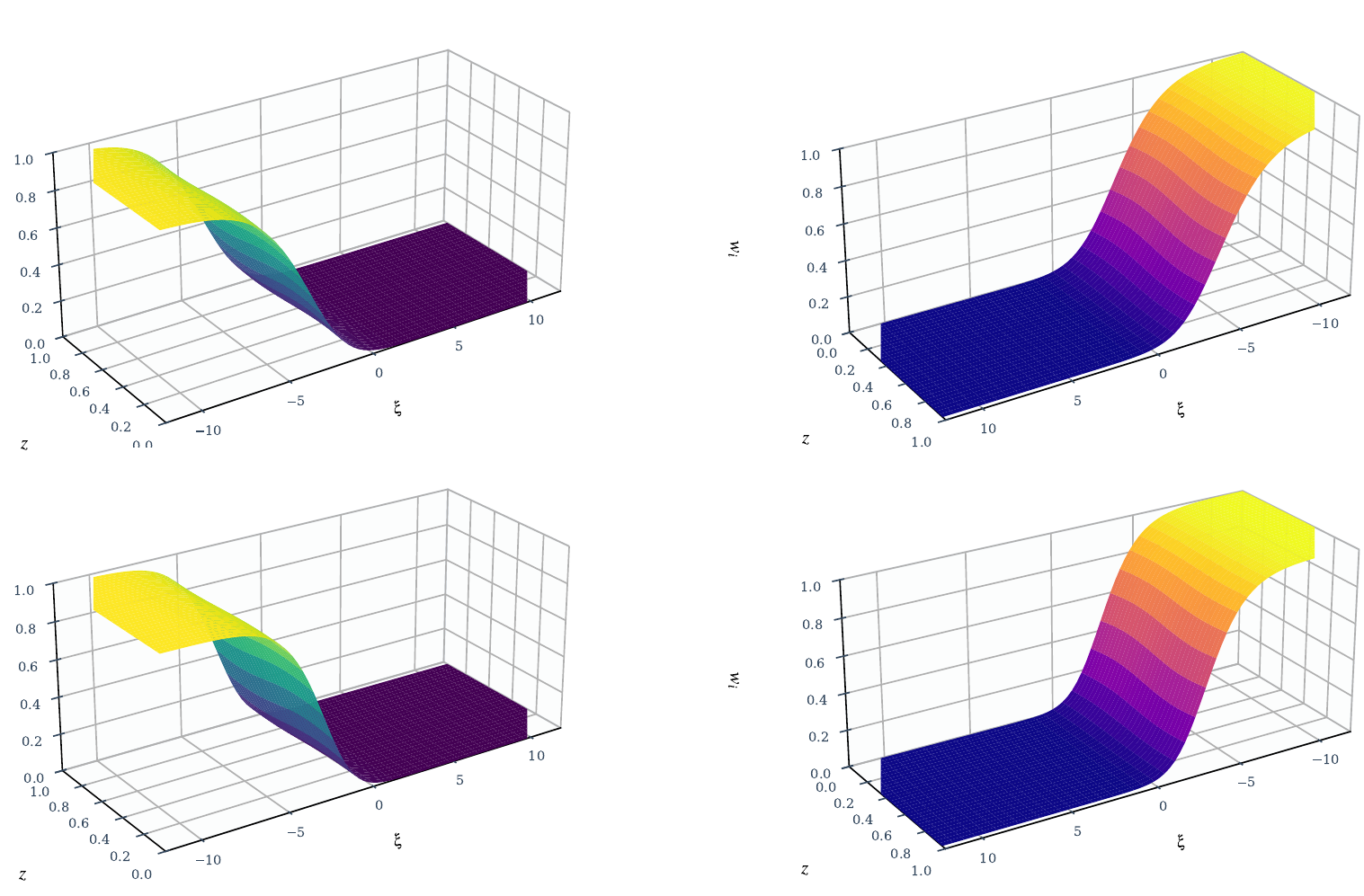}
\caption{Three-dimensional surfaces of the two solutions of \eqref{eq:app-moving-pde} at $t=10$ and $A=5$. The upper panels correspond to Case~1, the lower panels to Case~2; the left and right panels show $w_1$ and $w_2$, respectively. The coordinates are $\xi$ and periodic $z$. The two components are nearly equal at this observation time.}
\label{fig:app-pde-fronts}
\end{figure}

\medskip
\noindent\textbf{Interface diagnostics and numerical sensitivity.}
For each component, we locate the intersections with the level $w_i=1/2$ by linear interpolation in $\xi$, for every periodic node $z$. We then calculate the mean position and the range of these intersections in $z$. At $t=10$ the computed mean positions for the first component are $-4.097796$ in Case~1 and $-3.183715$ in Case~2. The respective transverse ranges for $w_1$ are $0.254014$ and $0.179427$. The final-time maximal nodal differences $\|w_1-w_2\|_\infty$ are approximately $1.8\times10^{-10}$ and $2.4\times10^{-15}$, respectively. This decay from the initially distinct profiles agrees with the symmetry of the coefficients and is not an imposed constraint.

\begin{figure}[!htbp]
\centering
\includegraphics[width=.76\textwidth]{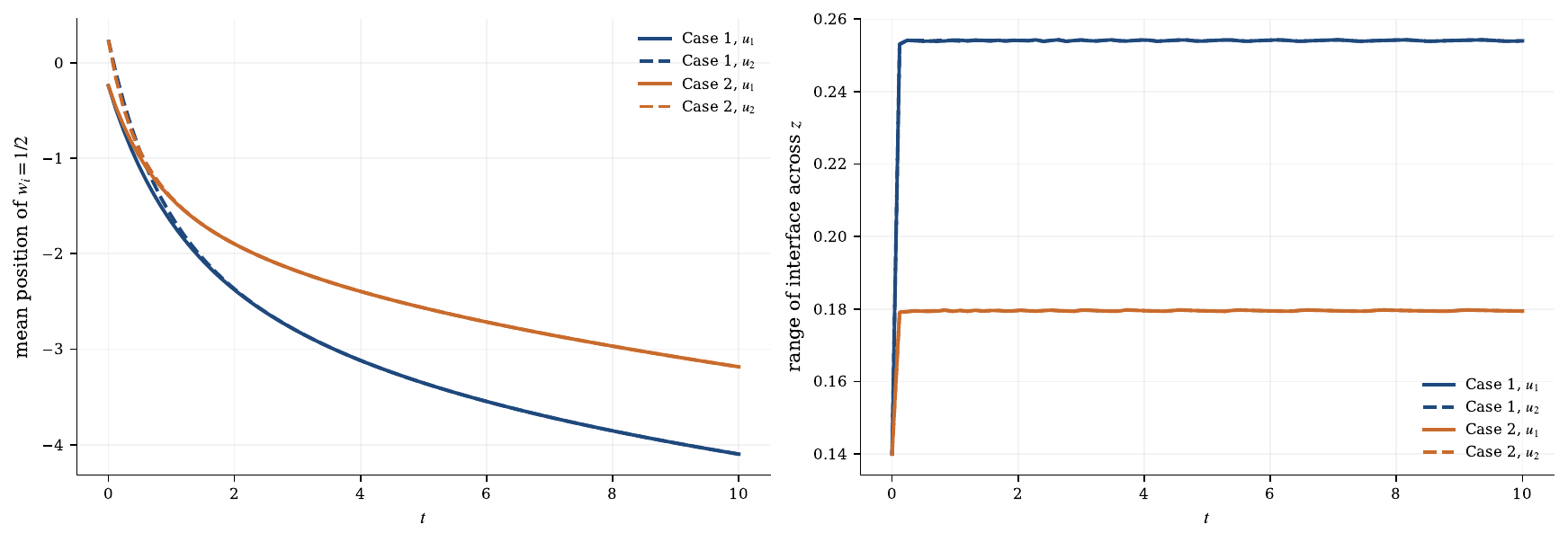}
\caption{Interface dynamics for the solutions of \eqref{eq:app-moving-pde}. Panel (a) shows the mean position of the level $w_i=1/2$ in the critical-speed frame; panel (b) shows the corresponding transverse range. Solid and dashed lines identify the first and second components. These observations complement the speed characterization in Theorem~\ref{minimalspeed} without establishing convergence to an entire critical front.}
\label{fig:app-pde-diagnostics}
\end{figure}

\begin{table}[!htbp]
\centering
\caption{Sensitivity of the mean $w_1=1/2$ position in the moving frame at $T=10$. The main and refined spatial grids use $\Delta t=0.01$; the time-refined computation uses the main grid with $\Delta t=0.005$. The values are direct numerical outputs, not certified error bounds.}
\label{tab:app-mesh}
\renewcommand{\arraystretch}{1.10}
\small
\begin{tabular}{@{}lrrrr@{}}\hline
& $321\times32$, $0.01$ & $361\times40$, $0.01$ & $321\times32$, $0.005$ & $|\Delta_t|$\\\hline
Case~1 & $-4.097796$ & $-4.068733$ & $-4.096106$ & $0.001689$\\
Case~2 & $-3.183715$ & $-3.113971$ & $-3.177467$ & $0.006248$\\\hline
\end{tabular}
\end{table}
The corresponding spatial differences are $0.029063$ and $0.069744$. Both temporal differences are smaller than in the first-order splitting, but the spatial sensitivity, especially in Case~2, remains visible. Therefore the results illustrate finite-time dynamics at the stated resolutions; they neither constitute a certified convergence study nor establish that the computed interfaces have reached an entire critical traveling front. The periodic spectral calculations and the exact solutions \eqref{eq:app-exact-wave} are independent of the nonlinear time-stepping method.

The three-dimensional time sequence in Figure~\ref{fig:app-multiple-fronts} uses the same computed solutions at their actual observation times. The numerical results illustrate the minimal-speed formula and its large-advection limit in two explicit periodic media. The nonlinear computations concern a bounded space--time region and do not enter the analytic proofs. In particular, the existence of entire critical transition fronts is established by the results of Section~\ref{sec:critical-proofs}, not by the stationarity of any plotted profile.
\FloatBarrier

\clearpage
\begin{figure}[p]
\centering
\vspace*{\fill}
\includegraphics[width=.99\textwidth,height=.84\textheight,keepaspectratio]{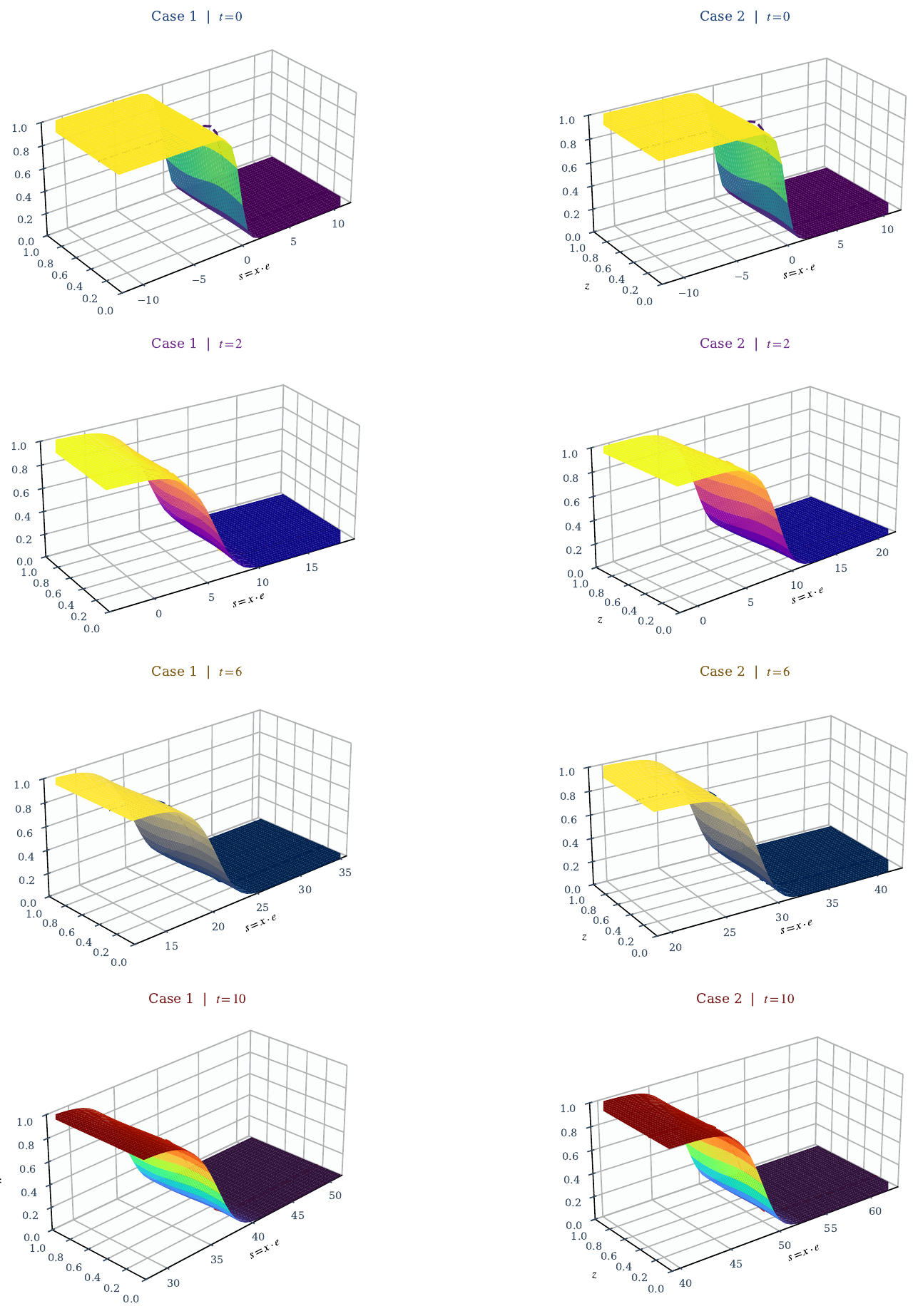}
\caption{Computed three-dimensional solutions of \eqref{eq:app-full-pde}. From top to bottom the rows show $t=0,2,6,10$; the left and right columns show Cases~1 and~2. Each colored surface is the independently evolved component $u_1$ in the physical coordinate $s=\xi+c_j^*t$; the dashed section at $z=1/4$ represents $u_2$.  The front positions are obtained directly from the computed nodal solutions, not from shifted copies of a fixed profile.}
\label{fig:app-multiple-fronts}
\vspace*{\fill}
\end{figure}
\clearpage

\begin{center}
\medskip\noindent\textbf{Declarations.}
\end{center}

\textbf{Conflict of interest statement:} The authors have no conflicts of interest to declare that are relevant
to the content of this article.

\textbf{Data availability statement:} No external data were used. The data in the numerical section were computed with the periodic finite-element eigenvalue scheme and the two-component nonlinear finite-element scheme described therein. The Python program and the numerical illustrations are supplied with the manuscript.

\endgroup

\begingroup\small\raggedbottom

\endgroup

\end{document}